\documentclass[a4paper,10pt]{amsart}
\usepackage[margin=3cm]{geometry}
\usepackage{a4wide,enumerate,xcolor,appendix,lscape}
\usepackage[colorlinks,linkcolor=blue,citecolor=blue]{hyperref}
\usepackage{amsmath,amssymb}
\usepackage{scrextend}
\usepackage{enumitem}
\usepackage{graphicx}
\allowdisplaybreaks
\numberwithin{equation}{section}

\let\pa\partial
\let\na\nabla
\let\eps\varepsilon
\newcommand{\N}{{\mathbb N}}
\newcommand{\R}{{\mathbb R}}
\newcommand{\diver}{\operatorname{div}}

\newcommand{\HS}{{\mathcal L}_{2}}
\renewcommand{\d}{\,\textnormal{d}}
\newcommand{\dt}{\,\textnormal{d}t}
\newcommand{\ds}{\,\textnormal{d}s}
\newcommand{\dr}{\,\textnormal{d}r}
\newcommand{\dx}{\,\textnormal{d}x}
\newcommand{\dy}{\,\textnormal{d}y}

\newcommand{\dd}{\mathrm{d}}
\newcommand{\E}{{\mathbb E}}
\newcommand{\F}{{\mathbb F}}
\newcommand{\Prob}{{\mathbb P}}

\newcommand{\br}[1]{\langle #1\rangle}
\newcommand{\Ufine}{\mathcal{U}}
\newcommand{\Hmid}{\mathcal{H}}
\newcommand{\Usharp}{\mathcal{U}^{\sharp}}
\newcommand{\Vcoarse}{\mathcal{V}}
\newcommand{\fluc}{V}
\newcommand{\gold}{\mathfrak g}

\newtheorem{theorem}{Theorem}[section]
\newtheorem{lemma}[theorem]{Lemma}
\newtheorem{proposition}[theorem]{Proposition}
\newtheorem{corollary}[theorem]{Corollary}

\theoremstyle{definition}
\newtheorem{definition}[theorem]{Definition}
\newtheorem{assumption}[theorem]{Assumption}
\newtheorem{remark}[theorem]{Remark}
\newtheorem{condition}[theorem]{Condition}

\providecommand{\Fcoarse}{\mathcal W}
\providecommand{\tr}{\operatorname{tr}}
\providecommand{\Q}{\mathbb Q}
\providecommand{\wt}{\langle y\rangle}
\providecommand{\Test}{\mathcal T}
\providecommand{\dz}{\,\mathrm dz}\providecommand{\dtau}{\,\mathrm d\tau}
\newif\ifmarkchanges
\markchangesfalse

\begin{document}
\title[Gaussian fluctuations of the SKT particle system]
{Gaussian fluctuations of the moderately interacting\\
Shigesada--Kawasaki--Teramoto particle system\\
around the mollified mean-field limit}
\author{Florian Huber}
\date{\today}
\maketitle

\begin{abstract}
We study, at the logarithmic mollification scale $\eta^{-2(d+1+\alpha)}\le\delta\log N$ of
\cite{chen_holzinger_2021_rigorous_SKT}, the fluctuation field $\fluc^N=\sqrt N(\mu^N-u_\eta)$ of
the moderately interacting multi-species particle system whose mean-field limit is the
Shigesada--Kawasaki--Teramoto (SKT) cross-diffusion system. The field is centred at the
solution $u_\eta$ of the mollified equation, which depends on $N$ through $\eta$.
On a scale of weighted Bessel-potential spaces of negative order we prove that the
laws of $(\fluc^N,M^N)$ are tight, where $M^N$ is the martingale part. Every limit
point of $M^N$ is a continuous Gaussian martingale with an explicit deterministic
bracket. Every limit point of $\fluc^N$ is a weak solution, and the mild solution,
of a linear Ornstein--Uhlenbeck SPDE whose drift is the linearisation of the SKT
operator at the limit $u$. The initial field and the martingale are independent, and
weak solutions are unique given both, so the whole sequence converges in law, jointly
with the martingale part, and the limit is unique in law among solutions of the martingale
problem. The
noise is $\dd M_i=\diver\big(\sqrt{2u_i\Theta_i}\,\dd W_i\big)$, with
$\Theta_i=\sigma_i+\sum_jf(a_{ij}u_j)$ and $W$ a cylindrical Wiener process on
$L^2(\R^d;\R^{nd})$. The hypotheses are those of
\cite{chen_holzinger_2021_rigorous_SKT}, including the smallness and the nonnegativity on $\R$ of the nonlinearity
$f$. We add radial kernels $B_{ij}$, the independence
of the initial data and the Brownian motions across particles and species, a smaller
$\delta$, a regularity budget on $u_0$, a strengthened smallness condition on $f$ that
makes the backward system diagonally dominant in weighted spaces, and an integer Sobolev
index $s$.
\end{abstract}

\tableofcontents

\section{Introduction}\label{sec:intro}

\subsection*{The system and the field}

Fix $n\ge1$ species. For each $i\le n$ the environmental potential is
$U_i(x)=-\frac12|x|^2$, as in \cite[(A2)]{chen_holzinger_2021_rigorous_SKT}. It is
dispersive: the drift $-\na U_i(x)=x$ pushes the particles outward. The mollifiers are
$B^\eta_{ij}(x)=\eta^{-d}B_{ij}(|x|/\eta)$. The particle system
\eqref{eq:particles} consists of $N$ diffusions per species. Species $j$ acts on
species $i$ through $f_\eta(B^\eta_{ij}*\mu^N_j)$, a nonlinear function of a
mollified empirical measure, which enters the diffusion coefficient. In the moderate
regime $\eta=\eta(N)\to0$ the empirical measures converge to the solution $u$ of the
SKT system \eqref{eq:skt_det} \cite{chen_holzinger_2021_rigorous_SKT}. We take the
hypotheses of that paper as given, together with its intermediate mollified system
\eqref{eq:nonloc}. This paper studies the next order,
\[
 \fluc^N_i\ :=\ \sqrt N\big(\mu^N_i-u_{\eta,i}\big).
\]

\subsection*{Why the centring is at $u_\eta$ and not at $u$}

The difference $\sqrt{N}(u_\eta-u)$ is deterministic. By
\cite[Thm.~3]{chen_holzinger_2021_rigorous_SKT},
$\|u_\eta-u\|_{L^\infty(0,T;L^2)}=O(\eta)$, see \eqref{eq:CHbias}. This bound gives
$\sqrt N\|u_\eta-u\|\to0$ only when $\eta=o(N^{-1/2})$, or in other words, when $\eta$ decays at least like a power of $N$. The scale used here, \ref{H:scaling}, is logarithmic,
$\Lambda:=\eta^{-2(d+1+\alpha)}\le\delta\log N$, and it excludes every power
$\eta=N^{-\varpi}$, $\varpi>0$. We therefore centre at $u_\eta$, as \cite{CHJ26} do. They
motivate this choice by \cite{Oel87}. We do not treat centring at $u$. A deterministic
correction of fixed order does not help at this scale, by the following heuristic. Suppose
that $u_\eta$ has an asymptotic expansion $u+\sum_{k\ge1}\eta^{2k}u_k$ in $L^2$, which we do not
prove. Fix $K\ge1$ and centre at $u+\sum_{1\le k<K}\eta^{2k}u_k$. If $K'\ge K$ is the first index
with $u_{K'}\ne0$, the centring leaves an error of exact order $\eta^{2K'}$. Under \ref{H:scaling},
$\eta\ge(\delta\log N)^{-1/(2(d+1+\alpha))}$, so $\eta^{2K'}$ decays only like a power of $\log N$
and $\sqrt N\eta^{2K'}\to\infty$.

\subsection*{Results}

All hypotheses are collected in Assumption~\ref{ass:main} of Section~\ref{sec:setting}.
The field and its martingale part live on the weighted Bessel-potential spaces
$\Ufine=H^{-s_1}$ and $\Hmid=H^{-s_1-3}_{-1}$ of negative order, with
$s_1=\lfloor\frac d2\rfloor+2$ (see \eqref{eq:Hsa} and \eqref{eq:spaces}). The mobilities are
$\Theta_i=\sigma_i+\sum_jf(a_{ij}u_j)$, and $\mathbf A_{ij}(u)=\delta_{ij}\Theta_i+a_{ij}u_if'(a_{ij}u_j)$ is
the diffusion matrix of the SKT system linearised at $u$.

\begin{theorem}[Main theorem]\label{thm:main}
Under Assumption~\ref{ass:main}, let $M^N$ be the martingale part of $\fluc^N$
\textup{(}\eqref{eq:MNval}\textup{)}. Then:
\begin{enumerate}[label=\normalfont(\roman*)]
\item the laws of $(\fluc^N,M^N)$ are tight on
$\mathcal Z_T=\big(C([0,T];\Hmid)\cap L^2_w(0,T;\Ufine)\big)^n\times C([0,T];\Hmid)^n$;
\item $(\fluc^N,M^N)\Rightarrow(\fluc,M)$ in $C([0,T];\Hmid^n)^2$. The limit $M$ is a continuous
centred Gaussian martingale with independent increments and
\[
 \big\langle M_i(\cdot,\varphi),M_k(\cdot,\psi)\big\rangle_t
 =\delta_{ik}\int_0^t\br{u_i,\ 2\Theta_i\,\na\varphi\cdot\na\psi}\dr
\]
for $\varphi,\psi\in\Ufine^*$. Moreover $M$ takes values in $\Ufine^n$. The initial value
$\fluc(0)=\fluc^0$ is the centred Gaussian field of Proposition~\ref{prop:initialCLT}, and
$\fluc^0$ and $M$ are independent. The limit $\fluc$ is the unique mild solution, and a weak solution, of the
Ornstein--Uhlenbeck equation
\[
 \d\fluc_i=\Big(\diver(\fluc_i\na U_i)+\Delta\sum_{j=1}^n\mathbf A_{ij}(u)\fluc_j\Big)\dt+\d M_i,
 \qquad \fluc_i(0)=\fluc^0_i ;
\]
\item the law of $\fluc$ is the unique law of a solution of the martingale problem of this
equation with initial law $\mathrm{Law}(\fluc^0)$ \textup{(}Definition~\ref{def:mgp}\textup{)};
\item on an extension of the probability space there is an
$L^2(\R^d;\R^{nd})$-cylindrical Wiener process $W$ with
$M_i(t)=\int_0^t\diver\big(\sqrt{2u_i\Theta_i}\,\d W_i\big)$, $i\le n$.
\end{enumerate}
\end{theorem}

Section~\ref{sec:assembly} gives the proof and lists the results it uses. The main steps are
the following. Proposition~\ref{prop:tight} gives tightness. Theorem~\ref{thm:noise} shows
that every limit point of the martingale part is a continuous Gaussian martingale with a
deterministic bracket. Theorem~\ref{thm:OU} shows that every limit point of $\fluc^N$ is a
weak solution, and the mild solution, of the limit equation, whose drift is the
linearisation of the SKT operator at $u$ (Theorem~\ref{thm:G}(ii)).
Theorem~\ref{thm:uniqueness} gives uniqueness of weak solutions with a given datum and a
given martingale. Corollary~\ref{cor:whole} shows that the datum and the martingale of a
limit point are independent, so that the whole sequence converges, jointly with $M^N$.
Corollary~\ref{cor:mgp} gives uniqueness in law. Theorem~\ref{thm:representation}
represents the noise. It is a direct consequence of Proposition~\ref{prop:rep}, which the
construction of the stochastic convolution in Section~\ref{sec:limitmart} also uses.
Lemma~\ref{H:W} is \cite[Prop.~4(i)]{chen_holzinger_2021_rigorous_SKT}:
existence and pathwise uniqueness for the particle system and for the intermediate
system. Lemma~\ref{H:U} is \cite[Thm.~2]{chen_holzinger_2021_rigorous_SKT}. In
Lemma~\ref{H:chaos}, the law of an intermediate particle at time $t$ is
$u_\eta(t)\dx$. We prove this at every $t$, in Lemma~\ref{lem:lawid}, by a duality
argument with the backward Kolmogorov equation. \cite[Thms.~7--8]{chen_holzinger_2021_rigorous_SKT}
state it for almost every $t$, with a proof given in outline, we do not use them. The existence and regularity of $u_\eta$ and $u$ are
\cite[Thms.~2--3]{chen_holzinger_2021_rigorous_SKT}, and the bias estimate
\eqref{eq:CHbias} is \cite[Thm.~3]{chen_holzinger_2021_rigorous_SKT}.

\emph{The logarithmic scale.} The scale enters through our own coupling estimates,
Proposition~\ref{prop:phi4} and \eqref{eq:phi2}: their Gr\"onwall factor $e^{C\Lambda T}$ is
polynomial in $N$ only when $\Lambda\le\delta\log N$. Remark~\ref{rem:betaceiling} lists the
other conditions the proof places on $\eta$. With
$\alpha=0$ and $\delta=10^{-2}$, the condition $\eta\le\frac12$ already requires
$\log N\ge1.6\cdot10^3$ in $d=1$ and $\log N\ge6.4\cdot10^3$ in $d=2$, so the result
is asymptotic in a strong sense.

\emph{Classical inputs.} Standard results are cited where they are used. A few results
are proved again because the available versions do not fit the present setting, the reason is
given at the point of use. The main ones are the Moser estimate (Lemma~\ref{lem:moser}), in a form with the
dependence on the range of the argument made explicit. The tightness criterion
(Proposition~\ref{prop:criterion}), for a triple of weighted Bessel spaces and the passage from weak
to mild solutions (Theorem~\ref{thm:weakmild}), for a time-dependent evolution family. The backward
evolution family is obtained from J.-L.~Lions' theorem after a conjugation that removes the
unbounded drift (Appendix~\ref{sec:appevol}). The remaining results are proved here. In particular the fixed-time fourth moment
\eqref{eq:K4-mild}, which Proposition~\ref{prop:close} uses, is Theorem~\ref{thm:K4}. Remark~\ref{rem:chain} lists the results the
proof of Theorem~\ref{thm:main} uses.

\subsection*{Related work, and what is new}

The closest work is \cite{CHJ26}, by three of the four authors of
\cite{chen_holzinger_2021_rigorous_SKT}. They prove a central limit theorem for a single population of
moderately interacting particles in $\R^d$, $d\ge3$, whose interaction potentials approximate
attractive sub-Coulomb Riesz potentials. Their normalisation is $\sqrt N$ and their
centring is the intermediate mollified solution, as here. They describe their method
as inspired by Oelschl\"ager's work on fluctuations, and motivate the centring by
\cite{Oel87}.

\emph{What is new compared with \cite{CHJ26}.} The two papers differ in the particle system, in
the form of the result and in the method.
\begin{itemize}
\item \emph{The model.} In \cite{CHJ26} the particles interact through the drift,
$\frac\kappa N\sum_l\na V^\eta(X^k-X^l)$, and the noise $\sqrt{2\sigma}\,\d W^k$ is additive. Here
$n$ species interact through the diffusion coefficient: the intensity $\sqrt{2\Theta^\eta_i}$ of
the noise of a particle of species $i$ depends on the mollified empirical measures of all species. As a consequence the
fluctuation martingale has a state-dependent bracket, the noise of the limit,
$\diver(\sqrt{2u_i\Theta_i}\,\d W_i)$, depends on $u$, and the principal part of the limit drift is
the non-diagonal diffusion matrix $\mathbf A(u)$ of the linearised SKT operator instead of a
constant multiple of $\Delta$. Its coercivity is not automatic. We obtain it from the diagonal
dominance \eqref{eq:epscond}. We also cover every dimension $d\ge1$.
\item \emph{The result.} \cite{CHJ26} prove convergence of conditional characteristic functions of
the fluctuation field tested against $\varphi$, for each test function $\varphi$ and each time $t$, and the Ornstein--Uhlenbeck
SPDE for the limit is stated there formally. We prove tightness of $(\fluc^N,M^N)$ on a path space,
convergence in law of the whole sequence jointly with its martingale part, and identify the limit as the
unique mild solution, which is also a weak solution, of the Ornstein--Uhlenbeck SPDE. We also prove uniqueness in law for the
martingale problem and a representation of the noise by a cylindrical Wiener process
(Theorem~\ref{thm:main}).
\item \emph{The method.} A process-level statement needs more than the dual problem for one test
function at a time. We construct the backward linearised system as an evolution family on weighted
Bessel spaces of negative order, uniformly in $\eta$, for a nonlocal system with the unbounded
dispersive drift $x\cdot\na$ (Section~\ref{sec:evolfam} and Appendix~\ref{sec:appevol}). We prove a
tightness criterion for a triple of weighted spaces (Appendix~\ref{sec:appaldous}), and we control the
remainders with a fourth-order version of the coupling of \cite{chen_holzinger_2021_rigorous_SKT}
(Appendix~\ref{sec:appk4}).
\item \emph{The price.} \cite{CHJ26} reach polynomial mollification $\eta=N^{-\varpi}$,
$0<\varpi<1/(8\lambda+12)$, for $0<\lambda<d-2$, with initial data in $W^{K,1}\cap W^{K,\infty}$,
$K>1/(2\varpi)$. This paper reaches only the logarithmic scale \ref{H:scaling}, with the fixed
regularity budget $u_0\in H^s$, $s\ge\lfloor\frac d2\rfloor+6$.
\end{itemize}
Neither result contains the other.

\emph{Relation to \cite{Hub25}.} In \cite{Hub25} the author studies a stochastic SKT system
with a multiplicative conservative noise of intensity $N^{-1/2}$ and an entropy-based correction term, on a
bounded domain and with coloured noise, and proves the existence of global nonnegative martingale
solutions. There the noise is motivated by the particle system of
\cite{chen_holzinger_2021_rigorous_SKT}, but it is not derived from it. Theorem~\ref{thm:main}
provides a derivation of this type at the level of Gaussian fluctuations, on $\R^d$ and with space-time
white noise. The limit noise $\diver(\sqrt{2u_i\Theta_i}\,\d W_i)$ has the structure of the noise
in \cite{Hub25}, with the mobility $\Theta_i=\sigma_i+\sum_jf(a_{ij}u_j)$ in place of the linear SKT
mobility used there. Formally the two agree for $f(s)=s$, which is not covered by Assumption~\ref{ass:main}
(Remark~\ref{rem:linearf}). The present result does not single out the correction term of \cite{Hub25}, which is of
order $N^{-1}$ and does not affect the Gaussian fluctuations.

Oelschl\"ager \cite{Oel87} proved a Gaussian fluctuation theorem for
moderately interacting diffusions. There a single species interacts through the drift,
$-\frac1{2N}\sum_l\na V_N(X^k-X^l)$, the diffusion coefficient is constant, and the
limit is the nonlinear diffusion equation
$\pa_tp=\frac12\na\cdot((1+p)\na p)$. The mollification is polynomial,
$V_N=\chi_N^dV_1(\chi_N\cdot)$ with $\chi_N=N^{\varpi/d}$ and
$\varpi\in(0,\frac d{2(d+2)})$ in the fluctuation theorem \cite[Thm.~2]{Oel87}. The field
is centred at the limit plus a deterministic correction $C_N$ with $C_N\to0$, and the
conclusion is the convergence of conditional characteristic functions of
$\br{Y_N(t),f}$ for each test function $f$. \cite[Remark~g]{Oel87} interprets it as an
asymptotic Ornstein--Uhlenbeck process in $\mathcal S'(\R^d)$. \cite{Oel89} is his
derivation of reaction--diffusion systems, a law of large numbers. Jourdain and M\'el\'eard \cite{JoMe98} treat the single-species
version of the present particle system, with the interaction in the diffusion
coefficient. \cite{chen_holzinger_2021_rigorous_SKT} recover their scaling at
$\alpha=0$. Their fluctuation result \cite[Thm.~3.14]{JoMe98} is in $d=1$, with a
normalisation smaller than $\sqrt N$, and its limit is deterministic, in
$C([0,T];W^{-4,1}_0)$. The weighted Hilbert-space approach to
fluctuations with an Ornstein--Uhlenbeck limit goes back to Fernandez and
M\'el\'eard \cite{FeMe97} for the McKean--Vlasov model, and the method here is of
that type.

The new points are the following. The setting is multi-species SKT cross-diffusion, with
the interaction in the diffusion coefficient, under the hypotheses of
\cite{chen_holzinger_2021_rigorous_SKT}. For this setting, the conclusion has two
properties that \cite{Oel87}, \cite{CHJ26} and \cite{JoMe98} do not have together: the
normalisation $\sqrt N$ with a Gaussian limit, and convergence of processes. The convergence is in law on $C([0,T];\Hmid^n)$ to the unique solution of an
Ornstein--Uhlenbeck SPDE, with a representation of its noise. In the setting of
\cite{Oel87}, the conditional characteristic functions \cite[(2.46)]{Oel87} already describe
the limit as an Ornstein--Uhlenbeck process in $\mathcal S'(\R^d)$ \cite[Remark~g]{Oel87}. The
technical contributions are three: the path-space tightness and the SPDE formulation,
a fourth-order version of the coupling of \cite{chen_holzinger_2021_rigorous_SKT}
(Appendix~\ref{sec:appk4}) and a backward evolution family on weighted spaces, uniform in
$\eta$, for a nonlocal system with the dispersive drift $x\cdot\na$
(Appendix~\ref{sec:appevol}). The scale is logarithmic, the centring depends on $N$, and
under the smallness hypothesis the limit equation is a small perturbation of $n$
decoupled Ornstein--Uhlenbeck equations. The coercivity of $\mathbf A(u)$ comes from diagonal dominance, a property of its
symmetric part: under the smallness condition \eqref{eq:epscond}, the off-diagonal terms in the
energy estimate \eqref{eq:S4} of the backward system are at most a quarter of the diagonal
dissipation. The pointwise nonlinearity $f$ is handled by the usual
Taylor expansion at $B^\eta_{ij}*u_{\eta,j}$.

\subsection*{Method, in one paragraph}

We prove tightness of the laws of $(\fluc^N,M^N)$ on a weighted negative Sobolev
scale. The Skorokhod--Jakubowski theorem gives almost-sure representatives on a new
probability space. We pass to the limit in the exact finite-$N$ identity
\eqref{eq:fluc_decomp} for $\br{\fluc^N_i,\varphi}$, identify the drift as the
linearised SKT operator and the bracket of the martingale part as the predicted
covariance, and represent the limit martingale. The scale is a chain
$\Ufine\hookrightarrow\Usharp\hookrightarrow\hookrightarrow\Hmid\hookrightarrow\hookrightarrow\Vcoarse$
of weighted Bessel spaces of negative order, \eqref{eq:chainA}. The field is bounded in
the two finest spaces, in $\Usharp$ uniformly in time and in $\Ufine$ in the time
integral. Its time increments are estimated in the coarsest space $\Vcoarse$, and those
of the martingale in $\Ufine$. The compact embeddings turn these bounds into tightness,
and convergence in law is stated in the middle space $\Hmid$.

\subsection*{Structure}

Section~\ref{sec:setting} fixes the hypotheses, the weighted scale, the particle system
and the elementary estimates. Section~\ref{sec:field} derives the exact finite-$N$
identity and the properties of $M^N$. Section~\ref{sec:evolfam} records the backward
evolution family, constructed in Appendix~\ref{sec:appevol}.
Section~\ref{sec:estimates} proves the uniform estimates, the convergence of the drift
coefficient, and that the three remainders vanish. Section~\ref{sec:tight} proves
tightness, with the criterion of Appendix~\ref{sec:appaldous}, and ends with the
Skorokhod representation. Section~\ref{sec:identify} identifies the drift and the
bracket, and with them the limit martingale (Theorem~\ref{thm:noise}).
Section~\ref{sec:limitmart} constructs the stochastic convolution against that
martingale and proves the initial CLT. Section~\ref{sec:limiteq} proves
that weak solutions are mild (Theorem~\ref{thm:weakmild}), and gives Theorem~\ref{thm:OU}, uniqueness,
convergence of the whole sequence, jointly with the martingale, and uniqueness in law.
Section~\ref{sec:assembly} identifies the noise and proves Theorem~\ref{thm:main}.
Appendix~\ref{sec:appaldous} gives the tightness criterion, Appendix~\ref{sec:appevol} the
evolution family for an abstract coefficient pair, and Appendix~\ref{sec:appk4} the
fixed-time fourth moment that Proposition~\ref{prop:close} uses.

\section{The setting and the particle system}\label{sec:setting}

\noindent
This section introduces the objects of the paper in the order in which they depend on one
another. Everything is stated on a fixed time horizon $[0,T]$ and for a fixed number $n$ of
species. We first fix the data, then the two deterministic systems, the SKT system and its
nonlocal approximation, then the three particle systems, two of which have these as their
one-particle laws, and the logarithmic scale that couples $N$ and $\eta$. The hypotheses are collected afterwards in
Assumption~\ref{ass:main}, and the section ends with the weighted scale of spaces and the
elementary lemmas used throughout.

\medskip
\noindent\emph{The data.} For each species $i\le n$ we are given a diffusivity $\sigma_i>0$ and the environmental potential
$U_i(x)=-\frac12|x|^2$, whose drift $-\na U_i(x)=x$ is dispersive. For each pair $(i,j)$ we are
given a radial mollifier $B_{ij}$ with mass $a_{ij}=\int_{\R^d}B_{ij}$, and its rescaling
$B^\eta_{ij}(x)=\eta^{-d}B_{ij}(|x|/\eta)$, which has the same mass. The interaction acts through
a nonlinearity $f\ge0$ and its cut-off $f_\eta\ge0$, which agrees with $f$ on $[-a_\eta,a_\eta]$
and has Lipschitz constant at most $\eta^{-\alpha}$. The random data are initial positions
$\xi^k_i$, $k\le N$, with density $u_{0,i}$, and Brownian motions $W^k_i$ in $\R^d$. The precise
requirements are \ref{H:CH} and \ref{H:cutoff} of Assumption~\ref{ass:main}.

\medskip
\noindent\emph{The two deterministic systems.} The SKT system is
\begin{equation}\label{eq:skt_det}
 \pa_t u_i=\diver(u_i\na U_i)+\Delta\Big(\sigma_iu_i+u_i\sum_{j=1}^nf(a_{ij}u_j)\Big),
 \qquad u_i(0)=u_{0,i}\quad\text{in }\R^d,\ t>0.
\end{equation}
\noindent Its nonlocal approximation, in which the local density $a_{ij}u_j$ in the argument of
$f$ is replaced by the mollified density $B^\eta_{ij}*u_{\eta,j}$ and $f$ by $f_\eta$, is
\begin{equation}\label{eq:nonloc}
 \pa_tu_{\eta,i}=\diver(u_{\eta,i}\na U_i)
 +\Delta\Big(\sigma_iu_{\eta,i}+u_{\eta,i}\sum_{j=1}^nf_\eta(B^\eta_{ij}*u_{\eta,j})\Big),
 \qquad u_{\eta,i}(0)=u_{0,i}.
\end{equation}
\noindent Under Assumption~\ref{ass:main} both systems have unique global solutions, bounded in
$H^s$ uniformly in $\eta$ by \ref{H:small}, and they are close: the bias estimate of
\cite[Thm.~3]{chen_holzinger_2021_rigorous_SKT} is
\begin{equation}\label{eq:CHbias}
 \|u_\eta-u\|_{L^\infty(0,T;L^2(\R^d))}
 +\|\na(u_\eta-u)\|_{L^2(0,T;L^2(\R^d))}\le C(T)\,\eta.
\end{equation}

\medskip
\noindent\emph{The particle system and the mobilities.} The particle system consists of $N$ diffusions per species. Species $j$ acts on species $i$
through $f_\eta(B^\eta_{ij}*\mu^N_j)$, evaluated along the empirical measures
$\mu^N_j:=\frac1N\sum_k\delta_{X^{N,\eta}_{k,j}}$, and this quantity enters the diffusion
coefficient. The dynamics are given by the It\^o SDEs
\begin{equation}\label{eq:particles}
 \d X_{k,i}^{N,\eta}=-\na U_i(X_{k,i}^{N,\eta})\dt
 +\Big(2\sigma_i+2\sum_{j=1}^n f_\eta\big(\tfrac1N\!\!\!\sum_{\substack{\ell=1\\(\ell,j)\neq(k,i)}}^N\!\!\!
 B_{ij}^\eta(X_{k,i}^{N,\eta}-X_{\ell,j}^{N,\eta})\big)\Big)^{1/2}\d W_i^k(t),
 \quad X_{k,i}^{N,\eta}(0)=\xi_i^k.
\end{equation}
A particle does not interact with itself: the inner sum omits $(\ell,j)=(k,i)$. The
mobilities of the particle system, of \eqref{eq:nonloc} and of \eqref{eq:skt_det} are
\begin{equation}\label{eq:mobilities}
\begin{aligned}
 \Theta_i^N(x)&:=\sigma_i+\sum_{j=1}^nf_\eta\big((B^\eta_{ij}*\mu_j^N)(x)\big),
 &\Theta_i^\eta(x)&:=\sigma_i+\sum_{j=1}^nf_\eta\big((B^\eta_{ij}*u_{\eta,j})(x)\big),\\
 \Theta_i(x)&:=\sigma_i+\sum_{j=1}^nf\big(a_{ij}u_j(x)\big),
\end{aligned}
\end{equation}
and we denote $w^\eta_{ij}:=B^\eta_{ij}*u_{\eta,j}$, so that $\Theta^\eta_i=\sigma_i+\sum_jf_\eta(w^\eta_{ij})$. Then \eqref{eq:nonloc} and \eqref{eq:skt_det} take the form
$\pa_tu_{\eta,i}=\diver(u_{\eta,i}\na U_i)+\Delta(\Theta^\eta_iu_{\eta,i})$ and
$\pa_tu_i=\diver(u_i\na U_i)+\Delta(\Theta_iu_i)$, and the first has the weak form
\begin{equation}\label{eq:weak_nonloc}
 \frac{\dd}{\dt}\br{u_{\eta,i},\varphi}
 =\br{u_{\eta,i},\,-\na U_i\cdot\na\varphi+\Theta^\eta_i\Delta\varphi}.
\end{equation}
The diffusion coefficient $\Sigma_{k,i}$ of the particle system \eqref{eq:particles} is
given by
\begin{equation}\label{eq:Sig}
 \Sigma^2_{k,i}=2\sigma_i+2\sum_{j=1}^n f_\eta\big(G^N_{k,ij}\big),
 \qquad G^N_{k,ij}:=\frac1N\!\!\sum_{(\ell,j)\ne(k,i)}\!\! B^\eta_{ij}\big(X^{N,\eta}_{k,i}-X^{N,\eta}_{\ell,j}\big).
\end{equation}
It differs from $2\Theta^N_i(X^{N,\eta}_{k,i})$ only in the term $j=i$, where the inner
sum omits $\ell=k$, so that $\Sigma^2_{k,i}=2\Theta^N_i(X^{N,\eta}_{k,i})+2\varrho^N_{k,i}$ with the
self-interaction defect
\begin{equation}\label{eq:self}
 \varrho_{k,i}^N:=f_\eta\big((B^\eta_{ii}*\mu_i^N)(X_{k,i})-\tfrac1N\eta^{-d}B_{ii}(0)\big)
 -f_\eta\big((B^\eta_{ii}*\mu_i^N)(X_{k,i})\big),
 \qquad |\varrho^N_{k,i}|\le \frac{C}{N\eta^{d+\alpha}}.
\end{equation}
The bound follows from $\mathrm{Lip}(f_\eta)\le\eta^{-\alpha}$, since the two arguments
differ by $\frac1N\eta^{-d}B_{ii}(0)$.

\medskip
\emph{The intermediate and the macroscopic particle systems.} Replacing the empirical interaction in \eqref{eq:particles} by the solution $u_\eta$ of
\eqref{eq:nonloc}, and keeping the initial data and the Brownian motions, gives the intermediate
system
\begin{equation}\label{eq:intermediate}
 \d\bar X_{k,i}^\eta=-\na U_i(\bar X_{k,i}^\eta)\dt
 +\big(2\sigma_i+2\textstyle\sum_j f_\eta(B_{ij}^\eta* u_{\eta,j}(t,\bar X_{k,i}^\eta))\big)^{1/2}\d W_i^k,
 \quad \bar X^\eta_{k,i}(0)=\xi^k_i.
\end{equation}
Its coefficients are deterministic, so its particles are independent, and its
diffusion coefficient is $\bar\Sigma_{k,i}$, with
\begin{equation}\label{eq:Sigbar}
 \bar\Sigma^2_{k,i}=2\sigma_i+2\sum_{j=1}^n f_\eta\big(H_{k,ij}\big),
 \qquad H_{k,ij}:=\big(B^\eta_{ij}*u_{\eta,j}\big)\big(\bar X^{\eta}_{k,i}\big),
\end{equation}
so that $\bar\Sigma^2_{k,i}=2\Theta^\eta_i(\bar X^\eta_{k,i})$. We denote by
$\bar\mu^N_j$ its empirical measures. By Lemma~\ref{lem:lawid}(ii), \eqref{eq:intermediate} has a unique strong solution, and the law of $\bar X^\eta_{k,i}(t)$ is
$u_{\eta,i}(t)\dx$ at every $t$. The same construction with the local coefficient
$f(a_{ij}u_j)$ in place of $f_\eta(B^\eta_{ij}*u_{\eta,j})$ gives the macroscopic system
\begin{equation}\label{eq:macroscopic}
 \d\widehat X_{k,i}=-\na U_i(\widehat X_{k,i})\dt
 +\big(2\sigma_i+2\textstyle\sum_j f(a_{ij}u_j(t,\widehat X_{k,i}))\big)^{1/2}\d W_i^k,
\quad \widehat X_{k,i}(0)=\xi^k_i.
\end{equation}
Its diffusion coefficient is $(2\Theta_i(\widehat X_{k,i}))^{1/2}$.

\medskip
\phantomsection\label{par:eta0}\emph{Convention at $\eta=0$.}Several statements cover the
nonlocal system and the SKT system at once, by allowing $\eta=0$. At $\eta=0$, the symbols
$u_\eta$, $\Theta^\eta_i$ and $\bar X^\eta_{k,i}$ stand for $u$, $\Theta_i$ and
$\widehat X_{k,i}$, and $B^0_{ij}*$ and $f_0$ stand for multiplication by $a_{ij}$ and for
$f$. With this convention, \eqref{eq:nonloc} and \eqref{eq:intermediate} at $\eta=0$ are
\eqref{eq:skt_det} and \eqref{eq:macroscopic}, and $\Theta^0_i=\Theta_i$. By
Lemma~\ref{lem:lawid}(ii) at $\eta=0$, \eqref{eq:macroscopic} has a unique strong solution, and the law of $\widehat X_{k,i}(t)$ is $u_i(t)\dx$.

The three systems are driven by the same data and are coupled particle by particle. The coupling
defect $Y_{k,i}:=X^{N,\eta}_{k,i}-\bar X^\eta_{k,i}$ between the first two solves, since
$-\na U_i(x)=x$,
\begin{equation}\label{eq:dY}
 \dd Y_{k,i}=Y_{k,i}\,\dd t+\Sigma^{\rm d}_{k,i}\,\dd W^k_i,
 \qquad \Sigma^{\rm d}_{k,i}:=\Sigma_{k,i}-\bar\Sigma_{k,i}.
\end{equation}
 Its initial value is $Y_{k,i}(0)=0$. Through this coupling the empirical measure
$\mu^N$ of \eqref{eq:particles} approaches the law $u_\eta$ of \eqref{eq:intermediate}, in the
moderate regime of \cite{chen_holzinger_2021_rigorous_SKT}, and $u_\eta$ approaches the law $u$ of
\eqref{eq:macroscopic} by \eqref{eq:CHbias}. The fluctuation field of this paper,
$\fluc^N=\sqrt N(\mu^N-u_\eta)$, measures the first of these two steps at the next order;
Section~\ref{sec:intro} explains why the centring is at $u_\eta$.

\medskip
\emph{The moderate scale.} The two parameters are coupled, $\eta=\eta(N)\to0$, at the logarithmic scale of \ref{H:scaling},
\begin{equation}\label{eq:H3}
 \Lambda:=\eta^{-2(d+1+\alpha)}\ \le\ \delta\log N,\qquad \eta\le\eta_0\le1,
 \qquad\alpha\in[0,1).
\end{equation}
where, for a fixed $\vartheta_0\in(0,\frac16]$, $\delta$ satisfies the smallness
condition
\begin{equation}\label{eq:smallness}
 \eps_\delta\ =\ C_\star\,C_2\,(T+1)\,\delta\ <\ 1,
 \qquad C_\star=(1+\vartheta_0)\Big(2+\frac4d\Big)\le7\ \ (\vartheta_0\le\tfrac16),
\end{equation}
with $C_2$ the constant of Proposition~\ref{prop:phi4}. Several estimates hold only for large $N$. We fix
\begin{equation}\label{eq:N0}
 N_0:=\min\big\{N_1\in\N:\ \eta(N)\le\eta_0,\ \ \eta(N)^{-2(d+1+\alpha)}\le\delta\log N\ \text{and}\ \eta(N)^{-d}\le N\ \text{for every }N\ge N_1\big\}.
\end{equation}
The minimum exists: $\eta(N)\to0$, the scale condition of \ref{H:scaling} holds for large $N$, and by it,
$\eta(N)^{-d}\le(\delta\log N)^{d/(2(d+1+\alpha))}\le N$ for large $N$.

\medskip
\noindent\emph{Constants and indices.} We denote
\[
 A_0:=\max_{i,j}\|B_{ij}\|_{L^1(\R^d)},\qquad A:=\max(1,A_0),\qquad
 \sigma_{\min}:=\min_i\sigma_i,\qquad \gold:=\tfrac{1+\sqrt5}2 .
\]
The rescaling $B^\eta_{ij}=\eta^{-d}B_{ij}(\cdot/\eta)$ preserves the $L^1$ norm, so
$\|B^\eta_{ij}\|_{L^1}\le A_0$ for every $\eta$. For the index $s$ of \ref{H:CH} below and the norm
$\|v\|_{H^s}=\|(\mathrm{Id}-\Delta)^{s/2}v\|_{L^2}$, we denote by $C_s$ the norm of the embedding
$H^s(\R^d)\hookrightarrow L^\infty(\R^d)$. The bound on the arguments of $f$ and the interval on
which $f$ is controlled are
\begin{equation}\label{eq:b0}
 \mathsf b_0:=A_0\,C_s\|u_0\|_{H^s},
 \qquad I:=[-2\mathsf b_0,2\mathsf b_0]\ \text{ as in \cite[eq.~(11)]{chen_holzinger_2021_rigorous_SKT}}.
\end{equation}
The indices of the weighted scale \eqref{eq:spaces} below are fixed once and for all as
\begin{equation}\label{eq:choice}
 s_1:=\big\lfloor\tfrac d2\big\rfloor+2,\quad a_1:=0,\qquad
 s_m:=s_1+3,\quad a_m:=1,\qquad s_2:=s_1+4,\quad a_2:=2.
\end{equation}

\begin{assumption}[Standing assumptions]\label{ass:main}\leavevmode
\begin{enumerate}[label=\normalfont(H\arabic*)]
\item \label{H:CH} \emph{The data.} Assumptions (A1)--(A4) of
\cite{chen_holzinger_2021_rigorous_SKT} hold, that is:
\begin{enumerate}[label=\normalfont(A\arabic*),leftmargin=*]
\item \label{H:A1} $\sigma_i>0$, and each $u_{0,i}$ is a probability density with finite
second moment;
\item \label{H:A2} $U_i(x)=-\frac12|x|^2$ for every $i$;
\item \label{H:A3} $B_{ij}\in C_0^\infty(\R^d)$ with $\operatorname{supp}B_{ij}\subset B_1(0)$,
for every $i,j$;
\item \label{H:A4} $s>\frac d2+1$, and $f\in W^{s+1,\infty}_{\rm loc}(\R;[0,\infty))$; for every
$\eta$, $f_\eta\in W^{s+1,\infty}(\R;[0,\infty))$ with $f_\eta=f$ on $[-a_\eta,a_\eta]$ and
$\mathrm{Lip}(f_\eta)\le\eta^{-\alpha}$, where $\alpha\in[0,1)$ and $a_\eta\to\infty$ as $\eta\to0$.
\end{enumerate}
In addition, each $B_{ij}$ is radial, $s\in\N$, and $u_0=(u_{0,1},\dots,u_{0,n})\in H^s$.
The random data are defined on one probability space $(\Omega,\mathcal F,\Prob)$ with a
filtration $(\mathcal F_t)$ satisfying the usual conditions, which we call the particle
filtration. The $\xi^k_i$ are $\mathcal F_0$-measurable, $(W^k_i)_{k\le N,\,i\le n}$ is an
$(\mathcal F_t)$-Brownian motion in $\R^{nNd}$, and
\begin{align}
 &\text{for each }i:\ \{(\xi^k_i,W^k_i)\}_{k\le N}\ \text{are i.i.d., with }
 \xi^k_i\sim u_{0,i}\ \text{independent of }W^k_i;
 \tag{W}\label{eq:dataw}\\
 &\text{the }n\ \text{arrays}\ \{(\xi^k_i,W^k_i)\}_{k\le N},\ i\le n,
 \ \text{are independent of one another.}
 \tag{X}\label{eq:datax}
\end{align}

\item \label{H:cutoff} \emph{The threshold $\eta_0$.} There is $\eta_0\in(0,1]$ with
\[
 \eta_0e^T\le1,\qquad \tfrac{3n}2\,\eta_0^{\,d+1+\alpha}\le\vartheta_0,\qquad
 a_\eta\ge2A\,C_s\|u_0\|_{H^s}\ \text{ for every }\eta\le\eta_0,
\]
where $\vartheta_0\in(0,\frac16]$ is the constant fixed in \eqref{eq:smallness}.

\item \label{H:small} \emph{Smallness of $f$.} Let $\eps=\eps(u_0)>0$ be the constant of
\cite[Thms.~2--3]{chen_holzinger_2021_rigorous_SKT}. Then
\[
 \|f\|_{C^{s+1}(I)}\le\eps
 \qquad\text{and}\qquad
 n\,A\,\gold^{\,a_1+1}\,C_s\|u_0\|_{H^s}\,\eps\ \le\ \tfrac14\,\sigma_{\min},
\]
where $\|\cdot\|_{C^{s+1}(I)}$ is the full norm, so that $\|f^{(m)}\|_{L^\infty(I)}\le\eps$ for
every $m\le s+1$.

\item \label{H:scaling} \emph{Moderate scaling.} $\eta=\eta(N)\to0$ as $N\to\infty$, with
$\eta^{-2(d+1+\alpha)}\le\delta\log N$, where $\delta>0$ satisfies \eqref{eq:smallness}.

\item \label{H:budget} \emph{Regularity.} The index $s$ of \ref{H:CH} satisfies $s\ge s_1+4$.
\end{enumerate}
\end{assumption}

\noindent\emph{Convention.} Assumption~\ref{ass:main} is in force throughout, and we do not
repeat it in the statements. All statements are for $\eta\le\eta_0$, and statements involving $N$ are for
$N\ge N_0$, with $N_0$ as in \eqref{eq:N0}. In particular $\sup_N$ denotes the supremum over $N\ge N_0$.
Finitely many $N$ do not affect tightness or limits in law. Appendix~\ref{sec:appevol}
works under Condition~\ref{ass:budget}, a condition on a level $(\sigma,a)$ and a coefficient
pair, which Assumption~\ref{ass:main} implies at the levels $(s_1,a_1)$ and $(s_1+1,a_1)$, the
latter with upper level, which is all Appendix~\ref{sec:appevol} needs
(Proposition~\ref{prop:levels}). Where fewer hypotheses suffice and this matters, the statement
says so.

\begin{remark}[On Assumption~\ref{ass:main}]\label{rem:hyp}
\emph{(i) Origin of the hypotheses.} Items \ref{H:CH}--\ref{H:scaling} are the assumptions of
\cite{chen_holzinger_2021_rigorous_SKT}, together with the regularity $u_0\in H^s$ of its
Theorems~2--3, and with five changes: the radial symmetry of $B_{ij}$ and the integrality of $s$ in \ref{H:CH}, the independence structure
\eqref{eq:dataw}--\eqref{eq:datax} of the driving data, the second condition of \ref{H:small},
the bounds $\eta_0e^T\le1$ and $\frac{3n}2\eta_0^{\,d+1+\alpha}\le\vartheta_0$ in
\ref{H:cutoff}, and the smaller $\delta$ of \eqref{eq:smallness} in \ref{H:scaling}.
Item \ref{H:budget} is the regularity this paper needs and  Proposition~\ref{prop:levels} lists where it is used. The index $s$ is an integer because
the Moser estimates \cite[Lemmas~12--14]{chen_holzinger_2021_rigorous_SKT}, used in the proofs
of \cite[Thms.~2--3]{chen_holzinger_2021_rigorous_SKT}, are stated for integer $s$. By
\eqref{eq:choice}, $s_1>\frac d2+1$ is an integer, and \ref{H:budget} is the condition
$s\ge\lfloor\frac d2\rfloor+6$.

\emph{(ii) The driving data.} Together, \eqref{eq:dataw} and \eqref{eq:datax} state that all the
$\xi^k_i$ and $W^k_i$, $k\le N$, $i\le n$, are mutually independent (for correlated species see
Remark~\ref{rem:X}). Consequently the augmented filtration generated by the $\xi^k_i$ and the
$W^k_i$ is a particle filtration.

\emph{(iii) The deterministic systems.} By \ref{H:CH}, \ref{H:small} and
\cite[Thms.~2--3]{chen_holzinger_2021_rigorous_SKT}, \eqref{eq:skt_det} and \eqref{eq:nonloc}
have unique global solutions, with
$\|u\|_{L^\infty(0,\infty;H^s)}^2+\sigma_*\|\na u\|_{L^2(0,\infty;H^s)}^2\le\|u_0\|^2_{H^s}$
for a constant $\sigma_*\in(0,\sigma_{\min})$ of \cite{chen_holzinger_2021_rigorous_SKT}, the same bound for $u_\eta$,
uniformly in $\eta$, and the bias estimate \eqref{eq:CHbias} \cite[Thm.~3]{chen_holzinger_2021_rigorous_SKT}.

\emph{(iv) The bound on $u_\eta$ and the interval $I$.} By (iii), $\|u_\eta(t)\|_{H^s}\le\|u_0\|_{H^s}$
for almost every $t$ and every $\eta\le\eta_0$. For the representative of $u_\eta$ that is weakly continuous into $H^s$
(Lemma~\ref{lem:lawid}(i)), which is the one used throughout, the bound holds for every $t$, because the norm is weakly lower
semicontinuous. So the embedding $H^s\hookrightarrow L^\infty$ gives
\begin{equation}\label{eq:Cu}
 C_u:=\sup_{\eta\le\eta_0}\ \sup_{t\le T}\ \max_j\|u_{\eta,j}(t)\|_{L^\infty}\ \le\ C_s\|u_0\|_{H^s}<\infty .
\end{equation}
By Young's inequality, $|B^\eta_{ij}*u_{\eta,j}|\le\|B^\eta_{ij}\|_{L^1}C_u\le\mathsf b_0$, so the
arguments of $f_\eta$ in $\Theta^\eta$ lie in $[-\mathsf b_0,\mathsf b_0]\subset I$. By \ref{H:cutoff},
$a_\eta\ge2AC_s\|u_0\|_{H^s}\ge2\mathsf b_0$, so $I\subseteq[-a_\eta,a_\eta]$, and $f_\eta=f$
on $I$ by \ref{H:A4}. Hence \ref{H:small} controls $f_\eta$ and its derivatives up to order
$s+1$ along $u_\eta$, uniformly in $\eta$.

\emph{(v) The second smallness condition.} The second condition of \ref{H:small} is
\eqref{eq:epscondup} at $a=a_1$, with $\sigma_{\min}=\min_k\sigma_k$ and
$\gamma=C_s\|u_0\|_{H^s}\eps$, the values of Lemma~\ref{lem:coeff}.
\end{remark}

\medskip
\emph{The weighted scale.}
We denote by $\langle x\rangle:=(1+|x|^2)^{1/2}$ the weight on $\R^d$ (also used with the
variables $y$ and $\xi$), and by $\langle\cdot\rangle^{a}$ the function $x\mapsto\langle x\rangle^{a}$.
It is smooth and positive, and for every $a\in\R$ all its derivatives are bounded by multiples of
$\langle x\rangle^{a}$, so multiplication by $\langle\cdot\rangle^{a}$ maps $\mathcal S'(\R^d)$ to
itself and $\mathcal S(\R^d)$ onto itself.
Throughout, for $s,a\in\R$ the weighted Bessel-potential space $H^s_a(\R^d)$ is given by
\begin{equation}\label{eq:Hsa}
 H^s_a(\R^d):=\big\{\varphi\in\mathcal S'(\R^d):
 \ \|\varphi\|_{H^s_a}:=\big\|\langle\cdot\rangle^{a}\varphi\big\|_{H^s(\R^d)}<\infty\big\}.
\end{equation}
We abbreviate the spaces we will work with by
\begin{equation}\label{eq:spaces}
  \Ufine=H^{-s_1}_{-a_1},\quad \Usharp=H^{-s_1-1}_{-a_1},\quad
  \Fcoarse:=H^{-s_1-3}_{-a_1-1},\quad
  \Hmid=H^{-s_m}_{-a_m},\quad \Vcoarse=H^{-s_2}_{-a_2},
  \qquad \Test:=H^{s_1+3}_{a_1+1},
\end{equation}
so that $\Fcoarse=\Test^*$. By Lemma~\ref{lem:spaces}, the spaces of \eqref{eq:spaces} embed as follows:
\begin{equation}\label{eq:chainA}
 \Ufine\hookrightarrow\Usharp\hookrightarrow\hookrightarrow\Hmid
 \hookrightarrow\hookrightarrow\Vcoarse,
 \qquad \Ufine\hookrightarrow\Fcoarse\hookrightarrow\Vcoarse ,
\end{equation}
where the double arrows are compact embeddings.

The indices \eqref{eq:choice} satisfy
\begin{equation}\label{eq:tightness_scale}
 s_1\in\N,\quad \tfrac d2+1<s_1,\qquad s_1+1<s_m<s_2,\qquad 0\le a_1<a_m<a_2,
 \qquad s_2\ge s_1+3,\quad a_2\ge a_1+1,
\end{equation}
and all arrows of \eqref{eq:chainA} follow from these relations and Lemma~\ref{lem:spaces}(ii). The
strict inequalities give the two compact embeddings. With $\Hmid=\Fcoarse$, we call $(\Ufine,\Hmid,\Vcoarse)$ the tightness triple:
$\Ufine$ is the fine space in which the uniform bounds are proved, $\Hmid$ the middle space in which
convergence in law is stated, and $\Vcoarse$ the coarse space in which the time-increment (Aldous)
estimates are proved. We denote by $c_{\rm emb}$ the norm of $\Ufine\hookrightarrow\Hmid$ and by
$c_\sharp$ the norm of $\Ufine\hookrightarrow\Usharp$; $c_\sharp\le1$, since the two spaces carry the
same weight.

With $\langle\nabla\rangle^\sigma:=(\mathrm{Id}-\Delta)^{\sigma/2}$,
\begin{equation}\label{eq:sob}
 \|v\|^2_{H^{\sigma+1}}=\|\langle\nabla\rangle^\sigma v\|^2_{L^2}+\|\na\langle\nabla\rangle^\sigma v\|^2_{L^2}.
\end{equation}
Ehrling's inequality \cite[Ch.~1, Lemme~5.1]{Lions69} states: if $X\hookrightarrow\hookrightarrow Y\hookrightarrow Z$ are Banach spaces, the second embedding injective, then for every $\eps>0$ there is $C_\eps$ with
\begin{equation}\label{eq:ehrling}
 \|u\|_Y\ \le\ \eps\,\|u\|_X+C_\eps\|u\|_Z\qquad\text{for all }u\in X.
\end{equation}

The section closes with the elementary facts used throughout: the duality,
embedding and compactness properties of the scale (Lemma~\ref{lem:spaces}), the reproducing kernel of the gradient
(Lemma~\ref{lem:kernel}), the composition estimate at integer order (Lemma~\ref{lem:moser}), the second moments of the smoothed fluctuation
(Lemma~\ref{lem:smoothed}), and the principle that under \ref{H:scaling} every
negative power of $\eta$ is polylogarithmic in $N$ (Lemma~\ref{lem:polylog}).
Lemma~\ref{lem:secondmoment} bounds the second moment of $u_\eta$ on $[0,T]$.
Lemma~\ref{lem:CH2} bounds the mobility defect in $L^2$. It is used in
Lemma~\ref{lem:L2field} and Corollaries~\ref{cor:trace} and \ref{cor:Maldous}, and
Corollary~\ref{cor:CH2} in Appendix~\ref{sec:appk4} proves it.
Lemma~\ref{lem:lawid} shows that $u_\eta(t)$ is the law of an intermediate particle at
every $t$.

\begin{lemma}[Duality, embedding, compactness]\label{lem:spaces}
Let $s,s',a,a'\in\R$ and recall $\|v\|_{H^s_a}=\|\br\cdot^av\|_{H^s}$, so that
$J_a:v\mapsto\br\cdot^av$ is an isometric isomorphism $H^s_a\to H^s$ for every
$s$. Throughout, $\br{\cdot,\cdot}$ denotes the extension by continuity of the
$L^2$ pairing
\begin{equation}\label{eq:wpair}
 \br{u,v}:=\int_{\R^d}u\,v\dx,\qquad u,v\in\mathcal S(\R^d),
\end{equation}
to $H^{-s}_{-a}\times H^{s}_{a}$.
\begin{enumerate}[label=\normalfont(\roman*)]
\item Under the $L^2$-pairing, $(H^s_a)^*=H^{-s}_{-a}$, isometrically.
\item If $s'\ge s$ and $a'\ge a$ then $H^{s'}_{a'}\hookrightarrow H^s_a$
continuously. If in addition $s'>s$ and $a'>a$, the embedding is compact.
\item If $s>\frac d2+m$, $m\in\N_0$, then for all multi-indices $|\nu|\le m$
\begin{equation}\label{eq:ptwisew}
 |D^\nu\varphi(x)|\ \le\ C(d,s,m,a)\,\br x^{-a}\,\|\varphi\|_{H^s_a}
 \qquad\text{for every }x,
\end{equation}
and consequently $H^s_a\hookrightarrow C^m_b$ when $a\ge0$.
\item For every $\phi\in C^\infty_b(\R^d)$, multiplication by $\phi$ is a bounded operator on
$H^s_a$.
\end{enumerate}
\end{lemma}

\begin{proof}
All four statements reduce to the unweighted case through the isometry $J_a$, and
$\mathcal S$ is dense in every $H^\sigma_b$ because $J_b$ maps $\mathcal S$ onto $\mathcal S$.
\emph{(i)} and the bound $|\br{u,v}|\le\|u\|_{H^{-s}_{-a}}\|v\|_{H^s_a}$ follow from the
$H^{-s}$--$H^s$ duality applied to $\br\cdot^{-a}u$ and $\br\cdot^av$.
\emph{(ii)} Put $\delta:=a'-a\ge0$ and $\phi:=\br\cdot^{-\delta}\in C^\infty_b$. Multiplication by a
$C^\infty_b$ function is bounded on $H^\sigma$ for every real $\sigma$ (Leibniz at integer
$\sigma\ge0$, interpolation, and duality). Since $\br\cdot^av=\phi\,\br\cdot^{a'}v$, the embedding is
$J_a^{-1}M_\phi J_{a'}$ composed with $H^{s'}\hookrightarrow H^s$, which gives continuity. If
$s'>s$ and $\delta>0$, split $\phi=\phi\chi_R+\phi(1-\chi_R)$ with a cut-off $\chi_R$ equal to
$1$ on $B_R$. The second multiplier has norm $O(R^{-\delta})$ on $H^s$, and
$M_{\phi\chi_R}:H^{s'}\to H^s$ is compact by Rellich's theorem for distributions supported
in the fixed ball $B_{2R}$. So $M_\phi:H^{s'}\to H^s$ is an operator-norm limit of compact
operators, hence compact.
\emph{(iii)} With $u:=\br\cdot^a\varphi\in H^s$, Leibniz gives
$D^\nu\varphi=\sum_{\nu'\le\nu}c_{\nu'}(D^{\nu-\nu'}\br\cdot^{-a})D^{\nu'} u$, where
$|D^{\nu-\nu'}\br x^{-a}|\le C\br x^{-a}$, and the Sobolev embedding $H^s\hookrightarrow C^m_b$
bounds $D^{\nu'} u$.
\emph{(iv)} Multiplication by $\phi$ commutes with $J_a$. So (iv) follows from the
boundedness of multiplication by $\phi$ on $H^s$, shown in the proof of (ii) above.

\end{proof}

\begin{lemma}[Polylogarithmic principle]\label{lem:polylog}
Assume \ref{H:scaling}. Then for every fixed $\kappa\ge0$,
\[
 N^{-1/2}\,\eta^{-\kappa}
 \ \le\ N^{-1/2}\,(\delta\log N)^{\frac{\kappa}{2(d+1+\alpha)}}
 \ \longrightarrow\ 0\qquad(N\to\infty).
\]
Consequently \emph{any} quantity of the form $N^{-1/2}\eta^{-\kappa}$, with
$\kappa$ independent of $N$, is negligible in the moderate regime of
\cite{chen_holzinger_2021_rigorous_SKT}.
\end{lemma}

\begin{proof}
\ref{H:scaling} gives $\eta^{-(d+1+\alpha)}\le(\delta\log N)^{1/2}$, hence
$\eta^{-1}\le(\delta\log N)^{1/(2(d+1+\alpha))}$, and $(\log N)^p=o(N^{1/2})$ for every $p>0$.
\end{proof}

\begin{lemma}[The law of an intermediate particle]\label{lem:lawid}
Let $\eta\in[0,\eta_0]$, with the \hyperref[par:eta0]{convention at $\eta=0$} of
Section~\ref{sec:setting}. Part~(i) uses only \ref{H:CH}--\ref{H:small}. Denote $g^\eta_{ij}:=u_{\eta,i}f'_\eta(B^\eta_{ij}*u_{\eta,j})$ and
$L_r:=x\cdot\na+\Theta^\eta_i(r)\Delta$.
\begin{enumerate}[label=\normalfont(\roman*)]
\item $u_\eta$ has a representative that is weakly continuous from $[0,T]$ into $H^s$,
with $u_\eta(0)=u_0$ and $\|u_\eta(t)\|_{H^s}\le\|u_0\|_{H^s}$ for every $t$. For it,
$u_{\eta,i}$, $B^\eta_{ij}*u_{\eta,j}$, $\Theta^\eta_i$ and $g^\eta_{ij}$ are bounded and
continuous on $[0,T]\times\R^d$; $|B^\eta_{ij}*u_{\eta,j}|\le\mathsf b_0$ everywhere;
$t\mapsto\Theta^\eta_i(t)v$ and $t\mapsto g^\eta_{ij}(t)v$ are continuous into $L^2$ for
every $v\in L^2$; and for every $g\in C^\infty_0(\R^d)$ and $0\le p\le q\le T$,
\begin{equation}\label{eq:lawweak}
 \br{u_{\eta,i}(q)-u_{\eta,i}(p),g}=\int_p^q\br{u_{\eta,i}(r),L_rg}\dd r .
\end{equation}
\item For every $k$ and $i$, \eqref{eq:intermediate} (at $\eta=0$,
\eqref{eq:macroscopic}) has a unique strong solution with continuous paths, adapted to
the particle filtration. For every $t\in[0,T]$,
every $k$ and every $i$, the law of $\bar X^\eta_{k,i}(t)$ is $u_{\eta,i}(t,x)\dx$. In
particular $u_{\eta,i}(t)\ge0$ and $\int u_{\eta,i}(t)=1$.
\end{enumerate}
\end{lemma}

\begin{proof}
Part~(i) builds the weakly continuous representative from the weak form of \eqref{eq:nonloc}
and a countable dense set of test functions. In part~(ii), the Lipschitz theory of SDEs gives the
strong solution. The law is identified by duality: both $\E\phi(\bar X^\eta_{k,i}(t))$ and
$\br{u_{\eta,i}(t),\phi}$ equal $\br{u_{0,i},\Phi(0)}$, where $\Phi$ solves the backward Kolmogorov
equation \eqref{eq:kolm}.

\emph{(i).} By \ref{H:small}, $u_\eta\in L^\infty(0,T;H^s)$ with norm at most
$\|u_0\|_{H^s}$, and it solves \eqref{eq:nonloc} (at $\eta=0$, \eqref{eq:skt_det}) in the
sense of distributions, with datum $u_0$. So, for $\phi\in C^\infty_0(\R^d)$, the map
$t\mapsto\br{u_{\eta,i}(t),\phi}$ agrees almost everywhere with a Lipschitz function
$F_\phi$ with $F_\phi(0)=\br{u_{0,i},\phi}$, whose derivative is
$\br{u_{\eta,i},x\cdot\na\phi+\Theta^\eta_i\Delta\phi}$, as in \eqref{eq:weak_nonloc}.
Fix a countable set $\{\phi_j\}\subset C^\infty_0$ dense in $L^2$. For $t$ in a set of full
measure, $\|u_\eta(t)\|_{H^s}\le\|u_0\|_{H^s}$ and $\br{u_\eta(t),\phi_j}=F_{\phi_j}(t)$ for
all $j$. For the remaining $t$, let $u_\eta(t)$ be the weak limit in $H^s$ of
$u_\eta(t_m)$ along $t_m\to t$ in that set. It exists and does not depend on the sequence,
because the sequence is bounded and its pairings with every $\phi_j$ converge to
$F_{\phi_j}(t)$. The same argument shows weak continuity at every $t$, and the bound
holds at every $t$ by weak lower semicontinuity of the norm. For $g\in C^\infty_0$, the function $F_g$ and
$t\mapsto\br{u_{\eta,i}(t),g}$ are both continuous and agree almost everywhere, so they
agree everywhere, and \eqref{eq:lawweak} is $F_g(q)-F_g(p)=\int_p^qF_g'$.

For each $x$, $u_{\eta,j}(t,x)$ and
$(B^\eta_{ij}*u_{\eta,j}(t))(x)=\br{u_{\eta,j}(t),B^\eta_{ij}(x-\cdot)}$ are continuous in
$t$, since point evaluation is continuous on $H^s$ ($s>\frac d2$) and $u_\eta$ is weakly
continuous in $H^s$. They are Lipschitz in $x$, uniformly in $t$, since
$H^s\hookrightarrow W^{1,\infty}$ ($s>\frac d2+1$) and $\|B^\eta_{ij}\|_{L^1}\le A_0$.
Hence they are continuous on $[0,T]\times\R^d$. By the bound at every $t$,
$|B^\eta_{ij}*u_{\eta,j}(t,x)|\le A_0C_s\|u_0\|_{H^s}=\mathsf b_0$ for every $(t,x)$, and
likewise $|a_{ij}u_j(t,x)|\le\mathsf b_0$. So the arguments of $f_\eta$ and $f'_\eta$ lie
in $I$, where $f_\eta=f$ by \ref{H:cutoff}. The functions $f$ and $f'$ are continuous on $I$, so
$\Theta^\eta_i$ and $g^\eta_{ij}$ are continuous and bounded, by $\sigma_i+n\eps$ and by
$C_s\|u_0\|_{H^s}\eps$. Continuity of $t\mapsto\Theta^\eta_i(t)v$ in $L^2$ follows by
dominated convergence, and likewise for $g^\eta_{ij}$.

\emph{(ii).}\emph{Existence.} The drift $x\mapsto x$ of \eqref{eq:intermediate} is linear. By (i),
$\Theta^\eta_i\ge\sigma_{\min}$, since $f_\eta\ge0$, and $\Theta^\eta_i$ is continuous and bounded. The function $\Theta^\eta_i$ is
Lipschitz in $x$, uniformly in $t$: the arguments $B^\eta_{ij}*u_{\eta,j}(t)$ are Lipschitz in $x$
with constant $A_0\sup_t\|\na u_{\eta,j}(t)\|_{L^\infty}<\infty$, since $H^s\hookrightarrow W^{1,\infty}$,
and they lie in $I$, where $|f'_\eta|=|f'|\le\eps$. As
$|\sqrt a-\sqrt b|\le|a-b|/(2\sqrt{\sigma_{\min}})$ for $a,b\ge\sigma_{\min}$, the diffusion
coefficient $(2\Theta^\eta_i)^{1/2}$ is bounded and Lipschitz in $x$, uniformly in $t$, and
measurable in $(t,x)$. The datum $\xi^k_i$ is $\mathcal F_0$-measurable and square integrable, and
$W^k_i$ is an $(\mathcal F_t)$-Brownian motion. Hence \cite[Thms.~5.2.5 and 5.2.9]{KaSh91} give a
strong solution with continuous paths, and pathwise uniqueness, so the solution is unique among
continuous processes adapted to any filtration for which $W^k_i$ is a Brownian motion. The same argument applies at $\eta=0$, with
$a_{ij}u_j$ in place of $B^\eta_{ij}*u_{\eta,j}$.

\emph{The law.} Fix $t\in(0,T]$, $k$, $i$ and $\phi\in C^\infty_0(\R^d)$. The case $t=0$ is
\eqref{eq:dataw}. Let $\Phi$ solve the scalar backward problem
\begin{equation}\label{eq:kolm}
 \pa_r\Phi+L_r\Phi=0\quad\text{on }[0,t],\qquad\Phi(t)=\phi.
\end{equation}
Problem \eqref{eq:kolm} is \eqref{eq:backward-prop66} with $n=1$ and $g=0$. By Lemma~\ref{lem:coeff}, the scalar
pair $(\Theta^\eta_i,0)$ satisfies Condition~\ref{ass:budget} with upper level at
$(s_1+1,a_1)$, with $\gamma=0$. Lemma~\ref{lem:coeff} depends on the present lemma only
through part~(i), proved above, so the argument is not circular. So
Propositions~\ref{prop:evolution-prop66} and \ref{prop:persist} apply. At the level $(s_1+1,a_1)$, which
Assumption~\ref{ass:main} covers by Proposition~\ref{prop:levels}, this gives
$\Phi\in C([0,t];\Test)$, and \eqref{eq:kolm} holds in $L^2(0,t;H^{s_1+2}_{a_1})$. So
\begin{equation}\label{eq:kolmint}
 \Phi(p)=\Phi(q)+\int_p^qh_r\dd r\quad(0\le p\le q\le t),\qquad h_r:=L_r\Phi(r),
\end{equation}
a Bochner integral in $H^{s_1+2}_{a_1}$. We use four facts. First,
$\Test\hookrightarrow C^2_b$ by Lemma~\ref{lem:spaces}(iii), since $s_1+3>\frac d2+2$.
Second, $L_r$ is bounded from $\Test$ into $H^{s_1+1}_{a_1}$, uniformly in $r$: $x\cdot\na$
maps $\Test$ into $H^{s_1+2}_{a_1}$ (Proposition~\ref{prop:persist}), and
$\Theta^\eta_i(r)-\kappa_i$ is bounded in $H^{s_1+2}$ by Lemma~\ref{lem:coeff}. Third,
$H^{s_1+1}_{a_1}\hookrightarrow C^1_b\cap L^2$, since $s_1+1>\frac d2+1$ and $a_1\ge0$.
Fourth, $r\mapsto h_r=x\cdot\na\Phi(r)+\Theta^\eta_i(r)\Delta\Phi(r)$ is continuous into
$L^2$, by the continuity of $\Phi$ in $\Test$ and by (i). Let $\pi:0=t_0<\dots<t_L=t$ be a
partition, with mesh $|\pi|:=\max_l(t_{l+1}-t_l)$.

\emph{The particle.} Put $X:=\bar X^\eta_{k,i}$, so $\dd X=X\dd r+(2\Theta^\eta_i(r,X))^{1/2}\dd W^k_i$.
For $g\in\Test$, It\^o's formula gives
$\E g(X(q))-\E g(X(p))=\E\int_p^q(L_rg)(X(r))\dd r$, since the stochastic integral has a
bounded integrand. With $g=\Phi(t_{l+1})$ and \eqref{eq:kolmint}, evaluated at the point
$X(t_l)$,
\[
 \E\phi(X(t))-\E\Phi(0,X(0))
 =\sum_l\E\int_{t_l}^{t_{l+1}}\Big[\big(L_r\Phi(t_{l+1})\big)(X(r))-h_r(X(t_l))\Big]\dd r.
\]
The integrand is at most $C\|\Phi(t_{l+1})-\Phi(r)\|_{\Test}+C\min\{1,|X(r)-X(t_l)|\}$,
by the second and third facts. It tends to $0$ as $|\pi|\to0$, by uniform continuity of
$\Phi$ and of the paths of $X$, and by dominated convergence. So
$\E\phi(X(t))=\E\Phi(0,\xi^k_i)=\br{u_{0,i},\Phi(0)}$.

\emph{The density.} For $g\in C^\infty_0$ and $0\le p\le q\le T$,
\eqref{eq:lawweak} holds by (i). Both sides
are continuous in $g\in\Test$, since $u_\eta$ is bounded in $L^2$ and $L_r$ is bounded from
$\Test$ into $L^2$; so the identity holds for $g\in\Test$. Hence
\[
 \br{u_{\eta,i}(t),\phi}-\br{u_{0,i},\Phi(0)}
 =\sum_l\int_{t_l}^{t_{l+1}}\Big[\br{u_{\eta,i}(r),L_r\Phi(t_{l+1})}-\br{u_{\eta,i}(t_l),h_r}\Big]\dd r.
\]
The integrand is at most $C\|\Phi(t_{l+1})-\Phi(r)\|_{\Test}+|\br{u_{\eta,i}(r)-u_{\eta,i}(t_l),h_r}|$.
The set $\{h_r:r\le t\}$ is compact in $L^2$ by the fourth fact, and $u_{\eta,i}$ is weakly
continuous and bounded in $L^2$. So the second term tends to $0$ uniformly in $r$ as
$|\pi|\to0$, and $\br{u_{\eta,i}(t),\phi}=\br{u_{0,i},\Phi(0)}$.

Both sides agree, so $\E\phi(X(t))=\br{u_{\eta,i}(t),\phi}$ for every
$\phi\in C^\infty_0(\R^d)$. The law of $X(t)$ is a probability measure, and it coincides
with $u_{\eta,i}(t)\dx$ as a distribution. So $u_{\eta,i}(t)\ge0$, the two measures are
equal, and $\int u_{\eta,i}(t)=1$.
\end{proof}

\begin{lemma}[Second moment on a finite horizon]\label{lem:secondmoment}
For every $\eta\in[0,\eta_0]$ and every $i$,
\begin{equation}\label{eq:secondmoment}
 \sup_{t\le T}\int_{\R^d}|x|^2u_{\eta,i}(t,x)\dd x
 \ \le\ e^{2T}\Big(\max_i\E\big|\xi^1_i\big|^2+d\max_i\big(\sigma_i+n\eps\big)\Big)\ =:\ \mathsf m_T,
\end{equation}
and consequently $\sup_{t\le T}\|\br\cdot u_{\eta,i}(t)\|_{L^2}\le C(T)$, both
uniformly in $\eta$ and in $i$.

\end{lemma}

\begin{proof}
We compute $\int|x|^2u_{\eta,i}(t)$ as the second moment of an intermediate particle, by
It\^o's formula localised with stopping times and by Gr\"onwall's inequality (Step~1). The
weighted bound follows from the second moment and the $L^\infty$ bound on $u_\eta$ (Step~2). The
bound $\Theta^\eta_i\le\sigma_i+n\eps$ used in Step~1 is proved at the end.

By Lemma~\ref{lem:lawid}, $u_{\eta,i}(t)$ is the density of $\bar X^\eta_{k,i}(t)$ for
every $t\le T$. So $\int|x|^2u_{\eta,i}(t)=\E|\bar X^\eta_{k,i}(t)|^2$, and we may argue on
the SDE. The function $u_\eta$ is given (Remark~\ref{rem:hyp}(iii)), so \eqref{eq:intermediate} is
not a McKean--Vlasov equation but an SDE with the fixed time-dependent coefficient
$\Theta^\eta_i(t,x)$. By Lemma~\ref{lem:lawid}(ii) this SDE has a
strong solution with continuous paths. So, with
$\tau_R:=\inf\{t:|\bar X^\eta_{k,i}(t)|\ge R\}$, a.s.\
$\sup_{[0,T]}|\bar X^\eta_{k,i}|<\infty$ and hence $\tau_R\uparrow\infty$.

\emph{Step 1: localised It\^o.} On $[\![0,\tau_R]\!]$ the integrand of the stochastic
term is bounded, so it is a true martingale, and
\[
 \E\big|\bar X^\eta_{k,i}(t\wedge\tau_R)\big|^2
 =\E|\xi^k_i|^2+\E\!\int_0^{t\wedge\tau_R}\!\!\Big(2|\bar X^\eta_{k,i}|^2
 +2d\,\Theta^\eta_i(s,\bar X^\eta_{k,i})\Big)\dd s,
\]
the $+2|\bar X|^2$ from the drift $+\bar X$ and the $2d\Theta^\eta_i$ from
$\tr\big(2\Theta^\eta_i\mathrm{Id}_d\big)$. Let $y_R(t):=\E|\bar X^\eta_{k,i}(t\wedge\tau_R)|^2$. Since
$|\bar X(t\wedge\tau_R)|\le\max(R,|\xi^k_i|)$, $y_R(t)\le R^2+\E|\xi^k_i|^2<\infty$ by \ref{H:CH}. On $\{s<\tau_R\}$ one has
$s\wedge\tau_R=s$, so $|\bar X(s)|\mathbb 1_{s<\tau_R}\le|\bar
X(s\wedge\tau_R)|$. With $\Theta^\eta_i\le\sigma_i+n\eps$ and Fubini's theorem,
\[
 y_R(t)\ \le\ \E|\xi^k_i|^2+2d(\sigma_i+n\eps)\,t+2\!\int_0^ty_R(r)\dd r.
\]
Set $z(t):=\E|\xi^k_i|^2+2d(\sigma_i+n\eps)t+2\int_0^ty_R$. Then
$z'=2d(\sigma_i+n\eps)+2y_R\le2d(\sigma_i+n\eps)+2z$, so
$y_R(t)\le z(t)\le e^{2t}\E|\xi^k_i|^2+d(\sigma_i+n\eps)(e^{2t}-1)$,
uniformly in $R$. Let $R\to\infty$: $\tau_R\uparrow\infty$ a.s., so Fatou gives
$\E|\bar X^\eta_{k,i}(t)|^2\le e^{2t}\E|\xi^k_i|^2+d(\sigma_i+n\eps)(e^{2t}-1)$
for \emph{every} $t\le T$ ($\E|\xi^k_i|^2<\infty$ because \ref{H:CH} makes
$\xi^k_i$ square-integrable), which is \eqref{eq:secondmoment} after the crude $(e^{2t}-1)\le e^{2T}$.

\emph{Step 2: the weighted bound.} By Lemma~\ref{lem:lawid}, $u_{\eta,i}(t)\ge0$ and
$\int u_{\eta,i}(t)=1$ for every $t$, and $\|u_{\eta,i}(t)\|_{L^\infty}\le C_s\|u_0\|_{H^s}$.
Hence
\[
 \|\br\cdot u_{\eta,i}(t)\|^2_{L^2}=\int(1+|x|^2)u_{\eta,i}(t)^2
 \le\|u_{\eta,i}(t)\|_{L^\infty}(1+\mathsf m_T)\le C_s\|u_0\|_{H^s}(1+\mathsf m_T)=:C(T)^2 .
\]

\emph{The bound $\Theta^\eta_i\le\sigma_i+n\eps$.} Young's inequality and
\ref{H:small} give
$|B^\eta_{ij}*u_{\eta,j}|\le\|B^\eta_{ij}\|_{L^1}\|u_{\eta,j}\|_{L^\infty}
\le A_0C_s\|u_0\|_{H^s}=\mathsf b_0$, so the argument lies in
$I=[-2\mathsf b_0,2\mathsf b_0]$. By \ref{H:cutoff},
$a_\eta\ge2\max(1,A_0)C_s\|u_0\|_{H^s}\ge2\mathsf b_0$, so
$I\subseteq[-a_\eta,a_\eta]$ and $f_\eta=f$ on $I$. Finally
$\|f\|_{L^\infty(I)}\le\|f\|_{C^{s+1}(I)}\le\eps$ by \ref{H:small}, the
$C^{s+1}$ norm being the full norm.

\end{proof}

\begin{lemma}[Second moments of the smoothed fluctuation]\label{lem:smoothed}
Let $\bar\fluc^N_j:=\sqrt N(\bar\mu^N_j-u_{\eta,j})$ be the fluctuation field of
the \emph{intermediate} system \eqref{eq:intermediate}, where
$\bar\mu^N_j=\frac1N\sum_k\delta_{\bar X^\eta_{k,j}}$. Then for all $t\le T$ and
$\sigma\ge0$,
\begin{equation}\label{eq:smoothedbound}
 \E\big\|B^\eta_{ij}*\bar\fluc^N_j(t)\big\|_{L^2(\R^d)}^2\le C\,\eta^{-d},
 \qquad
 \E\big\|B^\eta_{ij}*\bar\fluc^N_j(t)\big\|_{H^{\sigma}(\R^d)}^2\le C_\sigma\,\eta^{-d-2\sigma}.
\end{equation}

\end{lemma}

\begin{proof}
Both bounds follow from the variance of a normalised sum of i.i.d.\ centred variables. For the
$H^\sigma$ bound, a change of kernel in Fourier variables reduces it to the $L^2$ bound.

For fixed $t$, the $\bar X^\eta_{k,j}(t)$, $k=1,\dots,N$, are i.i.d.\ by Lemma~\ref{lem:crossspecies}(i), with
density $u_{\eta,j}(t)$ by Lemma~\ref{lem:lawid}(ii).
 Hence
$B^\eta_{ij}*\bar\fluc^N_j(x)=N^{-1/2}\sum_k\big(B^\eta_{ij}(x-\bar
X^\eta_{k,j})-(B^\eta_{ij}*u_{\eta,j})(x)\big)$ is a normalised sum of i.i.d.\
centred variables, so
\[
 \E\big|B^\eta_{ij}*\bar\fluc^N_j(x)\big|^2
 =\mathrm{Var}\big(B^\eta_{ij}(x-\bar X^\eta_{1,j})\big)
 \le\int_{\R^d}B^\eta_{ij}(x-y)^2u_{\eta,j}(y)\dy,
\]
and integrating in $x$ and using Fubini gives
$\|B^\eta_{ij}\|_{L^2}^2\int u_{\eta,j}=\|B^\eta_{ij}\|^2_{L^2}\le C\eta^{-d}$.
For the $H^\sigma$ bound, note $\widehat{B^\eta_{ij}}(\xi)=\widehat
B_{ij}(\eta\xi)$ and, for $\eta\le1$, $\br\xi^{\sigma}\le
C\eta^{-\sigma}\br{\eta\xi}^{\sigma}$, so that
$\br\xi^{\sigma}|\widehat B_{ij}(\eta\xi)|\le
C\eta^{-\sigma}|\widehat{\widetilde B}_{ij}(\eta\xi)|$ with
$\widehat{\widetilde B}_{ij}(\zeta):=\br\zeta^{\sigma}\widehat B_{ij}(\zeta)$ again a
Schwartz function. Denote $\widetilde B^\eta_{ij}:=\eta^{-d}\widetilde B_{ij}(\cdot/\eta)$, so that
$\widehat{\widetilde B^\eta_{ij}}(\xi)=\widehat{\widetilde B}_{ij}(\eta\xi)$. By Plancherel,
\[
 \big\|B^\eta_{ij}*\bar\fluc^N_j\big\|_{H^\sigma}\ \le\ C\eta^{-\sigma}\big\|\widetilde B^\eta_{ij}*\bar\fluc^N_j\big\|_{L^2},
\]
and the previous computation with $\widetilde B^\eta_{ij}$ in place of $B^\eta_{ij}$ gives
$\E\|\widetilde B^\eta_{ij}*\bar\fluc^N_j\|^2_{L^2}\le\|\widetilde B^\eta_{ij}\|^2_{L^2}=\eta^{-d}\|\widetilde B_{ij}\|^2_{L^2}$.
Squaring the Plancherel bound and taking expectations gives the $H^\sigma$ bound of
\eqref{eq:smoothedbound}, with $C_\sigma:=C\max_{ij}\|\widetilde B_{ij}\|^2_{L^2}$, which depends only on $B$, $\sigma$
and $d$.
\end{proof}

\begin{lemma}[Reproducing kernel of the gradient]\label{lem:kernel}
Let $s>\frac d2+1$, $a\ge0$, and let $(e_m)_{m\in\N}$ be any orthonormal basis of
$H^s_a(\R^d)$. Then
\begin{equation}\label{eq:kernel}
 K_{s,a}(x):=\sum_{m=1}^\infty|\na e_m(x)|^2
 =\sum_{\ell=1}^d\big\|\pa_\ell\delta_x\big\|_{(H^s_a)^*}^2
 \le C_{s,a}\,\langle x\rangle^{-2a}\le C_{s,a}
 \qquad\text{for all }x\in\R^d,
\end{equation}
with $C_{s,a}$ independent of the basis, and consequently
\begin{equation}\label{eq:kernelptwise}
 \big|\na\chi(x)\big|^2\ \le\ K_{s,a}(x)\,\|\chi\|^2_{H^s_a}
 \qquad\text{for every }\chi\in H^s_a\text{ and every }x\in\R^d.
\end{equation}
The same statements hold for $K^0_{s,a}(x):=\sum_m|e_m(x)|^2$, with
$|\chi(x)|^2$ in place of $|\na\chi(x)|^2$, and with the same constant $C_{s,a}$, which we
take to be the larger of the two suprema.
\end{lemma}

\begin{proof}
For a bounded linear functional $L$ on a Hilbert space $X$ with orthonormal
basis $(e_m)$, Parseval gives $\sum_m|L(e_m)|^2=\|L\|^2_{X^*}$. Apply this to
$L=\pa_\ell\delta_x$, which is bounded on $H^s_a$ with
$\|\pa_\ell\delta_x\|_{(H^s_a)^*}\le C\langle x\rangle^{-a}$ by
Lemma~\ref{lem:spaces}(iii). For \eqref{eq:kernelptwise}, $\pa_\ell\chi(x)
=\br{\chi,\pa_\ell\delta_x}$ gives
$|\pa_\ell\chi(x)|\le\|\pa_\ell\delta_x\|_{(H^s_a)^*}\|\chi\|_{H^s_a}$, and
summing over $\ell\le d$ is the displayed bound.
\end{proof}

\begin{lemma}[Composition at integer order]\label{lem:moser}
Let $m\ge1$ be an integer, let $J\subset\R$ be a compact interval containing $0$, and let $F\in C^m(J)$
with $F(0)=0$. There is $C=C(m,d)$ such that, for every $w\in H^m(\R^d)$ with values in
$J$ almost everywhere,
\[
 \|F(w)\|_{H^m}\ \le\ C\,\|F\|_{C^m(J)}\big(1+\|w\|_{L^\infty}\big)^{m-1}\|w\|_{H^m}.
\]
\end{lemma}

\begin{proof}
Lemma~\ref{lem:moser} is the Moser estimate of \cite[Prop.~2.1(C)]{Majda84}, in the form with
$(1+\|w\|_{L^\infty})^{m-1}$ and with the dependence on $F$ localised to $J$. We prove the bound for
a mollification $w_\eps$ of $w$, by the chain rule, H\"older's inequality and the
Gagliardo--Nirenberg inequality, and let $\eps\to0$ by weak lower semicontinuity of the $H^m$
norm. Let $w_\eps:=j_\eps*w$ with a mollifier
$j_\eps\ge0$ of mass $1$. Then $w_\eps$ is smooth, takes values in the convex set $J$,
and $\|w_\eps\|_{L^\infty}\le\|w\|_{L^\infty}$, $\|w_\eps\|_{H^k}\le\|w\|_{H^k}$. Since $F(0)=0$ and
$J$ is an interval containing $0$, $|F(w_\eps)|\le\|F'\|_{L^\infty(J)}|w_\eps|$. For
$1\le|\nu|=k\le m$, the chain rule expresses $D^\nu F(w_\eps)$ as a sum of terms
$c\,F^{(j)}(w_\eps)\prod_{l=1}^jD^{\nu_l}w_\eps$ with $1\le j\le k$, $|\nu_l|\ge1$ and
$\sum_l|\nu_l|=k$. By H\"older's inequality with exponents $2k/|\nu_l|$ and the
Gagliardo--Nirenberg inequality
$\|D^{\nu'} v\|_{L^{2k/|\nu'|}}\le C\|v\|_{L^\infty}^{1-|\nu'|/k}\|D^kv\|_{L^2}^{|\nu'|/k}$, each term
has $L^2$ norm at most $C\|F\|_{C^m(J)}\|w\|^{j-1}_{L^\infty}\|w\|_{H^k}$. So $F(w_\eps)$ is
bounded in $H^m$ by the right side of the claim. As $\eps\to0$, $F(w_\eps)\to F(w)$ in $L^2$,
since $F$ is Lipschitz on $J$, and weak lower semicontinuity of the $H^m$ norm gives the
claim.
\end{proof}

\begin{lemma}[Mobility defect in $L^2$]\label{lem:CH2}

Let $L^{\rm full}_{k,ij}:=f_\eta\big((B^\eta_{ij}*\mu^N_j)(X^{N,\eta}_{k,i})\big)
-f_\eta\big((B^\eta_{ij}*u_{\eta,j})(\bar X^{\eta}_{k,i})\big)$, the difference of the mobility terms of the
particle and of the intermediate system. Then
\begin{equation}\label{eq:CH2}
 \sup_{N}\ \sup_{k\le N}\ \E\int_0^T\sum_{j=1}^n\big|L^{\rm full}_{k,ij}(s)\big|^2\dd s\ \le\ C.
\end{equation}
\end{lemma}

\begin{proof}
This is Corollary~\ref{cor:CH2} in Appendix~\ref{sec:appk4}.
\end{proof}

\section{The fluctuation field and the exact identity}\label{sec:field}

\noindent
The argument rests on one exact identity.
Proposition~\ref{prop:ito} applies It\^o's formula to $\br{\fluc^N_i,\varphi}$, and
Theorem~\ref{thm:decomp} derives from it \eqref{eq:fluc_decomp}, a
semimartingale decomposition into a drift that is the linearised operator applied
to $\fluc^N$, three remainders, and a martingale $M^N_i$ whose bracket is
\eqref{eq:cov_finiteN}. The rest of the paper shows that the remainders vanish, that the drift passes to
the limit, and that the bracket converges.

\noindent
We fix two conventions about the notation of
\eqref{eq:fluc_decomp}. First, the drift operator \eqref{eq:GN} is denoted $\mathcal G^\eta$: it depends on $N$ only through
$\eta=\eta(N)$, since its coefficients are built from $u_\eta$ and not from $\mu^N$. Second, the quadratic remainder
\eqref{eq:Rquad} has two summands, and they are referred to separately: the
bilinear part $\mathcal R^{N,\rm bil}_i$ is the first, carrying
$f'_\eta(B^\eta_{ij}*u_{\eta,j})\,(B^\eta_{ij}*\fluc^N_j)$, and the Taylor
part $\mathcal R^{N,\rm Tay}_i$ is the second, carrying $\mathsf T^N_{ij}$ of
\eqref{eq:Taylor-outline2}.

\noindent
The section also establishes what is needed about $M^N$ itself before any limit is
taken: that it is a continuous $\Ufine$-valued square-integrable martingale
(Lemma~\ref{lem:MNval}), whose bracket has an absolutely
continuous trace with an explicit density (Corollary~\ref{cor:trace}) and which
has uniformly bounded fourth moments (Proposition~\ref{prop:M4}). The maximal inequality
for stochastic integrals is Lemma~\ref{lem:BDG}, and Lemma~\ref{lem:field0} records the initial fluctuation.
Two auxiliary facts are proved here as well, since the drift limit,
the evolution family and the bracket density use them:
Lemma~\ref{lem:coeffcont}, the time-continuity of the coefficients
$u_{\eta,i}$ and $u_{\eta,i}f'_\eta(B^\eta_{ij}*u_{\eta,j})$ of the backward
system, with a modulus uniform in $\eta$,
and Lemma~\ref{lem:L2field}, a bound on the square of the averaged mobility $\br{\mu^N_i,\Theta^N_i}$.

\emph{Convention (representative).} Throughout, for every $\eta\in[0,\eta_0]$, $u_\eta$ is the
weakly continuous representative of Lemma~\ref{lem:lawid}(i). So the bound of
\ref{H:small} holds for every $t\in[0,T]$. By weak continuity in $L^2$ and weak lower
semicontinuity of the norm, the first term of \eqref{eq:CHbias} is also a bound at every
$t$.

\begin{lemma}[Time-continuity of the coefficients]\label{lem:coeffcont}
Let $\eta\in[0,\eta_0]$, with the \hyperref[par:eta0]{convention at $\eta=0$} of Section~\ref{sec:setting},
and let $g^\eta_{ij}=u_{\eta,i}f'_\eta(w^\eta_{ij})$, the coefficient of the nonlocal term of
\eqref{eq:backward-mild}, be as in Lemma~\ref{lem:lawid}. Assume \ref{H:CH}--\ref{H:small}, and let
$C(T)$ be the bound $\|\br\cdot u_{\eta,i}(t)\|_{L^2}\le C(T)$ of Lemma~\ref{lem:secondmoment}.
Let $C_{s-1}$ be the norm of the embedding $H^{s-1}(\R^d)\hookrightarrow L^\infty(\R^d)$, and put
\begin{equation}\label{eq:Lu}
 M_0:=\|u_0\|_{H^s},\qquad
 L_u:=\|\diver\|_{L^2\to H^{-1}}\Big(C(T)+M_0\max_i\sigma_i+n\eps M_0+n\eps A_0C_{s-1}M_0^2\Big).
\end{equation}
Then $\|\pa_tu_{\eta,i}(t)\|_{H^{-1}}\le L_u$ for almost every $t\in[0,T]$, and:
\begin{enumerate}[label=\normalfont(\roman*)]
\item for every $s'\in[0,s)$, $u_{\eta,i}\in C\big([0,T];H^{s'}\big)$, with the
$\eta$-uniform H\"older modulus
\begin{equation}\label{eq:coeffholder}
 \|u_{\eta,i}(t)-u_{\eta,i}(t_0)\|_{H^{s'}}
 \ \le\ \big(L_u|t-t_0|\big)^{1-\vartheta}(2M_0)^{\vartheta},
 \qquad \vartheta=\frac{s'+1}{s+1}\in(0,1),
\end{equation}
with $M_0$ and $L_u$ as in \eqref{eq:Lu};
\item if $s\ge s_1+4$, then for every $s'\in[0,s-2]$,
\begin{equation}\label{eq:coeffcont}
 \Theta^\eta_i-\kappa_i\in C\big([0,T];H^{s'}\big),
 \qquad
 g^\eta_{ij}\in C\big([0,T];H^{s'}\big),
 \qquad \kappa_i:=\sigma_i+nf_\eta(0).
\end{equation}
\end{enumerate}
In particular, since $s\ge s_1+4$
and hence $s_1+2\le s-2$, the case $s'=s_1+2$ shows that the coefficients of
$\Delta$ in $\mathcal G^\eta(r)$, after subtraction of their limits at infinity,
are continuous in $r$ as multipliers on $H^{s_1+1}_{a_1}$.

\end{lemma}

\begin{proof}
Step~1 bounds $\pa_tu_{\eta,i}$ in $H^{-1}$ by $L_u$, from the divergence form of
\eqref{eq:nonloc}. Step~2 interpolates the resulting Lipschitz bound in $H^{-1}$ against the $H^s$
bound of \ref{H:small}, which gives (i). Step~3 proves (ii): a Moser composition estimate at an
integer index is interpolated against continuity in $L^2$.

\emph{Step 1: $\pa_tu_{\eta,i}\in L^\infty(0,T;H^{-1})$, $\eta$-uniformly.}
In divergence form, \eqref{eq:nonloc} becomes
$\pa_tu_{\eta,i}=\diver\big(u_{\eta,i}\na U_i
+\na[\sigma_iu_{\eta,i}+u_{\eta,i}\sum_jf_\eta(B^\eta_{ij}*u_{\eta,j})]\big)$,
whose bracket is $\Theta^\eta_iu_{\eta,i}$. Since $\diver:L^2\to H^{-1}$ is
bounded, it suffices that the argument lie in $L^\infty(0,T;L^2)$.
The first summand is $-xu_{\eta,i}$. Its coefficient $\na U_i=-x$ is
unbounded, and we use the bound from Lemma~\ref{lem:secondmoment} assumed in the
statement: $\|xu_{\eta,i}\|_{L^2}\le\|\br\cdot u_{\eta,i}\|_{L^2}\le C(T)$.
For the second, differentiate the product directly:
\[
 \Big\|\na\Big(u_{\eta,i}\sum_jf_\eta(w^\eta_{ij})\Big)\Big\|_{L^2}
 \le n\|f_\eta\|_{L^\infty(I)}\|\na u_{\eta,i}\|_{L^2}
 +n\|f'_\eta\|_{L^\infty(I)}\,A_0\,\max_j\|\na u_{\eta,j}\|_{L^\infty}\,\|u_{\eta,i}\|_{L^2},
\]
since each of the $n$ summands is bounded by its maximum over $j$, and
$|B^\eta_{ij}*\na u_{\eta,j}|\le\|B^\eta_{ij}\|_{L^1}\|\na u_{\eta,j}\|_{L^\infty}$ by Young's inequality.
The first summand of the bracket adds $\sigma_i\|\na u_{\eta,i}\|_{L^2}$.

Every factor is finite and $\eta$-uniform. First,
$\|\na u_{\eta,i}(t)\|_{L^2}\le\|u_{\eta,i}(t)\|_{H^s}\le M_0$ for a.e.\ $t$, since
$s\ge1$. Next,
$\|\na u_{\eta,j}\|_{L^\infty}\le C_{s-1}\|u_{\eta,j}\|_{H^{s}}$ by
$H^{s-1}\hookrightarrow L^\infty$ ($s>\frac d2+1$), and
$\|B^\eta_{ij}\|_{L^1}=\|B_{ij}\|_{L^1}\le A_0$. Finally,
$\|f_\eta\|_{L^\infty(I)},\|f'_\eta\|_{L^\infty(I)}\le\eps$ by \ref{H:small}
together with \ref{H:cutoff}, which gives $f_\eta=f$ on $I$.Hence $\|\pa_tu_{\eta,i}(t)\|_{H^{-1}}\le L_u$ for almost every $t$, with $L_u$ as in
\eqref{eq:Lu}.

\emph{Step 2: interpolation, with a modulus.} $u_{\eta,i}\in L^\infty(0,T;H^s)$
and $\pa_tu_{\eta,i}\in L^\infty(0,T;H^{-1})$ give a representative that is
$L_u$-Lipschitz $[0,T]\to H^{-1}$, namely the Bochner integral of the derivative. It
coincides with the representative of Lemma~\ref{lem:lawid}(i), since both are weakly
continuous in $H^{-1}$ and agree almost everywhere. So its $H^s$ norm is at most $M_0$ at
every $t$.
We work with this representative throughout. For $-1<s'<s$, H\"older on the Fourier side
with $\br\xi^{2s'}=(\br\xi^{-2})^{1-\vartheta}(\br\xi^{2s})^{\vartheta}$ requires
$s'=-(1-\vartheta)+s\vartheta$, i.e.\ $\vartheta=\frac{s'+1}{s+1}$, and yields
$\|v\|_{H^{s'}}\le\|v\|^{1-\vartheta}_{H^{-1}}\|v\|^{\vartheta}_{H^{s}}$. Apply
it to $v=u_{\eta,i}(t)-u_{\eta,i}(t_0)$, whose $H^s$ norm is $\le2M_0$ by
\ref{H:small}: this is \eqref{eq:coeffholder}, and it proves (i).

\emph{Step 3: composition, with continuity.} Fix $s'\in[0,s-2]$ and put
$m:=\lfloor s'\rfloor+1$, an integer with $s'<m\le s-1$. Convolution is a Fourier
multiplier, so $B^\eta_{ij}*$ is bounded on every $H^k$ with norm
$\le\|B_{ij}\|_{L^1}\le A_0$, uniformly in $\eta$, and equals $a_{ij}$ at
$\eta=0$. Recall $w^\eta_{ij}=B^\eta_{ij}*u_{\eta,j}$. By (i), $w^\eta_{ij}\in
L^\infty(0,T;H^{m})\cap C([0,T];L^2)$, and by Lemma~\ref{lem:lawid}(i) it takes values in
$[-\mathsf b_0,\mathsf b_0]\subset I$. Then:
\emph{(a) a uniform bound at the integer index $m$.} Lemma~\ref{lem:moser} with $J=I$,
applied to $f-f(0)$ and to $f'-f'(0)$, bounds $f_\eta(w^\eta_{ij})-f(0)$ and
$f'_\eta(w^\eta_{ij})-f'(0)$ in $H^m$, uniformly in $t$ and $\eta$. It applies because
$f_\eta=f$ on $I$ by \ref{H:cutoff}, $0\in I$, and $\|f\|_{C^{m+1}(I)}\le\eps$ by
\ref{H:small}, since $m+1\le s+1$.
\emph{(b) continuity at the bottom.} $|f_\eta(w^\eta_{ij}(t))-f_\eta(w^\eta_{ij}(t_0))|\le
\|f'_\eta\|_{L^\infty(I)}|w^\eta_{ij}(t)-w^\eta_{ij}(t_0)|$ pointwise, so $t\mapsto
f_\eta(w^\eta_{ij}(t))-f_\eta(0)$ is $L^2$-continuous, and likewise for $f'_\eta$
using $\|f''_\eta\|_{L^\infty(I)}$. Both norms are $\le\eps$, as in Step~1.
\emph{(c) Interpolate} (a) against (b) at $s'<m$, as in Step~2. Summing the
$n$ terms gives the first inclusion of \eqref{eq:coeffcont} with $\kappa_i$ as
stated. For the second inclusion it suffices to treat
$s'\in[s_1+2,s-2]$, a non-empty range since
$s\ge s_1+4$; for $s'<s_1+2$ the claim then follows from the case $s'=s_1+2$ by
$H^{s_1+2}\hookrightarrow H^{s'}$. On $[s_1+2,s-2]$, split
$g^\eta_{ij}=u_{\eta,i}\big(f'_\eta(w^\eta_{ij})-f'_\eta(0)\big)
+f'_\eta(0)\,u_{\eta,i}$: both terms are continuous into $H^{s'}$, the first by
the $H^{s'}$ algebra ($s'\ge s_1+2>\frac d2$) applied to two factors continuous
by (i) and (c), the second by (i).

\end{proof}

\begin{lemma}\label{lem:L2field}
$\sup_N\E\int_0^T\br{\mu^N_i,\Theta^N_i}^2\dd s\le C$.
\end{lemma}

\begin{proof}
Decompose $\Theta^N_i(X_{k,i})=\sigma_i+\sum_jf_\eta\big((B^\eta_{ij}*u_{\eta,j})(\bar X^\eta_{k,i})\big)+\sum_jL^{\rm full}_{k,ij}$
with $L^{\rm full}_{k,ij}$ as in Lemma~\ref{lem:CH2}. The comparison value
satisfies $f_\eta\big((B^\eta_{ij}*u_{\eta,j})(\bar X^\eta_{k,i})\big)\le\sup_If=:F$,
because $\|B^\eta_{ij}*u_{\eta,j}\|_\infty\le A_0C_s\|u_0\|_{H^s}=\mathsf b_0$, so the
argument lies in $I=[-2\mathsf b_0,2\mathsf b_0]$, where $f_\eta=f$ by \ref{H:cutoff}. Hence, with
$b_k:=\sum_j|L^{\rm full}_{k,ij}|$,
\[
 \br{\mu^N_i,\Theta^N_i}\le\sigma_i+nF+\frac1N\sum_kb_k,
 \qquad
 \br{\mu^N_i,\Theta^N_i}^2\le2(\sigma_i+nF)^2+\frac{2n}{N}\sum_k\sum_j|L^{\rm full}_{k,ij}|^2,
\]
by $(x+y)^2\le2x^2+2y^2$, Jensen in $k$ and
Cauchy--Schwarz in $j$. The first bound may be squared because
$\br{\mu^N_i,\Theta^N_i}\ge\sigma_i>0$, since $f_\eta\ge0$. Integrate, take expectations, and use \eqref{eq:CH2}:
$\frac{2n}{N}\sum_k\E\int\sum_j|L|^2\le 2nC$.
\end{proof}

\begin{lemma}[Maximal inequality for stochastic integrals]\label{lem:BDG}
Let $H$ and $\mathsf U$ be separable Hilbert spaces, $W$ a cylindrical Wiener process on $\mathsf U$
with respect to a normal filtration (for instance a standard Brownian motion in
$\mathsf U=\R^m$), and $\Upsilon$ a predictable process with values in the Hilbert--Schmidt operators
$\HS(\mathsf U;H)$ such that $\int_0^T\|\Upsilon\|^2_{\HS}\dd s<\infty$
almost surely. Put $M(t):=\int_0^t\Upsilon\dd W$. Then, for every $p>0$,
\begin{equation}\label{eq:BDG}
 \E\Big[\sup_{t\le T}\|M(t)\|^p_H\Big]\ \le\ C_p\,
 \E\Big[\Big(\int_0^T\|\Upsilon(s)\|^2_{\HS}\dd s\Big)^{p/2}\Big],
\end{equation}
with $C_p$ depending only on $p$. If $\E\int_0^T\|\Upsilon\|^2_{\HS}\dd s<\infty$, then $M$ is
a continuous square-integrable martingale, and
$\tr\br{\!\br M\!}_t:=\sum_m\br{\br{M,h_m}}_t=\int_0^t\|\Upsilon(s)\|^2_{\HS}\dd s$ for any
orthonormal basis $(h_m)$ of $H$. Let $H'$ be a separable Hilbert space and
$\mathsf J(r)\in\mathcal L(H;H')$ deterministic, such that $\mathsf J\Upsilon$ is predictable with
values in $\HS(\mathsf U;H')$ and $\E\int_0^T\|\mathsf J\Upsilon\|^2_{\HS}\dd r<\infty$. We define
$\int_0^t\mathsf J\dd M:=\int_0^t\mathsf J\Upsilon\dd W$. Then \eqref{eq:BDG} applies to it, and
$\|\mathsf J\Upsilon\|_{\HS}\le\|\mathsf J\|\,\|\Upsilon\|_{\HS}$. For elementary $\mathsf J$ this is the
Riemann sum $\sum_l\mathsf J_l(M(t_{l+1})-M(t_l))$, and for $h\in H'$ the scalar process
$\br{\int_0^t\mathsf J\dd M,h}$ is the real stochastic integral of
$\br{\Upsilon(\cdot)\,\cdot\,,\mathsf J^*h}$.
\end{lemma}

\begin{proof}
The inequality is \cite[Thm.~4.36]{da_prato_2014_SPDEs} with $Q=\mathrm{Id}$, so that
$L^0_2=\HS(\mathsf U;H)$. For a cylindrical $W$ the integral and its properties are those of
\cite[\S4.2.1]{da_prato_2014_SPDEs}. The martingale property and the quadratic variation
$\br{\!\br M\!}_t=\int_0^t\Upsilon\Upsilon^*\dd s$ are \cite[Thm.~4.27]{da_prato_2014_SPDEs}. The
trace of $\Upsilon\Upsilon^*$ is $\|\Upsilon\|^2_{\HS}$, and
$\br{\!\br M\!}_th\cdot h=\br{\br{M,h}}_t$ by the definition of the quadratic
variation \cite[\S3.4]{da_prato_2014_SPDEs}. In the last part, the integral is a
definition. A bounded deterministic operator $\mathsf J_0$ satisfies $\mathsf J_0\int\Upsilon\dd W=\int \mathsf J_0\Upsilon\dd W$: this
holds for elementary $\Upsilon$ and extends by the isometry. Applied on each interval of
constancy of an elementary $\mathsf J$, it gives the Riemann sum. Applied to the functional
$\br{\cdot,h}$, it gives the pairing, where $\mathsf J^*$ is the Hilbert adjoint.
\end{proof}

\begin{proposition}[Exact finite-$N$ identity]\label{prop:ito}
Let $\varphi\in H^{m}_{a}(\R^d)$ with $m>\frac d2+2$, $a\ge0$. Then, for every
$i$, $\Prob$-a.s.\ and for all $t\in[0,T]$,
\begin{equation}\label{eq:ito_emp}
 \d\br{\mu_i^N,\varphi}
 =\br{\mu_i^N,\ -\na U_i\cdot\na\varphi+\Theta_i^N\Delta\varphi}\dt
 +\br{\mu_i^N,\varrho^N_{\cdot,i}\Delta\varphi}\dt
 +\frac1{\sqrt N}\,\d M_i^N(\varphi),
\end{equation}
where $\Theta^N_i$ and $\varrho^N_{k,i}$ are as in
\eqref{eq:mobilities}--\eqref{eq:self}, the second term is shorthand for
$\frac1N\sum_k\varrho^N_{k,i}\Delta\varphi(X^{N,\eta}_{k,i})$, and
\begin{equation}\label{eq:MiN}
 M_i^N(t,\varphi):=\int_0^t\frac1{\sqrt N}\sum_{k=1}^N
 \big(2\Theta_i^N+2\varrho_{k,i}^N\big)^{1/2}(X_{k,i}^{N,\eta})\,
 \na\varphi(X_{k,i}^{N,\eta})\cdot\d W_i^k(s)
\end{equation}
is a continuous square-integrable real martingale with predictable bracket
\begin{equation}\label{eq:cov_finiteN}
 \br{M_i^N(\cdot,\varphi),M_i^N(\cdot,\psi)}_t
 =\int_0^t\br{\mu_i^N,\ \na\varphi\cdot\na\psi\,
 \big(2\Theta_i^N+2\varrho^N_{\cdot,i}\big)}\ds,
\end{equation}
and $\br{M_i^N(\cdot,\varphi),M_k^N(\cdot,\psi)}\equiv0$ for $i\ne k$. The definition
\eqref{eq:MiN}, the martingale property and the brackets hold for all $\varphi,\psi\in H^m_a$ with
$m>\frac d2+1$, since they involve only $\na\varphi$ and $\na\psi$, which are bounded by
Lemma~\ref{lem:spaces}(iii).
\end{proposition}

\begin{proof}
The generator of \eqref{eq:particles} acting on the $k$-th particle of species
$i$ is $b_i\cdot\na+A_{k,i}\Delta$ with $b_i=-\na U_i$ and
$A_{k,i}=\Theta^N_i(X_{k,i})+\varrho^N_{k,i}$: indeed the diffusion coefficient
in \eqref{eq:particles} is the scalar $(2A_{k,i})^{1/2}$ times the
identity on $\R^d$, so $\frac12\sigma\sigma^{\!\top}=A_{k,i}\mathrm{Id}$ and the
second-order part of the generator is $A_{k,i}\Delta$. It\^o's formula applied
to $\varphi(X^{N,\eta}_{k,i})$, summed over $k$ and divided by $N$, gives
\eqref{eq:ito_emp} with the martingale part
$\frac1N\sum_k(2A_{k,i})^{1/2}\na\varphi(X_{k,i})\cdot\d W^k_i$, which is
$N^{-1/2}$ times \eqref{eq:MiN}. The bracket follows from independence of the
$(W_i^k)_{k,i}$.
The stochastic integrand is bounded, since $\na\varphi$ is bounded by
Lemma~\ref{lem:spaces}(iii) with $m>\frac d2+2$ and the diffusion coefficient is bounded.
The drift is integrable along the paths: $\Delta\varphi$ is bounded,
$|x\cdot\na\varphi(x)|\le|x|\,\|\na\varphi\|_{L^\infty}$, and the paths are
square-integrable by \cite[Prop.~4(i)]{chen_holzinger_2021_rigorous_SKT}.
\end{proof}

\begin{theorem}[Finite-$N$ fluctuation dynamics]\label{thm:decomp}
Let $\varphi\in H^{s_1+2}_{a_1}(\R^d)$. Then, $\Prob$-a.s.\ for all $t\in[0,T]$,
\begin{equation}\label{eq:fluc_decomp}
 \d\br{\fluc_i^{N},\varphi}
 =\sum_{j=1}^n\br{\fluc_j^{N},\,(\mathcal G^\eta\varphi)_{ij}}\dt
 \;+\;\mathcal R_i^N(\varphi)\dt\;+\;\d M_i^N(\varphi),
\end{equation}
where the deterministic linear operator $\mathcal G^\eta$ is given by
\begin{equation}\label{eq:GN}
 (\mathcal G^\eta\varphi)_{ij}
:=\delta_{ij}\Big(-\na U_i\cdot\na\varphi+\Theta^\eta_i\,\Delta\varphi\Big)
 \;+\;B^\eta_{ij}*\Big(u_{\eta,i}\,f_\eta'\big(B^\eta_{ij}*u_{\eta,j}\big)\,\Delta\varphi\Big),
\end{equation}
and the remainder is $\mathcal R^N_i=\mathcal R^{N,\rm self}_i+\mathcal
R^{N,\rm quad}_i$ with
\begin{align}
 \mathcal R^{N,\rm self}_i(\varphi)
 &:=\sqrt N\,\br{\mu_i^N,\ \varrho^N_{\cdot,i}\,\Delta\varphi},
 \label{eq:Rself}\\
 \mathcal R^{N,\rm quad}_i(\varphi)
 &:=\frac1{\sqrt N}\sum_{j=1}^n\br{\fluc^N_i,\
 f_\eta'\big(B^\eta_{ij}*u_{\eta,j}\big)\,\big(B^\eta_{ij}*\fluc^N_j\big)\,\Delta\varphi}
 \;+\;\br{\mu^N_i,\ \sum_{j=1}^n\mathsf T^N_{ij}\,\Delta\varphi},
 \label{eq:Rquad}
\end{align}
where $\mathsf T^N_{ij}$ is the second-order Taylor remainder
\begin{equation}\label{eq:Taylor-outline2}
 \mathsf T^N_{ij}
:=\sqrt N\Big[f_\eta\big(B^\eta_{ij}*\mu^N_j\big)-f_\eta\big(B^\eta_{ij}*u_{\eta,j}\big)
 -f_\eta'\big(B^\eta_{ij}*u_{\eta,j}\big)\,\tfrac1{\sqrt N}\big(B^\eta_{ij}*\fluc^N_j\big)\Big],
 \qquad
 |\mathsf T^N_{ij}|\le\frac{\|f_\eta''\|_{L^\infty(\R)}}{2\sqrt N}\big|B^\eta_{ij}*\fluc^N_j\big|^2.
\end{equation}
The last bound is not uniform in $\eta$ and is used only at fixed $N$
(Lemma~\ref{lem:Rdet}). The uniform bound is Lemma~\ref{lem:ptwise}.
\end{theorem}

\begin{proof}
We subtract the weak form \eqref{eq:weak_nonloc} of the equation for $u_\eta$ from
\eqref{eq:ito_emp} and multiply by $\sqrt N$. The transport and diffusion terms give the diagonal part
of \eqref{eq:GN}. A second-order Taylor expansion of $f_\eta$ around $B^\eta_{ij}*u_{\eta,j}$ gives the
off-diagonal part and $\mathcal R^{N,\rm quad}_i$.

By \eqref{eq:lawweak}, \eqref{eq:weak_nonloc} holds for $\varphi\in C^\infty_0$. It extends to
$\varphi\in H^{s_1+2}_{a_1}$ by density, since $|\br{u_{\eta,i},x\cdot\na\varphi}|\le
\|\br\cdot u_{\eta,i}\|_{L^2}\|\na\varphi\|_{L^2}$ is controlled by Lemma~\ref{lem:secondmoment}.
Subtract \eqref{eq:weak_nonloc} from \eqref{eq:ito_emp} and multiply by
$\sqrt N$. The drift terms split as
\[
 \sqrt N\br{\mu^N_i-u_{\eta,i},-\na U_i\cdot\na\varphi}
 =\br{\fluc^N_i,-\na U_i\cdot\na\varphi},
\]
which is the first, diagonal, part of \eqref{eq:GN}, and
\[
 \sqrt N\Big[\br{\mu^N_i,\Theta^N_i\Delta\varphi}-\br{u_{\eta,i},\Theta^\eta_i\Delta\varphi}\Big]
 =\underbrace{\br{\fluc^N_i,\Theta^\eta_i\Delta\varphi}}_{=:\,\mathrm{I}}
 +\underbrace{\br{\mu^N_i,\sqrt N\big(\Theta^N_i-\Theta^\eta_i\big)\Delta\varphi}}_{=:\,\mathrm{II}},
\]
by adding and subtracting $\br{\mu^N_i,\Theta^\eta_i\Delta\varphi}$. Term
$\mathrm I$ is the second, diagonal, part of \eqref{eq:GN}. For $\mathrm{II}$,
Taylor's theorem applied to $f_\eta$ at the point $B^\eta_{ij}*u_{\eta,j}$ gives,
using $\sqrt N(B^\eta_{ij}*\mu^N_j-B^\eta_{ij}*u_{\eta,j})=B^\eta_{ij}*\fluc^N_j$,
\[
 \sqrt N\big(\Theta^N_i-\Theta^\eta_i\big)
 =\sum_{j=1}^n f_\eta'\big(B^\eta_{ij}*u_{\eta,j}\big)\big(B^\eta_{ij}*\fluc^N_j\big)
 +\sum_{j=1}^n\mathsf T^N_{ij},
\]
with $\mathsf T^N_{ij}$ as in \eqref{eq:Taylor-outline2}. Inserting this identity and substituting
$\mu^N_i=u_{\eta,i}+N^{-1/2}\fluc^N_i$ in the first sum,
\[
 \mathrm{II}
 =\sum_{j}\br{u_{\eta,i}\,f_\eta'\big(B^\eta_{ij}*u_{\eta,j}\big)\Delta\varphi,\
 B^\eta_{ij}*\fluc^N_j}
 +\mathcal R^{N,\rm quad}_i(\varphi).
\]
Moving the convolution onto the other factor of the pairing (by Fubini, since
$B^\eta_{ij}$ is even) turns the first
sum into $\sum_j\br{\fluc^N_j,B^\eta_{ij}*(u_{\eta,i}f'_\eta(B^\eta_{ij}*u_{\eta,j})\Delta\varphi)}$,
the off-diagonal part of \eqref{eq:GN}. The self-interaction term of
\eqref{eq:ito_emp}, multiplied by $\sqrt N$, is \eqref{eq:Rself}. Finally the
martingale part of \eqref{eq:ito_emp} carries a factor $N^{-1/2}$ which is
cancelled by the $\sqrt N$, leaving $M^N_i(\varphi)$ exactly.
\end{proof}

\begin{lemma}[Initial fluctuation]\label{lem:field0}
$\displaystyle\sup_N\E\|\fluc^N_i(0)\|^2_{\Ufine}\le C_{s_1,a_1}$, with no moment
assumption on $u_{0,i}$ beyond $u_{0,i}$ being a probability density.
\end{lemma}

\begin{proof}
$\fluc^N_i(0)=N^{-1/2}\sum_k(\delta_{\xi^k_i}-u_{0,i})$ with $\xi^k_i$ i.i.d.
Hence, for an ONB $(e_m)$ of $H^{s_1}_{a_1}$,
$\E\|\fluc^N_i(0)\|^2_{\Ufine}=\sum_m\mathrm{Var}\,e_m(\xi^1_i)
\le\sum_m\E[e_m(\xi^1_i)^2]=\E\big[K^0_{s_1,a_1}(\xi^1_i)\big]\le C_{s_1,a_1}$
by Lemma~\ref{lem:kernel}.
\end{proof}

\begin{lemma}[$M^N_i$ as a continuous $\Ufine$-valued martingale]\label{lem:MNval}
Fix $N$ and $\eta\in(0,\eta_0]$. For $x\in\R^d$ let $\na\delta_x\in\Ufine^{\,d}$ be the
vector of functionals $\varphi\mapsto\pa_\ell\varphi(x)$, and put
$\Gamma^N_{k,i}:=\big(2\Theta^N_i(X^{N,\eta}_{k,i})+2\varrho^N_{k,i}\big)^{1/2}$. Then
\begin{equation}\label{eq:MNval}
 M^N_i(t)\:=\ \frac1{\sqrt N}\sum_{k=1}^N\int_0^t
 \Gamma^N_{k,i}(s)\,\na\delta_{X^{N,\eta}_{k,i}(s)}\cdot\dd W^k_i(s)
\end{equation}
is a \emph{continuous, square-integrable, $\Ufine$-valued} martingale with
$M^N_i(0)=0$. Its pairing with $\varphi\in H^{s_1}_{a_1}$ is the real
martingale \eqref{eq:MiN} of Proposition~\ref{prop:ito}, the species cross-brackets vanish,
$M^N=(M^N_1,\dots,M^N_n)=\int_0^\cdot\Upsilon^N\dd W^N$ with $W^N:=(W^k_i)_{k\le N,i\le n}$, a
Brownian motion in $\R^{nNd}$, and the predictable integrand
$\Upsilon^N(s)\in\HS(\R^{nNd};\Ufine^n)$ that maps the unit vector $e_{k,i,\ell}$ to
$N^{-1/2}\Gamma^N_{k,i}(s)\,\pa_\ell\delta_{X^{N,\eta}_{k,i}(s)}$ in the $i$-th component,
and
\begin{equation}\label{eq:trMN}
 \tr\br{\!\br{M^N_i}\!}_t=\int_0^t\mathsf q^N_i\dd s,
 \qquad
 \mathsf q^N_i:=\br{\mu^N_i,\,K_{s_1,a_1}\big(2\Theta^N_i+2\varrho^N_{\cdot,i}\big)},
 \qquad\|\Upsilon^N(s)\|^2_{\HS}=\sum_i\mathsf q^N_i(s).
\end{equation}
Moreover $2\sigma_i\le(\Gamma^N_{k,i})^2\le c_\eta$ and $0\le\mathsf q^N_i\le c_\eta C_{s_1,a_1}$, where
$c_\eta:=2\max_i\sigma_i+2n\|f_\eta\|_{L^\infty(\R)}$ and $C_{s_1,a_1}$ is the constant of
Lemma~\ref{lem:kernel}; $c_\eta$ does not depend on $N$ or $i$.
\end{lemma}

\begin{proof}
We show that the integrand of \eqref{eq:MNval} is a bounded predictable process with
values in $\Ufine^{\,d}$, so that Lemma~\ref{lem:BDG} applies. The trace identity then follows from the
Hilbert--Schmidt norm of $\Upsilon^N$ and Lemma~\ref{lem:kernel}.

$\|\na\delta_x\|^2_{\Ufine^{\,d}}=\sum_{\ell\le d}\|\pa_\ell\delta_x\|^2_{(H^{s_1}_{a_1})^*}
=K_{s_1,a_1}(x)\le C_{s_1,a_1}$ for every $x$, by Lemma~\ref{lem:kernel} at
$(s,a)=(s_1,a_1)$. Next,
$2\Theta^N_i+2\varrho^N_{k,i}=2\sigma_i+2\sum_jf_\eta(\cdots)$ with the
self-term deleted, so with $f_\eta\ge0$ and $f_\eta\in L^\infty(\R)$ globally
by \ref{H:CH},
\[
 0\ <\ 2\sigma_i\ \le\ 2\Theta^N_i+2\varrho^N_{k,i}\ \le\
 2\max_i\sigma_i+2n\|f_\eta\|_{L^\infty(\R)}\ =\ c_\eta,
\]
with $c_\eta$ as in the statement. The integrand in \eqref{eq:MNval} is adapted with
$\Ufine^{\,d}$-norm $\le(c_\eta C_{s_1,a_1})^{1/2}$ uniformly in $s$ and $\omega$,
and is predictable: $s\mapsto\br{\na\delta_{X_{k,i}(s)},\varphi}
=\na\varphi(X_{k,i}(s))$ is continuous for each $\varphi\in H^{s_1+1}_{a_1}$,
which is dense in $\Ufine^*$. Each such scalar process is continuous and
adapted, hence predictable, the Borel $\sigma$-field of the separable
$\Ufine^{\,d}$ is generated by countably many of these functionals. Hence, the
integrand is predictable.
So the It\^o integral against the $N$ independent $d$-dimensional Brownian
motions $W^k_i$, $k\le N$, exists and is a
continuous square-integrable $\Ufine$-valued
martingale vanishing at $0$. Pairing with $\varphi$ and using
$\br{\na\delta_x,\varphi}=\na\varphi(x)$ returns \eqref{eq:MiN}. The
cross-brackets vanish by independence of the $(W^k_i)_{k,i}$ across $i$, as in
Proposition~\ref{prop:ito}. Finally,
$\|\Upsilon^N(s)\|^2_{\HS}=\sum_{k,i}N^{-1}(\Gamma^N_{k,i})^2\|\na\delta_{X_{k,i}}\|^2_{\Ufine^d}
=\sum_i\br{\mu^N_i,K_{s_1,a_1}(2\Theta^N_i+2\varrho^N_{\cdot,i})}=\sum_i\mathsf q^N_i$, and
\eqref{eq:trMN} is Lemma~\ref{lem:BDG}. $\mathsf q^N_i\ge0$ because
$K_{s_1,a_1}\ge0$, and $\mathsf q^N_i\le c_\eta C_{s_1,a_1}$ because $\mu^N_i$ is a probability measure.
 \end{proof}

\begin{corollary}\label{cor:trace}
The process $\rho^N_i(t):=\tr\br{\!\br{M^N_i}\!}_t$ is absolutely continuous,
$\dd\rho^N_i=\mathsf q^N_i\dd s$, with $\mathsf q^N_i$ as in \eqref{eq:trMN}, and
\begin{equation}\label{eq:psi}
 0\le\mathsf q^N_i\le2C_{s_1,a_1}\big(\br{\mu^N_i,\Theta^N_i}+\|\varrho^N\|_\infty\big),
\end{equation}
with $C_{s_1,a_1}$ the constant of Lemma~\ref{lem:kernel}, and
\begin{equation}\label{eq:psi2}
 \sup_N\E\int_0^T\big(\mathsf q^N_i\big)^2\dd s\le C,
 \qquad
 \sup_N\E\Big[\big(\rho^N_i(T)\big)^{2}\Big]\le C.
\end{equation}
\end{corollary}

\begin{proof}
By \eqref{eq:trMN}, $\rho^N_i$ is absolutely continuous with density $\mathsf q^N_i$, and
$\mathsf q^N_i\ge0$ by Lemma~\ref{lem:MNval}. The upper bound in \eqref{eq:psi} follows from the
definition of $\mathsf q^N_i$ in \eqref{eq:trMN}, from $K_{s_1,a_1}\le C_{s_1,a_1}$
(Lemma~\ref{lem:kernel}) and from $\Theta^N_i\ge0$, since $\mu^N_i$ is a probability measure.

\emph{\eqref{eq:psi2}.} Lemma~\ref{lem:L2field},
$\|\varrho^N\|_\infty\le CN^{-1}\eta^{-d-\alpha}\to0$ (Lemma~\ref{lem:polylog}), and $\E[(\int_0^T\mathsf q^N_i\dd s)^2]\le T\,\E\int_0^T(\mathsf q^N_i)^2\dd s$.
\end{proof}

\begin{proposition}\label{prop:M4}
$\sup_N\E\big[\sup_{t\le T}\|M^N_i(t)\|^4_{\Ufine}\big]\le C$.
\end{proposition}

\begin{proof}
Lemma~\ref{lem:BDG} at $p=4$ with \eqref{eq:trMN} and \eqref{eq:psi2}.
\end{proof}

\section{The backward evolution family}\label{sec:evolfam}

\noindent
The mild formulation needs a two-parameter family $(S^\eta(t,s))$ on $\Ufine$
whose adjoint solves the backward system dual to the drift operator. The
construction is deterministic parabolic theory on weighted Bessel spaces and is
carried out in Appendix~\ref{sec:appevol}. The standard existence theory for linear parabolic
equations, J.-L.~Lions' theorem \cite[Ch.~XVIII, \S3]{DaLi92}, does not apply to the backward system
directly, because the drift $x\cdot\na$ loses one power of the weight and has an unbounded
coefficient. In the appendix a time reversal and a dilation remove the drift, and a weight
conjugation moves the problem to unweighted spaces; Lions' theorem then applies, and our energy estimate gives the bounds uniformly in $\eta$. This section records the statements
the fluctuation argument uses. Proposition~\ref{prop:evol} gives the backward system,
its unique solution at each of the three levels used in the argument, and the family with its bound on
$\Ufine$ and on $\Usharp$. Lemma~\ref{lem:Gbound} bounds the drift operator on the test class $\Test$.
Lemma~\ref{lem:strongcont} gives the strong continuity in the evaluation time, and
Lemma~\ref{lem:fwdlip} the Lipschitz dependence on the terminal time, which gives the
continuity in the forward variable of Corollary~\ref{cor:fwdweak}. Lemma~\ref{lem:levelcoin} shows that the families at the
two levels used by the mild formulation agree where both act.

\begin{proposition}[Backward evolution family]\label{prop:evol}\label{prop:evolution-outline2}
Let $\eta\in[0,\eta_0]$, $t\in[0,T]$, let $(\sigma,a)$ be one of the levels $(s_1,a_1)$,
$(s_1+1,a_1)$ and $(s_1+3,a_1+1)$ of Proposition~\ref{prop:levels}, and let
$\varphi=(\varphi_1,\dots,\varphi_n)\in(H^\sigma_a)^n$. With
$g^\eta_{ij}=u_{\eta,i}f'_\eta(B^\eta_{ij}*u_{\eta,j})$ as in Lemma~\ref{lem:lawid}, consider the
backward system, for $0\le s\le t$,
\begin{equation}\label{eq:backward-mild}
 \pa_s\Phi_j+x\cdot\na\Phi_j+\Theta^\eta_j\Delta\Phi_j
 +\sum_{i=1}^nB^\eta_{ij}*\big(g^\eta_{ij}\Delta\Phi_i\big)=0,
 \qquad \Phi(t)=\varphi .
\end{equation}
In the notation \eqref{eq:GN}, its left side is $\pa_s\Phi_j+\sum_i\big(\mathcal G^\eta(s)\Phi_i\big)_{ij}$.
A solution is a $\Phi\in C([0,t];(H^\sigma_a)^n)\cap L^2(0,t;(H^{\sigma+1}_a)^n)$ with
$\Phi(t)=\varphi$ that satisfies \eqref{eq:backward-mild} in $L^2(0,t;(H^{\sigma-1}_{a-1})^n)$. This is
the class $\mathcal X_{\rm ph}$ of Appendix~\ref{sec:appevol}.
\begin{enumerate}[label=\normalfont(\roman*)]
\item \eqref{eq:backward-mild} has a unique solution $\Phi$, and
\begin{equation}\label{eq:evolbound-outline2}
 \sup_{0\le s\le t}\|\Phi(s)\|_{H^{\sigma}_{a}}
 \ \le\ C^{\rm ev}_{\sigma,a}(T)\,\|\varphi\|_{H^{\sigma}_{a}} .
\end{equation}
The constant $C^{\rm ev}_{\sigma,a}(T)$ depends only on $T$, $d$, $n$, the level, $\sigma_{\min}$, $A$,
$\|u_0\|_{H^s}$, $\eps$, $\max_j\sigma_j$ and $f(0)$, and not on $N$, $\eta$ or $t$. We denote
$C_1(T):=C^{\rm ev}_{s_1,a_1}(T)$, $C_1^\sharp(T):=C^{\rm ev}_{s_1+1,a_1}(T)$ and $C_1'(T):=C^{\rm ev}_{s_1+3,a_1+1}(T)$.
In particular
\begin{equation}\label{eq:evolcont}
 \Phi\in C\big([0,t];(H^{\sigma}_{a})^n\big)
 \qquad\text{at each of the three levels.}
\end{equation}
\item If $\varphi$ lies in a finer one of the three levels, the solutions at that level and
at a coarser level coincide.
\item Denote $\Phi(s)=:S^\eta(t,s)^*\varphi$. At $(\sigma,a)=(s_1,a_1)$, $S^\eta(t,s)$ is a family of
bounded operators on $\Ufine^n$, the Banach adjoint of $\varphi\mapsto\Phi(s)$ on
$(H^{s_1}_{a_1})^n=(\Ufine^n)^*$, with
\[
 S^\eta(t,r)S^\eta(r,s)=S^\eta(t,s)\quad(s\le r\le t),\qquad
 \sup_{s\le t\le T}\|S^\eta(t,s)\|_{\mathcal L(\Ufine^n)}\le C_1(T).
\]
At $(s_1+1,a_1)$ the same holds on $(\Usharp)^n$, with $C_1^\sharp(T)$ in place of $C_1(T)$.
\item For $\varphi\in\Test^n$, $\sup_{s\le t}\|\Phi(s)\|_{\Test}\le C_1'(T)\|\varphi\|_{\Test}$.
\end{enumerate}
\end{proposition}

\begin{proof}
Part~(i) follows from the existence theory of Appendix~\ref{sec:appevol}, applied at each
of the three levels. Parts (ii)--(iv) follow from (i) and uniqueness.

\emph{(i).} At $(s_1,a_1)$ and $(s_1+1,a_1)$, Condition~\ref{ass:budget} holds for the SKT
pair by Lemma~\ref{lem:coeff}, with constants independent of $\eta$, $N$ and $t$. At
$(s_1+3,a_1+1)$, Proposition~\ref{prop:persist} at $(s_1+1,a_1)$ applies, because
Lemma~\ref{lem:coeff} gives the condition with upper level there, since $s_1+4\le s$ by
\ref{H:budget}; Proposition~\ref{prop:levels} lists these requirements. So
Proposition~\ref{prop:evolution-prop66} applies at each level. It gives a unique solution whose
conjugated variable $\zeta=w\Xi$ lies in $\mathcal X_{\rm cj}$, together with
\eqref{eq:evolbound-outline2}. By Lemma~\ref{lem:bij}, solutions in $\mathcal X_{\rm cj}$ and in
$\mathcal X_{\rm ph}$ correspond to each other, so $\Phi$ is the unique solution in
$\mathcal X_{\rm ph}$. The class $\mathcal X_{\rm ph}$ is contained in
$C([0,t];(H^\sigma_a)^n)$, which is \eqref{eq:evolcont}.

\emph{(ii).} Let $(\sigma',a')$ be the finer and $(\sigma,a)$ the coarser level, so $\sigma'\ge\sigma$
and $a'\ge a$. By Lemma~\ref{lem:spaces}(ii), $H^{m'}_{b'}\hookrightarrow H^m_b$ whenever $m'\ge m$,
$b'\ge b$. Hence the solution at $(\sigma',a')$ lies in $\mathcal X_{\rm ph}$ at $(\sigma,a)$, and it
satisfies \eqref{eq:backward-mild} in $L^2(0,t;(H^{\sigma-1}_{a-1})^n)$. By uniqueness at
$(\sigma,a)$, the two solutions coincide.

\emph{(iii).} The last sentence of Proposition~\ref{prop:evolution-prop66}, at
$(s_1,a_1)$ and at $(s_1+1,a_1)$, gives (iii).

\emph{(iv).} By (ii), for $\varphi\in\Test^n$ the solution at the coarser levels is the solution at
$(s_1+3,a_1+1)$, and (i) at that level gives the bound.
\end{proof}

\begin{lemma}[The drift operator on the test class]\label{lem:Gbound}
For $\eta\in[0,\eta_0]$ and $r\in[0,T]$, the operator $\mathcal G^\eta(r)$ of \eqref{eq:GN} maps
$\Test$ into $H^{s_1+1}_{a_1}=(\Usharp)^*$, and there is a constant $C_{\mathcal G}$, independent
of $r$ and $\eta$, with
\begin{equation}\label{eq:Gbound}
 \big\|(\mathcal G^\eta(r)\varphi)_{ij}\big\|_{H^{s_1+1}_{a_1}}\ \le\ C_{\mathcal G}\,\|\varphi\|_{\Test}
 \qquad(\varphi\in\Test,\ i,j\le n).
\end{equation}
Moreover $r\mapsto\mathcal G^\eta(r)$ is continuous in the norm of $\mathcal L(\Test,H^{s_1+1}_{a_1})$.
\end{lemma}

\begin{proof}
The bound follows term by term from \eqref{eq:GN}, with the multiplier bounds of
Lemma~\ref{lem:coeff} for the coefficients. The continuity in $r$ follows from
Lemma~\ref{lem:coeffcont}(ii).

We treat the terms one by one, on
$\varphi\in\Test=H^{s_1+3}_{a_1+1}$: $\na$ maps $\Test\to H^{s_1+2}_{a_1+1}$ and
multiplication by $x$ maps $H^m_{a_1+1}\to H^m_{a_1}$ for every $m$,
losing no derivative, since $\br x^{a_1}x_k=\frac{x_k}{\br x}\br x^{a_1+1}$ with
$x/\br x\in C^\infty_b$, and multiplication by a $C^\infty_b$ function is bounded on every
$H^m$ (Lemma~\ref{lem:spaces}(iv)). Hence the transport
term lies in $H^{s_1+2}_{a_1}$. $\Delta$ maps
$\Test\to H^{s_1+1}_{a_1+1}\hookrightarrow H^{s_1+1}_{a_1}$, and the coefficients
$\Theta^\eta_j-\kappa_j$ and $g^\eta_{ij}$ are multipliers on
$H^{s_1+1}_{a_1}$: multiplication by an unweighted function commutes with
$\br\cdot^{a_1}$, so this is the unweighted statement
$H^{s'}\cdot H^{m}\subset H^{m}$ at $s'=s_1+2>\frac d2$ and $|m|=s_1+1\le s'$.
The constant $\kappa_j=\sigma_j+nf_\eta(0)$ is subtracted first and
treated separately, as multiplication by a scalar. The nonlocal
coefficient $g^\eta_{ij}$ vanishes at infinity and needs no
subtraction. Finally, $B^\eta_{ij}*$ is bounded on $H^{s_1+1}_{a_1}$, with norm at most $A_0$:
since $a_1=0$, it is a Fourier multiplier on $H^{s_1+1}$ whose symbol is bounded by
$\|B^\eta_{ij}\|_{L^1}\le A_0$. At $\eta=0$,
$B^0_{ij}*$ is multiplication by the scalar $a_{ij}$.
Hence $\mathcal G^\eta(r):\Test\to H^{s_1+1}_{a_1}$ is bounded, uniformly in
$r$ and in $\eta\in[0,\eta_0]$: the multiplier norms are bounded by the constants of
Lemma~\ref{lem:coeff} at level $s_1+1$, the convolution norm by $A_0$, and
none of these depends on $\eta$.

For the continuity in $r$, only the coefficients depend on $r$. By
Lemma~\ref{lem:coeffcont}(ii) at $s'=s_1+2$, which is admissible since $s_1+2\le s-2$ by
\ref{H:budget}, $\Theta^\eta_j-\kappa_j$ and $g^\eta_{ij}$ are continuous in $r$ with values in
$H^{s_1+2}$. The multiplier norm on $H^{s_1+1}_{a_1}$ of a function $h$ is at most
$C\|h\|_{H^{s_1+2}}$, so the multiplications are continuous in operator norm, and so is
$\mathcal G^\eta(r)$.
\end{proof}

\begin{lemma}[Time-continuity of $\Phi$]
\label{lem:strongcont}
Let $t\in(0,T]$, $\varphi\in\Test^n$ and
$\Phi(s):=S^\eta(t,s)^*\varphi$ as in Proposition~\ref{prop:evol}. Then
$\Phi\in C([0,t];\Test^n)$, and consequently
$\Phi\in C^1([0,t];(H^{s_1+1}_{a_1})^n)$ with
$\pa_s\Phi_j=-\sum_i(\mathcal G^\eta(s)\Phi_i(s))_{ij}$ as a Bochner identity.
\end{lemma}

\begin{proof}
By Proposition~\ref{prop:evol}(ii), $\Phi$ is the solution at the level $(s_1+3,a_1+1)$. So
$\Phi\in C([0,t];\Test^n)$ by \eqref{eq:evolcont}, and \eqref{eq:backward-mild} holds in
$L^2(0,t;(H^{s_1+2}_{a_1})^n)$. That is, the distributional derivative of $\Phi$ is
$\pa_s\Phi_j=-\sum_i(\mathcal G^\eta(s)\Phi_i(s))_{ij}\in L^2(0,t;H^{s_1+2}_{a_1})$, so
$\Phi\in H^1(0,t;(H^{s_1+2}_{a_1})^n)$ and
\begin{equation}\label{eq:intbackward}
 \Phi_j(s)\ =\ \varphi_j+\int_s^t\sum_i\big(\mathcal G^\eta(r)\Phi_i(r)\big)_{ij}\dd r
 \qquad\text{as a Bochner integral in }H^{s_1+1}_{a_1}.
\end{equation}
The integrand $r\mapsto\mathcal G^\eta(r)\Phi(r)$ is continuous into $(H^{s_1+1}_{a_1})^n$, by
Lemma~\ref{lem:Gbound} and the continuity of $\Phi$ in $\Test^n$:
$\|\mathcal G^\eta(r)\Phi(r)-\mathcal G^\eta(r_0)\Phi(r_0)\|\le
C_{\mathcal G}\|\Phi(r)-\Phi(r_0)\|_{\Test}+\|\mathcal G^\eta(r)-\mathcal G^\eta(r_0)\|\,\|\Phi(r_0)\|_{\Test}$.
By the fundamental theorem of calculus for a continuous integrand, \eqref{eq:intbackward}
gives $\Phi\in C^1([0,t];(H^{s_1+1}_{a_1})^n)$ with the stated derivative.
\end{proof}

\begin{lemma}[The two evolution families have the same adjoint, hence coincide on
$\Ufine$]\label{lem:levelcoin}
Denote by $S^\eta_\Ufine$ the family
on $\Ufine$ of Proposition~\ref{prop:evol}(iii), and by $S^\eta_\sharp$ the family on
$\Usharp$ of Proposition~\ref{prop:evol}(iii). Then, for all $s\le t$:
\begin{enumerate}[label=\normalfont(\roman*)]
\item \emph{(the adjoints coincide on $\Test$)} for every $\varphi\in\Test^n$,
\begin{equation}\label{eq:adjcoin}
 S^\eta_\sharp(t,s)^*\varphi\ =\ S^\eta_\Ufine(t,s)^*\varphi
 \qquad\text{in }\Test^n,
\end{equation}
so that ``$\Phi(s)=S^\eta(t,s)^*\varphi$'' is unambiguous;
\item \emph{(hence the families agree where both act)} for every
$\nu\in\Ufine^n$, $\ S^\eta_\sharp(t,s)\nu=S^\eta_\Ufine(t,s)\nu$.
\end{enumerate}

\end{lemma}

\begin{proof}
\textup{(i)} Fix $\varphi\in\Test^n$ and let $\Phi_\Ufine,\Phi_\sharp$ be the two
adjoint solutions, at the levels $(s_1,a_1)$ and $(s_1+1,a_1)$. By Proposition~\ref{prop:evol}(ii),
$\Phi_\sharp=\Phi_\Ufine$, and both are $\Test^n$-valued by Proposition~\ref{prop:evol}(iv). Hence
\eqref{eq:adjcoin} holds.
\textup{(ii)} Let $\nu\in\Ufine^n$. Then
$S^\eta_\Ufine(t,s)\nu\in\Ufine^n\hookrightarrow(\Usharp)^n$ and
$S^\eta_\sharp(t,s)\nu\in(\Usharp)^n$, and for every $\varphi\in\Test^n$
\[
 \br{S^\eta_\Ufine(t,s)\nu,\varphi}=\br{\nu,\Phi_\Ufine(s)}
 =\br{\nu,\Phi_\sharp(s)}=\br{S^\eta_\sharp(t,s)\nu,\varphi}.
\]
$\Test\supset\mathcal S$ is dense in $(\Usharp)^*=H^{s_1+1}_{a_1}$
\textup{(}Lemma~\ref{lem:spaces}\textup{)}, and both sides are elements of $(\Usharp)^n$,
so they are equal there.
\end{proof}

\begin{lemma}[The adjoint family is Lipschitz in the terminal time, with a
two-derivative loss]\label{lem:fwdlip}
There is $C_\flat=C_\flat(T)$, independent of
$\eta\in[0,\eta_0]$, with
\begin{equation}\label{eq:fwdlip}
 \big\|S^\eta(t',t)^*\varphi-\varphi\big\|_{(H^{s_1+1}_{a_1})^n}
 \ \le\ C_\flat\,(t'-t)\,\|\varphi\|_{\Test^n}
 \qquad\text{for all }0\le t\le t'\le T,\ \varphi\in\Test^n,
\end{equation}
and consequently, for $\mu\in(\Usharp)^n$,
\begin{equation}\label{eq:fwdlipdual}
 \big\|\big(S^\eta(t',t)-I\big)\mu\big\|_{(\Test^*)^n}
 \ \le\ C_\flat\,(t'-t)\,\|\mu\|_{(\Usharp)^n},
 \qquad \Test^*=H^{-(s_1+3)}_{-(a_1+1)}.
\end{equation}

\end{lemma}

\begin{proof}
We integrate the Bochner identity of Lemma~\ref{lem:strongcont} over $[t,t']$ and bound the
integrand with Lemma~\ref{lem:Gbound} and Proposition~\ref{prop:evol}(iv). The dual bound
\eqref{eq:fwdlipdual} then follows by duality.

Let $\Phi(r):=S^\eta(t',r)^*\varphi$, the backward solution with terminal time
$t'$. Lemma~\ref{lem:strongcont} gives
$\Phi\in C^1([0,t'];(H^{s_1+1}_{a_1})^n)$ with
$\pa_r\Phi_j=-\sum_i(\mathcal G^\eta(r)\Phi_i(r))_{ij}$ as a Bochner identity, so
$\Phi(t)-\varphi=\Phi(t)-\Phi(t')=\int_t^{t'}\sum_i(\mathcal G^\eta(r)\Phi_i(r))_{ij}\dd r$
in $(H^{s_1+1}_{a_1})^n$ and
$\|\Phi(t)-\varphi\|\le(t'-t)\sup_{r\le t'}\|\mathcal G^\eta(r)\Phi(r)\|_{(H^{s_1+1}_{a_1})^n}$.
Now $\mathcal G^\eta(r):\Test\to H^{s_1+1}_{a_1}$ boundedly, uniformly in $r$ and
in $\eta\in[0,\eta_0]$, by Lemma~\ref{lem:Gbound} and $\sup_{r\le t'}\|\Phi(r)\|_{\Test^n}\le C_1'(T)\|\varphi\|_{\Test^n}$
by Proposition~\ref{prop:evol}(iv), whose constant does not depend on $\eta$ or on the terminal time $t'\le T$. These two bounds give \eqref{eq:fwdlip}.

For \eqref{eq:fwdlipdual}: $(H^{s_1+1}_{a_1})^*=\Usharp$ and $\Test^*$ is as
in \eqref{eq:fwdlipdual}, both isometrically under the $L^2$ pairing by Lemma~\ref{lem:spaces}(i), so
for $\mu\in(\Usharp)^n$ (on which the family $S^\eta_\sharp$ of Lemma~\ref{lem:levelcoin} acts, by
Proposition~\ref{prop:evol}(iii); for $\varphi\in\Test^n$, $S^\eta_\sharp(t',t)^*\varphi$ is the
$\Phi(t)$ above, by Proposition~\ref{prop:evol}(ii)) the supremum of
$|\br{\mu,(S^\eta(t',t)^*-I)\varphi}|$ over $\|\varphi\|_{\Test^n}\le1$
is the $(\Test^*)^n$ norm of $(S^\eta(t',t)-I)\mu$.
\end{proof}

\begin{corollary}[Continuity in the forward variable]\label{cor:fwdweak}
Let $\eta\in[0,\eta_0]$.
\begin{enumerate}[label=\normalfont(\roman*)]
\item for every $s$ and $\nu\in\Ufine^n$, $t\mapsto S^\eta(t,s)\nu$ is
weakly continuous from $[s,T]$ to $\Ufine^n$, uniformly in
$\eta\in[0,\eta_0]$;
\item for every $s$ and $\varphi\in\Test^n$, $t\mapsto S^\eta(t,s)^*\varphi$ is
Lipschitz from $[s,T]$ to $(H^{s_1+1}_{a_1})^n$, with constant
$C^\sharp_1(T)C_\flat\|\varphi\|_{\Test^n}$;
\item for every $s$ and $\varphi\in\Test^n$, $t\mapsto S^\eta(t,s)^*\varphi$ is
$\frac12$-H\"older from $[s,T]$ to $(H^{s_1+2}_{a_1})^n$, uniformly in $\eta\in[0,\eta_0]$.
\end{enumerate}
\end{corollary}

\begin{proof}
Part (ii) combines the evolution law with \eqref{eq:fwdlip}. Part (i) applies
\eqref{eq:fwdlipdual} and a density argument. Part (iii) interpolates (ii) against the $\Test$ bound
of Proposition~\ref{prop:evol}(iv).

(ii): by Proposition~\ref{prop:evol}(ii), for $\varphi\in\Test^n$ the adjoints of the
$\Ufine$ family and of the $\Usharp$ family agree (Lemma~\ref{lem:levelcoin}(i)), so we may use the
latter. By the evolution law $S^\eta(t',s)=S^\eta(t',t)S^\eta(t,s)$, so
$S^\eta(t',s)^*\varphi-S^\eta(t,s)^*\varphi
=S^\eta(t,s)^*\big[(S^\eta(t',t)^*-I)\varphi\big]$, and $S^\eta(t,s)^*$ is
bounded on $H^{s_1+1}_{a_1}=(\Usharp)^*$, with norm $C_1^\sharp(T)$, by
Proposition~\ref{prop:evol}(iii) and apply \eqref{eq:fwdlip}.
(i): put $v_r:=S^\eta(r,s)\nu$, so $\|v_r\|_{\Ufine^n}\le C_1(T)\|\nu\|$ and
hence $\|v_r\|_{(\Usharp)^n}\le c_\sharp\,C_1(T)\|\nu\|$. For $t<t'$,
$v_{t'}-v_t=(S^\eta(t',t)-I)v_t$ and for $t'<t$,
$v_t-v_{t'}=(S^\eta(t,t')-I)v_{t'}$. On $\Ufine^n$ the two families agree
(Lemma~\ref{lem:levelcoin}(ii)), so in either case
\eqref{eq:fwdlipdual} gives, for $\varphi\in\Test^n$,
\[
 \big|\br{v_{t'}-v_t,\varphi}\big|\ \le\ C_\flat\,c_\sharp\,C_1(T)\,|t'-t|\,
 \|\nu\|_{\Ufine^n}\,\|\varphi\|_{\Test^n}.
\]
For general $\varphi\in\Ufine^{*n}$ take
$\varphi_k\in\Test^n$ with $\varphi_k\to\varphi$ in $\Ufine^{*n}$ ($\Test\supset
C^\infty_0$ is dense there). Then
$|\br{v_{t'}-v_t,\varphi}|\le|\br{v_{t'}-v_t,\varphi_k}|+2C_1(T)\|\nu\|_{\Ufine^n}\|\varphi-\varphi_k\|_{\Ufine^{*n}}$.
The second term is small for $k$ large, uniformly in $t,t'$ and $\eta$; for this $k$, the first term is small
for $|t'-t|$ small, uniformly in $\eta$.
(iii): put $h:=S^\eta(t',s)^*\varphi-S^\eta(t,s)^*\varphi$. By (ii),
$\|h\|_{H^{s_1+1}_{a_1}}\le C^\sharp_1(T)C_\flat|t'-t|\,\|\varphi\|_{\Test^n}$, and by
Proposition~\ref{prop:evol}(iv) and Lemma~\ref{lem:spaces}(ii),
$\|h\|_{H^{s_1+3}_{a_1}}\le c\|h\|_{\Test^n}\le2cC_1'(T)\|\varphi\|_{\Test^n}$.
Since $\|v\|_{H^m_a}=\|\br\cdot^av\|_{H^m}$, H\"older on the Fourier side gives
$\|h\|_{H^{s_1+2}_{a_1}}\le\|h\|^{1/2}_{H^{s_1+1}_{a_1}}\|h\|^{1/2}_{H^{s_1+3}_{a_1}}$.

\end{proof}

\begin{proposition}[The forward variable in the $\Hmid$ norm]\label{prop:fwdHmid}
Fix $\eta\in[0,\eta_0]$. For every $s\in[0,T]$ and every
$\nu\in\Ufine^n$ the map $[s,T]\ni t\mapsto S^\eta(t,s)\nu\in\Hmid^n$ is
continuous in norm.
\end{proposition}

\begin{proof}
Let $t_k\to t$ in $[s,T]$ and put $v_r:=S^\eta(r,s)\nu$. By
Corollary~\ref{cor:fwdweak}\textup{(i)}, $v_{t_k}\rightharpoonup v_t$ weakly in
$\Ufine^n$, and $\|v_r\|_{\Ufine^n}\le C_1(T)\|\nu\|_{\Ufine^n}$ by
Proposition~\ref{prop:evol}. The embedding
$\iota:\Ufine^n\hookrightarrow\Hmid^n$ is compact (Lemma~\ref{lem:spaces}(ii),
both $s_m>s_1$ and $a_m>a_1$ being strict), so $\{\iota v_{t_k}\}$ is
relatively norm-compact, while $\iota v_{t_k}\rightharpoonup\iota v_t$ because
$\iota$ is bounded. Every subsequence therefore has a norm-convergent
sub-subsequence, whose limit must be $\iota v_t$. Hence
$\iota v_{t_k}\to\iota v_t$ in $\Hmid^n$.
\end{proof}

\section{Uniform estimates and the remainders}\label{sec:estimates}

\noindent
Three remainders appear in \eqref{eq:fluc_decomp}: a self-interaction term, a
bilinear term quadratic in the fluctuation, and a second-order Taylor term. Each
must vanish in the limit, and the bilinear one is the difficult case: it is
quadratic in an object that is only bounded, so it needs a fourth moment. That
moment, \eqref{eq:K4-mild}, is the hypothesis of Proposition~\ref{prop:close}. It is
proved in Appendix~\ref{sec:appk4}, as is Lemma~\ref{lem:CH2}, which the proposition also uses.
Lemma~\ref{lem:remainders} estimates the three remainders. Its condition (a) and the bilinear bound \eqref{eq:cnought}
are pathwise and assume no field bound. Lemma~\ref{lem:Rdet} gives, at fixed $N$ and $\eta$ and with
no field bound, the deterministic regularity of $\mathcal R^N$ that the
limit passage needs. Corollary~\ref{cor:remvanish} states that the remainders
vanish in $L^2(\Omega)$, and Corollary~\ref{cor:remscalar} that their pairings with a fixed test function vanish in $L^1$.

Two other things are settled here because the remainder estimates need them.
The drift coefficient converges: Theorem~\ref{thm:G} identifies the limit of
$\mathcal G^\eta$ as the linearised operator, at rate $\eta$ along the full sequence,
and Lemma~\ref{lem:Gupgrade} upgrades the topology. The finite-$N$ mild formula
\eqref{eq:mildzero} is Proposition~\ref{prop:mildN}. its proof and that of Theorem~\ref{thm:weakmild} share the telescoping identity of Lemma~\ref{lem:telescope}. Lemma~\ref{lem:uetacont} gives the continuity and size of the field at fixed
$N$, and Lemma~\ref{lem:crossspecies} the independence across species of the
intermediate ensemble.

For each $t\in[0,T]$, Proposition~\ref{prop:mildN} gives, almost surely, the mild
form of \eqref{eq:fluc_decomp}:
\begin{equation}\label{eq:mildzero}
 \fluc^N(t)=\underbrace{S^\eta(t,0)\fluc^N(0)}_{=:P^N(t)}
 +\underbrace{\int_0^tS^\eta(t,s)\dd M^N(s)}_{=:Z^N(t)}
 +\underbrace{\int_0^tS^\eta(t,s)\mathcal R^N(s)\dd s}_{=:D^N(t)}.
\end{equation}
The stochastic convolution $Z^N(t)$ is the integral
\[
 Z^N(t)=\int_0^tS^\eta(t,r)\Upsilon^N(r)\dd W^N(r)
\]
of Lemma~\ref{lem:BDG} with $\mathsf J(r)=S^\eta(t,r)$, where $\Upsilon^N$ and $W^N$ are as in
Lemma~\ref{lem:MNval}. Proposition~\ref{prop:mildN} shows that this integral is well defined. The term $D^N(t)\in(\Test^*)^n$ is defined by
$\br{D^N(t),\varphi}:=\int_0^t\br{\mathcal R^N(s),S^\eta(t,s)^*\varphi}\dd s$ for $\varphi\in\Test^n$.
The integral exists by Lemma~\ref{lem:Rdet} and Proposition~\ref{prop:evol}(iv).

\begin{lemma}[Lipschitz bound, valid everywhere]\label{lem:lip}
Under \ref{H:A4}, let $i,j\le n$, let $w^\eta_{ij}$ be as in Section~\ref{sec:setting},
and put $h_{ij}:=N^{-1/2}B^\eta_{ij}*\fluc^N_j$. Then the Taylor remainder of
\eqref{eq:Taylor-outline2} satisfies
\begin{equation}\label{eq:Taylor-forced}
 \mathsf T^N_{ij}=\sqrt N\Big[f_\eta(w^\eta_{ij}+h_{ij})-f_\eta(w^\eta_{ij})
 -f'_\eta(w^\eta_{ij})h_{ij}\Big]
\end{equation}
and
\begin{equation}\label{eq:lip}
 \big|\mathsf T^N_{ij}(t)\big|\ \le\ 2\,\eta^{-\alpha}\,
 \big|B^\eta_{ij}*\fluc^N_j(t)\big|\qquad\text{pointwise, for every }\omega,t.
\end{equation}
\end{lemma}

\begin{proof}
Since $\fluc^N_j=\sqrt N(\mu^N_j-u_{\eta,j})$, we have $B^\eta_{ij}*\mu^N_j=w^\eta_{ij}+h_{ij}$, and
\eqref{eq:Taylor-outline2} becomes \eqref{eq:Taylor-forced}.
$\mathrm{Lip}(f_\eta)=\|f'_\eta\|_{L^\infty}\le\eta^{-\alpha}$, so
$|f_\eta(a+h)-f_\eta(a)|\le\eta^{-\alpha}|h|$ and
$|f'_\eta(a)h|\le\eta^{-\alpha}|h|$. Insert in \eqref{eq:Taylor-forced} and use
$\sqrt N|h_{ij}|=|B^\eta_{ij}*\fluc^N_j|$.
\end{proof}

\begin{lemma}[The pathwise pointwise bound]\label{lem:ptwise}
Assume $\mathsf b_0>0$. Then, with
$\mathsf C_\natural:=\max\{\tfrac\eps2,\tfrac2{\mathsf b_0}\}$,
\begin{equation}\label{eq:ptwise}
 \big|\mathsf T^N_{ij}(t,x)\big|\ \le\ \mathsf C_\natural\,
 N^{-1/2}\,\eta^{-\alpha}\,\big|B^\eta_{ij}*\fluc^N_j(t,x)\big|^2
 \qquad\text{for every }\omega,\ t\le T,\ x\in\R^d.
\end{equation}
\end{lemma}

If $\mathsf b_0=0$, then $w^\eta_{ij}\equiv0$, $B^\eta_{ij}*\fluc^N_j\equiv0$ and $\mathsf T^N_{ij}\equiv0$.

\begin{proof}
We split on the size of the increment $h_{ij}(x)$ of \eqref{eq:Taylor-forced}. If it is at most
$\mathsf b_0$, the Taylor formula with Lagrange remainder applies on an interval where $f_\eta=f$
and $|f''|\le\eps$. Otherwise the Lipschitz bound \eqref{eq:lip} is turned into a quadratic bound.

Fix $\omega,t,x$ and put $v:=|B^\eta_{ij}*\fluc^N_j(t,x)|$, so that the
increment in \eqref{eq:Taylor-forced} is $|h_{ij}(x)|=N^{-1/2}v$. Split on the
pointwise set $\{v\le\mathsf b_0\sqrt N\}$.

If $v\le\mathsf b_0\sqrt N$, then $|h_{ij}(x)|\le\mathsf b_0$, and
$|w^\eta_{ij}(x)|\le\mathsf b_0$ at every $x$, by Lemma~\ref{lem:lawid}(i). So the
segment joining $w^\eta_{ij}(x)$ to $w^\eta_{ij}(x)+h_{ij}(x)$ lies in
$[-2\mathsf b_0,2\mathsf b_0]=I\subset[-a_\eta,a_\eta]$, by \ref{H:cutoff}.
There $f_\eta=f$, and so $f''_\eta=f''$ on $I$, since both are continuous
($f_\eta\in W^{s+1,\infty}$ with $s\ge2$). By \ref{H:small},
$|f''|\le\eps$ on $I$. Taylor with Lagrange remainder gives
$|\mathsf T^N_{ij}(x)|\le\frac{\eps}{2}N^{-1/2}v^2$, which is at most the right
side of \eqref{eq:ptwise} since $\eta^{-\alpha}\ge1$.

\emph{If $v>\mathsf b_0\sqrt N$}, Lemma~\ref{lem:lip} gives
$|\mathsf T^N_{ij}(x)|\le2\eta^{-\alpha}v$, and dividing and multiplying by $v$,
\[
 2\eta^{-\alpha}v\ =\ 2\eta^{-\alpha}\frac{v^2}{v}
 \ <\ \frac{2}{\mathsf b_0}\,N^{-1/2}\,\eta^{-\alpha}\,v^2,
\]
the inequality being exactly the case hypothesis $v>\mathsf b_0\sqrt N$.
Since $\mathsf C_\natural$ is the larger of the two constants, \eqref{eq:ptwise}
holds in both cases.
\end{proof}

\begin{lemma}[The smoothed field and the Taylor remainder]\label{lem:taylor}
Let $i,j\le n$, let $\mathcal R^{N,\rm Tay}$ be the Taylor part of \eqref{eq:Rquad}, and denote
\begin{equation}\label{eq:K4def}
 \overline K_4:=\sup_{t\le T}\max_{i,j}\E\big\|B^\eta_{ij}*\fluc^N_j(t)\big\|^4_{H^{s_1}}.
\end{equation}
Then
\eqref{eq:cheap1} and \eqref{eq:Tayl2} hold for every $t\le T$, and \eqref{eq:cheap2}
holds for every $t\le T$ and every $N\ge N_0$, with $N_0$ as in \eqref{eq:N0}:
\begin{align}
 \|B^\eta_{ij}*\fluc^N_j\|_{H^\sigma}&\le2\sqrt N\,\|B^\eta_{ij}\|_{H^{\sigma}}\le2\|B_{ij}\|_{H^\sigma}\sqrt N\eta^{-(d/2+\sigma)}
 \qquad\text{a.s., for every }\sigma\ge0,\label{eq:cheap1}\\
 \E\|B^\eta_{ij}*\fluc^N_j\|^2_{H^{s_1}}&\le C\,N^{C_2(T+1)\delta}\,\eta^{-(2d+2s_1+2+2\alpha)}+C\eta^{-(d+2s_1)}
 \ =\ N^{C_2(T+1)\delta}\times\text{polylog}(N),
 \label{eq:cheap2}\\
 \E\int_0^T\big\|\mathcal R^{N,\rm Tay}(s)\big\|^2_{H^{-s_1-2}_{-a_1}}\dd s
 &\le C\,T\,\mathsf C_\natural^2\,N^{-1}\,\eta^{-2\alpha}\,\overline K_4,\label{eq:Tayl2}
\end{align}
with $\mathsf C_\natural$ the constant of Lemma~\ref{lem:ptwise}.
\end{lemma}

\begin{proof}
The three bounds are proved in turn. \eqref{eq:cheap1} follows from the triangle inequality and the
scaling of $B^\eta_{ij}$. For \eqref{eq:cheap2} we compare the particles with the intermediate ensemble.
\eqref{eq:Tayl2} follows from the pointwise bound \eqref{eq:ptwise}.

\eqref{eq:cheap1}. $B^\eta_{ij}*\mu^N_j=\frac1N\sum_\ell B^\eta_{ij}(\cdot-X_{\ell,j})$,
and translation is an isometry of $H^\sigma$. By Lemma~\ref{lem:lawid},
$u_{\eta,j}(t)$ is a probability density, so $|\widehat{u_{\eta,j}(t)}|\le1$ and
$\|B^\eta_{ij}*u_{\eta,j}(t)\|_{H^\sigma}\le\|B^\eta_{ij}\|_{H^\sigma}$. The second
inequality is Lemma~\ref{lem:scaling-fourth}.

\eqref{eq:cheap2}. Let $\bar\fluc^N_j:=\sqrt N(\bar\mu^N_j-u_{\eta,j})$ and
$Y_{\ell,j}:=X^{N,\eta}_{\ell,j}-\bar X^\eta_{\ell,j}$. Then
$B^\eta_{ij}*\fluc^N_j=B^\eta_{ij}*\bar\fluc^N_j+N^{-1/2}\sum_\ell\big(B^\eta_{ij}(\cdot-X^{N,\eta}_{\ell,j})
-B^\eta_{ij}(\cdot-\bar X^\eta_{\ell,j})\big)$. By the intermediate-system condition of
Lemma~\ref{lem:smoothed}, $\E\|B^\eta_{ij}*\bar\fluc^N_j\|^2_{H^{s_1}}\le C\eta^{-(d+2s_1)}$.
For the sum, $\|B^\eta_{ij}(\cdot-x)-B^\eta_{ij}(\cdot-y)\|_{H^{s_1}}\le|x-y|\,\|\na B^\eta_{ij}\|_{H^{s_1}}
\le C|x-y|\,\eta^{-(d/2+s_1+1)}$, by the fundamental theorem of calculus and
Lemma~\ref{lem:scaling-fourth} applied to $\na B_{ij}$. Jensen's inequality in
$\ell$ and \eqref{eq:phi2} give
\[
 \E\Big\|N^{-1/2}\sum_\ell\big(B^\eta_{ij}(\cdot-X^{N,\eta}_{\ell,j})-B^\eta_{ij}(\cdot-\bar X^\eta_{\ell,j})\big)\Big\|^2_{H^{s_1}}
 \le C\,N\,\eta^{-(d+2s_1+2)}\max_\ell\E|Y_{\ell,j}|^2
 \le C\,N^{C_2(T+1)\delta}\,\eta^{-(2d+2s_1+2+2\alpha)}.
\]
The last form of \eqref{eq:cheap2} is Lemma~\ref{lem:polylog}.

For \eqref{eq:Tayl2}, let $\varphi\in H^{s_1+2}_{a_1}$ with norm at most $1$. By
Lemma~\ref{lem:spaces}(iii), $\|\Delta\varphi\|_{L^\infty}\le C$, and since
$\mu^N_i$ is a probability measure,
$|\mathcal R^{N,\rm Tay}_i(\varphi)|\le C\sum_j\|\mathsf T^N_{ij}\|_{L^\infty}$. By
\eqref{eq:ptwise} and $H^{s_1}\hookrightarrow L^\infty$,
$\|\mathsf T^N_{ij}\|_{L^\infty}\le C\,\mathsf C_\natural N^{-1/2}\eta^{-\alpha}
\|B^\eta_{ij}*\fluc^N_j\|^2_{H^{s_1}}$. Taking the supremum over $\varphi$
(Lemma~\ref{lem:spaces}(i)), squaring, taking expectations and integrating over
$(0,T)$ gives \eqref{eq:Tayl2}.
\end{proof}

\begin{lemma}[Two elementary facts about $\fluc^N$ at fixed $N$]
\label{lem:uetacont}
Let $\eta\in(0,\eta_0]$.
\begin{enumerate}[label=\normalfont(\roman*)]
\item Deterministically, for every $t\le T$,
\begin{equation}\label{eq:flucNcrude}
 \|\fluc^N(t)\|_{\Ufine}\ \le\ \sqrt N\,\big(C(d,s_1,a_1)+\|u_0\|_{H^s}\big),
\end{equation}
and a fortiori the same bound in $\Usharp$.
\item $\Prob$-a.s., $\fluc^N\in C\big([0,T];\Ufine^n\big)$; a fortiori
$\fluc^N$ is continuous in $\Usharp$.
\end{enumerate}
\end{lemma}

\begin{proof}
By \eqref{eq:choice}, $a_1=0$ and $s_1>\frac d2+1$, and $\|v\|_{H^{-s_1}_{-a_1}}
=\|\br\cdot^{-a_1}v\|_{H^{-s_1}}\le\|\br\cdot^{-a_1}v\|_{L^2}\le\|v\|_{L^2}$,
so that $L^2\hookrightarrow\Ufine$ with norm $\le1$.

\emph{(i).} $\|\delta_x\|_{H^{-s_1}_{-a_1}}=\sup\{|\varphi(x)|:
\|\varphi\|_{H^{s_1}_{a_1}}\le1\}\le C\br x^{-a_1}\le C$ by
Lemma~\ref{lem:spaces}(iii) at $m=0$. Averaging over $k$ to get
$\|\mu^N_i(t)\|_{\Ufine}\le C$. For the other half,
$\|u_{\eta,i}(t)\|_{\Ufine}\le\|u_{\eta,i}(t)\|_{L^2}\le\|u_0\|_{H^s}$ at every
$t$, by Lemma~\ref{lem:lawid}(i).

\emph{(ii).} The empirical half is $\|\delta_x-\delta_y\|_{\Ufine}\le C|x-y|$ together with the continuity of the
particle paths. The other half is the continuity $u_{\eta,i}\in C([0,T];L^2)$ of
Lemma~\ref{lem:coeffcont}(i) at $s'=0$, with the $\eta$-uniform H\"older modulus
\eqref{eq:coeffholder} at $\vartheta=\frac1{s+1}$, composed with the embedding
$L^2\hookrightarrow\Ufine$.
\end{proof}

\begin{theorem}[The OU drift]\label{thm:G}
Define, for $i,j=1,\dots,n$,
\begin{equation}\label{eq:Amatrix}
 \mathbf A_{ij}(u):=\delta_{ij}\Big(\sigma_i+\sum_{k=1}^nf(a_{ik}u_k)\Big)
 +a_{ij}\,u_i\,f'(a_{ij}u_j)
 \;=\;\delta_{ij}\Theta_i+a_{ij}u_if'(a_{ij}u_j),
\end{equation}
and let $\mathcal G=\mathcal G(t)$ act on $\varphi\in\Test$ by
\begin{equation}\label{eq:Gweak}
 (\mathcal G\varphi)_{ij}:=\delta_{ij}\big(-\na U_i\cdot\na\varphi\big)
 +\mathbf A_{ij}(u)\,\Delta\varphi.
\end{equation}
Then:
\begin{enumerate}[label=\normalfont(\roman*)]
\item For every $\varphi\in C^\infty_0(\R^d)$ there is $C=C(\varphi,T)$ with
\begin{equation}\label{eq:Grate}
 \sup_{t\le T}\big\|(\mathcal G^\eta\varphi)_{ij}(t)-(\mathcal G\varphi)_{ij}(t)\big\|_{L^2(\R^d)}
 \ \le\ C\,\eta\qquad\text{for all }i,j\text{ and all }N,
\end{equation}
so $\mathcal G^\eta\varphi\to\mathcal G\varphi$ in $L^\infty(0,T;L^2(\R^d))$ along
the full sequence $\eta=\eta(N)$ of \ref{H:scaling}.

\item In the sense of distributions, the associated drift on the fluctuation
field is the linearisation of the SKT operator: for
$\xi=(\xi_1,\dots,\xi_n)$ and every $i$,
\begin{equation}\label{eq:Gstrong}
 \sum_j\br{\xi_j,(\mathcal G\varphi)_{ij}}
 =\br{\mathcal G_i\xi,\varphi},
 \qquad
 \mathcal G_i\xi=\diver\big(\xi_i\na U_i\big)+\Delta\Big(\sum_{j=1}^n\mathbf A_{ij}(u)\,\xi_j\Big),
\end{equation}
and $\mathcal G_i\xi=\frac{\dd}{\dd\tau}\big|_{\tau=0}\mathcal L_i(u+\tau\xi)$,
where $\mathcal L_i(v)=\diver(v_i\na U_i)+\Delta\big(\sigma_iv_i+v_i\sum_jf(a_{ij}v_j)\big)$
is the SKT operator of \eqref{eq:skt_det}.
\end{enumerate}
\end{theorem}

\begin{proof}
Part (ii) is a direct computation. For (i) we split $\mathcal G^\eta\varphi-\mathcal G\varphi$ into
three terms and bound each in $L^2$ by $C\eta$. The tools are that all arguments of $f$ lie in the
interval $I$, where $f_\eta=f$, and the mollification estimate \eqref{eq:mollerr}.

For (ii), $\br{\xi_j,\mathbf A_{ij}(u)\Delta\varphi}
=\br{\Delta(\mathbf A_{ij}(u)\xi_j),\varphi}$ and
$\br{\xi_i,-\na U_i\cdot\na\varphi}=\br{\diver(\xi_i\na U_i),\varphi}$, both in
$\mathcal D'$. For the linearisation statement, differentiate
$\mathcal L_i(v)$ at $v=u$ in the direction $\xi$: the term $\diver(v_i\na U_i)$
contributes $\diver(\xi_i\na U_i)$. The term $\Delta(\sigma_iv_i)$ contributes
$\Delta(\sigma_i\xi_i)$, $\Delta(v_i\sum_jf(a_{ij}v_j))$ contributes
\[
 \Delta\Big(\xi_i\sum_jf(a_{ij}u_j)\Big)
 +\Delta\Big(u_i\sum_ja_{ij}f'(a_{ij}u_j)\xi_j\Big),
\]
and collecting the coefficient of $\xi_j$ gives exactly $\mathbf A_{ij}(u)$ of
\eqref{eq:Amatrix}.

For (i), recall $w^\eta_{ij}=B^\eta_{ij}*u_{\eta,j}$, and put
$g_\eta:=u_{\eta,i}f'_\eta(w^\eta_{ij})\Delta\varphi$ and
$g:=u_if'(a_{ij}u_j)\Delta\varphi$. The first-order terms of \eqref{eq:GN} and
\eqref{eq:Gweak} coincide, so
\[
 (\mathcal G^\eta\varphi)_{ij}-(\mathcal G\varphi)_{ij}
 =\delta_{ij}\big(\Theta^\eta_i-\Theta_i\big)\Delta\varphi
 +B^\eta_{ij}*(g_\eta-g)+\big(B^\eta_{ij}*g-a_{ij}g\big).
\]

\emph{Step 1: the arguments of $f$ lie in $I$.} By \ref{H:small}, Lemma~\ref{lem:lawid}(i),
and $H^s\hookrightarrow L^\infty$, for every $t\le T$, $\|u_{\eta,j}(t)\|_{L^\infty}$ and
$\|u_j(t)\|_{L^\infty}$ are at most $C_s\|u_0\|_{H^s}$. Since
$\|B^\eta_{ij}\|_{L^1}\le A_0$ and $|a_{ij}|\le A_0$, both $w^\eta_{ij}$ and
$a_{ij}u_j$ take values in $[-\mathsf b_0,\mathsf b_0]\subset I$, and so does the
segment between them. On $I$ we have $f_\eta=f$ by \ref{H:cutoff}, and
$\|f\|_{C^2(I)}\le\eps$ by \ref{H:small}. Hence, pointwise,
\[
 \big|f_\eta(w^\eta_{ij})-f(a_{ij}u_j)\big|\le\eps\,\big|w^\eta_{ij}-a_{ij}u_j\big|,
 \qquad
 \big|f'_\eta(w^\eta_{ij})-f'(a_{ij}u_j)\big|\le\eps\,\big|w^\eta_{ij}-a_{ij}u_j\big|.
\]

\emph{Step 2: a mollification estimate.} Since $\operatorname{supp}B^\eta_{ij}\subset
B_\eta(0)$ and $\int B^\eta_{ij}=a_{ij}$, every $v\in H^1$ satisfies
$B^\eta_{ij}*v-a_{ij}v=\int B^\eta_{ij}(y)\big(v(\cdot-y)-v\big)\dy$, hence
\begin{equation}\label{eq:mollerr}
 \big\|B^\eta_{ij}*v-a_{ij}v\big\|_{L^2}
 \ \le\ \int|B^\eta_{ij}(y)|\,|y|\dy\;\|\na v\|_{L^2}
 \ \le\ A_0\,\eta\,\|\na v\|_{L^2}.
\end{equation}
With $w^\eta_{ij}-a_{ij}u_j=B^\eta_{ij}*(u_{\eta,j}-u_j)+(B^\eta_{ij}*u_j-a_{ij}u_j)$,
\eqref{eq:CHbias}, \eqref{eq:mollerr} and $\|\na u_j(t)\|_{L^2}\le\|u_0\|_{H^s}$
give $\sup_{t\le T}\|w^\eta_{ij}-a_{ij}u_j\|_{L^2}\le A_0\big(C(T)+\|u_0\|_{H^s}\big)\eta$.

\emph{Step 3: the three terms.} By Steps 1 and 2 the diagonal term has
$L^2$ norm at most $n\eps\|\Delta\varphi\|_{L^\infty}\max_k\|w^\eta_{ik}-a_{ik}u_k\|_{L^2}\le C\eta$.
Next, $g_\eta-g=(u_{\eta,i}-u_i)f'_\eta(w^\eta_{ij})\Delta\varphi
+u_i\big(f'_\eta(w^\eta_{ij})-f'(a_{ij}u_j)\big)\Delta\varphi$, so Step~1 and
\eqref{eq:CHbias} give $\|B^\eta_{ij}*(g_\eta-g)\|_{L^2}\le A_0\|g_\eta-g\|_{L^2}\le C\eta$.
Finally \eqref{eq:mollerr} with $v=g$ bounds the last term by $A_0\eta\|\na g\|_{L^2}$, and
$\na g=(\na u_i)f'(a_{ij}u_j)\Delta\varphi+u_ia_{ij}f''(a_{ij}u_j)(\na u_j)\Delta\varphi
+u_if'(a_{ij}u_j)\na\Delta\varphi$ has $L^2$ norm at most
$C(\|u_0\|_{H^s},\eps,\|\varphi\|_{C^3})$, uniformly in $t$. This is \eqref{eq:Grate}.
\end{proof}

\begin{remark}[The linear nonlinearity]\label{rem:linearf}
For $f(s)=cs$, \eqref{eq:Amatrix} becomes
$\mathbf A_{ij}(u)=\delta_{ij}\big(\sigma_i+c\sum_ka_{ik}u_k\big)+c\,a_{ij}u_i$, the diffusion matrix of the
classical SKT system with self-diffusion $\sigma_i$ and cross-diffusion coefficients $ca_{ij}$.
For $c\ne0$ this $f$ is not covered by Assumption~\ref{ass:main}, since \ref{H:A4} requires $f\ge0$
on all of $\R$. More generally, if $f\ge0$ on $\R$, $f(0)=0$ and $f$ is differentiable at $0$,
then $0$ is a minimum of $f$, so $f'(0)=0$. Assumption (A4) of
\cite{chen_holzinger_2021_rigorous_SKT} contains the same restriction. Whether nonnegativity of $f$
on $[0,\infty)$, together with $B_{ij}\ge0$, suffices for the results of this paper is not
addressed here.
\end{remark}

\begin{lemma}[Independence of the intermediate ensemble]\label{lem:crossspecies}
Assume \ref{H:CH}, with \eqref{eq:datax} possibly omitted, and let $u_\eta$ satisfy
\begin{equation}\label{eq:uetareg}
 u_{\eta,j}\in L^\infty(0,T;H^{s}(\R^d)),\quad s>\tfrac d2+1,\qquad
 \|u_\eta\|_{L^\infty(0,T;H^{s})}\le\|u_0\|_{H^{s}},
\end{equation}
as in \cite[Thm.~2]{chen_holzinger_2021_rigorous_SKT}. Then:
\begin{enumerate}[label=\normalfont(\roman*)]
\item under \eqref{eq:dataw} alone: for each species $j$, the processes
 $\bar X^\eta_{\ell,j}$, $\ell\le N$, are i.i.d.;
\item under \eqref{eq:dataw} and \eqref{eq:datax}: the processes
 $\bar X^\eta_{k,i}$, $(k,i)$ ranging over all pairs, are independent.
\end{enumerate}
\end{lemma}

\begin{proof}
We show that $\bar X^\eta_{k,i}$ is a fixed measurable function of its own data
$(\xi^k_i,W^k_i)$. Step~0 checks the hypotheses of the It\^o existence and uniqueness theorem, and
Step~1 deduces the measurable dependence. Independence of the data then gives (i) and (ii).

The intermediate system \eqref{eq:intermediate} is decoupled, and its
coefficients are deterministic: $\bar X^\eta_{k,i}$ solves
$\dd\bar X=-\na U_i(\bar X)\dd t+(2\Theta^\eta_i(t,\bar X))^{1/2}\dd W^k_i$ with
$\Theta^\eta_i(t,x)=\sigma_i+\sum_jf_\eta\big((B^\eta_{ij}*u_{\eta,j})(t,x)\big)$,
a function of $(t,x)$ alone. The other species enter only through the
deterministic densities $u_{\eta,j}$.
 \emph{Step 0: the coefficients.} The It\^o theory used below requires three properties of the
pair $\big(-\na U_i,(2\Theta^\eta_i)^{1/2}\big)$: Borel
measurability in $(t,x)$, a Lipschitz constant in $x$ uniform in $t$, and linear
growth. Linear growth holds, since $\na U_i$ is linear with $\na U_i(0)=0$
and $\Theta^\eta_i$ is bounded. For the measurability: $H^{s}$ is
separable, so \eqref{eq:uetareg} gives a jointly measurable representative of
$(t,y)\mapsto u_{\eta,j}(t,y)$. Then $(t,x,y)\mapsto B^\eta_{ij}(x-y)u_{\eta,j}(t,y)$
is jointly measurable and
$\int_0^T\!\!\int|B^\eta_{ij}(x-y)u_{\eta,j}(t,y)|\dd y\dd t
\le T\|B^\eta_{ij}\|_{L^2}\|u_0\|_{H^{s}}<\infty$, so the partial integral
$(t,x)\mapsto(B^\eta_{ij}*u_{\eta,j})(t,x)$ is jointly measurable by Tonelli, and
$\Theta^\eta_i$ is too, $f_\eta$ being continuous. For the Lipschitz bound, for a.e.\ $t$,
\[
 \|\na\big(B^\eta_{ij}*u_{\eta,j}(t)\big)\|_{L^\infty}
 \le\|\na B^\eta_{ij}\|_{L^2}\|u_{\eta,j}(t)\|_{L^2}
 =\eta^{-\frac d2-1}\|\na B_{ij}\|_{L^2}\|u_{\eta,j}(t)\|_{L^2}
 \le\eta^{-\frac d2-1}\|\na B_{ij}\|_{L^2}\|u_0\|_{H^{s}},
\]
by Cauchy--Schwarz and \eqref{eq:uetareg}, a bound independent of $t$. Modifying $\Theta^\eta_i$ on the exceptional null set of times does not changes a
solution. Composing, $\mathrm{Lip}(f_\eta)\le\eta^{-\alpha}$ by \ref{H:CH} and
$\tau\mapsto(2\sigma_i+2\tau)^{1/2}$ is Lipschitz on $[0,\infty)$ with constant
$(2\sigma_i)^{-1/2}$ because $f_\eta\ge0$ and $\sigma_i>0$, so that
$\mathrm{Lip}\big((2\Theta^\eta_i(t,\cdot))^{1/2}\big)
\le n(2\sigma_i)^{-1/2}\eta^{-\alpha}
\max_j\mathrm{Lip}(B^\eta_{ij}*u_{\eta,j}(t))$. By \ref{H:CH}, $\Theta^\eta_i\le\sigma_i+n\|f_\eta\|_{L^\infty(\R)}<\infty$.
And $\na U_i$ is linear by \ref{H:CH}.
\emph{Step 1: measurable dependence on the data.} The classical It\^o existence and uniqueness
theorem for coefficients Borel in $(t,x)$, Lipschitz in $x$ uniformly in $t$, and of
linear growth therefore applies, with the square-integrable datum $\xi^k_i$ of \ref{H:CH}.
Its Picard iterates are $\sigma(\xi^k_i,W^k_i)$-measurable.
Hence $\bar X^\eta_{k,i}$ agrees a.s.\ with $F_i(\xi^k_i,W^k_i)$ for a measurable
$F_i$ not depending on $k$, because the coefficients do not depend on $k$.
The clause $\xi^k_i\perp W^k_i$ of \eqref{eq:dataw} makes $\xi^k_i$ an admissible
$\mathcal F_0$-measurable datum for the Brownian motion $W^k_i$. Now \eqref{eq:dataw} makes the pairs
$(\xi^\ell_j,W^\ell_j)$, $\ell\le N$, i.i.d.\ at fixed $j$, which gives
\textup{(i)}; adding \eqref{eq:datax} makes them independent over all $(k,i)$,
which gives \textup{(ii)}. Independent inputs give independent outputs.
The uniqueness input is Lemma~\ref{H:W}. 
\end{proof}

\begin{lemma}[Topology upgrade for Theorem~\ref{thm:G}(i)]\label{lem:Gupgrade}
Let $\sigma:=s_1+1$, and let $N_0$ be as in \eqref{eq:N0}. Fix $\varphi\in C^\infty_0(\R^d)$ and put
$Q^N_{ij}:=(\mathcal G^\eta\varphi)_{ij}-(\mathcal G\varphi)_{ij}$. Then
\begin{equation}\label{eq:Gupgrade}
 \big\|Q^N_{ij}\big\|_{L^2(0,T;H^{s_1}_{a_1}(\R^d))}\ \longrightarrow\ 0
 \qquad(N\to\infty),
\end{equation}
which is the topology Lemma~\ref{lem:drift} requires.
\end{lemma}

\begin{proof}
We interpolate between the $L^2$ convergence of Theorem~\ref{thm:G}(i) and a bound in
$H^{s_1+1}$ that is uniform in $N$. Step~1 localises $Q^N_{ij}$ to a fixed compact set, so that the
weight plays no role. Step~2 gives the uniform bound, and Step~3 interpolates and applies H\"older's
inequality in time.

Let $K:=\operatorname{supp}\varphi$ and $K_1:=K+\overline{B_1(0)}$, and fix
$R$ with $K_1\subset B_R$. Throughout $\theta:=s_1/\sigma\in(0,1)$.

\emph{Step 1 (localisation: the transport term cancels, and everything lives in
a fixed compact).} Subtracting \eqref{eq:Gweak} from \eqref{eq:GN}, the
first-order terms $\delta_{ij}(-\na U_i\cdot\na\varphi)$ are identical
and cancel, leaving
\begin{equation}\label{eq:Ddiff}
 Q^N_{ij}=\delta_{ij}\big(\Theta^\eta_i-\Theta_i\big)\Delta\varphi
 \;+\;\Big[B^\eta_{ij}*\big(u_{\eta,i}f'_\eta(w^\eta_{ij})\Delta\varphi\big)
 -a_{ij}u_if'(a_{ij}u_j)\Delta\varphi\Big].
\end{equation}
Every term of \eqref{eq:Ddiff} carries a factor $\Delta\varphi$, supported in
$K$. The outer convolution spreads support by at most
$\operatorname{supp}B^\eta_{ij}\subset B_\eta$, so for $\eta\le1$
\[
 \operatorname{supp}Q^N_{ij}(t,\cdot)\ \subset\ K_1\qquad\text{for all }t\le T.
\]
Two consequences. \emph{(i) The weight is irrelevant:} choose
$\chi\in C^\infty_0(\R^d)$ with $\chi\equiv1$ on $K_1$. Then
$\br\cdot^{a_1}\chi\in C^\infty_0$, multiplication by it is bounded on
$H^{s_1}$, and for $v$ supported in $K_1$,
$\|v\|_{H^{s_1}_{a_1}}=\|\br\cdot^{a_1}\chi\,v\|_{H^{s_1}}\le C_{\rm mult}\|v\|_{H^{s_1}}$.
\emph{(ii)} Theorem~\ref{thm:G}(i) gives
\begin{equation}\label{eq:GupL2}
 \big\|Q^N_{ij}\big\|_{L^2(0,T;L^2(\R^d))}\ \longrightarrow\ 0.
\end{equation}

\emph{Step 2 (a uniform bound at a higher index).} We claim
\begin{equation}\label{eq:GupHs}
 M_Q:=\sup_{N\ge N_0}\ \operatorname*{ess\,sup}_{t\le T}\
 \big\|Q^N_{ij}(t)\big\|_{H^{\sigma}(\R^d)}\ <\ \infty.
\end{equation}
By \ref{H:CH} the scaling $B^\eta_{ij}=\eta^{-d}B_{ij}(\cdot/\eta)$ is
$L^1$-critical, so $\|B^\eta_{ij}\|_{L^1}=\|B_{ij}\|_{L^1}\le A_0$ uniformly in
$\eta$, and Young's inequality gives
$\|w^\eta_{ij}(t)\|_{H^{\sigma}}\le A_0\|u_{\eta,j}(t)\|_{H^{\sigma}}
\le A_0\|u_0\|_{H^s}$ by \ref{H:small} and $\sigma=s_1+1\le s$. Also
$\|w^\eta_{ij}\|_{L^\infty}\le\mathsf b_0$ (Lemma~\ref{lem:lawid}(i)), so the
argument of $f_\eta$ and $f'_\eta$ lies in the fixed compact
$I=[-2\mathsf b_0,2\mathsf b_0]$, on which $f_\eta=f$ and, by \ref{H:small},
$\|f\|_{C^{s+1}(I)}\le\eps$.

Here $\sigma=s_1+1$ is an integer, by \eqref{eq:choice}. Since
$w^\eta_{ij}(t,x)\in[-\mathsf b_0,\mathsf b_0]\subset I$ pointwise, $f_\eta=f$ on $I$, and
$\|f\|_{C^{\sigma+2}(I)}\le\eps$ by \ref{H:small} ($\sigma+2\le s+1$),
Lemma~\ref{lem:moser} at $m=\sigma$ with $J=I$, applied to $f-f(0)$ and $f'-f'(0)$, gives
\[
 \big\|f_\eta(w^\eta_{ij})-f(0)\big\|_{H^{\sigma}}\le C,
 \qquad
 \big\|f'_\eta(w^\eta_{ij})-f'(0)\big\|_{H^{\sigma}}\le C,
\]
uniformly in $\eta$ and $t$. Since $\sigma>\frac d2$, $H^{\sigma}$ is an
algebra, $\Delta\varphi\in C^\infty_0\subset H^{\sigma}$, and
$\|u_{\eta,i}\|_{H^{\sigma}}\le\|u_0\|_{H^s}$. Multiplying out and using
Young's inequality once more for the outer convolution, each term of
\eqref{eq:Ddiff} is bounded in $H^{\sigma}$ uniformly in $N$ and $t$. Two
book-keeping points. \emph{(i)} In the diagonal term the constants $\sigma_i$
and $f(0)$ cancel between $\Theta^\eta_i$ and $\Theta_i$, so only
$f(w^\eta_{ik})-f(a_{ik}u_k)$, summed over
$k$, survives there. The two bounds above are therefore needed for every $k$, not only for the fixed pair
$(i,j)$. \emph{(ii)} In the off-diagonal term the constant $f'(0)$ does
\emph{not} multiply $\Delta\varphi$ alone: its contribution is
$B^\eta_{ij}*(f'(0)u_{\eta,i}\Delta\varphi)$, bounded by
$A_0|f'(0)|\,C\|u_0\|_{H^s}\|\Delta\varphi\|_{H^\sigma}$ by the algebra property
and Young, and $|f'(0)|\le\eps$ since $0\in I$. The same bounds hold for the
limit terms, with $u,f,a_{ij}$ in place of $u_\eta,f_\eta,B^\eta_{ij}$; note
$|a_{ij}|=|\int B_{ij}|\le A_0$, so $\|a_{ij}u_j\|_{L^\infty}\le\mathsf b_0$ and
the same composition argument applies. This is \eqref{eq:GupHs}.

\emph{Step 3 (interpolation, then H\"older in time).} For $v\in H^{\sigma}$ the
Fourier-side H\"older inequality with exponents $\frac1{1-\theta},\frac1\theta$
gives $\|v\|_{H^{s_1}}\le\|v\|_{L^2}^{1-\theta}\|v\|_{H^{\sigma}}^{\theta}$
(with constant $1$, since $s_1=\theta\sigma$). Combining with Step~1(i) and
\eqref{eq:GupHs},
\[
 \big\|Q^N_{ij}(t)\big\|_{H^{s_1}_{a_1}}
 \ \le\ C_{\rm mult}\,M_Q^{\theta}\,\big\|Q^N_{ij}(t)\big\|_{L^2}^{1-\theta}.
\]
Square, integrate, and apply H\"older on $(0,T)$ with exponents
$\frac1{1-\theta},\frac1\theta$:
\[
 \int_0^T\big\|Q^N_{ij}\big\|^2_{H^{s_1}_{a_1}}\dt
 \ \le\ C_{\rm mult}^2M_Q^{2\theta}\int_0^T\big\|Q^N_{ij}\big\|_{L^2}^{2(1-\theta)}\dt
 \ \le\ C_{\rm mult}^2M_Q^{2\theta}\,T^{\theta}
 \Big(\int_0^T\big\|Q^N_{ij}\big\|^2_{L^2}\dt\Big)^{1-\theta},
\]
which tends to $0$ by \eqref{eq:GupL2}. This is \eqref{eq:Gupgrade}.
\end{proof}

\begin{lemma}[The remainders vanish]\label{lem:remainders}
Let $\varphi\in H^{s_1+2}_{a_1}$. 
Then, with $\mathcal R^N_i$ as in Theorem~\ref{thm:decomp},
\begin{enumerate}[label=\normalfont(\alph*)]
\item $\displaystyle\sup_{t\le T}\big|\mathcal R^{N,\rm self}_i(\varphi)\big|
 \le C\|\Delta\varphi\|_{L^\infty}\,N^{-1/2}\eta^{-(d+\alpha)}$;
\item $\displaystyle\E\int_0^T\Big|\br{\mu^N_i,\sum_j\mathsf T^N_{ij}\Delta\varphi}\Big|\ds
 \ \le\ C\,\|\Delta\varphi\|_{L^\infty}\,N^{-1/2}\eta^{-\alpha}
 \int_0^T\!\max_j\E\big\|B^\eta_{ij}*\fluc^N_j\big\|^2_{H^{s_1}}\ds;$
\item $\displaystyle\E\int_0^T\Big|\frac1{\sqrt N}\sum_j\br{\fluc^N_i,
 f'_\eta(B^\eta_{ij}*u_{\eta,j})\big(B^\eta_{ij}*\fluc^N_j\big)\Delta\varphi}\Big|\ds
 \ \le\ C\,\|\varphi\|_{H^{s_1+2}_{a_1}}\,N^{-1/2}
 \Big(\E\!\int_0^T\!\|\fluc^N_i\|^2_{\Ufine}\ds\Big)^{\!1/2}
 \Big(T\sup_{t\le T}\max_j\E\big\|B^\eta_{ij}*\fluc^N_j\big\|^2_{H^{s_1}}\Big)^{\!1/2}.$
\end{enumerate}
Moreover, let $\mathcal R^{N,\rm bil}_i$ and $\mathcal R^{N,\rm Tay}_i$ be the bilinear and
the Taylor part of \eqref{eq:Rquad}, as named in Section~\ref{sec:field}. The integrand in (b) is
$|\mathcal R^{N,\rm Tay}_i(\varphi)|$, and the left side of \eqref{eq:cnought} is
$|\mathcal R^{N,\rm bil}_i(\varphi)|$. The bound
\begin{equation}\label{eq:cnought}
 \Big|\frac1{\sqrt N}\sum_j\br{\fluc^N_i,
 f'_\eta(B^\eta_{ij}*u_{\eta,j})\big(B^\eta_{ij}*\fluc^N_j\big)\Delta\varphi}\Big|
 \ \le\ C\,\|\varphi\|_{H^{s_1+2}_{a_1}}\,N^{-1/2}\,
 \|\fluc^N_i\|_{\Ufine}\,\max_j\big\|B^\eta_{ij}*\fluc^N_j\big\|_{H^{s_1}}
\end{equation}
holds for every $\omega$ and $t$, under \ref{H:CH}, \ref{H:cutoff} and
\ref{H:small} alone.

\end{lemma}

\begin{proof}
Clauses (a) and (b) follow from pointwise bounds integrated against the probability measure
$\mu^N_i$. For (c) and \eqref{eq:cnought} we use the duality $\Ufine=(H^{s_1}_{a_1})^*$, a multiplier
bound for $f'_\eta(w^\eta_{ij})$, and Cauchy--Schwarz on $\Omega\times(0,T)$.

(a) is \eqref{eq:self} together with $|\br{\mu^N_i,\cdot}|\le\|\cdot\|_{L^\infty}$
($\mu^N_i$ a probability measure).

(b) By \eqref{eq:ptwise},
$|\mathsf T^N_{ij}|\le\mathsf C_\natural N^{-1/2}\eta^{-\alpha}
|B^\eta_{ij}*\fluc^N_j|^2$ pointwise for every $\omega,t,x$; since $\mu^N_i$ is
a probability measure and $s_1>\frac d2$,
\[
 \Big|\br{\mu^N_i,\textstyle\sum_j\mathsf T^N_{ij}\Delta\varphi}\Big|
 \le\|\Delta\varphi\|_{L^\infty}\sum_j\big\|\mathsf T^N_{ij}\big\|_{L^\infty}
 \le C\,\|\Delta\varphi\|_{L^\infty}N^{-1/2}\eta^{-\alpha}
 \sum_j\big\|B^\eta_{ij}*\fluc^N_j\big\|^2_{H^{s_1}},
\]
and one integrates, takes expectations, and bounds
$\sum_j\E\|\cdot\|^2\le n\max_j\E\|\cdot\|^2$. (c) By duality in $\Ufine=(H^{s_1}_{a_1})^*$ and the multiplication estimate of
\cite{behzadan2021multiplication},
\[
 \big|\br{\fluc^N_i,\,g\,(B^\eta_{ij}*\fluc^N_j)\,\Delta\varphi}\big|
 \le\|\fluc^N_i\|_{\Ufine}\,\big\|g\,\Delta\varphi\,(B^\eta_{ij}*\fluc^N_j)\big\|_{H^{s_1}_{a_1}}
 \le\|\fluc^N_i\|_{\Ufine}\,C\big\|g\Delta\varphi\big\|_{H^{s_1}_{a_1}}\,
 \big\|B^\eta_{ij}*\fluc^N_j\big\|_{H^{s_1}},
\]
the second step because $a_1\ge0$ and $H^{s_1}$ is an algebra, so multiplication
by an $H^{s_1}$ function is bounded on $H^{s_1}_{a_1}$,
with $g=f'_\eta(w^\eta_{ij})$. The function $w^\eta_{ij}\to0$
at infinity, so $g\to f'(0)$, which need not vanish. We subtract it and use
Lemma~\ref{lem:moser} at the integer $m=s_1$:
\[
 g=f'(0)+\big(f'(w^\eta_{ij})-f'(0)\big),\qquad
 \big\|f'(w^\eta_{ij})-f'(0)\big\|_{H^{s_1}}\le C
 \quad\text{uniformly in }\eta,t,
\]
by Lemma~\ref{lem:moser} with $J=I$ and $\|f\|_{C^{s_1+1}(I)}\le\eps$, $\|w^\eta_{ij}\|_{H^{s_1}}\le A_0\|u_0\|_{H^s}$
($L^1$-criticality and \ref{H:small}) and $w^\eta_{ij}(t,x)\in[-\mathsf b_0,\mathsf b_0]\subset I$
pointwise. Hence, since $a_1\ge0$ and $H^{s_1}$ is an algebra,
$\|g\Delta\varphi\|_{H^{s_1}_{a_1}}\le\big(|f'(0)|+C\big)C\|\Delta\varphi\|_{H^{s_1}_{a_1}}
\le C\|\varphi\|_{H^{s_1+2}_{a_1}}$, uniformly in $\eta$.

Now apply Cauchy--Schwarz on $\Omega\times(0,T)$:
\[
 \E\!\int_0^T\!\|\fluc^N_i\|_{\Ufine}\,\big\|B^\eta_{ij}*\fluc^N_j\big\|_{H^{s_1}}\ds
 \ \le\ \Big(\E\!\int_0^T\!\|\fluc^N_i\|^2_{\Ufine}\ds\Big)^{\!1/2}
 \Big(\E\!\int_0^T\!\big\|B^\eta_{ij}*\fluc^N_j\big\|^2_{H^{s_1}}\ds\Big)^{\!1/2},
\]
and bound the second factor by
$\big(T\sup_{t\le T}\E\|B^\eta_{ij}*\fluc^N_j\|^2_{H^{s_1}}\big)^{1/2}$. Summing
over $j$ (a factor $n$, absorbed into $C$) gives the stated form.
\end{proof}

\begin{lemma}[$\mathcal R^N$ is bounded, deterministically, at fixed $N$ and
$\eta$]\label{lem:Rdet}
Let $N\in\N$ and $\eta\in(0,\eta_0]$ be
fixed. Then, with $\mathcal R^N$ as in Theorem~\ref{thm:decomp},
\begin{equation}\label{eq:Rdet}
 \sup_{t\le T}\ \max_i\
 \big\|\mathcal R^N_i(t)\big\|_{H^{-s_1-2}_{-a_1}}
 \ \le\ \mathsf K_{N,\eta}
 \:=\ Cn\Big(\big\|f''_\eta\big\|_{L^\infty(\R)}\,\eta^{-2d}
 +\eta^{-(d+\alpha)}\Big)\,N^{1/2}\ <\ \infty,
\end{equation}
where $C=C(d,s_1,a_1,\max_{ij}\|B_{ij}\|_{L^\infty},A_0,C_s\|u_0\|_{H^s})$.
The bound \eqref{eq:Rdet} holds $\Prob$-a.s. In particular, for a.e.\ $\omega$,
$\mathcal R^N\in L^\infty\big(0,T;(H^{-s_1-2}_{-a_1})^n\big)
\subset L^1\big(0,T;(H^{-s_1-2}_{-a_1})^n\big)$.
\end{lemma}

\begin{proof}
At fixed $N$ and $\eta$ every term of $\mathcal R^N_i(\varphi)$ is a bounded function integrated
against a finite measure. Step~(1) bounds the field $B^\eta_{ij}*\fluc^N_j$ crudely in $L^\infty$,
Steps~(2)--(4) bound the three pieces of $\mathcal R^N$, and Step~(5) checks measurability in $t$.

Let $\varphi\in H^{s_1+2}_{a_1}$ with $\|\varphi\|_{H^{s_1+2}_{a_1}}\le1$.
By Lemma~\ref{lem:spaces}(iii) at $m=2$, $\|\Delta\varphi\|_{L^\infty}\le C$.
We use throughout, at every $t\le T$, that $\mu^N_i(t)$ is a
probability measure, that $u_{\eta,i}(t)\ge0$ with $\int u_{\eta,i}(t)=1$, and
that $\|u_{\eta,i}(t)\|_{L^\infty}\le C_u$, with $C_u\le C_s\|u_0\|_{H^s}$ as in \eqref{eq:Cu}. Hence
$\fluc^N_i(t)$ is a signed measure of total variation at most $2\sqrt N$.

\emph{(1) The field, at fixed $N$ and $\eta$.} By \ref{H:CH}, $B_{ij}\in
C^\infty_0$, so $\|B^\eta_{ij}\|_{L^\infty}=\eta^{-d}\|B_{ij}\|_{L^\infty}$ and
$\|B^\eta_{ij}\|_{L^1}=\|B_{ij}\|_{L^1}\le A_0$. Since
$B^\eta_{ij}*\fluc^N_j=\sqrt N\big(\frac1N\sum_kB^\eta_{ij}(\cdot-X^{N,\eta}_{k,j})
-B^\eta_{ij}*u_{\eta,j}\big)$,
\begin{equation}\label{eq:fieldcrude}
 \big\|B^\eta_{ij}*\fluc^N_j(t)\big\|_{L^\infty}
 \ \le\ \sqrt N\big(\eta^{-d}\|B_{ij}\|_{L^\infty}+A_0C_u\big)
 \ \le\ C\eta^{-d}\sqrt N
\end{equation}
for every $\omega$, $t$ and $x$.

\emph{(2) The self-interaction defect.} By \eqref{eq:self},
$|\varrho^N_{k,i}|\le CN^{-1}\eta^{-(d+\alpha)}$, whence
$|\mathcal R^{N,\rm self}_i(\varphi)|=\sqrt N|\br{\mu^N_i,\varrho^N_{\cdot,i}
\Delta\varphi}|\le CN^{-1/2}\eta^{-(d+\alpha)}$.

\emph{(3) The Taylor piece.} By \eqref{eq:Taylor-outline2} and \eqref{eq:fieldcrude},
$\|\mathsf T^N_{ij}(t)\|_{L^\infty}\le\frac{\|f''_\eta\|_{L^\infty(\R)}}
{2\sqrt N}\,C^2\eta^{-2d}N$, so
$|\mathcal R^{N,\rm Tay}_i(\varphi)|=|\br{\mu^N_i,\sum_j\mathsf T^N_{ij}
\Delta\varphi}|\le Cn\|f''_\eta\|_{L^\infty(\R)}\eta^{-2d}\sqrt N$.
Here $\|f''_\eta\|_{L^\infty(\R)}<\infty$ because \ref{H:CH} gives
$f_\eta\in W^{s+1,\infty}(\R)$ on all of $\R$, with $s+1>\frac d2+2\ge2$.

\emph{(4) The bilinear piece.} $\|f'_\eta\|_{L^\infty(\R)}\le\mathrm{Lip}(f_\eta)
\le\eta^{-\alpha}$ by \ref{H:CH}, so by \eqref{eq:fieldcrude} and the total
variation bound on $\fluc^N_i$,
\[
 \big|\mathcal R^{N,\rm bil}_i(\varphi)\big|
 =\Big|\tfrac1{\sqrt N}\sum_j\br{\fluc^N_i,f'_\eta(B^\eta_{ij}*u_{\eta,j})
 (B^\eta_{ij}*\fluc^N_j)\Delta\varphi}\Big|
 \le\tfrac1{\sqrt N}\cdot2\sqrt N\cdot n\,\eta^{-\alpha}\cdot C\eta^{-d}\sqrt N\cdot C
 \le Cn\,\eta^{-(d+\alpha)}\sqrt N.
\]
Every pairing in (2)--(4) is an integral of a bounded measurable function
against a finite measure, namely $\mu^N_i$, $u_{\eta,i}(t)\dd x$, or their
difference. These
integrals coincide with the $H^{-s_1-2}_{-a_1}$--$H^{s_1+2}_{a_1}$
duality pairing by Lemma~\ref{lem:spaces}(i) and the density of
$\mathcal S$ in $H^{s_1+2}_{a_1}$, which holds because $J_{a_1}$ maps $\mathcal S$ onto $\mathcal S$.

\emph{(5) Measurability.} For fixed $\varphi$, $t\mapsto\mathcal R^N_i(t,\varphi)$
is continuous: the particle paths are continuous, and $u_{\eta,i}\in
C([0,T];L^2)$ by Lemma~\ref{lem:coeffcont}(i) at $s'=0$. Weak measurability of
$t\mapsto\mathcal R^N(t)$ against $(H^{-s_1-2}_{-a_1})^*=H^{s_1+2}_{a_1}$
(Lemma~\ref{lem:spaces}(i)) follows; $H^{-s_1-2}_{-a_1}$ is separable, being
isometric to $H^{-s_1-2}$ by $J_{-a_1}$, so Pettis' theorem \cite[Ch.~II, \S1, Thm.~2]{DiUh77} gives strong
measurability of $t\mapsto\mathcal R^N(t)\in(H^{-s_1-2}_{-a_1})^n$.
 With the bound of (2)--(4)
this gives \eqref{eq:Rdet} and the stated inclusions.
\end{proof}

\begin{proposition}[The field bound under a fixed-time fourth moment]
\label{prop:close}
Assume that the fixed-time fourth moment $\overline K_4$ of \eqref{eq:K4def} satisfies
\begin{equation}\label{eq:K4-mild}
 \overline K_4\ =\ o\big(N\eta^{2\alpha}\big)\qquad(N\to\infty).
\end{equation}
Then $\ \sup_N\E\int_0^T\|\fluc^N(t)\|^2_{\Ufine}\dd t<\infty$, and consequently
$\ \sup_N\E\sup_{t\le T}\|\fluc^N(t)\|^2_{\Usharp}<\infty$.
\end{proposition}

\begin{proof}
Let $\mathsf I:=\int_0^T\|\fluc^N\|^2_{\Ufine}\dd t$. We bound $\E\mathsf I$ through the mild
formula \eqref{eq:mildzero}. In Step~1 the initial and martingale terms are bounded uniformly in $N$,
and the remainder term by $\int_0^T\|\mathcal R^N\|^2$. The bilinear part of the remainder is at most a
vanishing factor times $(\E\mathsf I)^{1/2}$ (Step~2), and the other two parts vanish (Step~3). Young's
inequality then absorbs the bilinear part (Step~4). Step~5 deduces the supremum bound in $\Usharp$
from the same formula.

\emph{Step 1 (a pathwise bound on $\mathsf I$ from the mild formula)}. With the notation of
\eqref{eq:mildzero}, put $\widetilde Z^N(t):=\fluc^N(t)-P^N(t)-D^N(t)$, an element of
$\Fcoarse^n=(\Test^*)^n$ defined for every $t$ and $\omega$. The three processes on the
right are weakly continuous in $t$ with values in $\Fcoarse^n$: $\fluc^N$ is continuous in
$\Ufine^n$ (Lemma~\ref{lem:uetacont}(ii)). $\br{P^N(t),\varphi}=\br{\fluc^N(0),S^\eta(t,0)^*\varphi}$
is continuous by Corollary~\ref{cor:fwdweak}(ii), and
$\br{D^N(t),\varphi}=\int_0^t\br{\mathcal R^N(s),S^\eta(t,s)^*\varphi}\dd s$ is continuous by
Corollary~\ref{cor:fwdweak}(iii), dominated convergence and Lemma~\ref{lem:Rdet}, and measurable
in $\omega$ for each $t$, since $\mathcal R^N_i(s,\varphi)$ is continuous in $s$ and measurable in
$\omega$. So
$\widetilde Z^N$ is jointly measurable in $(t,\omega)$, by Pettis' theorem, and
$\|\cdot\|_{\Ufine}$ is a Borel function on $\Fcoarse^n$ by Lemma~\ref{lem:lscsup}(i). By
Proposition~\ref{prop:mildN}, $\widetilde Z^N(t)=Z^N(t)$ almost surely for each $t$. By
Lemma~\ref{lem:BDG} and \eqref{eq:trMN},
$\E\|Z^N(t)\|^2_{\Ufine}\le C_1(T)^2\E\int_0^T\sum_i\mathsf q^N_i\dd r$.

Apply Corollary~\ref{cor:maxL2-forced} (a deterministic operator bound, applied
$\omega$ by $\omega$) to $F:=\mathcal R^N$, which lies in $L^2(0,T;(H^{-s_1-2}_{-a_1})^n)$ by
Lemma~\ref{lem:Rdet}. The element it bounds, defined by \eqref{eq:weakdef}, agrees with $D^N(r)$ in
$(\Test^*)^n$ for almost every $r$. Indeed, let $g\in C^1_c((0,T))$ be scalar, let
$\varphi\in\Test^n$, and take $G:=g\varphi$ in \eqref{eq:weakdef}. On $\{0\le s\le r\le T\}$ the function
$(s,r)\mapsto g(r)\br{\mathcal R^N(s),S^\eta(r,s)^*\varphi}$ is measurable in $s$ for fixed $r$, by
Lemma~\ref{lem:Rdet} and Lemma~\ref{lem:strongcont}, and continuous in $r$ for fixed $s$, by
Corollary~\ref{cor:fwdweak}(iii). Hence the function is jointly measurable, and it is bounded by
Lemma~\ref{lem:Rdet} and Proposition~\ref{prop:evol}(iv). Fubini's theorem therefore turns
\eqref{eq:weakdef} into $\int_0^Tg(r)\br{D^N(r),\varphi}\dd r$. Since $g$ is arbitrary and $\varphi$
may range over a countable dense subset of $\Test^n$, the two elements agree for almost every $r$.
Then, pathwise,
\[
 \mathsf I\ \le\ \underbrace{3TC_1(T)^2\|\fluc^N(0)\|^2_{\Ufine}+3\!\int_0^T\!\|\widetilde Z^N(t)\|^2_{\Ufine}\dd t}_{=:\mathsf A}
 \ +\ 3C_3(T)^2\!\int_0^T\!\big\|\mathcal R^N\big\|^2_{H^{-s_1-2}_{-a_1}}\dd s,
\]
and by Tonelli, Lemma~\ref{lem:field0}, \eqref{eq:psi2} and Cauchy--Schwarz in time,
$\E\mathsf A\le3TC_1(T)^2\big(\E\|\fluc^N(0)\|^2_{\Ufine}+\E\int_0^T\sum_i\mathsf q^N_i\big)\le\mathsf a$,
where the constant $\mathsf a$ does not depend on $N$.

\emph{Step 2 (the bilinear remainder).} Let
$\mathsf W(t):=\max_{i,j}\|B^\eta_{ij}*\fluc^N_j(t)\|_{H^{s_1}}$. The bound \eqref{eq:cnought} of
Lemma~\ref{lem:remainders} is pathwise, assumes no field bound, and is bilinear:
\begin{equation}\label{eq:bilin}
 \big\|\mathcal R^{N,\rm bil}_i(s)\big\|_{H^{-s_1-2}_{-a_1}}
 \le CN^{-1/2}\big\|\fluc^N_i(s)\big\|_{\Ufine}\,\mathsf W(s),
\end{equation}
and we combine it with the deterministic bound
\begin{equation}\label{eq:detbound}
 \big\|\fluc^N_i(t)\big\|_{\Ufine}\ \le\ 2C_\delta\sqrt N,
 \qquad C_\delta:=\sup_x\|\delta_x\|_{\Ufine}<\infty,
\end{equation}
(valid by Lemma~\ref{lem:spaces}(iii) with $s_1>\frac d2$, $a_1\ge0$, since both
$\mu^N_i$ and $u_{\eta,i}$ are probability measures) to lower the moment order,
at the cost of a factor $\sqrt N$:
\[
 \Big(\E\!\int_0^T\!\|\fluc^N\|^4_{\Ufine}\Big)^{1/2}
 \le\big(4nC_\delta^2N\big)^{1/2}\,(\E\mathsf I)^{1/2}.
\]
Hölder on $\Omega\times(0,T)$ then gives, with $\kappa:=2n^{3/2}C^2C_\delta T^{1/2}$,
\[
 \E\!\int_0^T\!\|\mathcal R^{N,\rm bil}\|^2
 \le C^2N^{-1}\Big(\E\!\int\!\|\fluc^N\|^4_{\Ufine}\Big)^{1/2}\Big(\E\!\int\!\mathsf W^4\Big)^{1/2}
 \le \kappa\,\overline K_4^{1/2}\,N^{-1/2}\,(\E\mathsf I)^{1/2}.
\]

\emph{Step 3 (the self-interaction and Taylor remainders).} By Lemma~\ref{lem:remainders}(a),
$\|\mathcal R^{N,\rm self}\|\le CN^{-1/2}\eta^{-(d+\alpha)}$ pathwise,
so $\E\int\|\mathcal R^{N,\rm self}\|^2\le CTN^{-1}\eta^{-2(d+\alpha)}\to0$ by
Lemma~\ref{lem:polylog}. For the Taylor piece we use a bound in expectation: Lemma~\ref{lem:taylor} gives
\[
 \E\int_0^T\big\|\mathcal R^{N,\rm Tay}\big\|^2_{H^{-s_1-2}_{-a_1}}\dd s
 \ \le\ C\,T\,N^{-1}\,\overline K_4\,\eta^{-2\alpha}\ \longrightarrow\ 0
\] by \eqref{eq:K4-mild}. We collect Steps~2 and~3, using that $\|\mathcal R^N\|^2$ is at most three times the sum of the
squared norms of $\mathcal R^{N,\rm self}$, $\mathcal R^{N,\rm Tay}$ and $\mathcal R^{N,\rm bil}$:
\begin{equation}\label{eq:Rcollect}
 \E\int_0^T\big\|\mathcal R^N\big\|^2_{H^{-s_1-2}_{-a_1}}\dd s
 \ \le\ C(T)\Big(N^{-1}\eta^{-2(d+\alpha)}+\overline K_4\,N^{-1}\eta^{-2\alpha}\Big)
 +3\kappa\,\overline K_4^{1/2}N^{-1/2}\big(\E\mathsf I\big)^{1/2}.
\end{equation}

\emph{Step 4 (absorption, and the bound on $\E\mathsf I$).} We take expectations in the
bound on $\mathsf I$ of Step~1 and insert \eqref{eq:Rcollect}. The coefficient in Young's inequality is $9C_3^2$: the bound on $\mathsf I$ above carries
$3C_3(T)^2$, and $\|\mathcal R^N\|^2\le3\sum$ of the three pieces
contributes a second factor $3$. So
$9C_3^2\kappa\overline K_4^{1/2}N^{-1/2}(\E\mathsf I)^{1/2}\le\frac12\E\mathsf I+
\frac{81}{2}C_3^4\kappa^2\overline K_4N^{-1}$.
Absorption of $\frac12\E\mathsf I$ into the left side is legitimate because $\E\mathsf I\le4nC_\delta^2NT<\infty$ by
\eqref{eq:detbound}. Since $\eta\le1$, the last term is at most $\frac{81}{2}C_3^4\kappa^2\overline K_4N^{-1}\eta^{-2\alpha}$. Hence
\[
 \E\int_0^T\|\fluc^N\|^2_{\Ufine}\dd t\ \le\ 2\E\mathsf A
 +C(T)\Big(N^{-1}\eta^{-2(d+\alpha)}+\overline K_4\,N^{-1}\eta^{-2\alpha}\Big),
\]
bounded uniformly in $N$ by $\E\mathsf A\le\mathsf a$, Lemma~\ref{lem:polylog} and \eqref{eq:K4-mild}.

\emph{Step 5 (the supremum in $\Usharp$)}. Fix a countable dense $\mathbb D\subset[0,T]$.
For each $t$, almost surely, \eqref{eq:mildN} holds, and Corollary~\ref{cor:maxreg} at
$\sigma=s_1+1$ bounds $\|D^N(t)\|_{\Usharp}$ by $C_{\rm mr}(T)\|\mathcal R^N\|_{L^2(0,T;H^{-s_1-2}_{-a_1})}$,
a right-hand side independent of $t$. So, off one null set, for all $t\in\mathbb D$,
\[
 \|\fluc^N(t)\|_{\Usharp}\le c_\sharp C_1(T)\|\fluc^N(0)\|_{\Ufine}
 +c_\sharp\|Z^N(t)\|_{\Ufine}+C_{\rm mr}(T)\|\mathcal R^N\|_{L^2(0,T;H^{-s_1-2}_{-a_1})}.
\]
$\fluc^N$ is continuous in $\Ufine$, hence in $\Usharp$, so the supremum over $[0,T]$ equals the
supremum over $\mathbb D$. Lemma~\ref{lem:field0}, Proposition~\ref{prop:conv} with Jensen,
\eqref{eq:Rcollect} and the bound on $\E\mathsf I$ of Step~4 give the claim.
\end{proof}

\begin{lemma}[Telescoping along dyadic partitions]\label{lem:telescope}
Let $\eta\in[0,\eta_0]$, $t\in(0,T]$ and $\varphi\in\Test^n$, and let
$\Phi(s):=S^\eta(t,s)^*\varphi$ for $s\in[0,t]$, as in Proposition~\ref{prop:evol}. For $k\in\N$ let
$\pi_k$ be the dyadic partition of $[0,t]$, with nodes $s_l:=lt2^{-k}$, $0\le l\le2^k$, and mesh
$|\pi_k|=t2^{-k}$. Let $v=(v_1,\dots,v_n):[0,t]\to(\Usharp)^n$ be bounded, and let
$r\mapsto\br{v_j(r),\psi}$ be continuous on $[0,t]$ for every $j\le n$ and every
$\psi\in(\Usharp)^*=H^{s_1+1}_{a_1}$. Let $v^0\in(\Test^*)^n$, and let real numbers
$\mathcal V_i(\tau,\chi)$ be given for $i\le n$, $\tau\in[0,t]$ and $\chi\in\Test$. Assume that
\begin{equation}\label{eq:teleweak}
 \br{v_i(\tau),\chi}=\br{v^0_i,\chi}
 +\int_0^\tau\sum_{j=1}^n\br{v_j(s),(\mathcal G^\eta(s)\chi)_{ij}}\dd s+\mathcal V_i(\tau,\chi)
\end{equation}
for every $i\le n$, every dyadic node $\tau$ of $[0,t]$, and every $\chi=\Phi_i(\varsigma)$
with $\varsigma$ a dyadic node of $[0,t]$. Then, for every $k\in\N$,
\begin{equation}\label{eq:telescope}
 \sum_i\br{v_i(t),\varphi_i}-\sum_i\br{v_i(0),\Phi_i(0)}=\mathcal E_{\pi_k}+\mathcal V_{\pi_k},
\end{equation}
where, for $\pi=\pi_k$ and with $l$ ranging over $0\le l<2^k$,
\begin{equation}\label{eq:Epi}
\begin{aligned}
 \mathcal E_\pi&:=\sum_l\int_{s_l}^{s_{l+1}}\sum_{i,j}
 \Big[\br{v_j(s),(\mathcal G^\eta(s)\Phi_i(s_{l+1}))_{ij}}
 -\br{v_j(s_l),(\mathcal G^\eta(s)\Phi_i(s))_{ij}}\Big]\dd s,\\
 \mathcal V_\pi&:=\sum_l\sum_i\Big[\mathcal V_i\big(s_{l+1},\Phi_i(s_{l+1})\big)
 -\mathcal V_i\big(s_l,\Phi_i(s_{l+1})\big)\Big].
\end{aligned}
\end{equation}
Moreover $\mathcal E_{\pi_k}\to0$ as $k\to\infty$.
\end{lemma}

\begin{proof}
The identity \eqref{eq:telescope} is exact. It follows from bilinearity, from
\eqref{eq:teleweak} at consecutive nodes, and from the integrated backward equation
\eqref{eq:intbackward}. The convergence of $\mathcal E_{\pi_k}$ uses the bound on $v$, the continuity
of $\Phi$ in $\Test^n$, and the compactness of the set $\mathcal K^\eta$ of \eqref{eq:compactN}.

\emph{Step 0 (tools).} Put $R_v:=\sup_{r\le t}\|v(r)\|_{(\Usharp)^n}<\infty$. By
Lemma~\ref{lem:strongcont}, $\Phi\in C([0,t];\Test^n)$. Let $\omega_\Phi$ be its modulus of
continuity in $\Test^n$. It tends to $0$ at $0$, by uniform continuity on the compact $[0,t]$. For
$s,s'\in[0,t]$, Lemma~\ref{lem:Gbound} gives
\[
 \big\|\mathcal G^\eta(s)\Phi_i(s)-\mathcal G^\eta(s')\Phi_i(s')\big\|
 \le C_{\mathcal G}\|\Phi_i(s)-\Phi_i(s')\|_{\Test}
 +\big\|\mathcal G^\eta(s)-\mathcal G^\eta(s')\big\|\,\|\Phi_i(s')\|_{\Test},
\]
with the operator norm of $\mathcal L(\Test,H^{s_1+1}_{a_1})$. So $s\mapsto(\mathcal G^\eta(s)\Phi_i(s))_{ij}$
is continuous into $H^{s_1+1}_{a_1}$, and
\begin{equation}\label{eq:compactN}
 \mathcal K^\eta:=\big\{(\mathcal G^\eta(s)\Phi_i(s))_{ij}:s\in[0,t],\,i,j\le n\big\}
 \ \text{is compact in }H^{s_1+1}_{a_1}=(\Usharp)^*,
\end{equation}
as a finite union of continuous images of $[0,t]$. Let $g:[0,t]\to H^{s_1+1}_{a_1}$ be continuous.
Then $s\mapsto\br{v_j(s),g(s)}$ is continuous, because
$|\br{v_j(s),g(s)}-\br{v_j(s'),g(s')}|\le R_v\|g(s)-g(s')\|+|\br{v_j(s)-v_j(s'),g(s')}|$. With
$g(s)=(\mathcal G^\eta(s)\chi)_{ij}$ and $g(s)=(\mathcal G^\eta(s)\Phi_i(s))_{ij}$, this shows that all
integrals in \eqref{eq:teleweak} and \eqref{eq:Epi} exist.

\emph{Step 1 (the exact identity).} Since $\Phi(t)=\varphi$, bilinearity gives, exactly,
\[
\sum_i\br{v_i(t),\varphi_i}-\sum_i\br{v_i(0),\Phi_i(0)}=\sum_l(A_l+B_l),
\]
where
\[
A_l:=\sum_i\br{v_i(s_{l+1})-v_i(s_l),\Phi_i(s_{l+1})},\qquad
 B_l:=\sum_j\br{v_j(s_l),\Phi_j(s_{l+1})-\Phi_j(s_l)}.
\]
Subtracting \eqref{eq:teleweak} at $\tau=s_l$ from \eqref{eq:teleweak} at $\tau=s_{l+1}$, with
the same test function $\chi=\Phi_i(s_{l+1})$, cancels $v^0$ and gives
\[
A_l=\int_{s_l}^{s_{l+1}}\sum_{i,j}\br{v_j(s),(\mathcal G^\eta(s)\Phi_i(s_{l+1}))_{ij}}\dd s
 +\sum_i\Big[\mathcal V_i\big(s_{l+1},\Phi_i(s_{l+1})\big)-\mathcal V_i\big(s_l,\Phi_i(s_{l+1})\big)\Big].
\]
By \eqref{eq:intbackward} at $s=s_l$ and $s=s_{l+1}$,
$\Phi_j(s_{l+1})-\Phi_j(s_l)=-\int_{s_l}^{s_{l+1}}\sum_i(\mathcal G^\eta(s)\Phi_i(s))_{ij}\dd s$, a Bochner
integral in $H^{s_1+1}_{a_1}=(\Usharp)^*$. Pairing with $v_j(s_l)\in\Usharp$ commutes with this
integral, so
\[
 B_l=-\int_{s_l}^{s_{l+1}}\sum_{i,j}\br{v_j(s_l),(\mathcal G^\eta(s)\Phi_i(s))_{ij}}\dd s.
\]
In both $A_l$ and $B_l$ the field carries the index $j$ and the test function the index $i$.
Summing over $l$ gives \eqref{eq:telescope}.

\emph{Step 2 ($\mathcal E_{\pi_k}\to0$).} For $s\in(s_l,s_{l+1}]$ we split the bracket in
\eqref{eq:Epi} as
\[
 \br{v_j(s),(\mathcal G^\eta(s)[\Phi_i(s_{l+1})-\Phi_i(s)])_{ij}}
 +\br{v_j(s)-v_j(s_l),(\mathcal G^\eta(s)\Phi_i(s))_{ij}}.
\]
By \eqref{eq:Gbound}, the first term is at most $R_vC_{\mathcal G}\,\omega_\Phi(|\pi|)$. For the
second, fix $\theta>0$ and cover $\mathcal K^\eta$ by finitely many balls $B(\psi^q,\theta)$, $q\le Q$,
of $H^{s_1+1}_{a_1}$ (a $\theta$-net). Then
\[
 \big|\br{v_j(s)-v_j(s_l),(\mathcal G^\eta(s)\Phi_i(s))_{ij}}\big|
 \le2\theta R_v+\max_{q\le Q}\big|\br{v_j(s)-v_j(s_l),\psi^q}\big|.
\]
Each $r\mapsto\br{v_j(r),\psi^q}$ is uniformly continuous on $[0,t]$, and $Q$ is finite, so the
maximum tends to $0$ as $|\pi|\to0$, uniformly in $s$ and $l$. Summing over the $n^2$ pairs $(i,j)$
and integrating over $[0,t]$ gives $\limsup_{k\to\infty}|\mathcal E_{\pi_k}|\le2n^2t\,\theta R_v$. Since
$\theta>0$ is arbitrary, $\mathcal E_{\pi_k}\to0$.
\end{proof}

\begin{proposition}[Mild formula at finite $N$]\label{prop:mildN}
Let $t\in[0,T]$. Then, $\Prob$-almost surely, the terms $P^N(t)$, $Z^N(t)$ and $D^N(t)$ of
\eqref{eq:mildzero} lie in $\Ufine^n$, and
\begin{equation}\label{eq:mildN}
 \fluc^N(t)=P^N(t)+Z^N(t)+D^N(t)\qquad\text{in }\Ufine^n.
\end{equation}
\end{proposition}

\begin{proof}
We fix $t$ and $\varphi\in\Test^n$ and apply Lemma~\ref{lem:telescope} to $\fluc^N$, with the
integrated remainder plus the martingale in the role of $\mathcal V_i$ (Step~1). The remainder term
converges because $\Phi$ is continuous in $H^{s_1+2}_{a_1}$ and $\mathcal R^N\in L^1(0,T;H^{-s_1-2}_{-a_1})$
by Lemma~\ref{lem:Rdet} (Step~2). The martingale term converges by the It\^o isometry (Step~3). Step~4
identifies the limit, first in $(\Test^*)^n$ and then in $\Ufine^n$.

Throughout, $N$ and $\eta=\eta(N)\in(0,\eta_0]$ are fixed.

\emph{Step 0: the operator is $\mathcal G^\eta$, and the weak form.}
The operator $\mathcal G^\eta$ of \eqref{eq:GN} is deterministic, and depends on $N$
only through $\eta=\eta(N)$: its coefficients are $\Theta^\eta_i$ and
$u_{\eta,i}f'_\eta(B^\eta_{ij}*u_{\eta,j})$, built from the solution $u_\eta$ of
\eqref{eq:nonloc} and not from the empirical measure. It is the operator of
Lemma~\ref{lem:Gbound} and of \eqref{eq:backward-mild}, whose backward family is $S^\eta$.

Integrating \eqref{eq:fluc_decomp} gives, for each $\varphi\in\Test$,
$\Prob$-a.s.\ for all $t\in[0,T]$,
\begin{equation}\label{eq:weakN}
 \br{\fluc^N_i(t),\varphi}=\br{\fluc^N_i(0),\varphi}
 +\int_0^t\sum_j\br{\fluc^N_j(s),(\mathcal G^\eta(s)\varphi)_{ij}}\dd s
 +\int_0^t\mathcal R^N_i(s,\varphi)\dd s+M^N_i(t,\varphi),
\end{equation}
which is \eqref{eq:weakform} with $\mathcal G^\eta$ in place of $\mathcal G$ and
the extra remainder. Theorem~\ref{thm:decomp} applies to $\varphi\in\Test$, since
$\Test\hookrightarrow H^{s_1+2}_{a_1}$ by Lemma~\ref{lem:spaces}(ii). Both integrands in \eqref{eq:weakN} are
absolutely integrable: the first because $\fluc^N_j(s)\in\Usharp$ with the
deterministic bound of Lemma~\ref{lem:uetacont}(i) while
$(\mathcal G^\eta(s)\varphi)_{ij}\in H^{s_1+1}_{a_1}=(\Usharp)^*$ uniformly in
$s$ (Lemma~\ref{lem:Gbound}). The second by Lemma~\ref{lem:Rdet}.

\emph{Step 1: telescoping.} Fix $t\in(0,T]$ and $\varphi\in\Test^n$ and set
$\Phi(s):=S^\eta(t,s)^*\varphi$. Let $(\pi_k)$ be the dyadic partitions of $[0,t]$, as in
Lemma~\ref{lem:telescope}. The triples formed by a dyadic node $\tau$ of $[0,t]$, an index $i\le n$
and a test function $\Phi_i(\varsigma)\in\Test$ with $\varsigma$ a dyadic node are countably many. So
there is one null set off which \eqref{eq:weakN} holds at all of them. Off this null set we apply
Lemma~\ref{lem:telescope}, $\omega$ by $\omega$, with $v:=\fluc^N$, $v^0:=\fluc^N(0)$ and
\[
 \mathcal V_i(\tau,\chi):=\int_0^\tau\mathcal R^N_i(s,\chi)\dd s+M^N_i(\tau,\chi).
\]
Its hypotheses hold. The bound $\sup_r\|\fluc^N(r)\|_{(\Usharp)^n}\le C\sqrt N$ is
Lemma~\ref{lem:uetacont}(i). The continuity of $r\mapsto\br{\fluc^N_j(r),\psi}$ for
$\psi\in(\Usharp)^*$ follows from Lemma~\ref{lem:uetacont}(ii). And \eqref{eq:weakN} is
\eqref{eq:teleweak}. The term $\mathcal V_\pi$ of \eqref{eq:Epi} splits as
$\mathcal V_\pi=\mathcal N_\pi+\mathcal M_\pi$, with
\[
 \mathcal N_\pi:=\sum_l\int_{s_l}^{s_{l+1}}\sum_i
 \mathcal R^N_i\big(s,\Phi_i(s_{l+1})\big)\dd s
 =\int_0^t\sum_i\br{\mathcal R^N_i(s),\Phi^\pi_i(s)}\dd s,
 \qquad\Phi^\pi:=\sum_l\Phi(s_{l+1})\mathbf 1_{(s_l,s_{l+1}]},
\]
\[
 \mathcal M_\pi:=\sum_l\sum_i\Big[M^N_i\big(s_{l+1},\Phi_i(s_{l+1})\big)-M^N_i\big(s_l,\Phi_i(s_{l+1})\big)\Big].
\]
So, almost surely, for every $k$,
\[
 \sum_i\br{\fluc^N_i(t),\varphi_i}-\sum_i\br{\fluc^N_i(0),\Phi_i(0)}
 =\mathcal E_{\pi_k}+\mathcal N_{\pi_k}+\mathcal M_{\pi_k},
 \qquad\mathcal E_{\pi_k}\to0\ (k\to\infty).
\]

\emph{Step 2: $\mathcal N_\pi\to\int_0^t\sum_i\br{\mathcal R^N_i(s),\Phi_i(s)}
\dd s$ a.s.} Let $\omega^\sharp_\Phi$ denote the modulus of continuity of $\Phi$
in $H^{s_1+2}_{a_1}$. It exists, and tends to $0$ at
$0$, by uniform continuity of $\Phi\in C([0,t];\Test^n)$ together with
$\Test\hookrightarrow H^{s_1+2}_{a_1}$ (Lemma~\ref{lem:spaces}(ii)). Since $\Phi^\pi$ is defined on $(0,t]$ and takes the right-endpoint
value, $\sup_{s\in(0,t]}\|\Phi^\pi(s)-\Phi(s)\|_{H^{s_1+2}_{a_1}}
\le\omega^\sharp_\Phi(|\pi|)$, and $\{0\}$ is Lebesgue-null, so
\begin{equation}\label{eq:Npi}
 \Big|\mathcal N_\pi-\int_0^t\sum_i\br{\mathcal R^N_i(s),\Phi_i(s)}\dd s\Big|
 \ \le\ \omega^\sharp_\Phi(|\pi|)\int_0^t\sum_i
 \big\|\mathcal R^N_i(s)\big\|_{H^{-s_1-2}_{-a_1}}\dd s
 \ \le\ n\,T\,\mathsf K_{N,\eta}\,\omega^\sharp_\Phi(|\pi|)\ \longrightarrow\ 0
\end{equation}
by Lemma~\ref{lem:Rdet}.

\emph{Step 3: $\mathcal M_\pi\to\sum_i\int_0^t\br{\dd M^N_i,\Phi_i}$ in
$L^2(\Omega)$.} $\Phi^\pi$ is left-continuous with deterministic
coefficients, hence predictable, so $\mathcal M_\pi=\sum_i\int_0^t
\br{\dd M^N_i,\Phi^\pi_i}$ exactly. The cross-brackets vanish for $i\ne k$ by
Proposition~\ref{prop:ito}, so with $\mathsf d_i:=\int_0^t\br{\dd M^N_i,\Phi^\pi_i-\Phi_i}$ one has
$\E[(\sum_i\mathsf d_i)^2]=\sum_i\E\br{\mathsf d_i}_t$, and, by the It\^o isometry for
$M^N_i$, which follows from \eqref{eq:MiN},
\[
 \E\Big|\mathcal M_\pi-\sum_i\int_0^t\br{\dd M^N_i,\Phi_i}\Big|^2
 =\sum_i\E\int_0^t\big\langle\mu^N_i,\big|\na(\Phi^\pi_i-\Phi_i)\big|^2
 \big(2\Theta^N_i+2\varrho^N_{\cdot,i}\big)\big\rangle\dd s
 \le\Big(\sum_i\E\int_0^T\mathsf q^N_i\dd s\Big)
 \sup_s\|\Phi^\pi(s)-\Phi(s)\|^2_{H^{s_1}_{a_1}}.
\]
The inequality is the pointwise reproducing-kernel bound
$|\na\chi(x)|^2\le K_{s_1,a_1}(x)\|\chi\|^2_{H^{s_1}_{a_1}}$ of Lemma~\ref{lem:kernel}
applied under the positive measure
$(2\Theta^N_i+2\varrho^N_{\cdot,i})\dd\mu^N_i$, together with the definition \eqref{eq:trMN} of
$\mathsf q^N_i$. This measure is positive because
$2\Theta^N_i(X^{N,\eta}_{k,i})+2\varrho^N_{k,i}=\Sigma^2_{k,i}\ge2\sigma_i>0$, by \eqref{eq:Sig}, the
identity stated after it, and $f_\eta\ge0$. The first factor is finite, by \eqref{eq:psi2} and
Cauchy--Schwarz, $\E\int_0^T\mathsf q^N_i\le(T\,\E\int_0^T(\mathsf q^N_i)^2)^{1/2}$. At fixed
$N$ and $\eta$ it is also bounded deterministically. Indeed \eqref{eq:psi} bounds
$\mathsf q^N_i\le2C_{s_1,a_1}(\br{\mu^N_i,\Theta^N_i}+\|\varrho^N\|_\infty)$ with
$C_{s_1,a_1}$ the bound of $K_{s_1,a_1}$ in Lemma~\ref{lem:kernel} ($s_1>\frac d2+1$,
$a_1\ge0$), $\Theta^N_i\le\sigma_i+n\|f_\eta\|_{L^\infty(\R)}$ pointwise,
because \ref{H:CH} puts $f_\eta$ in $W^{s+1,\infty}(\R;[0,\infty))$ on all of $\R$, and
$\|\varrho^N\|_\infty\le CN^{-1}\eta^{-(d+\alpha)}$ by \eqref{eq:self}. So
$\mathsf q^N_i$ is bounded by a deterministic constant $C(N,\eta)$ and
$\E\int_0^T\mathsf q^N_i\le TC(N,\eta)$. The second factor tends to
$0$ by uniform continuity of $\Phi$ in $\Test\hookrightarrow H^{s_1}_{a_1}$.

\emph{Step 4: assembling, and in which space the identity holds.} The left-hand
side does not depend on $\pi$, and $\mathcal E_\pi,\mathcal N_\pi$ converge
a.s.; hence $\mathcal M_\pi$ converges a.s.\ as well, and its a.s.\ limit
agrees with its $L^2$-limit. Using $\sum_i\br{\fluc^N_i(0),\Phi_i(0)}
=\br{S^\eta(t,0)\fluc^N(0),\varphi}$ and $\sum_i\int_0^t\br{\dd M^N_i,\Phi_i}
=\br{Z^N(t),\varphi}$ with $Z^N(t)$ as in \eqref{eq:mildzero},
\begin{equation}\label{eq:mildNtested}
 \br{\fluc^N(t),\varphi}=\br{S^\eta(t,0)\fluc^N(0),\varphi}+\br{Z^N(t),\varphi}
 +\int_0^t\br{\mathcal R^N(s),S^\eta(t,s)^*\varphi}\dd s
 \qquad\text{for every }\varphi\in\Test^n.
\end{equation}

By the definition of $D^N(t)$ after \eqref{eq:mildzero}, the last term of \eqref{eq:mildNtested} is
$\br{D^N(t),\varphi}$. Proposition~\ref{prop:evol}
gives $S^\eta$ on $\Ufine$ and on $\Usharp$, and $\mathcal R^N(s)$ lies in
neither: it lies two derivatives below $\Ufine$. But
$\varphi\mapsto\br{\mathcal R^N(s),S^\eta(t,s)^*\varphi}$ is a bounded
functional on $\Test^n$, by
$\|S^\eta(t,s)^*\varphi\|_{H^{s_1+2}_{a_1}}\le C\|S^\eta(t,s)^*\varphi\|_{\Test}
\le CC_1'(T)\|\varphi\|_{\Test}$, so it defines
$S^\eta(t,s)\mathcal R^N(s)\in(\Test^*)^n=(H^{-s_1-3}_{-a_1-1})^n$ with
$\|S^\eta(t,s)\mathcal R^N(s)\|_{\Test^*}\le CC_1'(T)
\|\mathcal R^N(s)\|_{H^{-s_1-2}_{-a_1}}$. The map $s\mapsto S^\eta(t,s)
\mathcal R^N(s)$ is weakly measurable ($\Phi\in C([0,t];\Test^n)$, and $\mathcal R^N$ is strongly
measurable by Lemma~\ref{lem:Rdet}) and $\Test^*$ is separable, hence it is strongly
measurable by Pettis and Bochner-integrable by Lemma~\ref{lem:Rdet}. Its Bochner integral over
$(0,t)$ is $D^N(t)$, since both have the same pairing with every $\varphi\in\Test^n$. With this,
\eqref{eq:mildNtested} gives \eqref{eq:mildN} as an identity in
$(\Test^*)^n$ holding simultaneously for all $\varphi\in\Test^n$. A single
exceptional set suffices because $\Test$ is separable and all four terms are
continuous in $\varphi$ on $\Test$.

Finally, three of the four terms of \eqref{eq:mildN} lie in $\Ufine^n$:
$\fluc^N(t)$ by Lemma~\ref{lem:uetacont}(i), $S^\eta(t,0)\fluc^N(0)$ by
Proposition~\ref{prop:evol}, and $Z^N(t)$ because it is an
$\Ufine$-valued It\^o integral in the sense of Lemma~\ref{lem:BDG}, as in \eqref{eq:mildzero}. We
check the hypotheses of Lemma~\ref{lem:BDG}. The operator $S^\eta(t,r)$ is deterministic, and $r\mapsto S^\eta(t,r)$ is strongly measurable into
$\mathcal L(\Ufine^n)$: for $\nu\in\Ufine^n$ and $\varphi\in\Test^n$,
$r\mapsto\br{S^\eta(t,r)\nu,\varphi}=\br{\nu,\Phi(r)}$ is continuous by
Lemma~\ref{lem:strongcont}, $\Test$ is dense in $\Ufine^*$, and
$\sup_r\|S^\eta(t,r)\|\le C_1(T)$, so $r\mapsto S^\eta(t,r)\nu$ is weakly and
hence (Pettis, $\Ufine$ separable) strongly measurable for each $\nu$. So
$S^\eta(t,\cdot)\Upsilon^N$ is predictable with values in $\HS(\R^{nNd};\Ufine^n)$: its entries are
$N^{-1/2}\Gamma^N_{k,i}(r)\,S^\eta(t,r)\pa_\ell\delta_x$ at $x=X^{N,\eta}_{k,i}(r)$, the map
$(r,x)\mapsto S^\eta(t,r)\pa_\ell\delta_x$ is strongly measurable in $r$ and continuous in $x$,
and $X^{N,\eta}_{k,i}$ and $\Gamma^N_{k,i}$ are continuous and adapted. And
$\E\int_0^t\|S^\eta(t,r)\Upsilon^N(r)\|^2_{\HS}\dd r
\le C_1(T)^2\,\E\int_0^t\sum_i\mathsf q^N_i\dd r<\infty$ by
$\|S\Upsilon^N\|_{\HS}\le\|S\|\,\|\Upsilon^N\|_{\HS}$, \eqref{eq:trMN} and Step~3's
deterministic bound. Pairing with $\varphi$ gives
$\br{Z^N(t),\varphi}=\sum_i\int_0^t\br{\dd M^N_i,\Phi_i}$, by the last part of
Lemma~\ref{lem:BDG}. So the
fourth term lies in $\Ufine^n$ too, by difference. The set $C^\infty_0\subset\Test$ is dense in
$\Ufine^*=H^{s_1}_{a_1}$, because $C^\infty_0$ is dense in $H^{s_1}$ and the isometry $J_{a_1}$ of
Lemma~\ref{lem:spaces} maps $C^\infty_0$ onto itself. Hence the identity \eqref{eq:mildN} holds in
$\Ufine^n$, for each fixed $t$, $\Prob$-a.s.
\end{proof}

\begin{corollary}[The remainder vanishes in $L^2(\Omega\times(0,T))$]
\label{cor:remvanish}
Under the hypotheses of Proposition~\ref{prop:close},
\begin{equation}\label{eq:remvanish}
 \E\int_0^T\big\|\mathcal R^N(s)\big\|^2_{H^{-s_1-2}_{-a_1}}\dd s
 \ \le\ C(T)\Big(N^{-1}\eta^{-2(d+\alpha)}
 \;+\;\overline K_4\,N^{-1}\eta^{-2\alpha}
 \;+\;\overline K_4^{1/2}N^{-1/2}\big(\E\mathsf I\big)^{1/2}\Big)
 \ \longrightarrow\ 0,
\end{equation}
where $\mathsf I=\int_0^T\|\fluc^N\|^2_{\Ufine}\dd t$ and $\E\mathsf I$ is bounded uniformly in
$N$ by Proposition~\ref{prop:close} itself.
\end{corollary}

\begin{proof}
The bound \eqref{eq:remvanish} is \eqref{eq:Rcollect}, which collects three bounds: the self-interaction piece by
Lemma~\ref{lem:remainders}(a) and Lemma~\ref{lem:polylog}, the Taylor piece by
\eqref{eq:Tayl2}, and the bilinear piece by \eqref{eq:bilin}, \eqref{eq:detbound} and
H\"older, which is where $\E\mathsf I$ enters. Each summand tends to $0$: the first by Lemma~\ref{lem:polylog},
the second by \eqref{eq:K4-mild}, since $\overline K_4N^{-1}\eta^{-2\alpha}=o(1)$
is exactly $\overline K_4=o(N\eta^{2\alpha})$, and the third because
\eqref{eq:K4-mild} gives $\overline K_4^{1/2}N^{-1/2}=o(\eta^{\alpha})=o(1)$,
$\alpha\ge0$ and $\eta\le1$.
\end{proof}

\begin{corollary}[The scalar remainders vanish]\label{cor:remscalar}
Under the hypotheses of Proposition~\ref{prop:close}, for every $\varphi\in H^{s_1+2}_{a_1}$
and every $i$, $\E\int_0^T|\mathcal R^N_i(\varphi)|\ds\to0$.
\end{corollary}

\begin{proof}
By Lemma~\ref{lem:remainders}, it suffices to show that the three bounds there tend to $0$. For (a) this follows
from Lemma~\ref{lem:polylog}. For (b), \eqref{eq:cheap2} gives
$\sup_{t\le T}\max_j\E\|B^\eta_{ij}*\fluc^N_j\|^2_{H^{s_1}}\le N^{C_2(T+1)\delta}\times$polylog,
so by Lemma~\ref{lem:polylog} the bound in (b) is
$O\big(N^{-\frac12(1-2C_2(T+1)\delta)}\big)\times$polylog.
This tends to $0$ because $C_2(T+1)\delta<\frac1{C_\star}<\frac12$ by \eqref{eq:smallness}.
For (c), let $\bar{\mathsf I}:=\sup_N\E\int_0^T\|\fluc^N\|^2_{\Ufine}\dd t$, which is finite by
Proposition~\ref{prop:close}. With \eqref{eq:cheap2}, the right side of (c) is at most
\[
 C\,N^{-1/2}\cdot \bar{\mathsf I}^{1/2}\cdot\big(CT\,N^{C_2(T+1)\delta}\big)^{1/2}
 \ =\ O\big(N^{-\frac12(1-C_2(T+1)\delta)}\big)\times\text{polylog}\ \to\ 0,
\]
since $C_2(T+1)\delta<1$.
\end{proof}

\section{Tightness}\label{sec:tight}

We prove tightness by the Aldous--Jakubowski method on the triple
\eqref{eq:chainA}. The criterion we use is Proposition~\ref{prop:criterion} of
Appendix~\ref{sec:appaldous}. It requires uniform bounds in the two finest
spaces, namely $\sup_t$ in $\Usharp$ and time-integrated in $\Ufine$ for the field,
and $\sup_t$ in $\Ufine$ for the martingale. It also requires the time-increment
\textup{(}Aldous\textup{)} condition in the coarsest space $\Vcoarse$ for
the field and in the finest space $\Ufine$ for the martingale. The compact
embeddings $\Usharp\hookrightarrow\hookrightarrow\Hmid
\hookrightarrow\hookrightarrow\Vcoarse$ turn these bounds into tightness on the
path space $\mathcal Z_T$. The continuous factor of $\mathcal Z_T$ is the middle space $\Hmid$, and
convergence in law is stated in $\Hmid$. The Aldous condition for the field follows
from the weak form \eqref{eq:weakN}, the uniform bound in $\Usharp$, and the Aldous
condition for the martingale. The proof of Proposition~\ref{prop:tight} matches each
ingredient to a hypothesis of the criterion,
and the section closes with Proposition~\ref{prop:skorokhod}: tightness together
with the separating family of Appendix~\ref{sec:appaldous} produces almost-sure
representatives $(\widetilde\fluc,\widetilde M)$ on a new probability space, and
everything from Section~\ref{sec:identify} onwards works with them.

\begin{lemma}[Martingale bounds]\label{lem:Mbounds}
For each $i$
the process $M^N_i$ is a square-integrable $\Ufine$-valued continuous martingale
(Lemma~\ref{lem:MNval}),
with
\begin{align}
 \sup_N\E\Big[\sup_{t\le T}\|M^N_i(t)\|_{\Ufine}^2\Big]&\le C,\label{eq:M2}\\
 \sup_N\ \sup_{\tau\le T}\ \E\Big[\big\|M^N_i(t_\vartheta)-M^N_i(\tau)\big\|_{\Ufine}^2\Big]
 &\le C\,\vartheta^{1/2},\quad t_\vartheta:=(\tau+\vartheta)\wedge T,\ 0<\vartheta\le1,\label{eq:Maldous}
\end{align}
the second for all $\F$-stopping times $\tau\le T$, with $C$ independent of
$N,\eta$.
\end{lemma}

\begin{proof}
We prove \eqref{eq:M2} by Doob's inequality \cite[Thm.~1.3.8(iv)]{KaSh91} for the coordinates of $M^N_i$ in an
orthonormal basis, and sum their brackets with the reproducing kernel of
Lemma~\ref{lem:kernel}.
Let $(e_m)_{m\in\N}$ be an orthonormal basis of $H^{s_1}_{a_1}$, so that
$\|\nu\|^2_{\Ufine}=\sum_m|\br{\nu,e_m}|^2$ for $\nu\in\Ufine$. By Doob's
inequality for the real martingales $M^N_i(\cdot,e_m)$ and monotone convergence,
\[
 \E\Big[\sup_{t\le T}\|M^N_i(t)\|^2_{\Ufine}\Big]
 \le\sum_{m}\E\Big[\sup_{t\le T}|M^N_i(t,e_m)|^2\Big]
 \le4\sum_m\E\br{M^N_i(\cdot,e_m)}_T.
\]
By \eqref{eq:cov_finiteN} and Tonelli,
\[
 \sum_m\E\br{M^N_i(\cdot,e_m)}_T
 =\E\int_0^T\br{\mu^N_i,\ \Big(\sum_m|\na e_m|^2\Big)
 \big(2\Theta^N_i+2\varrho^N_{\cdot,i}\big)}\ds
 =\E\int_0^T\br{\mu^N_i,\ K_{s_1,a_1}\big(2\Theta^N_i+2\varrho^N_{\cdot,i}\big)}\ds,
\]
where $K_{s_1,a_1}$ is the reproducing kernel of Lemma~\ref{lem:kernel}. Since
$s_1>\frac d2+1$, $K_{s_1,a_1}\le C_{s_1,a_1}$ pointwise, so
\[
 \sum_m\E\br{M^N_i(\cdot,e_m)}_T
 \le 2C_{s_1,a_1}\E\int_0^T\br{\mu^N_i,\Theta^N_i}\ds
 +2C_{s_1,a_1}\,T\,\sup_{k}\|\varrho^N_{k,i}\|_\infty,
\]
which is bounded uniformly in $N$ by Lemma~\ref{lem:L2field} and the Cauchy--Schwarz inequality in $(\omega,s)$,
and by \eqref{eq:self} together with Lemma~\ref{lem:polylog}. This proves
\eqref{eq:M2}. \eqref{eq:Maldous} is Corollary~\ref{cor:Maldous}, \eqref{eq:Maldousfix},
which uses only \eqref{eq:psi2}.

\end{proof}

\begin{corollary}[The Aldous clause of Lemma~\ref{lem:Mbounds}]
\label{cor:Maldous}
For every stopping time $\tau\le T$ and $\vartheta\in(0,1]$,
with $t_\vartheta:=(\tau+\vartheta)\wedge T$,
\begin{equation}\label{eq:Maldousfix}
 \sup_{N}\ \E\Big[\big\|M^N_i(t_\vartheta)-M^N_i(\tau)\big\|^2_{\Ufine}\Big]
 \ \le\ C\,\vartheta^{1/2}.
\end{equation}
\end{corollary}

\begin{proof}
Let $(e_m)$ be any orthonormal basis of $H^{s_1}_{a_1}$. Each real martingale
$M^N_i(\cdot,e_m)$ is square-integrable, since
$\sum_m\E\br{M^N_i(\cdot,e_m)}_T=\E\tr\br{\!\br{M^N_i}\!}_T=\E\int_0^T\mathsf q^N_i<\infty$,
the identity by Corollary~\ref{cor:trace}, the finiteness by \eqref{eq:psi2}. So its compensated square is a uniformly integrable martingale
on $[0,T]$ and Doob's optional sampling applies at the bounded stopping times
$\tau\le t_\vartheta\le T$. By that and monotone convergence, $\E\|M^N_i(t_\vartheta)-M^N_i(\tau)\|^2_{\Ufine}
=\E\int_\tau^{t_\vartheta}\mathsf q^N_i\dd s
=\E\int_0^T\mathsf q^N_i\mathbf 1_{[\tau,t_\vartheta]}\dd s$. Cauchy--Schwarz on
$\Omega\times(0,T)$ and \eqref{eq:psi2} give
$\le(\E\int_0^T(\mathsf q^N_i)^2)^{1/2}\,\vartheta^{1/2}\le C\vartheta^{1/2}$.
\end{proof}

\begin{theorem}[Aldous condition for $\fluc^N$]\label{thm:aldous}
There is $C=C(T)$ such that, with $N_0$ as in \eqref{eq:N0},
\begin{equation}\label{eq:aldousfull}
 \sup_{N\ge N_0}\ \sup_{\tau}\ \sup_{\vartheta\in(0,1]}\ \vartheta^{-1/2}\;
 \E\big\|\fluc^N(t_\vartheta)-\fluc^N(\tau)\big\|^2_{\Vcoarse}\ \le\ C,
\end{equation}
the middle supremum over stopping times $\tau\le T$, where $t_\vartheta:=(\tau+\vartheta)\wedge T$.
\end{theorem}

\begin{proof}
We subtract the weak form \eqref{eq:weakN} at the two times $\tau$ and $t_\vartheta$, bound
the drift, remainder and martingale increments against a dense set of test functions, and
take the supremum over the test functions to obtain the $\Fcoarse$ norm.
By Theorem~\ref{thm:K4}, \eqref{eq:K4-mild} holds. So
Proposition~\ref{prop:close} and Corollary~\ref{cor:remvanish} give, uniformly in
$N\ge N_0$,
\[
 \E\sup_{t\le T}\|\fluc^N(t)\|^2_{\Usharp}\le C,\qquad
 \E\int_0^T\|\mathcal R^N(r)\|^2_{H^{-s_1-2}_{-a_1}}\dd r\le C.
\]
Let $\{\varphi_k\}$ be a countable dense subset of the unit ball of $\Test$. Off one null
set, \eqref{eq:weakN} holds for every $\varphi_k$, every $i$ and every $t$. All its terms are continuous in $t$, so it may be
evaluated at $\tau$ and at $t_\vartheta$. Subtracting,
\[
 \br{\fluc^N_i(t_\vartheta)-\fluc^N_i(\tau),\varphi_k}
 =\int_\tau^{t_\vartheta}\sum_j\br{\fluc^N_j(r),(\mathcal G^\eta(r)\varphi_k)_{ij}}\dd r
 +\int_\tau^{t_\vartheta}\mathcal R^N_i(r,\varphi_k)\dd r
 +\br{M^N_i(t_\vartheta)-M^N_i(\tau),\varphi_k}.
\]
By Lemma~\ref{lem:Gbound}, $\mathcal G^\eta(r)$ maps $\Test$ into
$H^{s_1+1}_{a_1}=(\Usharp)^*$ with norm at most $C_{\mathcal G}$, uniformly in $r$ and
$\eta$. Moreover $\Test\hookrightarrow H^{s_1+2}_{a_1}\hookrightarrow H^{s_1}_{a_1}$. With
Cauchy--Schwarz in $r$ and $t_\vartheta-\tau\le\vartheta$, the three terms are bounded by
\[
 nC_{\mathcal G}\,\vartheta\sup_{r\le T}\|\fluc^N(r)\|_{\Usharp}
 +c\,\vartheta^{1/2}\Big(\int_0^T\|\mathcal R^N\|^2_{H^{-s_1-2}_{-a_1}}\dd r\Big)^{1/2}
 +c\,\big\|M^N_i(t_\vartheta)-M^N_i(\tau)\big\|_{\Ufine}.
\]
The supremum over $k$ is the norm of $\fluc^N_i(t_\vartheta)-\fluc^N_i(\tau)$ in
$\Test^*=\Fcoarse$, and $\Fcoarse\hookrightarrow\Vcoarse$ with norm $c_{\rm FV}$. The left
side is measurable, since $\fluc^N$ is continuous in $\Ufine$ and $\tau,t_\vartheta$ are
random variables. Square, take expectations, and use the two bounds above and
Corollary~\ref{cor:Maldous}:
\[
 \E\big\|\fluc^N(t_\vartheta)-\fluc^N(\tau)\big\|^2_{\Vcoarse}
 \le C\big(\vartheta^2+\vartheta+\vartheta^{1/2}\big)\le C\vartheta^{1/2}.
\]
\end{proof}

\begin{proposition}[Tightness]\label{prop:tight}
The laws of $(\fluc^N,M^N)$ are
tight on
\[
 \mathcal Z_T:=\Big(C([0,T];\Hmid)\cap L^2_w(0,T;\Ufine)\Big)^n
 \times C([0,T];\Hmid)^n,
\]
where $L^2_w$ denotes the weak topology.
\end{proposition}

\begin{proof}
We apply the second conclusion of Proposition~\ref{prop:criterion} to
$\xi^N:=\fluc^N$ and $m^N:=M^N$. Both have continuous $\Ufine^n$-valued paths, by
Lemma~\ref{lem:uetacont}(ii) and Lemma~\ref{lem:MNval}. In point~(a), the two
bounds on $\fluc^N$ are Proposition~\ref{prop:close}, whose hypothesis
\eqref{eq:K4-mild} is Theorem~\ref{thm:K4}, and the bound on $M^N$ is
\eqref{eq:M2} of Lemma~\ref{lem:Mbounds}. In point~(b), the bound on $\fluc^N$ is
Theorem~\ref{thm:aldous}, and the bound on $M^N$ is Corollary~\ref{cor:Maldous}.
\end{proof}

\begin{proposition}[Skorokhod representation]\label{prop:skorokhod}
There exist a subsequence (not relabelled), a
probability space $(\widetilde\Omega,\widetilde{\mathcal F},\widetilde\Prob)$ and
$\mathcal Z_T$-valued random variables
$(\widetilde\fluc^N,\widetilde M^N)$ and
$(\widetilde\fluc,\widetilde M)$ on it, with
$\mathrm{Law}(\widetilde\fluc^N,\widetilde M^N)
=\mathrm{Law}(\fluc^N,M^N)$, such that $\widetilde\Prob$-a.s.
\[
 \widetilde\fluc^N\to\widetilde\fluc,\quad
 \widetilde M^N\to\widetilde M\ \text{in }C([0,T];\Hmid)^n,\quad
 \widetilde\fluc^N\rightharpoonup\widetilde\fluc\ \text{in }L^2(0,T;\Ufine)^n.
\]
Here $\mathrm{Law}(\fluc^N,M^N)$ denotes the Radon Borel extension of
Lemma~\ref{lem:borel}. It agrees with the law of $(\fluc^N,M^N)$ on the Borel
$\sigma$-algebra of each $C([0,T];\Hmid)$ block.
Consequently, for every Borel function $F:\mathcal Z_T\to\R$, the random variables
$F(\fluc^N,M^N)$ and $F(\widetilde\fluc^N,\widetilde M^N)$ have the same law, after
completion of the probability spaces.
\end{proposition}

\begin{proof}
By Lemma~\ref{lem:separating}, $\mathcal Z_T$ carries a countable family $(f_j)$ of real
continuous functions separating its points. By the last conclusion of Proposition~\ref{prop:criterion}, the laws
$\mathrm{Law}(\fluc^N,M^N)$ are Radon Borel probability measures on $\mathcal Z_T$. With
tightness, Theorem~\ref{thm:SJ}, Jakubowski's version of the Skorokhod representation theorem \cite{Jak97} in the form of \cite[Thm.~C.1]{BrOn10}, applies. The Radon Borel
extension of Lemma~\ref{lem:borel} agrees with the law on $\sigma(f_j)$. Since
the $C([0,T];\Hmid)$ blocks are Polish, it therefore agrees with the law on each of
their Borel $\sigma$-algebras.
For the last claim, fix $N$ and put $\mathbf f:=(f_j)_{j\in\N}:\mathcal Z_T\to\R^{\N}$, a continuous
injective map. For a compact set $K\subset \mathcal Z_T$, the restriction $\mathbf f|_K$ is a continuous
bijection from a compact space onto a Hausdorff space, hence a homeomorphism onto the compact
set $\mathbf f(K)$. So, for every Borel set $B\subset \mathcal Z_T$, $\mathbf f(B\cap K)$ is Borel in
$\R^{\N}$, and $B\cap K=\mathbf f^{-1}\big(\mathbf f(B\cap K)\big)\in\sigma(f_j)$ by injectivity.
By Proposition~\ref{prop:tight}, there are compact sets $K_m\subset \mathcal Z_T$ with
$\Prob((\fluc^N,M^N)\in K_m)\ge1-\frac1m$. Put $K_\infty:=\bigcup_mK_m$. Then
$B\cap K_\infty\in\sigma(f_j)$ for every Borel set $B$. The set $K_\infty$ has full measure for
the law of $(\fluc^N,M^N)$, and hence also for the law of $(\widetilde\fluc^N,\widetilde M^N)$,
since the two laws agree on $\sigma(f_j)$. So $F^{-1}(A)\cap K_\infty\in\sigma(f_j)$ for every
Borel set $A\subset\R$, and the two laws of $F$ agree, after completion of the probability spaces.
\end{proof}

\section{Identification of the drift and of the bracket}\label{sec:identify}

\noindent
Two limits remain in \eqref{eq:fluc_decomp}. The drift term carries
$f'_\eta(B^\eta_{ij}*u_{\eta,j})$ against $\fluc^N_j$ and must converge to the
linearised coefficient $\mathbf A_{ij}(u)$. The bracket \eqref{eq:cov_finiteN} carries the
empirical measure and the empirical mobility and must converge to
\eqref{eq:cov_limit}. Both are limits of products in which one factor converges
only in law, and they are taken on different spaces. The bracket is identified on
the original probability space, where the coupling with the intermediate
system is available. The drift limit is taken
after the transfer, on the new space and along the subsequence of
Proposition~\ref{prop:skorokhod}. Lemma~\ref{lem:drift} is a convergence in
$L^1(\widetilde\Omega)$, because the product converges only for almost-sure
representatives of the field.
Theorem~\ref{thm:noise} uses the bracket identification to identify the limit martingale: a
continuous square-integrable $\Hmid$-valued martingale whose bracket is the
deterministic \eqref{eq:cov_limit}. Section~\ref{sec:limitmart} studies this limit
martingale rather than the particle system.

\begin{proposition}[Identification of the limiting bracket]\label{prop:bracketid}
Let $N\ge N_0$, with $N_0$ as in \eqref{eq:N0}, and let
$\varphi\in H^{s_m+1}_{a_m}$, $g:=|\na\varphi|^2$. Then, for every $t\le T$,
\begin{equation}\label{eq:bracketid}
 \E\Big|\int_0^t\br{\mu^N_i,\,2g\textstyle\sum_jf_\eta(B^\eta_{ij}*\mu^N_j)}\dr
 -\int_0^t\br{u_i,\,2g\textstyle\sum_jf(a_{ij}u_j)}\dr\Big|
 \ \longrightarrow\ 0.
\end{equation}
\end{proposition}

\begin{proof}
Insert the intermediate system and split into three differences. Step~1 replaces
$\mu^N$ by the intermediate empirical measure $\bar\mu^N$, using the coupling estimate
\eqref{eq:phi2}. Step~2 replaces $\bar\mu^N$ by $u_\eta$, by a law of large numbers. Step~3
replaces $u_\eta$ by $u$, using \eqref{eq:CHbias}. In each
difference the outer measure and the argument of $f_\eta$ are changed
separately. Throughout, $\|g\|_{L^\infty}+\|\na g\|_{L^\infty}<\infty$ by
Lemma~\ref{lem:spaces}(iii) and $s_m>\frac d2+1$, and $\|g\|_{L^2}<\infty$
because $H^{s_m+1}_{a_m}\hookrightarrow H^1$ by Lemma~\ref{lem:spaces}(ii)
($a_m\ge0$), so $\|g\|_{L^2}^2=\int|\na\varphi|^4\le\|\na\varphi\|_{L^\infty}^2
\|\na\varphi\|_{L^2}^2$; $\mu^N_i$, $\bar\mu^N_i$ are probability measures and $u_{\eta,i}(t)$ and
$u_i(t)$ are probability densities by Lemma~\ref{lem:lawid}(ii), at $\eta$ and at
$\eta=0$ and $t\le T$.

\emph{Step 1: $\mu^N\rightsquigarrow\bar\mu^N$, by the coupling.}
For the argument, $\mathrm{Lip}(f_\eta)\le\eta^{-\alpha}$ and the mean value
theorem give, pathwise,
\[
 \big|f_\eta(B^\eta_{ij}*\mu^N_j)(x)-f_\eta(B^\eta_{ij}*\bar\mu^N_j)(x)\big|
 \le\eta^{-\alpha}\|\na B^\eta_{ij}\|_{L^\infty}\tfrac1N\sum_\ell|Y_{\ell,j}|
 \le C_L\Lambda^{1/2}\tfrac1N\sum_\ell|Y_{\ell,j}|.
\]
For the outer measure, $|\br{\mu^N_i-\bar\mu^N_i,h}|\le\|\na
h\|_{L^\infty}\frac1N\sum_k|Y_{k,i}|$ with
$h:=2g\sum_jf_\eta(B^\eta_{ij}*\bar\mu^N_j)$, and $h$ is random, but its
gradient is bounded by a deterministic constant: $|f_\eta(z)|\le|f_\eta(0)|+\eta^{-\alpha}|z|$ with
$f_\eta(0)=f(0)$ for every $\eta$, since $0\in I\subseteq[-a_\eta,a_\eta]$
by \ref{H:cutoff}, and $\|B^\eta_{ij}*\bar\mu^N_j\|_{L^\infty}\le
\|B^\eta_{ij}\|_{L^\infty}\le C\eta^{-d}$ give $\|f_\eta(B^\eta_{ij}*\bar\mu^N_j)\|_{L^\infty}\le
C\eta^{-d-\alpha}$, while $\|f'_\eta\|_{L^\infty}\le\eta^{-\alpha}$ and
$\|\na(B^\eta_{ij}*\bar\mu^N_j)\|_{L^\infty}\le C\eta^{-d-1}$ give
$\|\na h\|_{L^\infty}\le C\eta^{-d-1-\alpha}=C\Lambda^{1/2}$.
Since $\E\frac1N\sum_k|Y_{k,i}|\le\max_k(\E|Y_{k,i}|^2)^{1/2}\le\sup_{r\le
T}\phi_2(r)^{1/2}$, both pieces are at most
$C\,T\,n\,\Lambda^{1/2}\sup_{r\le T}\phi_2(r)^{1/2}$, and by \eqref{eq:phi2},
\[
 \Lambda^{1/2}\phi_2^{1/2}
 \ \le\ C\,\eta^{-(d+1+\alpha)}\,N^{-\frac12+\frac{\eps_\delta}2}
 \eta^{-\frac{d+2\alpha}2}
 \ =\ N^{-\frac12+\frac{\eps_\delta}2}\times\text{polylog}\ \to\ 0,
\]
because the powers of $\eta^{-1}$ are powers of $\delta\log N$ by \eqref{eq:H3}, as in
Lemma~\ref{lem:polylog}, and $N^{-\frac12+\frac{\eps_\delta}2}$ decays polynomially, since
$\eps_\delta<1$ by \eqref{eq:smallness}.

\emph{Step 2: $\bar\mu^N\rightsquigarrow u_\eta$, a law of large numbers.}
By Lemma~\ref{H:chaos} the $\bar X^\eta_{\ell,j}(r)$ are i.i.d.\ with density
$u_{\eta,j}(r)$. For the argument, with
$\bar\fluc^N_j=\sqrt N(\bar\mu^N_j-u_{\eta,j})$ and $\sigma>\frac d2$ fixed,
Sobolev and Lemma~\ref{lem:smoothed} give
\[
 \eta^{-\alpha}\,\E\big\|B^\eta_{ij}*(\bar\mu^N_j-u_{\eta,j})\big\|_{L^\infty}
 \le C_\sigma\eta^{-\alpha}N^{-1/2}
 \big(\E\|B^\eta_{ij}*\bar\fluc^N_j\|^2_{H^\sigma}\big)^{1/2}
 \le C\,N^{-1/2}\eta^{-\alpha-\frac{d+2\sigma}2}\ \to\ 0,
\]
again by Lemma~\ref{lem:polylog}.
For the outer measure, $h_\eta:=2g\sum_jf_\eta(B^\eta_{ij}*u_{\eta,j})$ is
deterministic, and its argument lies in $I$ because
$\|B^\eta_{ij}*u_{\eta,j}\|_{L^\infty}\le\|B^\eta_{ij}\|_{L^1}\|u_{\eta,j}\|_{L^\infty}
\le\mathsf b_0$ by \ref{H:CH} and \ref{H:small}, so $f_\eta=f$ there by
\ref{H:cutoff} and $\|h_\eta\|_{L^\infty}\le2n\|g\|_{L^\infty}\sup_If$. Hence
$\E|\br{\bar\mu^N_i-u_{\eta,i},h_\eta}|
\le N^{-1/2}\|h_\eta\|_{L^\infty}\to0$.

\emph{Step 3: $u_\eta\rightsquigarrow u$.} By Lemma~\ref{lem:lawid}(i), at $\eta$ and at
$\eta=0$, $w^\eta_{ij}$ and $a_{ij}u_j$ take values in $[-\mathsf b_0,\mathsf b_0]\subset I$.
On $I$, $f_\eta=f$ by \ref{H:cutoff} and $\|f'\|_{L^\infty(I)}\le\eps$ by \ref{H:small}. Moreover
$w^\eta_{ij}-a_{ij}u_j=B^\eta_{ij}*(u_{\eta,j}-u_j)+(B^\eta_{ij}*u_j-a_{ij}u_j)$. The first
summand has $L^2$ norm at most $A_0\|u_\eta-u\|_{L^2}\le C\eta$, by \eqref{eq:CHbias}. The
second has $L^2$ norm at most $A_0\eta\|\na u_j\|_{L^2}\le A_0\eta\|u_0\|_{H^s}$, by
\eqref{eq:mollerr} and Lemma~\ref{lem:lawid}(i). So
$\sup_{r\le T}\|w^\eta_{ij}(r)-a_{ij}u_j(r)\|_{L^2}\le C\eta$. For the argument,
\[
 \Big|\br{u_{\eta,i},\,2g\textstyle\sum_j\big(f_\eta(w^\eta_{ij})-f(a_{ij}u_j)\big)}\Big|
 \le2n\eps\,\|g\|_{L^\infty}\|u_{\eta,i}\|_{L^2}\max_j\|w^\eta_{ij}-a_{ij}u_j\|_{L^2}
 \le C\eta,
\]
using $\|u_{\eta,i}\|_{L^2}\le\|u_\eta\|_{H^s}\le\|u_0\|_{H^s}$ (Lemma~\ref{lem:lawid}(i)). For the outer measure,
$h:=2g\sum_jf(a_{ij}u_j)$ is bounded and $g\in L^2$, so
$|\br{u_{\eta,i}-u_i,h}|\le2n\sup_If\,\|g\|_{L^2}\|u_\eta-u\|_{L^\infty(0,T;L^2)}\le C\eta$
by \eqref{eq:CHbias}. All three differences vanish, and \eqref{eq:bracketid} follows.
\end{proof}

\begin{lemma}[Drift limit]\label{lem:drift}
Let $(\widetilde\fluc,\widetilde M)$ and the subsequence be those of
Proposition~\ref{prop:skorokhod}.
For $\varphi\in C^\infty_0(\R^d)$ and $t\in[0,T]$ put
\begin{equation}\label{eq:Rtilde}
 \widetilde{\mathcal E}^N_i(t,\varphi):=\br{\widetilde\fluc^N_i(t),\varphi}-\br{\widetilde\fluc^N_i(0),\varphi}
 -\int_0^t\sum_j\br{\widetilde\fluc^N_j,(\mathcal G^\eta\varphi)_{ij}}\ds-\widetilde M^N_i(t,\varphi).
\end{equation}
Then, in $L^1(\widetilde\Omega)$, hence in probability, and $\widetilde\Prob$-a.s.\ along a
further subsequence,
\[
 \int_0^t\sum_j\br{\widetilde\fluc_j^{N},(\mathcal G^\eta\varphi)_{ij}}\ds
 +\widetilde{\mathcal E}^N_i(t,\varphi)
 \ \longrightarrow\
 \int_0^t\sum_j\br{\widetilde\fluc_j,(\mathcal G\varphi)_{ij}}\ds.
\]
\end{lemma}

\begin{proof}
We show that the remainder $\widetilde{\mathcal E}^N_i(t,\varphi)$ tends to $0$ in
$L^1(\widetilde\Omega)$, by transferring Corollary~\ref{cor:remscalar} to the new space through
the law identity of Proposition~\ref{prop:skorokhod}. Then we split the drift term into a term
that converges by the weak convergence of $\widetilde\fluc^N$ and a term that vanishes by
Lemma~\ref{lem:Gupgrade}.

\emph{The remainder.} Let $F:\mathcal Z_T\to\R$ send $(v,m)$ to
$\br{v_i(t),\varphi}-\br{v_i(0),\varphi}-\int_0^t\sum_j\br{v_j,(\mathcal G^\eta\varphi)_{ij}}\ds-m_i(t,\varphi)$.
It is continuous. The first, second and last terms are continuous on
$C([0,T];\Hmid)$, since $\varphi\in\Hmid^*$. The third pairs $v_j$ with
$\mathbf 1_{[0,t]}(\mathcal G^\eta\varphi)_{ij}$, an element of
$L^2(0,T;H^{s_1}_{a_1})=L^2(0,T;\Ufine)^*$ by Lemma~\ref{lem:coeffcont}(ii), since
$\varphi$ has compact support. So it is
continuous for the weak topology of $L^2(0,T;\Ufine)$. On the original space,
$F(\fluc^N,M^N)=\int_0^t\mathcal R^N_i(s,\varphi)\ds$ almost surely, by
\eqref{eq:fluc_decomp} and $\varphi\in H^{s_1+2}_{a_1}$. By
Proposition~\ref{prop:skorokhod}, $\widetilde{\mathcal E}^N_i(t,\varphi)=F(\widetilde\fluc^N,\widetilde M^N)$
has the same law. Hence
$\widetilde\E|\widetilde{\mathcal E}^N_i(t,\varphi)|\le\E\int_0^T|\mathcal R^N_i(\varphi)|\ds\to0$
by Corollary~\ref{cor:remscalar}, whose hypothesis \eqref{eq:K4-mild} holds by Theorem~\ref{thm:K4}. For the linear term, we split
\[
 \br{\widetilde\fluc^N_j,(\mathcal G^\eta\varphi)_{ij}}
 =\br{\widetilde\fluc^N_j,(\mathcal G\varphi)_{ij}}
 +\br{\widetilde\fluc^N_j,(\mathcal G^\eta\varphi)_{ij}-(\mathcal G\varphi)_{ij}}.
\]
The first term converges by the weak convergence
$\widetilde\fluc^N_j\rightharpoonup\widetilde\fluc_j$ in $L^2(0,T;\Ufine)$
tested against the \emph{fixed, deterministic} element
$(\mathcal G\varphi)_{ij}\in L^2(0,T;H^{s_1}_{a_1})$.
The second is bounded by
$\|\widetilde\fluc^N_j\|_{L^2(0,T;\Ufine)}\|(\mathcal G^\eta-\mathcal
G)\varphi\|_{L^2(0,T;H^{s_1}_{a_1})}$, the first factor bounded by
Proposition~\ref{prop:close} (the $L^2(0,T;\Ufine)$ bound, transferred to the new space through the law identity) and the second tending to $0$ by
Lemma~\ref{lem:Gupgrade}. The three terms converge as follows.
\begin{itemize}
\item \emph{Remainder}: $L^1(\widetilde\Omega)$, as shown above.
\item \emph{First (weak) term}: $\widetilde\Prob$-a.s., because
Proposition~\ref{prop:skorokhod} delivers
$\widetilde\fluc^N_j\rightharpoonup\widetilde\fluc_j$ in $L^2(0,T;\Ufine)$
a.s.; and also in $L^1(\widetilde\Omega)$, since
$\sup_N\E\|\widetilde\fluc^N\|^2_{L^2(0,T;\Ufine)}<\infty$ makes the sequence
$\int_0^t\br{\widetilde\fluc^N_j,(\mathcal G\varphi)_{ij}}\ds$ bounded in
$L^2(\widetilde\Omega)$, hence uniformly integrable.
\item \emph{Second term}: its second factor
$\|(\mathcal G^\eta-\mathcal G)\varphi\|_{L^2(0,T;H^{s_1}_{a_1})}$ is
deterministic and tends to $0$. The first
is bounded in $L^2(\widetilde\Omega)$. So the product tends to $0$ in
$L^2(\widetilde\Omega)$, hence in $L^1(\widetilde\Omega)$.
\end{itemize}
Adding the three, the convergence holds in $L^1(\widetilde\Omega)$, hence in
probability and, along a further subsequence, $\widetilde\Prob$-a.s.
\end{proof}

\begin{lemma}[Continuous martingales with deterministic brackets]\label{lem:gaussmart}
Let $(\mathcal F_t)_{t\in[0,T]}$ be a filtration, and let $\mathsf X=(\mathsf X^1,\dots,\mathsf X^p)$ be a continuous
$\R^p$-valued $(\mathcal F_t)$-martingale whose bracket matrix
$c(t):=(\br{\mathsf X^q,\mathsf X^{q'}}_t)_{q,q'\le p}$ is deterministic. Then, for $0\le r\le t\le T$, the
increment $\mathsf X(t)-\mathsf X(r)$ is a centred Gaussian vector with covariance $c(t)-c(r)$, independent of
$\mathcal F_r$. Consequently, for $0\le t_0<t_1<\dots<t_L\le T$, the increments
$\mathsf X(t_l)-\mathsf X(t_{l-1})$, $l=1,\dots,L$, are independent of each other and of $\mathcal F_{t_0}$.
\end{lemma}

\begin{proof}
The process $\mathsf X$ remains a continuous martingale, with the same brackets, for the
filtration $(\mathcal F_t)$ augmented by the null sets and made right-continuous, since its
paths are continuous. We apply It\^o's formula in that filtration. The bracket matrix $c$ is
continuous, as the bracket of a continuous martingale. For $\lambda\in\R^p$, It\^o's formula
shows that $\exp\big(\mathrm i\lambda\cdot \mathsf X(t)+\frac12\lambda^{\!\top}c(t)\lambda\big)$ is a
continuous local martingale. It is bounded by $\exp(\frac12\lambda^{\!\top}c(T)\lambda)$, hence a
martingale. So
\[
 \E\big[e^{\mathrm i\lambda\cdot(\mathsf X(t)-\mathsf X(r))}\,\big|\,\mathcal F_r\big]
 =e^{-\frac12\lambda^{\!\top}(c(t)-c(r))\lambda}
 \qquad(0\le r\le t\le T).
\]
The conditional characteristic function is deterministic. Hence $\mathsf X(t)-\mathsf X(r)$ is centred Gaussian
with covariance $c(t)-c(r)$ and independent of $\mathcal F_r$. For the consequence, let
$A\in\mathcal F_{t_0}$ and $\lambda_1,\dots,\lambda_L\in\R^p$. Conditioning successively on
$\mathcal F_{t_{L-1}},\dots,\mathcal F_{t_0}$ gives
\[
 \E\Big[\mathbf 1_A\,e^{\mathrm i\sum_l\lambda_l\cdot(\mathsf X(t_l)-\mathsf X(t_{l-1}))}\Big]
 =\Prob(A)\prod_{l=1}^Le^{-\frac12\lambda_l^{\!\top}(c(t_l)-c(t_{l-1}))\lambda_l},
\]
which is the asserted independence.
\end{proof}

\begin{theorem}[Limiting martingale and its covariance]\label{thm:noise}
Along the
subsequence of Proposition~\ref{prop:skorokhod},
$\widetilde M_i$ is a continuous, square-integrable, $\Hmid$-valued martingale
with $\widetilde M_i(0)=0$, with respect to
the filtration generated by $(\widetilde\fluc,\widetilde M)$, and its scalar predictable brackets
$\br{\widetilde M_i(\cdot,\varphi),\widetilde M_k(\cdot,\psi)}_t$ are
deterministic,
\begin{equation}\label{eq:cov_limit}
 \big\langle\widetilde M_i(\cdot,\varphi),\widetilde M_k(\cdot,\psi)\big\rangle_t
 =\delta_{ik}\int_0^{t}\br{u_i,\ \na\varphi\cdot\na\psi\,
 \Big(2\sigma_i+2\sum_{j=1}^nf(a_{ij}u_j)\Big)}\dr,
 \qquad\varphi,\psi\in H^{s_m+1}_{a_m},
\end{equation}

and consequently $\widetilde M=(\widetilde M_1,\dots,\widetilde M_n)$ is a
\emph{Gaussian} martingale with independent increments, with the
species-diagonal covariance \eqref{eq:cov_limit}. Moreover, $\widetilde M$ is independent of
$\widetilde{\mathcal F}_0$, in particular of $\widetilde\fluc(0)$.
\end{theorem}

\begin{proof}
Step~1 identifies the limit of the finite-$N$ brackets, by
Proposition~\ref{prop:bracketid}. Step~2 passes to the limit in the martingale identities of
the scalar pairings, on the space of Proposition~\ref{prop:skorokhod}. Step~3 lifts the scalar
martingales to an $\Hmid$-valued martingale, and Step~4 deduces Gaussianity from
Lemma~\ref{lem:gaussmart}.

\emph{Step 1: identification of the limit of \eqref{eq:cov_finiteN}.} Fix $i$ and
$\varphi\in H^{s_m+1}_{a_m}\subset\Hmid^*$. Denote by $b^N(t):=\br{M^N_i(\cdot,\varphi)}_t$ the
bracket \eqref{eq:cov_finiteN}, and by $b(t)$ the right side of \eqref{eq:cov_limit} at $k=i$,
$\psi=\varphi$. We show that $\E|b^N(t)-b(t)|\to0$ for every $t\le T$. We split
\[
 \br{M^N_i(\cdot,\varphi)}_t
 =\underbrace{\int_0^t\br{\mu^N_i,|\na\varphi|^2\,2\sigma_i}\dr}_{=:J^N_1}
 +\underbrace{\int_0^t\br{\mu^N_i,|\na\varphi|^2\,2\textstyle\sum_jf_\eta(B^\eta_{ij}*\mu^N_j)}\dr}_{=:J^N_2}
 +\underbrace{\int_0^t\br{\mu^N_i,|\na\varphi|^2\,2\varrho^N_{\cdot,i}}\dr}_{=:J^N_3}.
\]
$|J^N_3|\le CTN^{-1}\eta^{-d-\alpha}\|\na\varphi\|^2_{L^\infty}\to0$ by \eqref{eq:self}
and Lemma~\ref{lem:polylog}. For $J^N_1$, let $g:=2\sigma_i|\na\varphi|^2$. Then $g\in H^{s_m}_{a_m}$ by the
algebra property, since $\varphi\in H^{s_m+1}_{a_m}$ and $s_m>\frac d2+1$, so
$g\in H^{s_1}_{a_1}=\Ufine^*$ and $g\in L^2$. Since $\mu^N_i=u_{\eta,i}+N^{-1/2}\fluc^N_i$,
\[
 J^N_1-\int_0^t\br{u_i,g}\dr=N^{-1/2}\int_0^t\br{\fluc^N_i,g}\dr+\int_0^t\br{u_{\eta,i}-u_i,g}\dr .
\]
By Cauchy--Schwarz, the first term has expectation at most
$N^{-1/2}\|g\|_{\Ufine^*}\big(T\,\E\int_0^T\|\fluc^N_i\|^2_{\Ufine}\dr\big)^{1/2}\le CN^{-1/2}$, by
Proposition~\ref{prop:close}, whose hypothesis holds by Theorem~\ref{thm:K4}. The second is at most
$T\|g\|_{L^2}\sup_{r\le T}\|u_\eta(r)-u(r)\|_{L^2}\le C\eta$, by \eqref{eq:CHbias}. So
$\E|J^N_1-\int_0^t\br{u_i,g}\dr|\to0$.
For $J^N_2$, Proposition~\ref{prop:bracketid} gives
\[
 \E\Big|J^N_2-\int_0^t\br{u_i,|\na\varphi|^2\,2\sum_jf(a_{ij}u_j)}\dr\Big|
 \ \longrightarrow\ 0.
\]
Hence $\E|b^N(t)-b(t)|\to0$.

\emph{Step 2: the limit is a martingale with the limiting bracket.} With $i$,
$\varphi$, $b^N$ and $b$ as in Step~1, put $m^N(t):=M^N_i(t,\varphi)$. The
particle filtration $(\mathcal F_t)$ contains the filtration generated by
$(\fluc^N,M^N)$. Let $0\le r<t\le T$, let $r_1,\dots,r_p\le r$, let $h$ be a bounded
continuous function on $\Hmid^{2np}$, and put
$h^N:=h\big((\fluc^N(r_q),M^N(r_q))_{q\le p}\big)$. The processes $m^N$ and
$(m^N)^2-b^N$ are $(\mathcal F_t)$-martingales, and $h^N$ is $\mathcal F_r$-measurable, so
\begin{equation}\label{eq:martid}
 \E\big[(m^N(t)-m^N(r))\,h^N\big]=0,\qquad
 \E\big[\big(m^N(t)^2-m^N(r)^2\big)\,h^N\big]=\E\big[\big(b^N(t)-b^N(r)\big)\,h^N\big].
\end{equation}
By Step~1, $\E|b^N(t')-b(t')|\to0$ for every $t'\le T$. Since $|h^N|\le\|h\|_\infty$, the
right side of the second identity equals $(b(t)-b(r))\,\E h^N+o(1)$. The left sides of
\eqref{eq:martid} and $\E h^N$ are expectations of continuous functions of
$(\fluc^N,M^N)\in \mathcal Z_T$, integrable by Proposition~\ref{prop:M4}. By
Proposition~\ref{prop:skorokhod}, they do not change when $(\fluc^N,M^N)$ is replaced
by $(\widetilde\fluc^N,\widetilde M^N)$. Denote by $\widetilde m^N:=\widetilde M^N_i(\cdot,\varphi)$,
$\widetilde m:=\widetilde M_i(\cdot,\varphi)$, $\widetilde h^N$ and $\widetilde h$ the
corresponding functions of the new variables. Along the subsequence of
Proposition~\ref{prop:skorokhod}, $\widetilde m^N(t')\to\widetilde m(t')$ for every $t'$ and
$\widetilde h^N\to\widetilde h$ almost surely, by convergence in $C([0,T];\Hmid^n)$. By
Proposition~\ref{prop:M4}, $\sup_N\E\sup_{t\le T}|m^N(t)|^4<\infty$, so
$(\widetilde m^N(t'))^2$ is uniformly integrable. Passing to the limit,
\[
 \widetilde\E\big[(\widetilde m(t)-\widetilde m(r))\,\widetilde h\big]=0,\qquad
 \widetilde\E\big[\big(\widetilde m(t)^2-\widetilde m(r)^2-b(t)+b(r)\big)\,\widetilde h\big]=0 .
\]
A monotone class argument extends both to every bounded $\widetilde h$ measurable with
respect to $\widetilde{\mathcal F}_r:=\sigma\big(\widetilde\fluc(q),\widetilde M(q):q\le r\big)$.
So $\widetilde m$ is a continuous square-integrable $(\widetilde{\mathcal F}_t)$-martingale
with $\br{\widetilde m}_t=b(t)$. Polarisation gives the brackets of pairs
$(\varphi,\psi)$. For $i\ne k$ the finite-$N$ cross-bracket vanishes
(Proposition~\ref{prop:ito}), and the same argument applied to
$M^N_i(\cdot,\varphi)M^N_k(\cdot,\psi)$ shows that the limit cross-bracket vanishes.

\emph{Step 3: the $\Hmid$-valued martingale.} The map $v\mapsto\sup_t\|v(t)\|_{\Hmid}$ is
continuous on $C([0,T];\Hmid)$, so Fatou's lemma, the almost sure convergence and
\eqref{eq:M2} give
$\widetilde\E\sup_t\|\widetilde M_i(t)\|^2_{\Hmid}\le c_{\rm emb}^2\sup_N\E\sup_t\|M^N_i(t)\|^2_{\Ufine}<\infty$.
So $\widetilde M_i(t)$ is Bochner integrable in $\Hmid$, and so is its conditional expectation
given $\widetilde{\mathcal F}_r$. By Step~2, $\widetilde\E[\widetilde M_i(t)-\widetilde M_i(r)\,|\,\widetilde{\mathcal F}_r]$
pairs to $0$ with every $\varphi$ in $H^{s_m+1}_{a_m}$, which is dense in $\Hmid^*$. Hence it is
$0$, and $\widetilde M_i$ is a continuous square-integrable $\Hmid$-valued martingale.

\emph{Step 4: Gaussianity.} By Steps 1--2 the bracket of $\widetilde M_i$ is
the deterministic function of $t$ given by \eqref{eq:cov_limit} with
$\varphi=\psi$. Polarisation gives the mixed brackets, and the vanishing of
$\br{M^N_i,M^N_k}$ for $i\ne k$ (Proposition~\ref{prop:ito}) gives
$\delta_{ik}$. Fix finitely many pairs $(\varphi_q,i_q)$, $q\le p$, with
$\varphi_q\in H^{s_m+1}_{a_m}$, and put $\mathsf X^q(t):=\widetilde M_{i_q}(t,\varphi_q)$. Then
$\mathsf X=(\mathsf X^1,\dots,\mathsf X^p)$ is a continuous $\R^p$-valued $(\widetilde{\mathcal F}_t)$-martingale with
$\mathsf X(0)=0$, and its bracket matrix is deterministic, by Step~2. By Lemma~\ref{lem:gaussmart}, for
$0=t_0<t_1<\dots<t_L\le T$ the increments $\mathsf X(t_l)-\mathsf X(t_{l-1})$ are centred Gaussian, independent
of each other and of $\widetilde{\mathcal F}_0$. So the finite-dimensional distributions of
$(\widetilde M_i(\cdot,\varphi))_{i,\varphi}$ are centred Gaussian and independent of
$\widetilde{\mathcal F}_0$. The pairings with
$\varphi\in H^{s_m+1}_{a_m}$, a dense subset of $\Hmid^*$, generate the Borel
$\sigma$-field of $C([0,T];\Hmid^n)$. So $\widetilde M$ is a Gaussian process with independent
increments, and it is independent of $\widetilde{\mathcal F}_0$, in particular of
$\widetilde\fluc(0)$.

\end{proof}

\begin{remark}[Correlated species]\label{rem:X}
We discuss only correlated Brownian motions across species. Condition \eqref{eq:datax} then
enters in three places. It makes the finite-$N$ cross-brackets
of $M^N_i$ and $M^N_k$, $i\ne k$, vanish (Proposition~\ref{prop:ito}), and so the limit bracket
\eqref{eq:cov_limit} is species-diagonal. It also makes the intermediate ensemble
independent across species (Lemma~\ref{lem:crossspecies}(ii)), which Lemma~\ref{lem:D} uses. And
Lemma~\ref{lem:MNval} uses that $W^N$ is a standard Brownian motion in $\R^{nNd}$. If the
Brownian motions of different species were correlated, the first use would add
off-diagonal terms to the limit bracket. We have not checked whether the fourth-moment
estimate of Appendix~\ref{sec:appk4} survives without the second use.
\end{remark}

\section{The limit martingale and the stochastic convolution}\label{sec:limitmart}

\noindent
This section studies the limit object rather than the particle system. It opens with
the maximal bound for the stochastic convolution at finite $N$ (Proposition~\ref{prop:conv}),
which Proposition~\ref{prop:close} uses. Definition~\ref{def:weak} fixes the class of driving
martingales: continuous, $\Ufine^n$-valued, with species-diagonal brackets that have a
bounded density. Proposition~\ref{prop:rep} represents every such martingale as a stochastic
integral against a cylindrical Wiener process, on an extension of the probability space.
With it, the stochastic convolution $\int_0^tS^0(t,s)\dd M(s)$ is a stochastic integral
with a deterministic integrand (Proposition~\ref{prop:convint}), and its maximal bound
(Theorem~\ref{thm:maxlimit}) follows as at finite $N$. Proposition~\ref{prop:initialCLT} is
the classical central limit theorem for the initial fluctuation, and
Lemma~\ref{lem:initmom} places its limit in $\Ufine^n$.

\begin{lemma}[Measurability of the evolution family]\label{lem:Smeas}
Let $\eta\in[0,\eta_0]$ and let $S^\eta$ be the family on $\Ufine^n$ of
Proposition~\ref{prop:evol}\textup{(iii)}.
\begin{enumerate}[label=\normalfont(\roman*)]
\item For $t\in(0,T]$ and $\chi\in(\Ufine^*)^n$, the map $s\mapsto S^\eta(t,s)^*\chi$ is continuous
from $[0,t]$ into $(\Ufine^*)^n$. Hence, for $w\in\Ufine^n$, the map $s\mapsto S^\eta(t,s)w$ is
weakly continuous from $[0,t]$ into $\Ufine^n$ and norm continuous into $\Hmid^n$.
\item For $w\in\Ufine^n$, the map $(r,s)\mapsto S^\eta(s,r)w$ is Borel from
$\mathbb T:=\{(r,s):0\le r<s\le T\}$ into $\Ufine^n$.
\item Let $\mathsf U'$ be a separable Hilbert space, $\Upsilon$ a predictable process with values in
$\HS(\mathsf U';\Ufine^n)$, and $s\in(0,T]$. Then $r\mapsto S^\eta(s,r)\Upsilon(r)$ is predictable on
$[0,s]$ with values in $\HS(\mathsf U';\Ufine^n)$, and
$\|S^\eta(s,r)\Upsilon(r)\|_{\HS}\le C_1(T)\|\Upsilon(r)\|_{\HS}$.
\end{enumerate}
\end{lemma}

\begin{proof}
\emph{(i).} For $\chi\in\Test^n$, the map $s\mapsto S^\eta(t,s)^*\chi$ is continuous into
$\Test^n\hookrightarrow(\Ufine^*)^n$ by Lemma~\ref{lem:strongcont}. For general $\chi$, take
$\chi_k\in\Test^n$ with $\chi_k\to\chi$ in $(\Ufine^*)^n$, which is possible since
$C^\infty_0\subset\Test$ is dense in $\Ufine^*=H^{s_1}_{a_1}$. By
Proposition~\ref{prop:evol}\textup{(iii)}, $\|S^\eta(t,s)^*(\chi-\chi_k)\|\le C_1(T)\|\chi-\chi_k\|$
uniformly in $s$. So $s\mapsto S^\eta(t,s)^*\chi$ is a uniform limit of continuous maps. Weak
continuity follows from $\br{S^\eta(t,s)w,\chi}=\br{w,S^\eta(t,s)^*\chi}$. A sequence that converges
weakly in $\Ufine^n$ converges in norm in $\Hmid^n$, because the embedding
$\Ufine\hookrightarrow\Hmid$ of \eqref{eq:chainA} is compact.

\emph{(ii).} Fix $\chi\in(\Ufine^*)^n$ and put $F(r,s):=\br{S^\eta(s,r)w,\chi}$ on $\mathbb T$. It is
continuous in $r$ for fixed $s$, by (i), and continuous in $s$ for fixed $r$, by
Corollary~\ref{cor:fwdweak}\textup{(i)}. For $k\in\N$ put $s_k:=\min(\lceil ks\rceil/k,T)$. Then
$(r,s_k)\in\mathbb T$, and $F_k(r,s):=F(r,s_k)$ is Borel, since on each of the countably many sets
$\{s_k=c\}$ it is a continuous function of $r$. As $k\to\infty$, $s_k\to s$, so $F_k\to F$
pointwise, and $F$ is Borel. Hence $(r,s)\mapsto S^\eta(s,r)w$ is weakly measurable, and it is
Borel by Pettis' theorem, since $\Ufine^n$ is separable.

\emph{(iii).} For $e\in \mathsf U'$ and $\chi\in(\Ufine^*)^n$,
$\br{S^\eta(s,r)\Upsilon(r)e,\chi}=\br{\Upsilon(r)e,S^\eta(s,r)^*\chi}$ pairs the predictable
$\Ufine^n$-valued process $\Upsilon(\cdot)e$ with a deterministic continuous map, by (i). So it is
predictable. Let $e$ run through an orthonormal basis of $\mathsf U'$ and $\chi$ through a countable dense
subset of $(\Ufine^*)^n$. The functionals $Q\mapsto\br{Qe,\chi}$ obtained in this way are continuous
and separate the points of the separable Hilbert space $\HS(\mathsf U';\Ufine^n)$. So they generate its
Borel $\sigma$-field, and $S^\eta(s,\cdot)\Upsilon$ is predictable. The bound is
$\|SQ\|_{\HS}\le\|S\|\,\|Q\|_{\HS}$ with Proposition~\ref{prop:evol}\textup{(iii)}.
\end{proof}

\begin{proposition}[Maximal bound for the stochastic convolution]\label{prop:conv}
With $\Upsilon^N$ and $W^N$ as in
Lemma~\ref{lem:MNval}, let $\beta\in(\frac14,\frac12)$ and let
$\mathsf Y^N(s):=\int_0^s(s-r)^{-\beta}S^\eta(s,r)\Upsilon^N(r)\dd W^N(r)$, in the sense of
Lemma~\ref{lem:BDG} and in a jointly measurable version. Then $c_\beta R_\beta \mathsf Y^N$, where
$(R_\beta g)(t):=\int_0^t(t-s)^{\beta-1}S^\eta(t,s)g(s)\dd s$ and
$c_\beta=\frac{\sin\pi\beta}{\pi}$, is for each $t$ a version of
$Z^N(t):=\int_0^tS^\eta(t,s)\Upsilon^N(s)\dd W^N(s)$, and for any countable dense
$\mathbb D\subset[0,T]$,
\begin{equation}\label{eq:conv}
 \sup_N\ \E\Big[\sup_{t\in\mathbb D}\big\|Z^N(t)\big\|^4_{\Ufine}\Big]\ \le\ C(T,\beta).
\end{equation}
\end{proposition}

\begin{proof}
We use the factorisation method. By the
stochastic Fubini theorem, $Z^N$ is the image under the deterministic operator $c_\beta R_\beta$ of
the auxiliary process $\mathsf Y^N$ (Step~1). The operator $R_\beta$ is bounded from
$L^4(0,T;\Ufine)$ into bounded functions (Step~2), and $\mathsf Y^N$ is bounded in
$L^4(\Omega\times(0,T);\Ufine)$ by Lemma~\ref{lem:BDG} (Step~3).

\emph{Step 0 (measurability).} $\mathsf Y^N$ is defined by a stochastic integral, hence
only up to a null set for each $s$. The integrand
$(r,s)\mapsto(s-r)^{-\beta}\mathbf 1_{r<s}S^\eta(s,r)$ is strongly measurable, by
Lemma~\ref{lem:Smeas}\textup{(ii)}, and the proof of
\cite[Thm.~4.33]{da_prato_2014_SPDEs} gives a version of $\mathsf Y^N$ that is jointly
measurable in $(s,\omega)$. The
Bochner integral in $R_\beta$ and the Tonelli step of Step~3 use this version.

\emph{Step 1 (factorisation).} At fixed $N$ and $\eta$, $\mathsf q^N_i$ is bounded by the deterministic
constant $c_\eta C_{s_1,a_1}$, where $c_\eta:=2\max_i\sigma_i+2n\|f_\eta\|_{L^\infty(\R)}$. Indeed,
$\mathsf q^N_i=\br{\mu^N_i,K_{s_1,a_1}(2\Theta^N_i+2\varrho^N_{\cdot,i})}$ by \eqref{eq:trMN},
$K_{s_1,a_1}\le C_{s_1,a_1}$ by Lemma~\ref{lem:kernel}, $2\Theta^N_i+2\varrho^N_{k,i}\le c_\eta$ since
the self-term is deleted and $f_\eta$ is bounded by \ref{H:CH}, and $\mu^N_i$ is a probability measure. With $\|S\Upsilon^N\|_{\HS}\le\|S\|\,\|\Upsilon^N\|_{\HS}$,
\eqref{eq:trMN} and Proposition~\ref{prop:evol}, this gives
\begin{equation}\label{eq:Yfinite}
 \E\int_0^s(s-r)^{-2\beta}
 \big\|S^\eta(s,r)\Upsilon^N(r)\big\|^2_{\HS}\dd r
 \ \le\ C_1(T)^2\,\E\int_0^s(s-r)^{-2\beta}\sum_i\mathsf q^N_i\dd r
 \ \le\ C_1(T)^2\,n\,c_\eta C_{s_1,a_1}\,\frac{s^{1-2\beta}}{1-2\beta}\ <\ \infty,
\end{equation}
since $2\beta<1$. So $\mathsf Y^N(s)$ is well defined for every $s$, by Lemma~\ref{lem:BDG}. The stochastic Fubini theorem
\cite[Thm.~4.33]{da_prato_2014_SPDEs} applies, its integrability condition being satisfied:
by Cauchy--Schwarz against $(t-s)^{\beta-1}\dd s$ and Tonelli, using
$\int_r^t(t-s)^{\beta-1}(s-r)^{-2\beta}\dd s=\mathrm B(\beta,1-2\beta)(t-r)^{-\beta}$

\[
 \int_0^t(t-s)^{\beta-1}\Big(\int_0^s(s-r)^{-2\beta}\E\mathsf q^N\dd r\Big)^{1/2}\dd s
 \le\Big(\frac{t^\beta}{\beta}\Big)^{1/2}
 \Big(\mathrm B(\beta,1-2\beta)\!\int_0^t(t-r)^{-\beta}\E\mathsf q^N\dd r\Big)^{1/2}<\infty
\]
because $\beta<\frac12$ and $\E\mathsf q^N\in L^2(0,T)$ by \eqref{eq:psi2}.
 With the
evolution law $S^\eta(t,s)S^\eta(s,r)=S^\eta(t,r)$ and, substituting
$s=r+(t-r)u$,
$\int_r^t(t-s)^{\beta-1}(s-r)^{-\beta}\dd s=\mathrm B(\beta,1-\beta)=\pi/\sin\pi\beta$

one gets $Z^N(t)=c_\beta(R_\beta \mathsf Y^N)(t)$ for each fixed $t$, a.s.

\emph{Step 2 ($R_\beta$).} By Proposition~\ref{prop:evol} and H\"older with
exponents $(\frac43,4)$,
\[
 \|(R_\beta g)(t)\|_{\Ufine}\le C_1(T)\Big(\int_0^Tr^{\frac43(\beta-1)}\dd r\Big)^{3/4}
 \|g\|_{L^4(0,T;\Ufine)},
\]
finite if and only if $\frac43(\beta-1)>-1$, i.e.\ $\beta>\frac14$. The bound holds for every
$t$. The identity $Z^N(t)=c_\beta(R_\beta \mathsf Y^N)(t)$ of Step~1, however, holds only almost surely
for each fixed $t$, so \eqref{eq:conv} is stated with a supremum over a countable set.

\emph{Step 3 (fourth moment of $\mathsf Y^N$).} The integrand
$r\mapsto(s-r)^{-\beta}S^\eta(s,r)\Upsilon^N(r)$ is predictable, by Lemma~\ref{lem:Smeas}\textup{(iii)} with $\mathsf U'=\R^{nNd}$ and the
predictable integrand $\Upsilon^N$ of Lemma~\ref{lem:MNval}. It is square-integrable for every $s$, by \eqref{eq:Yfinite}. Lemma~\ref{lem:BDG}
with $p=4$ and the same Hilbert--Schmidt bound give
\[
 \E\|\mathsf Y^N(s)\|^4_{\Ufine}\le C\,C_1(T)^4\,
 \E\Big[\Big(\int_0^s(s-r)^{-2\beta}\textstyle\sum_{i}\mathsf q^N_i(r)\dd r\Big)^{2}\Big].
\]
Cauchy--Schwarz against the measure $(s-r)^{-2\beta}\dd r$, of mass
$\frac{s^{1-2\beta}}{1-2\beta}\le m_\beta:=\frac{T^{1-2\beta}}{1-2\beta}$,
then Tonelli with $\int_r^T(s-r)^{-2\beta}\dd s\le m_\beta$ (the same
constant), give
$\int_0^T\E\|\mathsf Y^N(s)\|^4\dd s\le C\,C_1(T)^4m_\beta^2n\sum_i\E\int_0^T(\mathsf q^N_i)^2\dd r\le C$
by \eqref{eq:psi2}.

\emph{Step 4 (conclusion).} Pathwise $\sup_{t}\|c_\beta R_\beta \mathsf Y^N(t)\|\le
c_\beta\|R_\beta\|\,\|\mathsf Y^N\|_{L^4(0,T;\Ufine)}$, with $\|R_\beta\|$ the operator norm of Step~2.
Off one null set, $Z^N(t)=c_\beta(R_\beta \mathsf Y^N)(t)$ for every $t$ in the countable set $\mathbb D$, by
Step~1. Hence
\[
 \E\sup_{t\in\mathbb D}\|Z^N(t)\|^4_{\Ufine}\le c_\beta^4\|R_\beta\|^4\int_0^T\E\|\mathsf Y^N(s)\|^4_{\Ufine}\dd s\le C,
\]
by Step~3, uniformly in $N$.
\end{proof}

\begin{definition}[Weak solution]\label{def:weak}
Let $(\Omega,\mathcal F,(\mathcal F_t),\Prob)$ be a filtered probability space satisfying the
usual conditions. Let $\fluc^0$ be $\Ufine^n$-valued and $\mathcal F_0$-measurable, and let $M$ be
a continuous square-integrable $\Ufine^n$-valued martingale, $M(0)=0$, such that
$\varphi\mapsto M_i(t,\varphi)$ is linear and continuous on $\Ufine^*$,
$t\mapsto M_i(t,\varphi)$ is continuous, and the brackets are species-diagonal
and absolutely continuous,
\begin{equation}\label{eq:limbracket}
 \br{M_i(\cdot,\varphi),M_k(\cdot,\psi)}_t=\delta_{ik}\int_0^t
 \br{u_i,\na\varphi\cdot\na\psi\,\Gamma_i}\dd r,
 \qquad 0\le\Gamma_i\in L^\infty\big((0,T)\times\R^d\big),
\end{equation}
for each pair $\varphi,\psi\in\Ufine^*=H^{s_1}_{a_1}\supset\Test$.

We denote $\mathsf q_i(r):=\br{u_i,K_{s_1,a_1}\Gamma_i}$. Then
$\br{M_i(\cdot,\chi)}_t\le\int_0^t\mathsf q_i\dd r\cdot\|\chi\|^2_{H^{s_1}_{a_1}}$ by
the pointwise bound $|\na\chi(x)|^2\le K_{s_1,a_1}(x)\|\chi\|^2_{H^{s_1}_{a_1}}$
of Lemma~\ref{lem:kernel}. Since $u_i(r)$ is a probability density \textup{(}Lemma~\ref{lem:lawid}(ii) at
$\eta=0$\textup{)} and $K_{s_1,a_1}\le C_{s_1,a_1}$, $\mathsf q_i$ is bounded by
$C_{s_1,a_1}\|\Gamma_i\|_\infty$, so $\mathsf q_i\in L^2(0,T)$ and
$\int_0^T\mathsf q_i\dd r<\infty$.

An adapted $\Ufine^n$-valued process $\fluc$ is a \emph{weak solution} if
$\Prob$-a.s.\ it is weakly continuous in $(\Usharp)^n$ with
$\sup_{t\le T}\|\fluc(t)\|_{\Usharp}<\infty$, and for every $\varphi\in\Test$ and
$t\in[0,T]$, $\Prob$-a.s.,
\begin{equation}\label{eq:weakform}
 \br{\fluc_i(t),\varphi}=\br{\fluc^0_i,\varphi}
 +\int_0^t\sum_j\br{\fluc_j(s),(\mathcal G(s)\varphi)_{ij}}\dd s+M_i(t,\varphi).
\end{equation}
\end{definition}

\begin{corollary}\label{cor:supupgrade}
Set $\|w\|_{\Ufine^n}:=+\infty$ for $w\in\Hmid^n\setminus\Ufine^n$. For any process $w$ with
$\Prob$-a.s.\ $\Hmid^n$-continuous paths and any countable dense
$\mathbb D\subset[0,T]$,
\begin{equation}\label{eq:supupgrade}
 \sup_{t\in[0,T]}\|w(t)\|_{\Ufine^n}\ =\ \sup_{t\in\mathbb D}\|w(t)\|_{\Ufine^n}
 \qquad\text{in }[0,\infty],\ \Prob\text{-a.s.}
\end{equation}
In particular, by Lemma~\ref{lem:flucHmid}, \eqref{eq:supupgrade} holds for any
weak solution $\fluc$.
\end{corollary}

\begin{proof}
$\Ufine^n$ is a Hilbert space, hence
reflexive, and $\Ufine^n\hookrightarrow\Hmid^n$ is continuous and injective by
\eqref{eq:chainA}, since both are spaces of tempered distributions and the
embedding is the identity. These are the hypotheses of
Lemma~\ref{lem:lscsup}\textup{(i)} at $(X,Y)=(\Ufine^n,\Hmid^n)$, so
$\|\cdot\|_{\Ufine^n}:\Hmid^n\to[0,\infty]$ is lower semicontinuous.
Composing with a continuous path makes
$t\mapsto\|w(t)\|_{\Ufine^n}$ lower semicontinuous on $[0,T]$, and a lower
semicontinuous $[0,\infty]$-valued function has the same supremum over a dense set
as over the interval.
\end{proof}

\begin{proposition}[Initial CLT]\label{prop:initialCLT}
Let $\xi^1_i,\xi^2_i,\dots$ be i.i.d.\ with density $u_{0,i}$, independent across
$i$. Then, with $\fluc^N_i(0)=N^{-1/2}\sum_k(\delta_{\xi^k_i}-u_{0,i})$,
the vector $\fluc^N(0)=(\fluc^N_1(0),\dots,\fluc^N_n(0))$ converges in law to $\fluc^0$ in $\Hmid^n$, where $(\fluc^0_1,\dots,\fluc^0_n)$ is a
centred Gaussian field, independent across species, with covariance
\[
 \E\big[\br{\fluc^0_i,\varphi}\br{\fluc^0_i,\psi}\big]
 =\int_{\R^d}\varphi\psi\,u_{0,i}\dx-\Big(\int\varphi u_{0,i}\Big)\Big(\int\psi u_{0,i}\Big),
\]
for all $\varphi,\psi\in\Hmid^*$. That is, the covariance of $\br{\fluc^0_i,\varphi}$ and
$\br{\fluc^0_i,\psi}$ is the covariance of $\varphi(\xi^1_i)$ and $\psi(\xi^1_i)$.
\end{proposition}

\begin{proof}
Finite-dimensional distributions: the classical multivariate CLT, applied to the i.i.d.\
vectors $\big(\varphi_q(\xi^k_{i_q})\big)_{q\le p}$, $k\le N$, for finitely many pairs
$(\varphi_q,i_q)$. The limit covariance vanishes between different species, because the
$\xi^k_i$ are independent across $i$. Tightness in
$\Hmid^n$: Lemma~\ref{lem:field0} gives a uniform bound in $\Ufine^n$ and
$\Ufine^n\hookrightarrow\hookrightarrow\Hmid^n$. A tight sequence whose finite-dimensional
distributions converge, against a set of $\varphi$ that separates points of $\Hmid^n$, converges
in law.
\end{proof}

\begin{lemma}[The initial datum lies in $\Ufine^n$]\label{lem:initmom}
Let $C_{s_1,a_1}$ be the constant of
Lemma~\ref{lem:kernel}. Then:
\begin{enumerate}[label=\normalfont(\roman*)]
\item for $\fluc^0$ the initial fluctuation of Proposition~\ref{prop:initialCLT},
\begin{equation}\label{eq:init2}
 \E\big\|\fluc^0\big\|^2_{\Ufine^n}\ \le\ n\,C_{s_1,a_1},\qquad
 \E\big\|\fluc^0\big\|^4_{\Ufine^n}\ \le\ 3\,n^2C_{s_1,a_1}^2,
\end{equation}
so in particular $\fluc^0\in\Ufine^n$ almost surely, and $\fluc^0$ is a Borel map into
$\Ufine^n$ on this full-measure set;
\item for $\widetilde\fluc$ a limit point of Proposition~\ref{prop:skorokhod}, the same
bounds hold for $\widetilde\fluc(0)$, and $\widetilde\fluc(0)$ is a Borel map into
$\Ufine^n$.
\end{enumerate}
\end{lemma}

\begin{proof}
Both parts pass to the limit in finite-$N$ moment bounds, using that the $\Ufine^n$ norm is
lower semicontinuous on $\Hmid^n$: through the portmanteau inequality in (i), and through Fatou's
lemma in (ii).
We use $\|\nu\|^2_{\Ufine^n}=\sum_i\|\nu_i\|^2_{\Ufine}$, extended by $+\infty$ off
$\Ufine^n$. By Lemma~\ref{lem:lscsup}\textup{(i)} at $(X,Y)=(\Ufine^n,\Hmid^n)$, the map
$\nu\mapsto\|\nu\|_{\Ufine^n}$ is lower semicontinuous on $\Hmid^n$, hence Borel.

\emph{(i).} Fix $i$ and put $\zeta_k:=N^{-1/2}(\delta_{\xi^k_i}-u_{0,i})$, so that
$\fluc^N_i(0)=\sum_{k\le N}\zeta_k$ with independent centred $\zeta_k\in\Ufine$. By
Lemma~\ref{lem:kernel}, $\|\delta_x\|^2_{\Ufine}=K^0_{s_1,a_1}(x)\le C_{s_1,a_1}$, so
$\|\zeta_k\|_{\Ufine}\le2N^{-1/2}C_{s_1,a_1}^{1/2}$. Lemma~\ref{lem:sum4} and
Lemma~\ref{lem:field0} give
$\E\|\fluc^N_i(0)\|^4_{\Ufine}\le16C_{s_1,a_1}^2N^{-1}+3C_{s_1,a_1}^2$. The laws of
$\fluc^N_i(0)$ converge in $\Hmid$ (Proposition~\ref{prop:initialCLT}), and the
portmanteau inequality for the nonnegative lower semicontinuous functions
$\|\cdot\|^2_{\Ufine}$ and $\|\cdot\|^4_{\Ufine}$ gives $\E\|\fluc^0_i\|^2_{\Ufine}\le
C_{s_1,a_1}$ and $\E\|\fluc^0_i\|^4_{\Ufine}\le3C_{s_1,a_1}^2$. Sum over $i$, using
$\big(\sum_i a_i\big)^2\le n\sum_ia_i^2$ for the fourth moment. Borel measurability into
$\Ufine^n$ follows as in (ii) below.

\emph{(ii).} By Proposition~\ref{prop:skorokhod}, $\widetilde\fluc^N(0)\to\widetilde\fluc(0)$
almost surely in $\Hmid^n$, and $\widetilde\fluc^N(0)$ has the law of $\fluc^N(0)$. The
bounds of (i) at finite $N$ and Fatou's lemma with lower semicontinuity give the bounds
for $\widetilde\fluc(0)$. On the full-measure set $\{\|\widetilde\fluc(0)\|_{\Ufine^n}<\infty\}$,
$\widetilde\fluc(0)$ is Borel into $\Ufine^n$ by the Lusin--Souslin theorem \cite[Thm.~15.1]{Kec95}, since
$\Ufine^n\hookrightarrow\Hmid^n$ is a continuous injection between Polish spaces.
\end{proof}

\begin{proposition}[Representation of the driving martingale]\label{prop:rep}
Let $H$ be $\Ufine^n$ or $\Hmid^n$, let $(\mathcal F_t)$
satisfy the usual conditions, and let $M$ be a continuous square-integrable $H$-valued
$(\mathcal F_t)$-martingale with $M(0)=0$. Let $0\le\Gamma_i\in L^\infty((0,T)\times\R^d)$,
and assume that \eqref{eq:limbracket} holds for all $i,k$ and all $\varphi,\psi$ in a
dense subset of $H^*$ \textup{(}componentwise\textup{)}. Put $\mathsf U:=L^2(\R^d;\R^{nd})$, let
$\mathsf q_i(r)=\br{u_i(r),K_{s_1,a_1}\Gamma_i(r)}$ be as in Definition~\ref{def:weak}, and put
\[
 \big(\Upsilon(r)v\big)_i:=\diver\Big(\big(u_i(r)\Gamma_i(r)\big)^{1/2}v_i\Big),
 \qquad v=(v_1,\dots,v_n)\in\mathsf U.
\]
Then:
\begin{enumerate}[label=\normalfont(\roman*)]
\item $\Upsilon(r)\in\HS(\mathsf U;\Ufine^n)$, $r\mapsto\Upsilon(r)$ is strongly measurable, and
\begin{equation}\label{eq:psibdd}
 \|\Upsilon(r)\|^2_{\HS(\mathsf U;\Ufine^n)}=\sum_i\mathsf q_i(r),\qquad
 0\le\mathsf q_i(r)\le C_{s_1,a_1}\|\Gamma_i\|_{L^\infty},\qquad
 C_{\mathsf q}:=\sum_iC_{s_1,a_1}\|\Gamma_i\|_{L^\infty};
\end{equation}
\item on the extension $\widehat\Omega:=\Omega\times\Omega'$ of
\cite[Thm.~8.2]{da_prato_2014_SPDEs}, with the augmented product filtration, there is
a cylindrical Wiener process $W$ on $\mathsf U$ such that $M(t)=\int_0^t\Upsilon\dd W$ for every
$t$, almost surely. The integral is an $\Ufine^n$-valued continuous martingale, and the
identity holds in $H$;
\item $M$ is a continuous square-integrable $\Ufine^n$-valued $(\mathcal F_t)$-martingale,
\eqref{eq:limbracket} holds for all $\varphi,\psi\in\Ufine^*$, and
$\E\|M(t)-M(s)\|^2_{\Ufine^n}=\int_s^t\sum_i\mathsf q_i\dd r$ for $0\le s\le t\le T$.
\end{enumerate}
\end{proposition}

\begin{proof}
Part (i) is a computation with the reproducing kernel of Lemma~\ref{lem:kernel}. For (ii) we
identify the operator quadratic variation of $M$ from its scalar brackets and apply the
martingale representation theorem \cite[Thm.~8.2]{da_prato_2014_SPDEs}. The proof of (ii) only checks the hypotheses of that theorem. Part (iii) transfers the
$\Ufine^n$-valued properties of the stochastic integral back to $M$.

\emph{(i).} By Lemma~\ref{lem:lawid}(ii) at $\eta=0$, $u_i(r)\ge0$ and $\int u_i(r)=1$, so
$\mathsf q_i(r)\le C_{s_1,a_1}\|\Gamma_i\|_{L^\infty}$ by Lemma~\ref{lem:kernel} at
$(s,a)=(s_1,a_1)$. For $\chi\in\Ufine^*$ and $v\in\mathsf U$,
$\br{(\Upsilon(r)v)_i,\chi}=-\int(u_i\Gamma_i)^{1/2}v_i\cdot\na\chi\dx$, so the duality
adjoint is $\Upsilon(r)^*\chi=-\big((u_i\Gamma_i)^{1/2}\na\chi_i\big)_i$. Let $(e_m)$ be an
orthonormal basis of $H^{s_1}_{a_1}$, and let $y_m\in\Ufine$ be the Riesz
representatives of the functionals $\br{\cdot,e_m}$. They form an orthonormal basis of
$\Ufine$ by Lemma~\ref{lem:spaces}(i). By monotone convergence and \eqref{eq:kernel},
\[
 \|\Upsilon(r)\|^2_{\HS(\mathsf U;\Ufine^n)}=\sum_i\sum_m\big\|(u_i\Gamma_i)^{1/2}\na e_m\big\|^2_{L^2}
 =\sum_i\int u_i\Gamma_iK_{s_1,a_1}\dx=\sum_i\mathsf q_i(r).
\]
For fixed $v$ and $\chi$, $r\mapsto\br{\Upsilon(r)v,\chi}$ is measurable by Fubini's
theorem. So $r\mapsto(\Upsilon(r),F)_{\HS}$ is measurable for finite-rank $F$, and, by the
uniform bound and density of finite-rank operators, for every $F\in\HS(\mathsf U;\Ufine^n)$. The space
$\HS(\mathsf U;\Ufine^n)$ is separable, so Pettis' theorem gives strong measurability.

\emph{(ii).} Denote by $\iota$ the embedding $\Ufine^n\hookrightarrow H$ ($\iota=\mathrm{Id}$
if $H=\Ufine^n$), and let $R_H:H^*\to H$ be the Riesz map. By (i), $\iota\Upsilon$ is a
deterministic, strongly measurable, bounded $\HS(\mathsf U;H)$-valued map, and for $a=R_H\varphi$,
$b=R_H\psi$ with $\varphi,\psi\in H^*$,
\[
 \int_0^t\big((\iota\Upsilon)^{*_H}a,(\iota\Upsilon)^{*_H}b\big)_{\mathsf U}\dd r
 =\sum_i\int_0^t\br{u_i,\na\varphi_i\cdot\na\psi_i\,\Gamma_i}\dd r ,
\]
since $(\iota\Upsilon)^{*_H}R_H\varphi=\Upsilon^*\varphi$. Denote the left side by $\mathsf A_t(a,b)$. The form $\mathsf A_t$ is
deterministic, bilinear and continuous on $H\times H$, since $\iota\Upsilon$ is bounded in
$\HS(\mathsf U;H)$. By \cite[Prop.~3.13]{da_prato_2014_SPDEs}, $M$ has an operator quadratic variation
$\br{\!\br M\!}_t$, a non-negative nuclear operator on $H$, continuous in $t$, such that
$(\br{\!\br M\!}_ta,b)_H$ is the bracket $\br{\br{M,a},\br{M,b}}_t$ for all $a,b\in H$, almost surely.
The dense subset of $H^*$ in the hypothesis contains a countable set $\mathbb D_0$ that is dense in
$H^*$, because $H^*$ is separable. By hypothesis and the identity for $\mathsf A_t$ above, off one null set,
$(\br{\!\br M\!}_ta,b)_H=\mathsf A_t(a,b)$ for all $a,b\in R_H\mathbb D_0$ and all rational $t$, hence
for all $t$, since both sides are continuous in $t$. For fixed $\omega$ and $t$ both sides are
continuous bilinear forms on $H\times H$, and $R_H\mathbb D_0$ is dense in $H$. So, off that null
set, they agree for all $a,b\in H$ and all $t$, that is,
$\br{\!\br M\!}_t=\int_0^t(\iota\Upsilon)(\iota\Upsilon)^{*_H}\dd r$. This is the hypothesis of
\cite[Thm.~8.2]{da_prato_2014_SPDEs} with $\mathsf U$ as above and $Q=\mathrm{Id}$: the theorem
requires $Q$ only to be bounded, symmetric and non-negative, and its proof treats
$Q=\mathrm{Id}$ with a cylindrical Wiener process. It gives $W$ with
$M=\int_0^\cdot\iota\Upsilon\dd W$ in $H$. There $M$ is a martingale and $W$ a Wiener process for the
product filtration; both stay so after augmentation by the null sets and passage to the
right-continuous filtration, since their paths are continuous. The integral $\int_0^\cdot\Upsilon\dd W$ exists as a
continuous square-integrable $\Ufine^n$-valued martingale by
\cite[Thm.~4.27]{da_prato_2014_SPDEs}, since $\int_0^T\|\Upsilon\|^2_{\HS}\dd r\le TC_{\mathsf q}$,
and $\iota$ commutes with it.

\emph{(iii).} By (ii), almost surely on $\widehat\Omega$ the path of $M$ lies in
$\iota\big(C([0,T];\Ufine^n)\big)$. This set is Borel in $C([0,T];H)$, by the
Lusin--Souslin theorem, so the same holds $\Prob$-almost surely on $\Omega$. The
isometry of \cite[Thm.~4.27]{da_prato_2014_SPDEs} gives the second moment of the
increments, and Doob's inequality gives $\E\sup_t\|M(t)\|^2_{\Ufine^n}\le4TC_{\mathsf q}$. For
$r\le t$ and $A\in\mathcal F_r$, $\iota\E[\mathbf 1_A(M(t)-M(r))]=0$ by the martingale
property in $H$, and $\iota$ is injective, so $M$ is an $\Ufine^n$-valued martingale.
Finally $\br{M(t),\chi}=\int_0^t\Upsilon^*\chi\dd W$ for $\chi\in(\Ufine^*)^n$, whose bracket
is $\sum_i\int_0^t\br{u_i,|\na\chi_i|^2\Gamma_i}\dd r$; polarisation gives
\eqref{eq:limbracket} for all $\varphi,\psi\in\Ufine^*$.
\end{proof}

\begin{lemma}[The scalar integral against $M$]\label{lem:scalint}
Let $M$ be as in Definition~\ref{def:weak}. For a
measurable $g:[0,T]\to(\Ufine^*)^n$ put
\begin{equation}\label{eq:Nnorm}
 \mathcal N_t(g)^2:=\sum_{i=1}^n\int_0^t\br{u_i(r),|\na g_i(r)|^2\,\Gamma_i(r)}\dd r
 \ \le\ \max_i\|\mathsf q_i\|_{L^\infty}\int_0^t\|g(r)\|^2_{\Ufine^{*n}}\dd r .
\end{equation}
With $W$ and $\Upsilon$ of Proposition~\ref{prop:rep} and $g\in L^2(0,T;(\Ufine^*)^n)$, put
$\mathcal I_t(g):=\int_0^t\Upsilon(r)^*g(r)\dd W(r)$, a real stochastic integral. Then:
\begin{enumerate}[label=\normalfont(\roman*)]
\item for a step function $g=\sum_lg_l\mathbf 1_{(s_l,s_{l+1}]}$,
$\mathcal I_t(g)=\sum_l\br{M(t\wedge s_{l+1})-M(t\wedge s_l),g_l}$;
\item $\E|\mathcal I_t(g)|^2=\mathcal N_t(g)^2$, and $\mathcal I_t(g)$ has a version measurable with
respect to the completion of $\sigma(M(r):r\le t)$;
\item $t\mapsto \mathcal I_t(g)$ is a continuous square-integrable martingale with deterministic
bracket $\mathcal N_t(g)^2$.
\end{enumerate}
We also denote $\mathcal I_t(g)$ by $\sum_i\int_0^t\br{\dd M_i,g_i}$.
\end{lemma}

\begin{proof}
The bound in \eqref{eq:Nnorm} is \eqref{eq:kernelptwise} with \eqref{eq:psibdd}. The
integrand $r\mapsto\Upsilon(r)^*g(r)\in\mathsf U$ is deterministic with
$\|\Upsilon(r)^*g(r)\|^2_{\mathsf U}=\sum_i\br{u_i,|\na g_i|^2\Gamma_i}$. So (ii)'s isometry and (iii) are
\cite[Thm.~4.27]{da_prato_2014_SPDEs}. For (i), $\br{\int_a^b\Upsilon\dd W,\chi}=
\int_a^b\Upsilon^*\chi\dd W$ for fixed $\chi$, and Proposition~\ref{prop:rep}(ii) applies.
Step functions are dense in $L^2(0,T;(\Ufine^*)^n)$, so the isometry makes $\mathcal I_t(g)$ an
$L^2$-limit of the sums in (i), which gives the measurability in (ii).
\end{proof}

\begin{proposition}[The stochastic convolution]\label{prop:convint}
Let $M$ be as in Definition~\ref{def:weak}, fix
$t\in(0,T]$, and let $S^0(t,s)$ be the evolution family of Proposition~\ref{prop:evol}
at $\eta=0$. For a partition $\pi:0=s_0<\dots<s_L=t$ set
\[
 J_\pi:=\sum_{l=0}^{L-1}S^0(t,s_{l+1})\big(M(s_{l+1})-M(s_l)\big)\ \in\ \Ufine^n .
\]
Then:
\begin{enumerate}[label=\normalfont(\roman*)]
\item $\sup_\pi\E\|J_\pi\|^2_{\Ufine^n}\le C_1(T)^2\int_0^t\sum_i\mathsf q_i\dd r$;
\item $J_\pi$ converges in $L^2(\Omega;\Hmid^n)$ as $|\pi|\to0$, to a limit $J(t)$ that
does not depend on the sequence of partitions;
\item $J(t)\in\Ufine^n$ almost surely, and
\begin{equation}\label{eq:convbound}
 \E\|J(t)\|^2_{\Ufine^n}\ \le\ C_1(T)^2\int_0^t\sum_i\mathsf q_i\dd r ;
\end{equation}
\item for each $\varphi\in(\Ufine^*)^n$,
\begin{equation}\label{eq:convint}
 \br{J(t),\varphi}=\sum_{i=1}^n\int_0^t\br{\dd M_i(s),\big(S^0(t,s)^*\varphi\big)_i}
 \qquad\Prob\text{-a.s.},
\end{equation}
and $J(t)$ is the unique element of $L^2(\Omega;\Ufine^n)$ with this property.
\end{enumerate}
On the extension of Proposition~\ref{prop:rep}, $J(t)=\int_0^tS^0(t,s)\Upsilon(s)\dd W(s)$
almost surely. We denote $J(t)$ by $\int_0^tS^0(t,s)\dd M(s)$, and set $J(0):=0$.
\end{proposition}

\begin{proof}
We work on the extension of Proposition~\ref{prop:rep}, where $M=\int_0^\cdot\Upsilon\dd W$.
There $J_\pi$ is a stochastic integral with a deterministic piecewise constant integrand, and the
It\^o isometry gives (i) and the convergence (ii). The limit is then transferred back to $\Omega$.

By Lemma~\ref{lem:Smeas}\textup{(i)} at $\eta=0$, $s\mapsto S^0(t,s)w$ is weakly
continuous in $\Ufine^n$ and norm continuous in $\Hmid^n$ for every $w\in\Ufine^n$. By
Lemma~\ref{lem:Smeas}\textup{(iii)}, applied to the deterministic integrand $\Upsilon$ of
Proposition~\ref{prop:rep}\textup{(i)}, $s\mapsto S^0(t,s)\Upsilon(s)$ is strongly measurable into
$\HS(\mathsf U;\Ufine^n)$, with norm at most $C_1(T)\|\Upsilon(s)\|_{\HS}$. Let
$J_W(t):=\int_0^tS^0(t,s)\Upsilon(s)\dd W(s)$ on the extension.

By Lemma~\ref{lem:BDG}, $J_\pi=\int_0^tS^\pi(s)\Upsilon(s)\dd W(s)$ with
$S^\pi(s):=S^0(t,s_{l+1})$ on $(s_l,s_{l+1}]$. The isometry and \eqref{eq:psibdd} give
(i) and \eqref{eq:convbound} for $J_W(t)$. Moreover
\[
 \E\|J_\pi-J_W(t)\|^2_{\Hmid^n}=\int_0^t\big\|\big(S^\pi(s)-S^0(t,s)\big)\Upsilon(s)\big\|^2_{\HS(\mathsf U;\Hmid^n)}\dd s .
\]
For each $s$ and each vector $e$, $(S^\pi(s)-S^0(t,s))\Upsilon(s)e\to0$ in $\Hmid^n$ by the
continuity above, and it is bounded by $2c_{\rm emb}C_1(T)\|\Upsilon(s)e\|$. Dominated
convergence, first over an orthonormal basis of $\mathsf U$ and then in $s$, gives
$J_\pi\to J_W(t)$ in $L^2(\Hmid^n)$. The $J_\pi$ are functions of $M$, so the limit $J(t)$
has a version on $\Omega$, and $J(t)=J_W(t)$ almost surely. This gives (ii). The event
$\{J(t)\in\Ufine^n\}$ is Borel in $\Hmid^n$ by Lusin--Souslin, and has full measure on the
extension, which gives (iii). For (iv),
$\br{J_W(t),\varphi}=\int_0^t\Upsilon(s)^*S^0(t,s)^*\varphi\dd W(s)=\mathcal I_t(S^0(t,\cdot)^*\varphi)$.
An element of $L^2(\Omega;\Ufine^n)$ is determined by its pairings with a countable set
dense in $(\Ufine^*)^n$.
\end{proof}

\begin{theorem}[The $\Ufine$ maximal bound at the limit]\label{thm:maxlimit}
Let $M$ be as in Definition~\ref{def:weak}, and let
$\beta\in(\frac14,\frac12)$. Put $\kappa_\beta:=\big(\int_0^Tr^{\frac43(\beta-1)}\dd r\big)^{3/4}$,
$m_\beta:=\frac{T^{1-2\beta}}{1-2\beta}$, $c_\beta:=\frac{\sin\pi\beta}\pi$, and let
$C_4$ be the constant of Lemma~\ref{lem:BDG} at $p=4$. Then for every countable
$\mathbb D\subset[0,T]$,
\begin{equation}\label{eq:maxD}
 \E\Big[\sup_{t\in\mathbb D}\|J(t)\|^4_{\Ufine^n}\Big]
 \ \le\ c_\beta^4\,C_1(T)^8\,\kappa_\beta^4\,C_4\,C_{\mathsf q}^2\,m_\beta^2\,T\ =:\ C_\beta(T).
\end{equation}
\end{theorem}

\begin{proof}
We work on the extension of Proposition~\ref{prop:rep}; by
Proposition~\ref{prop:convint}, $J(t)=J_W(t)$ almost surely for each $t$, so the countable
supremum is the same for $J$ and $J_W$. We repeat the factorisation argument of
Proposition~\ref{prop:conv}, that is, the stochastic Fubini theorem, the bound of the operator
$R_\beta$, and Lemma~\ref{lem:BDG}, with $M^N$ replaced by $M=\int_0^\cdot\Upsilon\dd W$, with $\eta=0$, and
with $\mathsf q^N_i$ replaced by the deterministic $\mathsf q_i\le C_{s_1,a_1}\|\Gamma_i\|_{L^\infty}$. Put
$\mathsf Y(s):=\int_0^s(s-r)^{-\beta}S^0(s,r)\Upsilon(r)\dd W(r)$. The map
$(r,s)\mapsto(s-r)^{-\beta}\mathbf 1_{r<s}S^0(s,r)\Upsilon(r)$ is measurable into
$\HS(\mathsf U;\Ufine^n)$, by Lemma~\ref{lem:Smeas}\textup{(ii)} at $\eta=0$ and
Proposition~\ref{prop:rep}\textup{(i)}. So the proof of \cite[Thm.~4.33]{da_prato_2014_SPDEs} gives a
version of $\mathsf Y$ that is jointly measurable in $(s,\omega)$, and the stochastic Fubini
theorem gives $J_W(t)=c_\beta\int_0^t(t-s)^{\beta-1}S^0(t,s)\mathsf Y(s)\dd s$ almost surely for
each $t$. Lemma~\ref{lem:BDG} gives
$\E\|\mathsf Y(s)\|^4_{\Ufine^n}\le C_4\big(C_1(T)^2C_{\mathsf q} m_\beta\big)^2$, and H\"older's inequality
with exponents $(\frac43,4)$ gives
$\|J_W(t)\|_{\Ufine^n}\le c_\beta C_1(T)\kappa_\beta\|\mathsf Y\|_{L^4(0,T;\Ufine^n)}$ for every
$t$. Take the fourth power and the expectation.
\end{proof}

\section{The limit equation}\label{sec:limiteq}

\noindent
The limit equation is
\begin{equation}\label{eq:OU_SPDE}
 \d\fluc_i=\Big(\diver\big(\fluc_i\na U_i\big)
 +\Delta\sum_{j=1}^n\mathbf A_{ij}(u)\fluc_j\Big)\dt+\d M_i,
 \qquad \fluc_i(0)=\fluc^0_i,\qquad i=1,\dots,n,
\end{equation}
with $\mathbf A_{ij}(u)$ as in \eqref{eq:Amatrix}. By Theorem~\ref{thm:G}(ii), its drift is the
adjoint of the operator $\mathcal G$ in the weak form \eqref{eq:weakform}.
Definitions~\ref{def:weak} and~\ref{def:mild} give the two
formulations of \eqref{eq:OU_SPDE}. For a linear equation with a well-behaved evolution
family, a weak solution is mild by variation of constants. For the generator of a $C_0$-semigroup this is classical, see \cite[Ch.~5]{da_prato_2014_SPDEs}. That result does not apply here, since the drift depends on time and only the evolution family of Proposition~\ref{prop:evol} is available, not a generator. On the scale used here
it is not immediate. The mild formulation asserts an identity in $\Ufine^n$
simultaneously in $t$, whereas the weak formulation gives an identity against each
test function at each fixed $t$, almost surely. Theorem~\ref{thm:weakmild} shows that weak solutions are mild,
and Theorem~\ref{thm:uniqueness} gives uniqueness given the pair $(\fluc^0,M)$.
Corollary~\ref{cor:whole} passes from convergence of subsequences to convergence of
the whole sequence. For this it uses that the datum and the martingale of a limit
point are independent, with fixed Gaussian laws, so that their joint law is fixed.

\begin{lemma}[A weak solution has $\Hmid^n$-norm-continuous paths]\label{lem:flucHmid}
Let $\fluc$ be a weak solution \textup{(}Definition~\ref{def:weak}\textup{)}. Then $\Prob$-a.s.\ $t\mapsto\fluc(t)$ is continuous
$[0,T]\to\Hmid^n$ in norm.
\end{lemma}

\begin{proof}
Fix $\omega$ off the null set on which Definition~\ref{def:weak}'s two path conditions
fail, and let $t_k\to t$ in $[0,T]$. Weak continuity in $(\Usharp)^n$ says
$\br{\fluc(t_k),\psi}\to\br{\fluc(t),\psi}$ for every
$\psi\in(\Usharp)^{*n}=(H^{s_1+1}_{a_1})^n$, i.e.\ $\fluc(t_k)\rightharpoonup\fluc(t)$
weakly in the Hilbert space $(\Usharp)^n$. The embedding
$\Usharp\hookrightarrow\hookrightarrow\Hmid$ of \eqref{eq:chainA} is compact,
by $s_m>s_1+1$ and $a_m>a_1$. A compact operator carries weakly convergent
sequences to norm convergent ones, so $\fluc(t_k)\to\fluc(t)$ in $\Hmid^n$.
\end{proof}

\begin{definition}[Mild solution]\label{def:mild}
Let $(S^0(t,s))_{0\le s\le t\le T}$ be the evolution family on $\Ufine^n$ of
Proposition~\ref{prop:evol}\textup{(iii)} at $\eta=0$. It is dual to the backward system
\eqref{eq:backward-mild} with $\eta=0$, which, by \eqref{eq:Amatrix}, is
\begin{equation}\label{eq:backward0}
 \pa_s\Phi_j-\na U_j\cdot\na\Phi_j+\sum_{i=1}^n\mathbf A_{ij}(u)\Delta\Phi_i=0,\qquad\Phi(t)=\varphi.
\end{equation}
A \emph{mild solution} of \eqref{eq:OU_SPDE} with initial datum $\fluc^0$ is an
$\Ufine^n$-valued adapted process with continuous $\Vcoarse^n$-valued paths for which there is a
process $J$ with continuous $\Vcoarse^n$-valued paths such that, for each $t$, $J(t)$ equals
almost surely the stochastic convolution $\int_0^tS^0(t,s)\dd M(s)$ of
Proposition~\ref{prop:convint}, and, $\Prob$-a.s.\ for all $t\in[0,T]$,
\[
 \fluc(t)=S^0(t,0)\fluc^0+J(t).
\]
\end{definition}

\begin{theorem}[Weak $\Rightarrow$ mild]\label{thm:weakmild}
Let $M$ be as in Definition~\ref{def:weak}'s preamble.
Let $\fluc$ be an adapted process such that, $\Prob$-a.s.,
\begin{enumerate}[label=\normalfont(M\arabic*),leftmargin=*]
\item\label{it:M1} $\fluc(t)\in(\Usharp)^n$ for every $t$ and
$R:=\sup_{t\le T}\|\fluc(t)\|_{(\Usharp)^n}<\infty$;
\item\label{it:M2} $t\mapsto\br{\fluc_j(t),\psi}$ is continuous on $[0,T]$ for
every $j$ and every $\psi\in(\Usharp)^*=H^{s_1+1}_{a_1}$;
\item\label{it:M3} \eqref{eq:weakform} holds for every $\varphi\in\Test$ and every
$t\in[0,T]$, $\Prob$-a.s.\ \textup{(}the per-pair quantifier of
Definition~\ref{def:weak}\textup{)}, for some $\Fcoarse^n$-valued
$\mathcal F_0$-measurable $\fluc^0$.
\end{enumerate}
Then:
\begin{enumerate}[label=\normalfont(\roman*)]
\item $\fluc^0=\fluc(0)\in(\Usharp)^n$ a.s., and for each fixed $t\in[0,T]$, $\Prob$-a.s.,
\begin{equation}\label{eq:sharpid}
 \fluc(t)\ =\ S^0_\sharp(t,0)\fluc^0+J(t)\qquad\text{in }(\Usharp)^n,
\end{equation}
hence in $\Hmid^n$ and in $\Vcoarse^n$; here $S^0_\sharp$ is the $\Usharp$ family of
Proposition~\ref{prop:evol} at $\eta=0$ and $J(t)$ is Proposition~\ref{prop:convint}'s object. If
moreover $\fluc^0\in\Ufine^n$, then $S^0_\sharp(t,0)\fluc^0=S^0(t,0)\fluc^0$ by
Lemma~\ref{lem:levelcoin}\textup{(ii)}.
\item A weak solution satisfies \ref{it:M1}--\ref{it:M3}, and for every $t\in[0,T]$,
$\Prob$-a.s.,
\begin{equation}\label{eq:mild-mild}
 \fluc(t)=S^0(t,0)\fluc^0+\int_0^tS^0(t,s)\dd M(s)\qquad\text{in }\Ufine^n.
\end{equation}
\item For a weak solution, $\widehat J:=\fluc-S^0(\cdot,0)\fluc^0$ has $\Prob$-a.s.\
$\Hmid^n$-continuous paths and is the $\Hmid^n$-path-continuous version of $J$, unique up
to indistinguishability. Off a single null set,
\begin{equation}\label{eq:mildsimult}
 \fluc(t)\ =\ S^0(t,0)\fluc^0+\widehat J(t)\qquad\text{in }\Ufine^n,
 \qquad\text{for \textbf{every} }t\in[0,T],
\end{equation}
and $\widehat J(t)\in\Ufine^n$ for every $t$. So $\fluc$ is a mild solution in the sense of
Definition~\ref{def:mild}.
\end{enumerate}
Consequently two weak solutions with the same $\fluc^0$ and the same $M$ are
indistinguishable.
\end{theorem}

\begin{proof}
Step~0 identifies the datum with $\fluc(0)$. Step~1 tests $\fluc$ against the backward
solution $\Phi(s)=S^0(t,s)^*\varphi$ along the dyadic partitions of $[0,t]$. A telescoping sum
splits the increment into a drift error, which tends to $0$ almost surely, and a martingale sum,
which converges in $L^2(\Omega)$ to the scalar integral of Lemma~\ref{lem:scalint}. Step~2
identifies the limit in $(\Usharp)^n$, which is (i), and Steps~3 and~4 prove (ii) and (iii).

\emph{Step 0: $\fluc^0=\fluc(0)$.} Taking $t=0$ in \eqref{eq:weakform} gives
$\br{\fluc_i(0)-\fluc^0_i,\varphi}=0$ a.s.\ for each $\varphi\in\Test$. Take a
countable $\{\varphi_p\}\subset C^\infty_0$ dense in $\Test$. Here $\fluc^0_i\in\Test^*$ only. Such a set exists, since $C^\infty_0$ is dense in every $H^s$ and $\br\cdot^{a}$
maps $H^s_a$ isometrically onto $H^s$ and $C^\infty_0$ onto itself. The density
gives $\fluc_i(0)=\fluc^0_i$ in $\Test^{*}$. Since $\fluc(0)\in(\Usharp)^n$ by
\ref{it:M1}, and $(\Usharp)^n\hookrightarrow\Test^{*n}$ is injective because both
are spaces of tempered distributions with the identity as embedding, the equality
holds in $(\Usharp)^n$ and $\fluc^0\in(\Usharp)^n$.

\emph{Step 1: the decomposition.} Fix $t\in(0,T]$, $\varphi\in\Test^n$, and set
$\Phi(s):=S^0(t,s)^*\varphi$. By Lemma~\ref{lem:levelcoin}\textup{(i)} this is also the
adjoint of the $\Usharp$ family applied to $\varphi$, and every pairing below is of an element of
$(\Usharp)^n$ against an element of $(\Usharp)^*$. By
Proposition~\ref{prop:evol} and Lemma~\ref{lem:strongcont},
$\Phi\in C([0,t];\Test^n)$ solves \eqref{eq:backward-mild} with $\eta=0$.
Let $(\pi_k)$ be the dyadic partitions of $[0,t]$. Off one null set, \eqref{eq:weakform} holds
for the countably many pairs $(\tau,\Phi_i(\varsigma))$ with $\tau,\varsigma$ dyadic nodes of $[0,t]$
and $i\le n$. On it, Lemma~\ref{lem:telescope} at $\eta=0$ applies $\omega$ by $\omega$, with
$v:=\fluc$, $v^0:=\fluc^0$ and $\mathcal V_i(\tau,\chi):=M_i(\tau,\chi)$: its hypotheses are
\ref{it:M1}--\ref{it:M3}, and $\mathcal G^0=\mathcal G$ by the convention at $\eta=0$. It gives
\[
 \sum_i\br{\fluc_i(t),\varphi_i}-\sum_i\br{\fluc_i(0),\Phi_i(0)}=\mathcal E_{\pi_k}+\mathcal M_{\pi_k},
 \qquad \mathcal E_{\pi_k}\to0\ \text{a.s.},
\]
with $\mathcal E_\pi$ as in \eqref{eq:Epi} and
$\mathcal M_\pi:=\mathcal V_\pi=\sum_l\sum_i\big(M_i(s_{l+1},\Phi_i(s_{l+1}))-M_i(s_l,\Phi_i(s_{l+1}))\big)$,
where $s_l$ are the nodes of $\pi$.

\emph{$\mathcal M_\pi\to\sum_i\int_0^t\br{\dd M_i,\Phi_i}$ in $L^2(\Omega)$.}
$\Phi^\pi:=\sum_l\Phi(s_{l+1})\mathbf 1_{(s_l,s_{l+1}]}$ is left-continuous with
deterministic coefficients, hence predictable, and by linearity of
$\varphi\mapsto M_i(t,\varphi)$ one has $\mathcal M_\pi=\sum_i\int_0^t
\br{\dd M_i,\Phi^\pi_i}$ exactly. The equality below is the isometry of Lemma~\ref{lem:scalint}\textup{(ii)}
applied to $\Phi^\pi-\Phi$. The cross terms
$i\ne k$ vanish by the factor $\delta_{ik}$ in \eqref{eq:limbracket}. So
\[
 \E\Big|\mathcal M_\pi-\sum_i\int_0^t\br{\dd M_i,\Phi_i}\Big|^2
 =\sum_i\E\int_0^t\br{u_i,\big|\na(\Phi^\pi_i-\Phi_i)\big|^2\Gamma_i}\dd s
 \le\Big(\sum_i\int_0^t\mathsf q_i\dd s\Big)\sup_s\|\Phi^\pi(s)-\Phi(s)\|^2_{H^{s_1}_{a_1}},
\]
the inequality by \eqref{eq:kernelptwise} and the definition of $\mathsf q_i$. The last
factor tends to $0$ by uniform continuity of $\Phi$ in
$\Test\subset H^{s_1}_{a_1}$ and the limit is $\mathcal I_t(\Phi)$ of Lemma~\ref{lem:scalint}.

\emph{Step 2: the identification at $\Usharp$, which is (i).} The left-hand side of Step~1 is
$\pi$-free, so $\mathcal M_{\pi_k}$ converges a.s.\ as well and the two limits
agree, giving
\[
 \sum_i\br{\fluc_i(t),\varphi_i}
 =\sum_i\br{\fluc_i(0),\Phi_i(0)}+\sum_i\int_0^t\br{\dd M_i,\Phi_i}.
\]
For the first term, $\fluc(0)=\fluc^0\in(\Usharp)^n$ by Step~0 and
$\Phi(0)=S^0_\sharp(t,0)^*\varphi$. So
$\sum_i\br{\fluc_i(0),\Phi_i(0)}=\br{S^0_\sharp(t,0)\fluc^0,\varphi}$ by the
$\Usharp$ clause of Proposition~\ref{prop:evol}, since $\Phi(0)$ is the adjoint of
the $\Usharp$ family by Lemma~\ref{lem:levelcoin}\textup{(i)}. For the second,
$\sum_i\int_0^t\br{\dd M_i,\Phi_i}=\br{J(t),\varphi}$ is \eqref{eq:convint},
available for $\varphi\in\Ufine^{*n}\supset(\Usharp)^{*n}\supset\Test^n$.

Now take $\{\varphi_p\}$ countable and dense in $(\Usharp)^{*n}$ with
$\varphi_p\in(C^\infty_0)^n$. Off a null set, the two identities above hold for every
$p$, and $\fluc(t)$ and
$S^0_\sharp(t,0)\fluc^0+J(t)$ are two elements of $(\Usharp)^n$ that agree
against a dense subset of $(\Usharp)^*$, hence are equal. The first lies in
$(\Usharp)^n$ by \ref{it:M1}, the second because $S^0_\sharp(t,0)\fluc^0\in(\Usharp)^n$ and
$J(t)\in\Ufine^n\subset(\Usharp)^n$ by
Proposition~\ref{prop:convint}\textup{(iii)} and \eqref{eq:chainA}. This is
\eqref{eq:sharpid}. The remaining inclusions are \eqref{eq:chainA}.

\emph{Step 3: (ii).} A weak solution satisfies \ref{it:M1}, since $\fluc(t)\in\Ufine^n\subset(\Usharp)^n$
with $\sup_t\|\fluc(t)\|_{\Usharp}<\infty$, \ref{it:M2}, which is its weak continuity in
$(\Usharp)^n$ and \ref{it:M3}, with $\fluc^0\in\Ufine^n\subset\Fcoarse^n$. By (i), \eqref{eq:sharpid}
holds with $S^0_\sharp(t,0)\fluc^0=S^0(t,0)\fluc^0$. Here $\fluc(t)\in\Ufine^n$ by
Definition~\ref{def:weak}, $S^0(t,0)\fluc^0\in\Ufine^n$ by Proposition~\ref{prop:evol}, and
$J(t)\in\Ufine^n$ by Proposition~\ref{prop:convint}\textup{(iii)}. The embedding
$\Ufine^n\hookrightarrow(\Usharp)^n$ is injective, so the identity holds in $\Ufine^n$.

\emph{Step 4: (iii).} $t\mapsto S^0(t,0)\fluc^0$ is $\Hmid^n$-norm continuous by
Proposition~\ref{prop:fwdHmid} at $s=0$, $\eta=0$, with
$\nu=\fluc^0(\omega)\in\Ufine^n$.
With Lemma~\ref{lem:flucHmid}, $\widehat J$ has $\Hmid^n$-continuous paths. By
(ii) it is a version of $J$. Two versions with continuous paths
agree on a countable dense set of times off one null set, hence everywhere, so it is
unique up to indistinguishability. \eqref{eq:mildsimult} is the definition of
$\widehat J$, rearranged, so it holds for every $\omega$ and $t$. Both $\fluc(t)$ and $S^0(t,0)\fluc^0$ lie
in $\Ufine^n$ at every $t$, by Definition~\ref{def:weak} and
Proposition~\ref{prop:evol}. So $\widehat J(t)\in\Ufine^n$, and the identity holds
between $\Ufine^n$-valued random variables at every $t$. 

Finally, Definition~\ref{def:mild} requires an $\Ufine^n$-valued adapted
process with continuous $\Vcoarse^n$ paths satisfying its identity $\Prob$-a.s.\ for
all $t$, with the convolution of Proposition~\ref{prop:convint}. $\Ufine^n$-valuedness
and adaptedness hold by Definition~\ref{def:weak}. The paths are continuous by
Lemma~\ref{lem:flucHmid} and $\Hmid\hookrightarrow\Vcoarse$. The identity for all
$t$ is \eqref{eq:mildsimult}, and $\widehat J$ is that convolution: by
the first part of this step it is its unique path-continuous version, and it agrees with it a.s.\
at each fixed $t$, which is how Definition~\ref{def:mild} identifies the
convolution.

The last claim follows from (ii) and Lemma~\ref{lem:flucHmid}: the two processes are both continuous in $\Hmid^n$, and
by \eqref{eq:mild-mild} they agree at each fixed $t$ almost surely.
\end{proof}

\begin{corollary}[What the $C^\infty_0$ weak identity gives]\label{cor:aclosed}
Let $\fluc$ be an adapted process with $\Prob$-a.s.\ $\Hmid^n$-continuous paths
and $M$ as in Definition~\ref{def:weak}, satisfying \ref{it:M1} and
\begin{enumerate}[label=\normalfont(M3$'$),leftmargin=*]
\item\label{it:M3p} \eqref{eq:weakform} holds for every $\varphi\in C^\infty_0$
and every $t\in[0,T]$, $\Prob$-a.s., for some $\Fcoarse^n$-valued
$\mathcal F_0$-measurable $\fluc^0$.
\end{enumerate}
In \textup{(iii)} let
$\beta\in(\frac14,\frac12)$.
Then:
\begin{enumerate}[label=\normalfont(\roman*)]
\item \ref{it:M3} holds, together with \textup{(ii)}, Theorem~\ref{thm:weakmild}(i) applies;
\item \ref{it:M2} holds, so \eqref{eq:sharpid} holds for every
fixed $t$;
\item if moreover $\fluc^0\in\Ufine^n$ and $\E\|\fluc^0\|^4_{\Ufine^n}<\infty$, then, off one
null set, $\fluc(t)\in\Ufine^n$ for every $t$, and
\begin{equation}\label{eq:LinftyUfine}
 \E\Big[\sup_{t\le T}\|\fluc(t)\|^4_{\Ufine^n}\Big]
 \ \le\ 8\,C_1(T)^4\,\E\|\fluc^0\|^4_{\Ufine^n}+8\,C_\beta(T),
\end{equation}
with $C_\beta(T)$ of \eqref{eq:maxD}.
\end{enumerate}
\end{corollary}

\begin{proof}
Part (i) extends the weak identity from $C^\infty_0$ to $\Test$ by density, and part (ii)
upgrades $\Hmid^n$-continuity to weak continuity in $(\Usharp)^n$ with the bound \ref{it:M1}. Part
(iii) combines Theorem~\ref{thm:weakmild}(i) with the maximal bound of Theorem~\ref{thm:maxlimit}
on a countable set of times, and Corollary~\ref{cor:supupgrade} extends the bound to all times.

\textup{(i)} Fix $t\in[0,T]$ and $\varphi\in\Test$, and take $\varphi_k\in
C^\infty_0$ with $\varphi_k\to\varphi$ in $\Test$, possible because $C^\infty_0$
is dense in $\Test=H^{s_1+3}_{a_1+1}$: $C^\infty_0$ is dense in $H^{s_1+3}$, and
$\br\cdot^{a_1+1}$ maps $H^{s_1+3}_{a_1+1}$ isometrically onto $H^{s_1+3}$ and $C^\infty_0$ onto
itself. Each term of
\eqref{eq:weakform} converges. For the two pairings:
$\br{\fluc_i(t),\varphi_k}\to\br{\fluc_i(t),\varphi}$ because
$\fluc(t)\in(\Usharp)^n$ by \ref{it:M1} and
$\Test\hookrightarrow(\Usharp)^*$; and
$\br{\fluc^0_i,\varphi_k}\to\br{\fluc^0_i,\varphi}$ because $\fluc^0_i\in\Test^*$
by \ref{it:M3p}.
For the drift, the integrand is bounded by
$n\,R\,C_{\mathcal G}\|\varphi_k-\varphi\|_{\Test}$
with $R$ from \ref{it:M1} and $C_{\mathcal G}$ of Lemma~\ref{lem:Gbound} at $\eta=0$, uniformly in $s$. And $M_i(t,\varphi_k)\to M_i(t,\varphi)$ by
continuity of $\varphi\mapsto M_i(t,\varphi)$ on $\Ufine^*\supset\Test$
\textup{(}Definition~\ref{def:weak}'s preamble, the inclusion by
Lemma~\ref{lem:spaces}\textup{)}.

\textup{(ii)} For $\psi\in(\Usharp)^*$ and $\chi\in\Hmid^*$,
$|\br{\fluc_j(s)-\fluc_j(s_0),\psi}|\le|\br{\fluc_j(s)-\fluc_j(s_0),\chi}|
+2R\|\psi-\chi\|_{(\Usharp)^*}$ with $R$ from \ref{it:M1}; $\Hmid^*$ is dense in
$(\Usharp)^*$ and the first
term tends to $0$ by $\Hmid^n$-continuity. Choosing first $\chi$ close to $\psi$ and then $s$
close to $s_0$ gives \ref{it:M2}.

\textup{(iii)} By \textup{(i)} and \textup{(ii)},
Theorem~\ref{thm:weakmild}(i) applies. Since $\fluc^0\in\Ufine^n$, it gives
$\fluc(t)=S^0(t,0)\fluc^0+J(t)$ in $\Hmid^n$ for each $t$, almost surely. Put
$\widehat J:=\fluc-S^0(\cdot,0)\fluc^0$. It has continuous $\Hmid^n$-paths, by hypothesis and
Proposition~\ref{prop:fwdHmid}. Fix a countable dense $\mathbb D\subset[0,T]$.
Then $\widehat J(t)=J(t)$ for all $t\in\mathbb D$ off one null set, so
Theorem~\ref{thm:maxlimit} and Corollary~\ref{cor:supupgrade} give
$\E\sup_{t\le T}\|\widehat J(t)\|^4_{\Ufine^n}\le C_\beta(T)$, and $\widehat J(t)\in\Ufine^n$
for every $t$ off one null set. Moreover $\|S^0(t,0)\fluc^0\|_{\Ufine^n}\le
C_1(T)\|\fluc^0\|_{\Ufine^n}$ by Proposition~\ref{prop:evol}. Then
$(a+b)^4\le8(a^4+b^4)$ gives \eqref{eq:LinftyUfine}.
\end{proof}

\begin{theorem}[Fluctuation SPDE]\label{thm:OU}
Let $(\widetilde\fluc,\widetilde M)$ be a limit
point of Proposition~\ref{prop:skorokhod}, with the filtration
$(\widetilde{\mathcal F}_t)$ generated by $(\widetilde\fluc,\widetilde M)$, augmented by the
$\widetilde\Prob$-null sets and made right-continuous. Then:
\begin{enumerate}[label=\normalfont(\roman*)]
\item $\widetilde M$ is the Gaussian martingale of Theorem~\ref{thm:noise}.
$\widetilde\fluc(0)$ is $\widetilde{\mathcal F}_0$-measurable, lies in $\Ufine^n$, and
has the law of the Gaussian field $\fluc^0$ of Proposition~\ref{prop:initialCLT};
\item $\widetilde\fluc$ is a weak solution, in the sense of Definition~\ref{def:weak}
with datum $\widetilde\fluc(0)$ and martingale $\widetilde M$, of \eqref{eq:OU_SPDE};
\item off one null set, $\widetilde\fluc(t)\in\Ufine^n$ for every $t\in[0,T]$, and
\begin{equation}\label{eq:Linftyexplicit}
 \E\Big[\sup_{t\in[0,T]}\big\|\widetilde\fluc(t)\big\|^4_{\Ufine^n}\Big]
 \ \le\ 24\,C_1(T)^4n^2C_{s_1,a_1}^2+8\,C_\beta(T)\ <\ \infty
\end{equation}
for every $\beta\in(\frac14,\frac12)$, with $C_\beta(T)$ of \eqref{eq:maxD};
\item $\widetilde\fluc$ is the mild solution of Definition~\ref{def:mild} with datum
$\widetilde\fluc(0)$ and martingale $\widetilde M$.
\end{enumerate}
\end{theorem}

\begin{proof}
We verify the hypotheses of Corollary~\ref{cor:aclosed} for the limit point: the bound
\ref{it:M1} by lower semicontinuity (Step~0), and the weak identity on $C^\infty_0$ by passing to the
limit in \eqref{eq:Rtilde} with Lemma~\ref{lem:drift} (Step~1). Step~2 identifies the datum and the
martingale. Corollary~\ref{cor:aclosed} and Theorem~\ref{thm:weakmild} then give (iii), (ii) and
(iv), in this order.
Fix $\beta\in(\frac14,\frac12)$.

\emph{Step 0: the bound \ref{it:M1}.} By Theorem~\ref{thm:K4}, \eqref{eq:K4-mild}
holds, so Proposition~\ref{prop:close} gives
$\sup_N\E\sup_{t\le T}\|\fluc^N(t)\|^2_{\Usharp}<\infty$. The map
$v\mapsto\sup_{t\le T}\|v(t)\|_{\Usharp}$ is lower semicontinuous on
$C([0,T];\Hmid)$ (Lemma~\ref{lem:lscsup}(ii)), so Fatou's lemma and
$\widetilde\fluc^N\to\widetilde\fluc$ in $C([0,T];\Hmid^n)$ give
$\E\sup_{t\le T}\|\widetilde\fluc(t)\|^2_{\Usharp}<\infty$. This is \ref{it:M1}.

\emph{Step 1: the weak identity for $\varphi\in C^\infty_0$.} Fix
$\varphi\in C^\infty_0(\R^d)$ and $t\le T$. By the definition \eqref{eq:Rtilde},
\[
 \br{\widetilde\fluc^N_i(t),\varphi}=\br{\widetilde\fluc^N_i(0),\varphi}
 +\int_0^t\sum_j\br{\widetilde\fluc^N_j,(\mathcal G^\eta\varphi)_{ij}}\ds
 +\widetilde{\mathcal E}^N_i(t,\varphi)+\widetilde M^N_i(t,\varphi).
\]
Along the subsequence of Proposition~\ref{prop:skorokhod}, the first, second and last
terms converge almost surely to the corresponding limits, and Lemma~\ref{lem:drift}
gives the convergence of the sum of the other two in $L^1(\widetilde\Omega)$. Hence
\eqref{eq:weakform} holds for $\varphi\in C^\infty_0$, with $\fluc^0:=\widetilde\fluc(0)$
and $M:=\widetilde M$. By Lemma~\ref{lem:initmom}(ii), $\widetilde\fluc(0)$ lies in
$\Ufine^n\subset\Fcoarse^n$ almost surely, by \eqref{eq:chainA}, and it is Borel into
$\Ufine^n$, hence into $\Fcoarse^n$. We set it to $0$ on the exceptional null set. This
changes no almost-sure statement. So $\fluc^0$ is an $\Fcoarse^n$-valued
$\widetilde{\mathcal F}_0$-measurable datum, and this is \ref{it:M3p}.

\emph{Step 2: (i).} By Theorem~\ref{thm:noise}, $\widetilde M$ is a continuous Gaussian
martingale with the bracket \eqref{eq:cov_limit}. It remains one for the augmented
right-continuous filtration, with the same brackets, since its paths are continuous. By
Fatou's lemma, lower semicontinuity and \eqref{eq:M2},
$\E\sup_t\|\widetilde M(t)\|^2_{\Hmid^n}\le c_{\rm emb}^2\liminf_N\E\sup_t\|M^N(t)\|^2_{\Ufine^n}<\infty$.
Also $\widetilde M(0)=0$, since $M^N(0)=0$.
So Proposition~\ref{prop:rep} applies with $H=\Hmid^n$, $\Gamma_i=2\Theta_i$, and the dense
set $H^{s_m+1}_{a_m}\subset\Hmid^*$ of \eqref{eq:cov_limit}. By its part \textup{(iii)},
$\widetilde M$ satisfies the requirements on $M$ in Definition~\ref{def:weak}. $\widetilde\fluc(0)$ is
$\widetilde{\mathcal F}_0$-measurable by definition of the filtration, and lies in
$\Ufine^n$ by Lemma~\ref{lem:initmom}(ii). Finally
$\mathrm{Law}(\widetilde\fluc^N(0))=\mathrm{Law}(\fluc^N(0))$, which converges to
$\mathrm{Law}(\fluc^0)$ by Proposition~\ref{prop:initialCLT}, while
$\widetilde\fluc^N(0)\to\widetilde\fluc(0)$ almost surely, so
$\mathrm{Law}(\widetilde\fluc(0))=\mathrm{Law}(\fluc^0)$.

\emph{Step 3: (iii).} Corollary~\ref{cor:aclosed}(iii) applies: its hypotheses are
\ref{it:M1}, \ref{it:M3p}, adaptedness, the $\Hmid^n$-continuous paths,
$\beta\in(\frac14,\frac12)$, and the datum
$\widetilde\fluc(0)\in\Ufine^n$ with a finite fourth moment. By Lemma~\ref{lem:initmom}, applied through the
law identity of Step~2, $\E\|\widetilde\fluc(0)\|^4_{\Ufine^n}\le3n^2C^2_{s_1,a_1}$,
so \eqref{eq:LinftyUfine} gives \eqref{eq:Linftyexplicit}.

\emph{Step 4: (ii).} Proposition~\ref{prop:skorokhod} gives $\Hmid^n$-continuous
paths, so Corollary~\ref{cor:aclosed}(ii) gives \ref{it:M2}. With \ref{it:M1} this is
the weak continuity in $(\Usharp)^n$ and the $\Usharp$ bound that
Definition~\ref{def:weak} requires. Corollary~\ref{cor:aclosed}(i) extends the identity
of Step~1 from $C^\infty_0$ to $\Test$. The process is adapted to
$(\widetilde{\mathcal F}_t)$ by construction. Its $\Ufine^n$-valuedness at every $t$ is
(iii), proved in Step~3. We redefine $\widetilde\fluc$ to be $0$ on the
null set of (iii). The filtration is augmented, so the modified process is still adapted. Measurability as an
$\Ufine^n$-valued map follows from the Lusin--Souslin theorem, since $\Ufine^n\hookrightarrow\Hmid^n$ is a
continuous injection between Polish spaces.

\emph{Step 5: (iv).} By (ii) and Theorem~\ref{thm:weakmild}, $\widetilde\fluc$ satisfies
\eqref{eq:mild-mild} for each fixed $t$, almost surely. By
Theorem~\ref{thm:weakmild}(iii) the identity holds simultaneously in $t$, which is the
requirement of Definition~\ref{def:mild}.
\end{proof}

\begin{theorem}[Uniqueness]\label{thm:uniqueness}
Two weak solutions of \eqref{eq:OU_SPDE}, in the sense
of Definition~\ref{def:weak}, with the same datum $\fluc^0$ and the same martingale $M$
are indistinguishable as $\Hmid^n$-valued processes. Two mild solutions with the same
$(\fluc^0,M)$ are indistinguishable as $\Vcoarse^n$-valued processes.
\end{theorem}

\begin{proof}
By Theorem~\ref{thm:weakmild}, for each $t$ both solutions equal
$S^0(t,0)\fluc^0+\int_0^tS^0(t,s)\dd M(s)$ almost surely. Both have
$\Hmid^n$-continuous paths by Lemma~\ref{lem:flucHmid}, so they agree for all $t$ off
one null set. Two mild solutions agree almost surely at each fixed $t$, by
Definition~\ref{def:mild}, and their paths are continuous in $\Vcoarse^n$.
\end{proof}

\begin{lemma}[The law of a mild solution]\label{lem:mildlaw}
Let $M$ be as in Definition~\ref{def:weak}, let $\fluc^0$ be an $\Ufine^n$-valued random variable,
and let $\fluc$ be a process with continuous $\Hmid^n$-valued paths such that, for each $t\in[0,T]$,
$\fluc(t)=S^0(t,0)\fluc^0+J(t)$ almost surely, with $J(t)$ of Proposition~\ref{prop:convint}. Then
the law of $(\fluc,M)$ on $C([0,T];\Hmid^n)^2$ is determined by the law of $(\fluc^0,M)$ on
$\Hmid^n\times C([0,T];\Hmid^n)$. If the latter is centred Gaussian, so is the former.
\end{lemma}

\begin{proof}
We approximate the finite-dimensional distributions of $(\fluc,M)$ by those of the Riemann sums
of Proposition~\ref{prop:convint}, which are functions of $(\fluc^0,M)$. For $t\in(0,T]$, let
$(\pi_k)$ be partitions $0=s^k_0<\dots<s^k_{L_k}=t$ of $[0,t]$ with $|\pi_k|\to0$. By
Proposition~\ref{prop:convint}\textup{(ii)}, $J(t)$ is the limit in $L^2(\Omega;\Hmid^n)$ of
$J_{\pi_k}(t)=\sum_lS^0(t,s^k_{l+1})\big(M(s^k_{l+1})-M(s^k_l)\big)$. For $t=0$ put
$J_{\pi_k}(0):=0=J(0)$. Fix $0\le t_1<\dots<t_p\le T$. The values $M(s)$ and $\fluc^0$ lie in
$\Ufine^n$, and $S^0(t,s)$ is bounded on $\Ufine^n$, so each $S^0(t_q,0)\fluc^0+J_{\pi_k}(t_q)$ is a
Borel function of finitely many of these values. Since $\Ufine^n\hookrightarrow\Hmid^n$ is
continuous and injective, the Borel sets of $\Ufine^n$ are Borel in $\Hmid^n$ (Lusin--Souslin), so
the $\Ufine^n$-valued joint law of these values is determined by $\mathrm{Law}(\fluc^0,M)$. So the
law of $\big(S^0(t_q,0)\fluc^0+J_{\pi_k}(t_q)\big)_{q\le p}$ is determined by
$\mathrm{Law}(\fluc^0,M)$, and it converges, as $k\to\infty$, to the law of $(\fluc(t_q))_{q\le p}$.
The same holds jointly with $(M(t_q))_{q\le p}$, since the approximants and these values are Borel
functions of the same finitely many values. A law on $C([0,T];\Hmid^n)^2$ is determined by its
finite-dimensional distributions, which proves the first claim. If $(\fluc^0,M)$ is centred
Gaussian, the finite-dimensional distributions of $(\fluc,M)$ are limits in $L^2$ of linear images
of $(\fluc^0,M)$, hence centred Gaussian.
\end{proof}

\begin{corollary}[Convergence of the whole sequence]\label{cor:whole}
$(\fluc^N,M^N)\Rightarrow(\fluc,M)$ in $C([0,T];\Hmid^n)^2$, where $\fluc$ is the mild solution of \eqref{eq:OU_SPDE} whose
datum $\fluc^0$ has the law of Proposition~\ref{prop:initialCLT}, whose martingale $M$
is a continuous centred Gaussian martingale with the bracket \eqref{eq:cov_limit}, and
where $\fluc^0$ and $M$ are independent. The law of $\fluc$ is centred Gaussian.
\end{corollary}

\begin{proof}
We show that all limit points have the same law. Step~1 determines the joint law of the
datum and the martingale of a limit point, by their independence (Theorem~\ref{thm:noise}).
Step~2 shows, by Lemma~\ref{lem:mildlaw}, that this joint law determines the law of the limit point.
By Proposition~\ref{prop:tight} the laws of $(\fluc^N,M^N)$ are tight in
$C([0,T];\Hmid^n)^2$, a Polish space. It suffices to show that every limit point has the same law. Let a
subsequence converge in law in $C([0,T];\Hmid^n)^2$. Proposition~\ref{prop:tight} and
Proposition~\ref{prop:skorokhod}, applied to this subsequence, give a further subsequence and a
limit point $(\widetilde\fluc,\widetilde M)$ whose law on $C([0,T];\Hmid^n)^2$ is the given limit.
We show that this law does not depend on the subsequence.

\emph{Step 1: the joint law of $(\widetilde\fluc(0),\widetilde M)$.} By
Theorem~\ref{thm:noise}, $\widetilde M$ is a continuous martingale for the filtration
$(\widetilde{\mathcal F}_t)$ generated by $(\widetilde\fluc,\widetilde M)$, with
deterministic brackets, and $\widetilde M$ is
independent of $\widetilde{\mathcal F}_0$, which contains $\sigma(\widetilde\fluc(0))$.
The law of $\widetilde M$ on $C([0,T];\Hmid^n)$ is determined by its finite-dimensional
distributions, which are centred Gaussian with covariances given by
\eqref{eq:cov_limit}. The law of $\widetilde\fluc(0)$ is that of $\fluc^0$, by
Theorem~\ref{thm:OU}(i). Hence
$\mathrm{Law}(\widetilde\fluc(0),\widetilde M)=\mathrm{Law}(\fluc^0)\otimes\mathrm{Law}(M)$,
the same for every limit point.

\emph{Step 2: the law of $(\widetilde\fluc,\widetilde M)$.} By Theorem~\ref{thm:OU}(iv), for each
$t$, almost surely, $\widetilde\fluc(t)=S^0(t,0)\widetilde\fluc(0)+J(t)$, and $\widetilde\fluc$ has
continuous $\Hmid^n$-valued paths. By Theorem~\ref{thm:OU}(i) and (ii), $\widetilde M$ is as in
Definition~\ref{def:weak} and $\widetilde\fluc(0)$ is $\Ufine^n$-valued. So
Lemma~\ref{lem:mildlaw} applies: $\mathrm{Law}(\widetilde\fluc,\widetilde M)$ is determined by
$\mathrm{Law}(\widetilde\fluc(0),\widetilde M)$, which does not depend on the limit point, by Step~1.
The pair $(\widetilde\fluc(0),\widetilde M)$ is centred Gaussian, since its components are centred
Gaussian and independent. By Lemma~\ref{lem:mildlaw}, the law of $(\widetilde\fluc,\widetilde M)$
is centred Gaussian.
\end{proof}

\begin{definition}[Martingale problem]\label{def:mgp}
Let $\mu_0$ be a Borel probability measure on $\Ufine^n$. A solution of the martingale
problem for \eqref{eq:OU_SPDE} with initial law $\mu_0$ is a filtered probability space
satisfying the usual conditions and an adapted process $\fluc$ with $\Prob$-a.s.\ continuous
$\Hmid^n$-valued paths such that: $\fluc(0)$ is $\Ufine^n$-valued with law $\mu_0$,
$\Prob$-a.s., $\fluc(t)\in(\Usharp)^n$ for every $t$ and $\sup_{t\le T}\|\fluc(t)\|_{\Usharp}<\infty$
and for every $i$ and $\varphi\in\Test$ the process
\[
 M_i(t,\varphi):=\br{\fluc_i(t),\varphi}-\br{\fluc_i(0),\varphi}
 -\int_0^t\sum_j\br{\fluc_j(s),(\mathcal G(s)\varphi)_{ij}}\dd s
\]
is a continuous martingale with the brackets \eqref{eq:cov_limit} for $\varphi,\psi\in\Test$.
\end{definition}

\begin{corollary}[Uniqueness in law]\label{cor:mgp}
All solutions of the martingale problem with the same initial law $\mu_0$ have the same
law on $C([0,T];\Hmid^n)$. The limit $\fluc$ of Corollary~\ref{cor:whole} is a solution with
$\mu_0=\mathrm{Law}(\fluc^0)$.
\end{corollary}

\begin{proof}
Given a solution, we build from the scalar martingales $M_i(\cdot,\varphi)$ an $\Ufine^n$-valued
martingale (Step~1), check the hypotheses of Theorem~\ref{thm:weakmild} (Step~2), and conclude
with Lemmas~\ref{lem:gaussmart} and~\ref{lem:mildlaw} (Step~3).

Let $\fluc$ be a solution. \emph{Step 1: a $\Ufine^n$-valued martingale.} Choose an orthonormal basis $(e_m)$ of $H^{s_1}_{a_1}$ made of
Schwartz functions whose span is dense in $\Test$: apply Gram--Schmidt in $H^{s_1}_{a_1}$ to a
sequence in $\mathcal S$ that is dense in $\Test$. By Doob's inequality and
\eqref{eq:kernel},
$\sum_m\E\sup_t|M_i(t,e_m)|^2\le4\sum_m\E\br{M_i(\cdot,e_m)}_T=4\int_0^T\mathsf q_i\dd r<\infty$, with
$\mathsf q_i$ as in Definition~\ref{def:weak} for $\Gamma_i=2\Theta_i$. So, almost surely, the series
$\overline M_i(t):=\sum_mM_i(t,e_m)y_m$, with $(y_m)$ the dual basis of $\Ufine$, converges in
$\Ufine$ uniformly in $t$. It is a continuous square-integrable $\Ufine$-valued martingale, and
$\br{\overline M_i(t),\varphi}=M_i(t,\varphi)$ for $\varphi$ in the span of $(e_m)$, hence, by
$L^2(\Prob)$-continuity of both sides in $\varphi\in\Test$, for every $\varphi\in\Test$, almost surely.

\emph{Step 2: the hypotheses of Theorem~\ref{thm:weakmild}.} By Proposition~\ref{prop:rep}(iii) with $H=\Ufine^n$ and the dense
subset $\Test$, $\overline M$ satisfies the requirements on $M$ in Definition~\ref{def:weak}.
By Step~1 and the definition of $M_i(t,\varphi)$, \eqref{eq:weakform} holds for every
$\varphi\in\Test$ and every $t$, almost surely, with datum $\fluc(0)\in\Ufine^n\subset\Fcoarse^n$ and
martingale $\overline M$. This is \ref{it:M3}, and it implies \ref{it:M3p}. The paths of $\fluc$ are
continuous in $\Hmid^n$, and \ref{it:M1} holds, by Definition~\ref{def:mgp}. So
Corollary~\ref{cor:aclosed}(ii) applies and gives \ref{it:M2}. So $\fluc$ satisfies \ref{it:M1}--\ref{it:M3} of
Theorem~\ref{thm:weakmild}, with datum $\fluc(0)\in\Ufine^n$ and martingale $\overline M$.

\emph{Step 3: the law.} For $\varphi\in\Test$, the processes
$\br{\overline M_i(\cdot),\varphi}=M_i(\cdot,\varphi)$ are continuous martingales with $M_i(0,\varphi)=0$
and with the deterministic brackets \eqref{eq:cov_limit}. By Lemma~\ref{lem:gaussmart}, applied to
finitely many of them, their finite-dimensional distributions are centred Gaussian, determined by
\eqref{eq:cov_limit}, and independent of $\mathcal F_0$. These pairings generate the Borel
$\sigma$-field of $C([0,T];\Hmid^n)$, since $\Test$ is dense in $\Hmid^*$. So $\overline M$ is a
centred Gaussian process with independent increments, independent of $\mathcal F_0$, whose law is
determined by \eqref{eq:cov_limit}. So $\mathrm{Law}(\fluc(0),\overline M)=\mu_0\otimes
\mathrm{Law}(\overline M)$ is determined. By Theorem~\ref{thm:weakmild}(i) and
Lemma~\ref{lem:levelcoin}(ii), $\fluc(t)=S^0(t,0)\fluc(0)+J(t)$ almost surely for each $t$, with
$J(t)$ of Proposition~\ref{prop:convint}. So Lemma~\ref{lem:mildlaw} applies, with
$(\fluc^0,M)=(\fluc(0),\overline M)$, and shows that the law of $\fluc$ is determined by this joint law.

The limit $\fluc$ is a solution: it is a weak solution by Theorem~\ref{thm:OU}(ii), so the
processes $M_i(\cdot,\varphi)$ are the pairings of its martingale, which have the brackets
\eqref{eq:cov_limit} on $\Test$ by Theorem~\ref{thm:noise} and Proposition~\ref{prop:rep}(iii), and it satisfies \ref{it:M1} by
Theorem~\ref{thm:OU}.
\end{proof}

\section[The noise, and the proof of the main theorem]{The noise, and the proof of the main theorem}\label{sec:assembly}

Theorem~\ref{thm:noise} describes the limit martingale through its bracket. The
proofs of Theorems~\ref{thm:OU} and \ref{thm:uniqueness} use, through
Section~\ref{sec:limitmart}, the representation of Proposition~\ref{prop:rep} on an extension
of the probability space.
We first state that representation for the limit martingale:
$M$ is a stochastic integral against a cylindrical Wiener process on
the unweighted $L^2(\R^d;\R^{nd})$, so that \eqref{eq:OU_SPDE} is driven by white
noise in divergence form. The noise space is $L^2(\R^d;\R^{nd})$ and not
$L^2(\R^d;\R^{d})$. It carries one independent cylindrical process per species, as
the vanishing of the cross-brackets in \eqref{eq:cov_limit} requires.

\begin{theorem}[Noise representation]\label{thm:representation}
Let $M$ be a continuous square-integrable $\Hmid^n$-valued martingale with $M(0)=0$, on a
probability space $(\Omega,\mathcal F,\Prob)$ with a filtration satisfying the usual conditions,
whose brackets are \eqref{eq:cov_limit} \textup{(}with $M$ in place of $\widetilde M$\textup{)} for
$\varphi,\psi\in H^{s_m+1}_{a_m}$. By Theorem~\ref{thm:noise}, the limit martingale
$\widetilde M$ with the usual augmentation of the filtration generated by
$(\widetilde\fluc,\widetilde M)$ is an example. Set
\begin{equation}\label{eq:Psi}
 \Psi_i(u):=u_i\Big(2\sigma_i+2\sum_{j=1}^nf(a_{ij}u_j)\Big)=2u_i\Theta_i.
\end{equation}
There is a filtered extension of $(\Omega,\mathcal F,\Prob)$ carrying an $L^2(\R^d;\R^{nd})$-cylindrical Wiener process
$W=(W_1,\dots,W_n)$ (equivalently, $n$ independent
$L^2(\R^d;\R^{d})$-cylindrical processes) such that
\[
 M_i(t)=\int_0^t\Upsilon_i(u(s))\,\d W(s),\qquad
 \Upsilon_i(u)v=\diver\big(\sqrt{\Psi_i(u)}\,v_i\big),\quad v=(v_1,\dots,v_n),
\]
with $\Upsilon_i(u)\in\HS\big(L^2(\R^d;\R^{nd});\Hmid\big)$.
\end{theorem}

\begin{proof}
The brackets \eqref{eq:cov_limit} are \eqref{eq:limbracket} with $\Gamma_i=2\Theta_i$, on
$H^{s_m+1}_{a_m}$, a dense subset of $\Hmid^*=H^{s_m}_{a_m}$. Moreover
$0\le2\Theta_i\le2(\sigma_i+n\eps)$ by Lemma~\ref{lem:lawid}(i) at $\eta=0$. So $M$ satisfies the
hypotheses of Proposition~\ref{prop:rep} with $H=\Hmid^n$ and $\Gamma_i=2\Theta_i$, and
$u_i\Gamma_i=\Psi_i(u)$. For the example, a continuous martingale remains a martingale, with the
same brackets, for the usual augmentation of its filtration. Proposition~\ref{prop:rep}\textup{(ii)} gives $W$ with
$M=\int_0^\cdot\Upsilon(u)\dd W$. By Proposition~\ref{prop:rep}\textup{(i)},
$\Upsilon_i(u)\in\HS(L^2(\R^d;\R^{nd});\Ufine)$, and composition with the embedding
$\Ufine\hookrightarrow\Hmid$ keeps it Hilbert--Schmidt.
\end{proof}

\begin{proof}[Proof of Theorem~\ref{thm:main}]
(i) is Proposition~\ref{prop:tight}. (ii) is Corollary~\ref{cor:whole} together with
Theorems~\ref{thm:noise}, \ref{thm:OU} and \ref{thm:uniqueness}. The bracket is
\eqref{eq:cov_limit}, extended to $\Ufine^*$ by Proposition~\ref{prop:rep}(iii). (iii) is
Corollary~\ref{cor:mgp}. (iv): by (ii), $M$ has the law of the martingale $\widetilde M$ of
Theorem~\ref{thm:noise}. Being a continuous square-integrable martingale with $M(0)=0$ and with
the brackets \eqref{eq:cov_limit}, with respect to the filtration generated by the process itself,
is a property of the law: it is expressed by the martingale property of $M(\cdot,\varphi)$ and of
$M_i(\cdot,\varphi)M_k(\cdot,\psi)-\br{M_i(\cdot,\varphi),M_k(\cdot,\psi)}$ for the generated
filtration. The process $\widetilde M$ has these properties for the larger filtration generated
by $(\widetilde\fluc,\widetilde M)$, hence for its own. So $M$ has them for the filtration it
generates, and also for its usual augmentation, since $M$ is continuous. Theorem~\ref{thm:representation}, applied to $M$, gives $W$ on an
extension of the probability space of $M$.
\end{proof}

\begin{remark}[The results the proof uses]\label{rem:chain}
\emph{Tightness}: Propositions~\ref{prop:mildN} and \ref{prop:conv},
Corollaries~\ref{cor:maxL2-forced} and \ref{cor:remvanish}, Lemmas~\ref{lem:taylor} and
\ref{lem:Rdet}, Proposition~\ref{prop:close} with Theorem~\ref{thm:K4} (which uses
Propositions~\ref{prop:iid} and \ref{prop:phi4} and Lemma~\ref{lem:transfer}),
Lemma~\ref{lem:Mbounds}, Theorem~\ref{thm:aldous}, Corollary~\ref{cor:Maldous} and
Proposition~\ref{prop:criterion}. \emph{Identification}: Theorem~\ref{thm:decomp},
Lemma~\ref{lem:drift} (with Theorem~\ref{thm:G}, Lemma~\ref{lem:Gupgrade} and
Lemma~\ref{lem:remainders}), Proposition~\ref{prop:bracketid} and
Theorem~\ref{thm:noise} (with Lemma~\ref{lem:gaussmart}). \emph{Limit equation}: the stochastic convolution of
Section~\ref{sec:limitmart}, Theorem~\ref{thm:weakmild}, Corollary~\ref{cor:aclosed} and Lemma~\ref{lem:mildlaw}. \emph{Evolution family}:
Section~\ref{sec:evolfam} and Appendix~\ref{sec:appevol}. \emph{Hypotheses and data}:
Lemma~\ref{lem:inputs} (with Lemmas~\ref{lem:crossspecies} and \ref{lem:lawid}),
Lemma~\ref{lem:CH2} with Corollary~\ref{cor:CH2}, Proposition~\ref{prop:initialCLT} and Lemma~\ref{lem:initmom}.
\emph{Uniqueness in law}: Corollary~\ref{cor:mgp}. \emph{Noise}: Theorem~\ref{thm:representation}. \emph{Imported from
\cite{chen_holzinger_2021_rigorous_SKT}}: Proposition~4(i) and Theorems~2 and 3; from
\cite{Majda84}, the Moser estimate behind Lemma~\ref{lem:moser}. From
Theorem~3 the proof uses the existence of $u$, with the bound of \ref{H:small}, and the
bias estimate \eqref{eq:CHbias}.
\end{remark}

\appendix
\section{The tightness criterion}\label{sec:appaldous}

\noindent
The criterion used in Section~\ref{sec:tight} is a version of the
Aldous--Jakubowski condition adapted to a triple of weighted Bessel spaces. We
give it with its proof because the version used in Section~\ref{sec:tight} differs
from the usual statement\ (see, e.g., \cite[Cor.~3.9]{BrMo13} for criteria of this type on a Gelfand triple with the weak topology of the finest space). It assumes uniform bounds in the two finest spaces of the
triple, the time-increment condition in the coarsest space for the field and in the
finest space for the martingale, and compact embeddings between the spaces. Its
conclusion holds in the intersection topology of Definition~\ref{def:ZThat}.

\begin{lemma}[The finer ball is closed, and two corollaries]\label{lem:lscsup}
Let $X$ be a reflexive Banach space, $Y$ a normed space and
$\iota:X\to Y$ a continuous injective embedding. Identify $X$ with
$\iota(X)\subseteq Y$, and set
$\|w\|_X:=+\infty$ for $w\in Y\setminus X$. Then, for every $R<\infty$:
\begin{enumerate}[label=\normalfont(\roman*)]
\item $B_R:=\{w\in Y:\|w\|_X\le R\}$
is closed in $Y$. Equivalently, $\|\cdot\|_X:Y\to[0,\infty]$ is lower
semicontinuous;
\item with $F_X(v):=\sup_{t\le T}\|v(t)\|_X\in[0,\infty]$ for
$v\in C([0,T];Y)$, the set $\{F_X\le R\}$ is closed in $C([0,T];Y)$. Equivalently
$F_X$ is lower semicontinuous, hence Borel, so that $F_X(\xi)$ is a random
variable and $\E F_X(\xi)^2\in[0,\infty]$ is well defined for every random
element $\xi$ of $C([0,T];Y)$;
\item if $v\in C([0,T];Y)$ has $\|v(t)\|_X\le R$ for \emph{a.e.} $t$, then
$v(t)\in X$ and $\|v(t)\|_X\le R$ for \emph{every} $t$.
\end{enumerate}
\end{lemma}

\begin{proof}
\emph{(i).} The closed $X$-ball is weakly compact in $X$, $X$ being reflexive
(Kakutani, \cite[Thm.~3.17]{Bre11}). $\iota$ is linear and continuous, hence weak-to-weak continuous, so
$B_R=\iota(\{\|x\|_X\le R\})$ is weakly compact in $Y$, hence weakly closed, and
being convex it is closed in the norm topology.

\emph{(ii) and (iii).} \emph{(i)}$\Rightarrow$\emph{(ii)}: $\{F_X\le R\}
=\bigcap_{t\le T}\mathrm{ev}_t^{-1}(B_R)$ is an intersection of preimages of a
closed set under continuous evaluations.
\emph{(i)}$\Rightarrow$\emph{(iii)}: the good times are dense, so
$v(t)=\lim v(t_k)$ in $Y$ with $v(t_k)\in B_R$. And
\emph{(iii)}$\Rightarrow$\emph{(i)}: given $w_k\to w$ in $Y$ with $w_k\in B_R$,
set $v(0):=w$, $v(1/k):=w_k$ and interpolate $Y$-linearly. $v$ is
$Y$-continuous at $0$ because on $[\frac1{k+1},\frac1k]$ one has
$\|v(t)-w\|_Y\le\max(\|w_k-w\|_Y,\|w_{k+1}-w\|_Y)$, and $\|v(t)\|_X\le R$ for
every $t>0$ because $B_R$ is convex. \emph{(iii)} then gives
$w=v(0)\in B_R$.
\end{proof}

\begin{lemma}[The $L^2_w$ ingredient]\label{lem:l2weak}
Let $\Ufine$ be a separable Hilbert space and let $(\xi^N)$ be
$L^2(0,T;\Ufine)$-valued random elements with
$\mathsf L:=\sup_N\E\int_0^T\|\xi^N(t)\|^2_{\Ufine}\dd t<\infty$. Then the laws of
$(\xi^N)$ are tight in $L^2_w(0,T;\Ufine)$.
\end{lemma}

\begin{proof}
$L^2(0,T;\Ufine)\cong L^2(0,T)\otimes\Ufine$ is a separable Hilbert space, so it
is reflexive and its closed balls $B_R$ are weakly compact (Kakutani, \cite[Thm.~3.17]{Bre11}). The weak
topology on $B_R$ is metrizable because the dual is separable \cite[\S3.6]{Bre11}. By Markov,
$\Prob(\|\xi^N\|_{L^2(0,T;\Ufine)}>R)\le \mathsf L/R^2\le\eps$ for $R$ large, uniformly
in $N$. On a separable Banach space the weak and
norm Borel $\sigma$-algebras coincide. So a strongly measurable $\xi^N$ is an
$L^2_w$-valued random element.
\end{proof}

\begin{lemma}[Compactness for the intersection topology]\label{lem:joincompact}
Let $E:=C([0,T];\Hmid)\cap L^2(0,T;\Ufine)$ carry the supremum $\tau_1\vee\tau_2$
of the uniform and the weak-$L^2$ topologies. If $K_1\subset C([0,T];\Hmid)$ is
$\tau_1$-compact and $K_2\subset L^2_w(0,T;\Ufine)$ is $\tau_2$-compact, then
$K_1\cap K_2$ is $\tau_1\vee\tau_2$-compact.
\end{lemma}

\begin{proof}
Let $\tau_0$ be the topology on $L^1(0,T;\Hmid)$ generated by
$v\mapsto\int_0^T\br{v(t),\phi(t)}\dd t$, $\phi\in C^\infty_c((0,T)\times\R^d)$. It
is Hausdorff, because an element of $L^1(0,T;\Hmid)$ is determined by its space-time
distributional action. Both $C([0,T];\Hmid)$ and $L^2(0,T;\Ufine)$ embed in
$L^1(0,T;\Hmid)$, the first injectively because a continuous path is determined by
its values almost everywhere, and both embeddings are continuous from $\tau_1$,
respectively $\tau_2$, to $\tau_0$. Hence the set
$\Delta=\{(x,y)\in K_1\times K_2:x=y\text{ in }L^1(0,T;\Hmid)\}$ is closed in the
compact product $K_1\times K_2$, so $\Delta$ is compact. The topology
$\tau_1\vee\tau_2$ on $K_1\cap K_2$ is the one induced by $x\mapsto(x,x)\in\Delta$,
so $K_1\cap K_2$ is compact.

\end{proof}

\begin{definition}[The path space of the criterion]\label{def:ZThat}
Let $U_0$ be a \emph{separable} Hilbert space and set
\[
 \mathcal Z_T:=\Big(C([0,T];\Hmid)\cap L^2_w(0,T;\Ufine)\Big)^n\times C([0,T];\Hmid)^n,
 \qquad
 \widehat{\mathcal Z}_T:=\mathcal Z_T\times C([0,T];U_0),
\]
$L^2_w$ carrying the weak topology.
\end{definition}

Proposition~\ref{prop:tight} is stated on $\mathcal Z_T$, and
Proposition~\ref{prop:criterion} on $\widehat{\mathcal Z}_T$.

\begin{lemma}[A countable separating family on $\widehat{\mathcal Z}_T$]\label{lem:separating}
Let $\widehat{\mathcal Z}_T$ be as in Definition~\ref{def:ZThat}. There is a countable
family of $\widehat{\mathcal Z}_T$-continuous real functions separating the points of
$\widehat{\mathcal Z}_T$. The same holds for $\mathcal Z_T$, which is the case $U_0=\{0\}$.
\end{lemma}

\begin{proof}
Let $(e_j)$ be an orthonormal basis of $\Hmid$, $(f_j)$ one of $U_0$, and
$(\varphi_k)$ a countable dense subset of $L^2(0,T;\Ufine)$. On a factor
$C([0,T];\Hmid)$ take $v\mapsto\br{v(q),e_j}_{\Hmid}$ for
$q\in\Q\cap[0,T]$ and $j\in\N$: each is continuous for the uniform topology,
and they separate, since a continuous path is determined by its values on a
dense subset of $[0,T]$ and a vector by its coefficients. On the weak-$L^2$
ingredient take $v\mapsto\int_0^T\br{v(t),\varphi_k(t)}_{\Ufine}\dd t$:
continuous by definition of the weak topology, and separating by density.
On $C([0,T];U_0)$ proceed as in the first case with $(f_j)$. The union of these
families, composed with the continuous coordinate projections, is countable,
continuous and separating. The intersection factor
$C([0,T];\Hmid)\cap L^2_w(0,T;\Ufine)$ carries the supremum of the two topologies.
This topology is finer than each, so functions continuous for either are continuous for
it, and the first subfamily alone already separates points.

\end{proof}

\begin{lemma}[From a separating family to a Radon Borel law]\label{lem:borel}
Let $X$ be a topological space carrying a countable family $(f_j)$ of real
continuous functions separating its points.
\begin{enumerate}[label=\normalfont(\roman*)]
\item Every compact $K\subseteq X$ belongs to $\sigma(f_j:j\in\N)$ and is metrizable.
\item Let $(\Omega,\mathcal F,\Prob)$ be a
probability space and let $\Xi:\Omega\to X$ be $\sigma(f_j:j\in\N)$-measurable.
Suppose $\Xi$ is tight: for every $\eps>0$ there is a compact $K\subseteq
X$ with $\Prob(\Xi\in K)\ge1-\eps$, an event by \textup{(i)}. Then there is a
$\sigma$-compact set $S\subseteq X$ with $\Prob(\Xi\in S)=1$ and $B\cap S\in\sigma(f_j)$ for every
Borel set $B\subseteq X$, and the law of $\Xi$ extends uniquely from
$\sigma(f_j)$ to a Radon Borel probability measure on $X$.
\end{enumerate}
\end{lemma}

\begin{proof}
We embed $X$ into the Polish space $\R^{\N}$ by the family $(f_j)$. The embedding is a
homeomorphism on compacts (Step~1), so on a $\sigma$-compact set of full measure every Borel
set is $\sigma(f_j)$-measurable (Step~2), and the law extends (Step~3).
Put $F:=(f_j)_{j\in\N}:X\to\R^{\N}$, continuous because each $f_j$ is and
injective because the family separates points, and give $\R^{\N}$ the product
topology, in which it is Polish; $\mathcal B(\R^{\N})$ is the product
$\sigma$-algebra, so $\sigma(F)=\sigma(f_j)$.

\emph{Step 1: compacts are $\sigma(f_j)$-measurable.} Let $K\subseteq X$ be
compact. Then $F(K)$ is compact in $\R^{\N}$, hence Borel, and $F$ being
injective, $K=F^{-1}(F(K))\in\sigma(F)=\sigma(f_j)$. Moreover $F|_K$ is a
homeomorphism onto $F(K)$: it is a continuous bijection from a compact space to
a Hausdorff one. In particular $K$ is compact metrizable, carrying the
metric of $\R^{\N}$. Hence (i) holds.

\emph{Step 2: on a $\sigma$-compact set the two $\sigma$-algebras agree.} Choose
compacts $K'_m$ with $\Prob(\Xi\in K'_m)\ge1-1/m$, put
$K_m:=K'_1\cup\dots\cup K'_m$ (compact, and increasing, which Step 3
needs) and $S:=\bigcup_mK_m$, so that $\Prob(\Xi\in S)=1$. Fix
$B\in\mathcal B(X)$ and $m$. The Borel $\sigma$-algebra of the subspace $K_m$ is
the trace $\mathcal B(X)\cap K_m$, so $B\cap K_m\in\mathcal B(K_m)$; by Step 1,
$F|_{K_m}$ is a homeomorphism onto $F(K_m)$, so
$F|_{K_m}(B\cap K_m)$ is Borel in $F(K_m)$ and hence Borel in $\R^{\N}$, the set
$F(K_m)$ being compact. Therefore
$B\cap K_m=F^{-1}\big(F|_{K_m}(B\cap K_m)\big)\in\sigma(f_j)$ (global
injectivity of $F$ making the intersection with $K_m$ superfluous) and
$B\cap S=\bigcup_m(B\cap K_m)\in\sigma(f_j)$. So $S$ has the properties stated in (ii).

\emph{Step 3: extension, and inner regularity.} Define
$\mu(B):=\Prob(\Xi\in B\cap S)$ for $B\in\mathcal B(X)$; Step 2 makes the
right-hand side meaningful, and $\mu$ is a Borel probability measure agreeing
with the law of $\Xi$ on $\sigma(f_j)$, because $\Prob(\Xi\in S)=1$.
Let $B\in\mathcal B(X)$ and $\eps>0$. Since $K_m\uparrow S$, we have
$B\cap K_m\uparrow B\cap S$, so $\mu(B\cap K_m)>\mu(B)-\eps/2$ for some $m$; and
$B\cap K_m$ is a Borel subset of the compact \emph{metric} space $K_m$, on which
the finite Borel measure $\Prob(\Xi\in\cdot\,\cap K_m)$ is regular, so there is
a compact $C\subseteq B\cap K_m$ with $\mu(C)>\mu(B\cap K_m)-\eps/2$. Hence
$\mu$ is inner regular by compacts, i.e.\ Radon. Uniqueness: two Radon measures
agreeing on the compacts agree, and compacts lie in $\sigma(f_j)$ by Step 1, on
which any two such extensions agree with the law of $\Xi$.
\end{proof}

\begin{theorem}[Skorokhod--Jakubowski {\cite{Jak97}, \cite[Thm.~C.1]{BrOn10}}]
\label{thm:SJ}
Let $X$ be a topological space with a sequence of real continuous functions that
separates its points, and let $(\nu_j)$ be a tight sequence of Borel probability measures
on $X$. Then there are a subsequence $(j_k)$ and Borel maps $\theta,\theta_k:[0,1]\to X$
with $\sigma$-compact ranges such that $\nu_{j_k}$ is the law of $\theta_k$ under Lebesgue
measure and $\theta_k(t)\to\theta(t)$ in $X$ for every $t\in[0,1]$.
\end{theorem}

\begin{lemma}[Upgrading the path space from $\Vcoarse$ to $\Hmid$]\label{lem:upgrade}
Let $\xi$ be a random element of $C([0,T];\Vcoarse)$. Then:
\begin{enumerate}[label=\normalfont(\roman*)]
\item $G_R:=\{v\in C([0,T];\Vcoarse):\ \sup_{t\le T}\|v(t)\|_{\Usharp}\le R\}$
 is
closed, so $\omega\mapsto\sup_t\|\xi(\omega,t)\|_{\Usharp}\in[0,\infty]$ is a
random variable and $\mathsf Q:=\E\sup_t\|\xi(t)\|^2_{\Usharp}$ is well defined;
\item every $v\in G_R$ lies in $C([0,T];\Hmid)$;
\item if $(\xi^N)$ are such random elements whose laws are tight in
$C([0,T];\Vcoarse)$ with $\sup_N\mathsf Q^N<\infty$, then the laws are tight in
$C([0,T];\Hmid)$.
\end{enumerate}
\end{lemma}

\begin{proof}
Part (i) follows from Lemma~\ref{lem:lscsup}(ii). Parts (ii) and (iii) follow from the
inequality \eqref{eq:ehrling}, which bounds the $\Hmid$-norm by a small multiple of the
$\Usharp$-norm plus a multiple of the $\Vcoarse$-norm.

\emph{(i).} Lemma~\ref{lem:lscsup}(ii) applies at $(X,Y)=(\Usharp,\Vcoarse)$,
because $\Usharp$ is Hilbert, hence reflexive, and
$\Usharp\hookrightarrow\Vcoarse$ is the identity, continuous by
\eqref{eq:chainA}, and injective, both being spaces of tempered distributions:
$G_R=\{F_{\Usharp}\le R\}$ is closed, and
$\{\sup_t\|\xi(t)\|_{\Usharp}\le R\}=\xi^{-1}(G_R)$, so the supremum is
measurable.

\emph{(ii).} The inequality \eqref{eq:ehrling} applies with $(X,Y,Z)=(\Usharp,\Hmid,\Vcoarse)$,
because $\Usharp\hookrightarrow\hookrightarrow\Hmid$ by \eqref{eq:chainA} and
$\Hmid\hookrightarrow\Vcoarse$ is continuous and injective; all these spaces are spaces of
tempered distributions, with the identity as embedding. Apply it to the \emph{vector}
$u=v(t)-v(s)$:
\[
 \|v(t)-v(s)\|_{\Hmid}\ \le\ \eps'\|v(t)-v(s)\|_{\Usharp}
 +C_{\eps'}\|v(t)-v(s)\|_{\Vcoarse}
 \ \le\ 2R\eps'+C_{\eps'}\|v(t)-v(s)\|_{\Vcoarse},
\]
so $\limsup_{s\to t}\|v(t)-v(s)\|_{\Hmid}\le2R\eps'$ for every $\eps'>0$, i.e.\
$v$ is $\Hmid$-continuous.

\emph{(iii).} Fix $\eps>0$, pick $R$ with $\sup_N\mathsf Q^N/R^2\le\eps/2$, so
$\Prob(\xi^N\notin G_R)\le\eps/2$ by Markov and (i), and pick
$K\subset C([0,T];\Vcoarse)$ compact with $\sup_N\Prob(\xi^N\notin K)\le\eps/2$.
Put $K_R:=K\cap G_R$, compact by (i), and contained in $C([0,T];\Hmid)$ by (ii).
For $v,w\in K_R$, \eqref{eq:ehrling} applied to $v(t)-w(t)$ and then supremised
in $t$ gives
\begin{equation}\label{eq:upgrade}
 \sup_t\|v(t)-w(t)\|_{\Hmid}
 \ \le\ 2R\eps'+C_{\eps'}\|v-w\|_{C_t\Vcoarse},
\end{equation}
so the identity map $(K_R,\|\cdot\|_{C_t\Vcoarse})\to C([0,T];\Hmid)$ is
uniformly continuous: given $\delta>0$ take $\eps'=\delta/(4R)$, making the
first term $\delta/2$, and then $\|v-w\|_{C_t\Vcoarse}<\delta/(2C_{\eps'})$.
Being a continuous image of a compact set, $K_R$ is compact in
$C([0,T];\Hmid)$, and $\Prob(\xi^N\notin K_R)\le\eps$.
\end{proof}

\begin{proposition}[The criterion]\label{prop:criterion}
Let $\xi^N,m^N$ be $n$-tuples of processes with continuous $\Ufine$-valued
paths, each adapted to a filtration $(\mathcal F^N_t)_{t\le T}$, and let $W^N$
be processes with continuous $U_0$-valued paths, $U_0$ and $\widehat{\mathcal Z}_T$ as in
Definition~\ref{def:ZThat}. Fix an integer $N_\ast\ge1$ and put
$t_\vartheta:=(\tau+\vartheta)\wedge T$. Assume
\begin{enumerate}[label=\normalfont(\alph*)]
\item $\sup_N\E\sup_{t\le T}\|\xi^N(t)\|^2_{\Usharp}<\infty$,
 $\sup_N\E\int_0^T\|\xi^N\|^2_{\Ufine}\dd t<\infty$, and
 $\sup_N\E\sup_{t\le T}\|m^N(t)\|^2_{\Ufine}<\infty$;
\item $\displaystyle\sup_{N\ge N_\ast}\sup_{\tau}\sup_{\vartheta\in(0,1]}
 \vartheta^{-1/2}\,\E\|\xi^N(t_\vartheta)-\xi^N(\tau)\|^2_{\Vcoarse}<\infty$
 and the same for $m^N$ in $\Ufine$, the middle supremum over
 $(\mathcal F^N_t)$-stopping times $\tau\le T$;
\item the laws of $W^N$ are tight in $C([0,T];U_0)$.
\end{enumerate}
Then the laws of $(\xi^N,m^N,W^N)$ are tight on $\widehat{\mathcal Z}_T$. Moreover
\textup{(a)} and \textup{(b)} alone give the laws of $(\xi^N,m^N)$ tight on
$\mathcal Z_T$. In both cases these laws are Radon Borel probability measures.
\end{proposition}

\begin{proof}
We prove tightness block by block. Jakubowski's criterion with Aldous's condition gives
tightness in $C([0,T];\Vcoarse)$ (Step~1), Lemma~\ref{lem:upgrade} upgrades it to
$C([0,T];\Hmid)$ (Step~2), Lemma~\ref{lem:l2weak} gives the weak-$L^2$ component (Step~3), and
Lemma~\ref{lem:joincompact} and Tychonoff's theorem assemble the blocks (Step~4). Step~5 shows
that the laws are Radon Borel measures.

The finitely many $N<N_\ast$ contribute finitely many laws, each individually tight
(Ulam's theorem, \cite[Thm.~1.3]{Bil99}, the spaces below being Polish or, in the $L^2_w$ ingredient, handled by
Lemma~\ref{lem:l2weak} which is stated for all $N$), so we may argue for
$N\ge N_\ast$.

\emph{Step 1: from (a) and (b) to tightness in $C([0,T];\Vcoarse)$.} Extend each
process to $[0,\infty)$ by $\xi^N(t):=\xi^N(T)$ for $t>T$. We use Jakubowski's
criterion \cite[Thm.~3.1]{Jak86}: a sequence of processes with paths in
$D([0,\infty);\Vcoarse)$ is tight there if (1) for every $\eps>0$ there is a
compact $K_\eps\subset\Vcoarse$ with $\Prob\big(\xi^N(t)\in K_\eps\text{ for all }t\big)\ge1-\eps$
for all $N$, and (2) for every $f$ in a family $\mathbb F$ of continuous functions
on $\Vcoarse$ that separates points and is closed under addition, the real
processes $f(\xi^N)$ are tight in $D([0,\infty);\R)$.

For (1), let $K_R$ be the closed $\Usharp$-ball of radius $R$. It is compact in
$\Vcoarse$: relatively compact by $\Usharp\hookrightarrow\hookrightarrow\Hmid\hookrightarrow\Vcoarse$,
and closed by Lemma~\ref{lem:lscsup}(i). By (a) and Chebyshev,
$\Prob\big(\sup_{t\le T}\|\xi^N(t)\|_{\Usharp}>R\big)\le \mathsf Q/R^2$, with $\mathsf Q$ the first
bound of (a).

For (2), let $(e_m)$ be an orthonormal basis of $\Vcoarse$ and let $\mathbb F$ be the
set of maps $v\mapsto(v,e)_{\Vcoarse}$ with $e$ a finite rational combination of the
$e_m$. It separates points and is closed under addition. Fix $f\in\mathbb F$. The
real process $f(\xi^N)$ is continuous and adapted to $(\mathcal F^N_t)$, and
$|f(\xi^N(t))|\le\|e\|_{\Vcoarse}\|\xi^N(t)\|_{\Vcoarse}\le C\|e\|_{\Vcoarse}\|\xi^N(t)\|_{\Usharp}$,
so its one-dimensional marginals are tight by (a). By (b) and Chebyshev, for every
$(\mathcal F^N_t)$-stopping time $\tau\le T$ and every $\vartheta\in(0,1]$,
$\Prob\big(|f(\xi^N(t_\vartheta))-f(\xi^N(\tau))|\ge\rho\big)\le C\|e\|^2_{\Vcoarse}\vartheta^{1/2}\rho^{-2}$,
uniformly in $N\ge N_\ast$. Stopping times of the filtration generated by $f(\xi^N)$
are $(\mathcal F^N_t)$-stopping times, so Aldous's criterion
\cite{aldous1978stopping} gives tightness of $f(\xi^N)$ in $D([0,\infty);\R)$.

Hence the laws of $\xi^N$ are tight in $D([0,\infty);\Vcoarse)$. The paths are
continuous and constant after $T$. $C([0,T];\Vcoarse)$ is closed in the Skorokhod
$J_1$ topology, and on continuous paths that topology coincides with the uniform one \cite[\S12]{Bil99},
so the laws are tight in $C([0,T];\Vcoarse)$. The same argument applies to $m^N$:
since $\Ufine\hookrightarrow\Vcoarse$ continuously, (b)'s second clause gives the
stopping-time bound in $\Vcoarse$, and (a)'s third bound gives (1) through
$\Ufine\hookrightarrow\Usharp$.

\emph{Step 2: upgrade to $C([0,T];\Hmid)$.} Lemma~\ref{lem:upgrade}(iii) with
Step~1 and the first bound of (a). For $m^N$, Step~1 gives the tightness and
(a)'s third bound gives the moment hypothesis, since
$\Ufine\hookrightarrow\Usharp$ gives $\|\cdot\|_{\Usharp}\le C\|\cdot\|_{\Ufine}$.

\emph{Step 3: the weak-$L^2$ ingredient.} Lemma~\ref{lem:l2weak} with the second
bound of (a).

\emph{Step 4: assembling.} For the first factor, intersect the Step-2 and
Step-3 compacts and apply Lemma~\ref{lem:joincompact}. The two exceptional sets
add, so an $\eps/2+\eps/2$ split gives $\eps$. A finite product of tight
families is tight: the product of the compacts is compact by Tychonoff and the
complements add. Steps 1--3 use only \textup{(a)} and \textup{(b)}, and
concern only the first two blocks. Assembling these two blocks gives the second
conclusion, the laws of $(\xi^N,m^N)$ tight on $\mathcal Z_T$. Adjoining (c), which is
tightness in the Polish space $C([0,T];U_0)$, gives tightness of the laws of
$(\xi^N,m^N,W^N)$ on $\widehat{\mathcal Z}_T$.

\emph{Step 5: the laws are Radon Borel measures.} Let $(f_j)$ be the countable
separating family of Lemma~\ref{lem:separating}. Each $f_j$ is one of
$v\mapsto\br{v(q),e_i}_{\Hmid}$ with $q\in\Q\cap[0,T]$, or
$v\mapsto\int_0^T\br{v(t),\varphi_k(t)}_{\Ufine}\dd t$, or the corresponding
functions on the other blocks; composed with $(\xi^N,m^N,W^N)$ each is a real
random variable: the first because $\omega\mapsto\xi^N(q,\omega)$ is
$\Ufine$-valued measurable and $\Ufine\hookrightarrow\Hmid$, the second because
the paths are continuous, hence the map is jointly measurable in $(t,\omega)$ and
Fubini applies. So $(\xi^N,m^N,W^N)$ is $\sigma(f_j)$-measurable, and it is tight
by Steps 1--4; Lemma~\ref{lem:borel}(ii) then gives a unique Radon Borel extension of
its law, and likewise for $(\xi^N,m^N)$ on $\mathcal Z_T$.
\end{proof}

\section{The evolution family and the forced problem}\label{sec:appevol}

This appendix constructs the backward evolution family of Section~\ref{sec:evolfam} and proves the estimates used in the fluctuation argument. Everything in it is deterministic. Existence and uniqueness for linear parabolic systems in variational form are classical:
J.-L.~Lions' theorem \cite[Ch.~XVIII, \S3]{DaLi92}, \cite[p.~112]{Sho97} solves the abstract Cauchy
problem for a measurable, bounded and quasi-coercive form
on a Gelfand triple, and evolution families for parabolic problems are constructed in
\cite[Ch.~5]{Pazy83} and \cite[Ch.~5]{Tan79}. These results do not apply directly to
\eqref{eq:backward-mild} on the weighted scale. The drift $x\cdot\na$ maps $H^{\sigma+1}_a$ only into
$H^\sigma_{a-1}$, not into $H^{\sigma-1}_a$, so the form of the backward system is not bounded on
$H^{\sigma+1}_a$, and the coefficient $x$ is unbounded. The route is therefore a time reversal and a
dilation, which together remove the drift, followed by the conjugation with the weight, which moves the problem
to unweighted spaces. For the conjugated problem we prove a G\aa rding type inequality whose constants
are uniform in $\eta$, and then apply Lions' theorem. Persistence of regularity and the forced
problem follow, the two-derivative maximal bound of the forced problem is needed for the Duhamel
term of the mild formulation. The results are stated for an abstract
coefficient pair satisfying Condition~\ref{ass:budget}. Lemma~\ref{lem:coeff} shows that the
SKT pair satisfies it.

\begin{condition}[A level and a coefficient pair]\label{ass:budget}
Fix a level $(\sigma,a)$ with $\sigma\in\N$, $\sigma>\frac d2+1$ and $a\ge0$, and
$\eta_0\in(0,1]$ with $\eta_0e^T\le1$. For $\eta\in[0,\eta_0]$ let $B^\eta_{ij}$ be as in
\ref{H:CH}, with $B^0_{ij}*$ the multiplication by $a_{ij}$, and
let $A=\max(1,A_0)$ be the constant defined before Assumption~\ref{ass:main}. A coefficient pair is a family of measurable
functions $\Theta^\eta_j$, $g^\eta_{ij}$ on $[0,T]\times\R^d$, $i,j\le n$, such that:
\begin{enumerate}[label=\normalfont(\alph*)]
\item there are constants $\sigma_{\min}>0$, $\kappa_j$, $C_\Theta$, $C_\sharp$ and $\gamma\ge0$,
independent of $\eta$ and $t$, with
\begin{equation}\label{eq:C}
 \Theta^\eta_j\ge\sigma_{\min},\quad
 \|\Theta^\eta_j-\kappa_j\|_{H^{\sigma+1}}+\|\na\Theta^\eta_j\|_{L^\infty}\le C_\Theta,\quad
 \|g^\eta_{ij}\|_{H^{\sigma+1}}\le C_\sharp,\quad
 \|g^\eta_{ij}\|_{L^\infty}\le\gamma;
\end{equation}
\item for every $v\in L^2$, $t\mapsto\Theta^\eta_j(t)v$ and $t\mapsto g^\eta_{ij}(t)v$ are continuous
into $L^2$;
\item the smallness condition
\begin{equation}\label{eq:epscond}
 n\,A\,\gold^{\,a}\,\gamma\ \le\ \tfrac14\,\sigma_{\min},
 \qquad \gold:=\tfrac{1+\sqrt5}2,
\end{equation}
holds, where $\gold^{\,a}$ is the constant of Lemma~\ref{lem:kernelop} at weight $a$.
\end{enumerate}
We say that the condition holds \emph{with upper level} if in addition (a) holds with
$H^{\sigma+3}$ in place of $H^{\sigma+1}$, and (c) holds at the weight $a+1$,
\begin{equation}\label{eq:epscondup}
 n\,A\,\gold^{\,a+1}\,\gamma\ \le\ \tfrac14\,\sigma_{\min}.
\end{equation}
Then the condition holds at the level $(\sigma+2,a+1)$ as well.
\end{condition}

Lemma~\ref{lem:coeff} below shows that the SKT pair satisfies Condition~\ref{ass:budget} at the
levels used. The scalar pair of the Kolmogorov equation in Lemma~\ref{lem:lawid}, with $n=1$
and $g=0$, satisfies it with $\gamma=0$.

\begin{lemma}[Dilation of Sobolev norms]\label{lem:scaling-prop66}
Let $m\ge0$, $a\ge0$, $V\in H^m_a(\R^d)$ and $\mu>0$. Put $V_\mu:=V(\mu\,\cdot)$.
Then
\begin{equation}\label{eq:scaling}
 \|V_\mu\|_{H^m}\le\mu^{-d/2}\max(1,\mu)^{m}\,\|V\|_{H^m},
 \qquad
 \|V_\mu\|_{H^m_a}\le C(a,m,d)\,\mu^{-d/2}\max(1,\mu)^{m}\max(1,\mu^{-1})^{a}\,\|V\|_{H^m_a}.
\end{equation}
The unweighted exponent is sharp. In particular
$\|V_\mu\|_{H^m_a}\le C\mu^{-d/2-a}\|V\|_{H^m_a}$ for $\mu\le1$, and
$\|V_\mu\|_{H^m_a}\le C\mu^{m-d/2}\|V\|_{H^m_a}$ for $\mu\ge1$.
\end{lemma}

\begin{proof}
The unweighted bound is a computation on the Fourier side. The weighted bound follows for
integer $m$ from the Leibniz rule, and for real $m$ by interpolation.
$\widehat{V_\mu}(\xi)=\mu^{-d}\widehat V(\xi/\mu)$, so
$\|V_\mu\|^2_{H^m}=\mu^{-d}\int\langle\mu\zeta\rangle^{2m}|\widehat V(\zeta)|^2$,
and $\langle\mu\zeta\rangle\le\max(1,\mu)\langle\zeta\rangle$. Sharpness is seen
by concentrating $\widehat V$ at $|\zeta|\gg1$ resp.\ $|\zeta|\ll1$. For the
weighted statement take $m$ a nonnegative integer (the general case follows by
interpolation) and expand, for $|\nu|\le m$,
$\pa^\nu[\langle x\rangle^aV(\mu x)]
=\sum_{\nu'\le\nu}c_{\nu'}\,\pa^{\nu-\nu'}\langle x\rangle^a\,
\mu^{|\nu'|}(\pa^{\nu'} V)(\mu x)$. Since
$|\pa^{\nu-\nu'}\langle x\rangle^a|\le C\langle x\rangle^a\le
C\max(1,\mu^{-1})^{a}\langle\mu x\rangle^a$, each term is bounded in $L^2$ by
$C\max(1,\mu^{-1})^{a}\max(1,\mu)^{|\nu'|}\mu^{-d/2}\|\langle\cdot\rangle^a\pa^{\nu'} V\|_{L^2}$,
and $\|\langle\cdot\rangle^a\pa^{\nu'} V\|_{L^2}\le C\|V\|_{H^{|\nu'|}_a}$
because $\pa^{\nu'}(\langle\cdot\rangle^aV)-\langle\cdot\rangle^a\pa^{\nu'} V$ is a
combination of lower-order derivatives of $V$ with $C^\infty_b$-multiples of
$\langle\cdot\rangle^{a}$ as coefficients. For non-integer $m$ interpolate: $v\mapsto\langle\cdot\rangle^av$
is an isometric isomorphism $H^m_a\to H^m$ for every $m$, so
$[H^{m_0}_a,H^{m_1}_a]_\vartheta=H^{m}_a$ and the bound is log-linear in $m$.
\end{proof}

\subsection*{Notation for this appendix}
Fix $t\in(0,T]$, a level $(\sigma,a)$ and a coefficient pair $(\Theta^\eta,g^\eta)$ as in
Condition~\ref{ass:budget}. By \ref{H:A2}, $U_j(x)=-\frac12|x|^2$, so $-\na U_j\cdot\na=x\cdot\na$. For the
SKT pair, $g^\eta_{ij}=u_{\eta,i}f'_\eta(B^\eta_{ij}*u_{\eta,j})$ is the coefficient of
Lemma~\ref{lem:coeffcont}. In the backward variable $s$ the
coefficients are evaluated at time $s$. The backward system of
Proposition~\ref{prop:evolution-outline2} is
\begin{equation}\label{eq:backward-prop66}
 \pa_s\Phi_j+x\cdot\na\Phi_j+\Theta^\eta_j\Delta\Phi_j
 +\sum_iB^\eta_{ij}*\big(g^\eta_{ij}\Delta\Phi_i\big)=0,\qquad\Phi(t)=\varphi.
\end{equation}
With $\mathcal A(s)\Phi_j:=x\cdot\na\Phi_j+\Theta^\eta_j\Delta\Phi_j
+\sum_iB^\eta_{ij}*(g^\eta_{ij}\Delta\Phi_i)$ and a forcing term
$G\in L^2(0,t;H^{\sigma-1}_a)$, the backward problem is
\begin{equation}\label{eq:backF}
 -\pa_s\Phi=\mathcal A(s)\Phi+G(s),\qquad \Phi(t)=\varphi.
\end{equation}
For $\Phi\in H^{\sigma+1}_a$, the terms $\Theta^\eta_j\Delta\Phi_j$ and
$B^\eta_{ij}*(g^\eta_{ij}\Delta\Phi_i)$ lie in $H^{\sigma-1}_a$, by Condition~\ref{ass:budget}(a)
and Lemma~\ref{lem:kernelop}, and $x\cdot\na\Phi_j$ lies in $H^\sigma_{a-1}$. So
$\mathcal A(s)\Phi\in H^{\sigma-1}_{a-1}$. The solution class $\mathcal X_{\rm ph}$ is the
space of $\Phi\in C([0,t];H^\sigma_a)\cap L^2(0,t;H^{\sigma+1}_a)$, with norm
$\|\Phi\|_{\mathcal X_{\rm ph}}:=\sup_{s\le t}\|\Phi(s)\|_{H^\sigma_a}+\|\Phi\|_{L^2(0,t;H^{\sigma+1}_a)}$.
A solution of \eqref{eq:backF} in $\mathcal X_{\rm ph}$ is a $\Phi\in\mathcal X_{\rm ph}$
with $\Phi(t)=\varphi$ that satisfies \eqref{eq:backF} in $L^2(0,t;H^{\sigma-1}_{a-1})$.

The time reversal $\Psi(\tau):=\Phi(t-\tau)$ turns \eqref{eq:backward-prop66} into the
forward system, with coefficients evaluated at time $t-\tau$,
\begin{equation}\label{eq:forward}
 \pa_\tau\Psi_j=x\cdot\na\Psi_j+\Theta^\eta_j\Delta\Psi_j
 +\sum_iB^\eta_{ij}*\big(g^\eta_{ij}\Delta\Psi_i\big),\qquad\Psi(0)=\varphi.
\end{equation}
Let $w(y):=\langle y\rangle^a$ and $\mathcal Dh:=w\,\Delta(w^{-1}h)$. By the Leibniz rule and
$w\na(w^{-1})=-\na\log w$,
\[
 \mathcal Dh=\Delta h+\mathbf b_w\cdot\na h+c_wh,\qquad
 \mathbf b_w:=-2\na\log w,\qquad c_w:=w\,\Delta(w^{-1}).
\]
Both coefficients lie in $C^\infty_b$, since $\na\log w(y)=a\,y\br y^{-2}$ and
$c_w(y)=a(a+2)|y|^2\br y^{-4}-ad\br y^{-2}$. For $\lambda>0$ we denote by $\mathrm{Dil}_\lambda$ the
dilation $(\mathrm{Dil}_\lambda v)(x):=v(\lambda x)$. The dilation of
Lemma~\ref{lem:dilation}, followed by the conjugation $\zeta:=w\,\Xi$, gives the
system
\begin{equation}\label{eq:conj}
 \pa_\tau\zeta_j=\widehat\Theta_j\,\mathcal D\zeta_j
 +\sum_i\mathcal K^\tau_{ij}\big(\widehat g_{ij}\,\mathcal D\zeta_i\big),
 \qquad
 \mathcal K^\tau_{ij}h:=w\,\big(B^{\eta_\tau}_{ij}*(w^{-1}h)\big),
\end{equation}
where $\widehat\Theta_j$, $\widehat g_{ij}$ and $\eta_\tau$ are those of
Lemma~\ref{lem:dilation}. We denote its right side by $\mathcal F(\tau)\zeta$. With a
forcing term $\mathsf F\in L^2(0,t;H^{\sigma-1})$, given by \eqref{eq:Fdef} below, it becomes
\begin{equation}\label{eq:conjF}
 \pa_\tau\zeta=\mathcal F(\tau)\zeta+\mathsf F,\qquad \zeta(0)=\zeta_0.
\end{equation}
The solution class $\mathcal X_{\rm cj}$ is $\mathcal X_{\rm ph}$ at $a=0$, in the
variable $\tau$, and \eqref{eq:conjF} is required in $L^2(0,t;H^{\sigma-1})$.
Multiplication by $w$ is an isometric isomorphism $H^m_a\to H^m$ for every $m$.
In both equations, the terms other than the time derivative lie in the stated
spaces: for $\mathcal X_{\rm ph}$ by the computation above, and for
$\mathcal X_{\rm cj}$ by Lemma~\ref{lem:ident}. So each equation holds in its space if
and only if it holds in the sense of distributions, and the change of variables of
Lemma~\ref{lem:bij} needs no estimate at the weight $a-1$.

\begin{lemma}[Dilation conjugation]\label{lem:dilation}
Let $\Xi_j(\tau,y):=\Psi_j(\tau,e^{-\tau}y)$. Then $\Psi$ solves
\eqref{eq:forward} if and only if
\begin{equation}\label{eq:dilated}
 \pa_\tau\Xi_j=\widehat\Theta_j\Delta\Xi_j
 +\sum_iB^{\eta_\tau}_{ij}*\big(\widehat g_{ij}\Delta\Xi_i\big),
 \qquad\Xi(0)=\varphi,\qquad \eta_\tau:=\eta e^{\tau},
\end{equation}
with $\widehat\Theta_j(\tau,y):=e^{2\tau}\Theta^\eta_j(t-\tau,e^{-\tau}y)$ and
$\widehat g_{ij}(\tau,y):=e^{2\tau}g^\eta_{ij}(t-\tau,e^{-\tau}y)$. There is no
transport term in \eqref{eq:dilated}.
\end{lemma}

\begin{proof}
$\pa_\tau\Psi_j(\tau,x)=(\pa_\tau\Xi_j)(\tau,e^\tau x)+e^\tau
x\cdot(\na\Xi_j)(\tau,e^\tau x)$ and
$x\cdot\na_x\Psi_j(\tau,x)=e^\tau x\cdot(\na\Xi_j)(\tau,e^\tau x)$: the transport
contributions cancel exactly. Also
$\Delta_x\Psi_j=e^{2\tau}(\Delta\Xi_j)(\tau,e^\tau x)$. In the nonlocal term
substitute $y'=e^\tau x'$ ($\mathrm dx'=e^{-d\tau}\mathrm dy'$,
$x-x'=e^{-\tau}(y-y')$) and use
\[
 e^{-d\tau}B^\eta_{ij}\big(e^{-\tau}u\big)
 =e^{-d\tau}\eta^{-d}B_{ij}\Big(\frac{|u|}{\eta e^{\tau}}\Big)
 =(\eta e^{\tau})^{-d}B_{ij}\Big(\frac{|u|}{\eta e^{\tau}}\Big)
 =B^{\eta_\tau}_{ij}(u). \qedhere
\]
\end{proof}

\begin{lemma}[The double conjugation is a bijection of solution classes]
\label{lem:bij}
The map $\mathcal C:\Phi\mapsto\zeta$,
$\zeta(\tau,y):=w(y)\,\Phi\big(t-\tau,e^{-\tau}y\big)$, is a linear bijection
$\mathcal X_{\rm ph}\to\mathcal X_{\rm cj}$ with
\[
 c(T)^{-1}\|\Phi\|_{\mathcal X_{\rm ph}}\le\|\mathcal C\Phi\|_{\mathcal X_{\rm cj}}
 \le c(T)\|\Phi\|_{\mathcal X_{\rm ph}},\qquad
 c(T):=C(a,\sigma,d)\,e^{(\sigma+1+a)T},
\]
and it carries solutions of \eqref{eq:backF} with forcing $G$ and datum
$\varphi$ onto solutions of \eqref{eq:conjF} with $\zeta_0=w\varphi$ and
\begin{equation}\label{eq:Fdef}
 \mathsf F_j(\tau,y)=w(y)\,G_j\big(t-\tau,e^{-\tau}y\big).
\end{equation}
In particular the forcing carries no factor $e^{2\tau}$.
\end{lemma}

\begin{proof}
We compare the norms of $\zeta$ and $\Phi$ with Lemma~\ref{lem:scaling-prop66}, derive the
continuity in time from the strong continuity of the dilation group, and then transform the
equations.

\emph{Norms.} $\|\zeta(\tau)\|_{H^m}=\|\Xi(\tau)\|_{H^m_a}$ and
$\Xi(\tau,\cdot)=\Psi(\tau,e^{-\tau}\cdot)$, $\Psi(\tau,\cdot)=\Xi(\tau,e^{\tau}\cdot)$.
Applying Lemma~\ref{lem:scaling-prop66} at $\mu=e^{-\tau}\le1$ and at $\mu=e^{\tau}\ge1$ gives,
for $m\in\{\sigma,\sigma+1\}$ and $0\le\tau\le T$,
\[
 C^{-1}e^{-(m-\frac d2)\tau}\|\Psi(\tau)\|_{H^m_a}
 \ \le\ \|\Xi(\tau)\|_{H^m_a}\ \le\ C\,e^{(\frac d2+a)\tau}\|\Psi(\tau)\|_{H^m_a},
\]
two-sided and uniform on $[0,T]$. $c(T)=C e^{(\sigma+1+a)T}$ dominates both
exponents at $m=\sigma+1$, since $\frac d2<\sigma+1$. For continuity in $\tau$ we show that the dilations $\mathrm{Dil}_\lambda$ form a strongly continuous group on $H^m_a$, $m,a\ge0$. It is
locally uniformly bounded by Lemma~\ref{lem:scaling-prop66}, both branches. On the
dense subspace $C^\infty_0$ (dense because $\br\cdot^aC^\infty_0=C^\infty_0$),
$\mathrm{Dil}_\lambda v\to v$ in $H^m_a$ as $\lambda\to1$. Indeed, for $v\in C^\infty_0$ and
$\lambda\in[\frac12,2]$, $\br\cdot^a(\mathrm{Dil}_\lambda v-v)$ is supported in a fixed compact set, and each
derivative $\pa^\nu[\br x^a(v(\lambda x)-v(x))]$ is a finite sum of terms
$(\pa^{\nu_1}\br x^a)\big(\lambda^{|\nu_2|}(\pa^{\nu_2}v)(\lambda x)-(\pa^{\nu_2}v)(x)\big)$, which tend
to $0$ uniformly by uniform continuity of $\pa^{\nu_2}v$. So the convergence holds in every $H^k$,
$k\in\N$, hence in $H^m$ for real $m\ge0$ by
$\|h\|^2_{H^m}\le\|h\|^{2(1-\vartheta)}_{H^{\lfloor m\rfloor}}\|h\|^{2\vartheta}_{H^{\lceil m\rceil}}$,
$\vartheta:=m-\lfloor m\rfloor$. An $\eps/3$ argument with the local uniform bound extends the
convergence to every $v\in H^m_a$. Continuity at an arbitrary $\lambda_0>0$ follows from the
group law $\mathrm{Dil}_\lambda=\mathrm{Dil}_{\lambda_0}\mathrm{Dil}_{\lambda/\lambda_0}$:
$\|\mathrm{Dil}_\lambda v-\mathrm{Dil}_{\lambda_0}v\|_{H^m_a}\le\|\mathrm{Dil}_{\lambda_0}\|\,\|\mathrm{Dil}_{\lambda/\lambda_0}v-v\|_{H^m_a}\to0$ as
$\lambda\to\lambda_0$. Continuity of the composition $\tau\mapsto \mathrm{Dil}_{e^{-\tau}}\Psi(\tau)$
then follows from
$\|\mathrm{Dil}_{\lambda}\Psi(\tau)-\mathrm{Dil}_{\lambda_0}\Psi(\tau_0)\|_{H^m_a}
\le(\sup_{\lambda\in[e^{-T},1]}\|\mathrm{Dil}_\lambda\|)\|\Psi(\tau)-\Psi(\tau_0)\|_{H^m_a}
+\|(\mathrm{Dil}_{\lambda}-\mathrm{Dil}_{\lambda_0})\Psi(\tau_0)\|_{H^m_a}$, the second term using
strong continuity at the fixed vector $\Psi(\tau_0)$. Hence
$\tau\mapsto\Xi(\tau)$ is $H^\sigma_a$-continuous if and only if $s\mapsto\Phi(s)$ is.
Multiplication by $w$ is an isometric isomorphism $H^m_a\to H^m$. So $\mathcal C$
is a bijection with the stated equivalence.

\emph{Equations.} The time reversal $s=t-\tau$ turns $-\pa_s$ into $\pa_\tau$ and
transports $G$ by composition. The dilation removes the transport term exactly
(Lemma~\ref{lem:dilation}) and multiplies $\widehat\Theta,\widehat g$ by $e^{2\tau}$,
that factor arising from $\Delta_x\Psi=e^{2\tau}(\Delta\Xi)$. The zeroth-order
term $G$ is transported by composition alone and picks up no such factor. The
weight conjugation multiplies by $w$ and replaces $\Delta$ by $\mathcal D$.
\end{proof}

\begin{lemma}[Uniform coefficient bounds]\label{lem:coeff}
Let $\sigma\in\N$ with $\frac d2+1<\sigma\le s-1$, and $a=a_1$. The SKT pair, namely
$\Theta^\eta_j$ of \eqref{eq:mobilities} and $g^\eta_{ij}=u_{\eta,i}f'_\eta(B^\eta_{ij}*u_{\eta,j})$,
satisfies Condition~\ref{ass:budget} at $(\sigma,a)$, with $\kappa_j:=\sigma_j+nf(0)$,
$\sigma_{\min}=\min_k\sigma_k$ and $\gamma:=C_s\|u_0\|_{H^s}\eps$. If $\sigma+3\le s$, it satisfies
it with upper level. The same holds for the scalar pair $(\Theta^\eta_i,0)$.
\end{lemma}

\begin{proof}
We check the three parts of Condition~\ref{ass:budget} in turn. The bounds in (a) follow
from the range of the arguments of $f_\eta$ (Lemma~\ref{lem:lawid}(i)), from Moser's estimate
(Lemma~\ref{lem:moser}) and from the algebra property of $H^{\sigma+1}$.

$\Theta^\eta_j\ge\sigma_{\min}$ because $f_\eta\ge0$ by \ref{H:A4}. Convolution is a Fourier multiplier, so
$\|B^\eta_{jk}*v\|_{H^m}\le\|B_{jk}\|_{L^1}\|v\|_{H^m}\le A\|v\|_{H^m}$ for every
real $m$, and $\|B^\eta_{jk}*u_{\eta,k}\|_{L^\infty}\le\mathsf b_0$ by Lemma~\ref{lem:lawid}(i). So the argument of
$f_\eta$ lies in $I=[-2\mathsf b_0,2\mathsf b_0]\subset[-a_\eta,a_\eta]$ by \ref{H:cutoff}, where
$f_\eta=f$, and $\|f^{(m)}\|_{L^\infty(I)}\le\eps$ for $m\le s+1$ by \ref{H:small}.

Since $B^\eta_{jk}*u_{\eta,k}\in H^{\sigma+1}\hookrightarrow C_0$, each summand
of $\Theta^\eta_j$ tends to $f(0)$ at infinity, so $\kappa_j=\sigma_j+nf(0)$ is
the limit, and $\Theta^\eta_j-\kappa_j=\sum_k\big(f(w^\eta_{jk})-f(0)\big)$ is the function to bound in $H^{\sigma+1}$. Moser's composition
estimate applied to $f-f(0)$ (which vanishes at $0$) and $\sigma+1\le s$ give the
second bound. Here Moser's estimate is Lemma~\ref{lem:moser} at the integer order
$\sigma+1$, with $J=I$. It requires $f-f(0)\in C^{\sigma+1}(I)$, and for $g^\eta_{ij}$ below
$f'-f'(0)\in C^{\sigma+1}(I)$, and both follow from $\|f\|_{C^{s+1}(I)}\le\eps$ since
$\sigma+2\le s+1$. $\|\na\Theta^\eta_j\|_{L^\infty}\le C\|\Theta^\eta_j-\kappa_j\|_{H^{\sigma+1}}$
since $\sigma>\frac d2$. For $g^\eta_{ij}$ split
$g^\eta_{ij}=f'(0)u_{\eta,i}+u_{\eta,i}\big(f'(w^\eta_{ij})-f'(0)\big)$: the first term
is in $H^{\sigma+1}$ because $u_{\eta,i}$ is, the second is a product of two
$H^{\sigma+1}$ functions ($H^{\sigma+1}$ is an algebra). The $L^\infty$ bound is
immediate. The same argument with $\sigma+2$ in place of $\sigma$ gives the upper level when
$\sigma+3\le s$. Part~(b) of the condition is Lemma~\ref{lem:lawid}(i). Part~(c) with the weight $a+1=a_1+1$, that is \eqref{eq:epscondup}, is the second
half of \ref{H:small}, and it implies (c) at the weight $a_1$, since $\gold>1$. For the scalar pair,
(a) and (b) are the statements on $\Theta^\eta_i$, and (c) holds with $\gamma=0$.
\end{proof}

\begin{lemma}[Cost of the dilation]\label{lem:dilcost}
For $0\le\tau\le t\le T$, $\eta\in[0,\eta_0]$, with
$\theta(\tau):=e^{2\tau}\sigma_{\min}$:
\begin{equation}\label{eq:dilbounds}
\begin{gathered}
 \widehat\Theta_j\ \ge\ \theta(\tau),\qquad
 \|\widehat g_{ij}\|_{L^\infty}\ \le\ e^{2\tau}\,\gamma
 \ =\ \frac{\gamma}{\sigma_{\min}}\,\theta(\tau),\\
 \|\widehat\Theta_j-e^{2\tau}\kappa_j\|_{H^{\sigma+1}}\le e^{(2+\frac d2)T}C_\Theta,
 \quad
 \|\na\widehat\Theta_j\|_{L^\infty}\le e^{T}C_\Theta,
 \quad
 \|\widehat g_{ij}\|_{H^{\sigma+1}}\le e^{(2+\frac d2)T}C_\sharp,
\end{gathered}
\end{equation}
and $\|B^{\eta_\tau}_{ij}\|_{L^1}\le A$,
$\operatorname{supp}B^{\eta_\tau}_{ij}\subset B_{\eta e^T}(0)\subset B_1(0)$.
\end{lemma}

\begin{proof}
The first two lines are Condition~\ref{ass:budget}(a) multiplied by $e^{2\tau}$;
\emph{both} the dissipation floor and the coupling size carry that factor. The Sobolev bounds use Lemma~\ref{lem:scaling-prop66} with
$\mu=e^{-\tau}\le1$, giving $\mu^{-d/2}=e^{d\tau/2}$;
$\na_y\widehat\Theta_j=e^{\tau}(\na\Theta^\eta_j)(t-\tau,e^{-\tau}y)$ gives the
$L^\infty$ bound. The last statement is the $L^1$-criticality of the scaling
together with $\eta e^T\le1$.
\end{proof}

\begin{lemma}[The conjugated mollifier]\label{lem:kernelop}
Put $\rho_z(y):=w(y)/w(y-z)$ and
$C_\rho(a):=\sup_{|z|\le1}\|\rho_z\|_{L^\infty}$. Then
\begin{equation}\label{eq:Crho}
 C_\rho(a)=\gold^{\,a},\qquad\gold=\tfrac{1+\sqrt5}2,
\end{equation}
$\rho_z\in C^\infty_b$ with $\|\rho_z\|_{C^k_b}\le C_k(a,d)$ uniformly in
$|z|\le1$, and
\begin{equation}\label{eq:krep}
 \mathcal K^\tau_{ij}h=\int_{\R^d}B^{\eta_\tau}_{ij}(z)\,\rho_z\cdot(\mathrm T_zh)\dz,
 \qquad(\mathrm T_zh)(y):=h(y-z),
\end{equation}
the integral converging in $H^m$ for every $m$ with $\|\mathcal
K^\tau_{ij}h\|_{H^m}\le C_m A\|h\|_{H^m}$, uniformly in $\eta\in[0,\eta_0]$ and
$\tau\in[0,T]$.
\end{lemma}

\begin{proof}
\eqref{eq:krep} is the substitution $z=y-y'$; only $|z|\le1$ occurs by
Lemma~\ref{lem:dilcost}. For \eqref{eq:Crho}, $\rho_z(y)^{2/a}
=\frac{1+|y|^2}{1+|y-z|^2}$ is maximised at $|z|=1$ and $y=\lambda z$ with
$\lambda$ maximising $(1+\lambda^2)/(1+(\lambda-1)^2)$, i.e.\
$\lambda^2-\lambda-1=0$, $\lambda=\gold$, with value $\gold^2$, hence
$C_\rho(a)=\gold^{a}$.
Uniform $C^k_b$ bounds: $\na\log\rho_z=a\big(y\wt^{-2}-(y-z)\langle
y-z\rangle^{-2}\big)$ is bounded by $a$ since $|y|/(1+|y|^2)\le\frac12$, and,
since $\partial^\beta\log\langle y\rangle$ is bounded on $\R^d$ for every
$|\beta|\ge1$, all derivatives of $\log\rho_z$ of order $\ge1$ are bounded by
constants $C_\beta(a,d)$, uniformly in $z$. By the chain rule,
$\partial^\alpha\rho_z=\rho_z\,P_\alpha\big(\partial^\beta\log\rho_z:1\le|\beta|\le|\alpha|\big)$
with a universal polynomial $P_\alpha$ (e.g.\ $\na\rho_z=\rho_z\na\log\rho_z$).
Together with $\rho_z\le\gold^a$ for $|z|\le1$, this gives
$\|\partial^\alpha\rho_z\|_{L^\infty}\le\gold^aC_\alpha(a,d)$. The map $z\mapsto\rho_z\mathrm T_zh$ is continuous
$H^m$-valued ($\mathrm T_z$ is strongly continuous on $H^m$ and $z\mapsto\rho_z$ is
Lipschitz in $C^k_b$), so Minkowski's integral inequality applies. Multiplication
by a $C^{k}_b$ function with $k>|m|$ is bounded on $H^m$ (for $m<0$ by duality,
multiplication by a real function being self-adjoint for the $L^2$ pairing), and
translation is an isometry.
\end{proof}

\emph{The energy and its terms.} For $\tau\in[0,t]$ and $\zeta\in H^{\sigma+1}$ we denote
\begin{equation}\label{eq:energy}
\begin{gathered}
 v_j:=\langle\nabla\rangle^\sigma\zeta_j,\qquad
 E:=\sum_j\|v_j\|^2_{L^2},\qquad
 D:=\sum_j\|\na v_j\|^2_{L^2},\\
 D_j:=\big\langle\langle\nabla\rangle^\sigma(\widehat\Theta_j\mathcal D\zeta_j),v_j\big\rangle,\qquad
 C_{ij}:=\big\langle\langle\nabla\rangle^\sigma\mathcal K^\tau_{ij}(\widehat g_{ij}\mathcal D\zeta_i),v_j\big\rangle,
\end{gathered}
\end{equation}
with $\widehat\Theta_j$ and $\widehat g_{ij}$ evaluated at $\tau$, and with the pairing of $H^{-1}$
and $H^1$. For a function $\zeta$ of $\tau$, the same symbols denote these quantities at
$\zeta(\tau)$. So $E=\|\zeta\|^2_{H^\sigma}$ and $E+D=\|\zeta\|^2_{H^{\sigma+1}}$. The term $D_j$ is
the signed diffusion term, it is not $\|\na v_j\|^2_{L^2}$.

\begin{lemma}[Admissibility of the energy identity]\label{lem:ident}
There is a constant $C$, independent of $\tau\in[0,t]$ and $\eta\in[0,\eta_0]$, with
\begin{equation}\label{eq:Fbound}
 \|\mathcal F(\tau)\zeta\|_{H^{\sigma-1}}\le C\|\zeta\|_{H^{\sigma+1}},\qquad\zeta\in H^{\sigma+1}.
\end{equation}
If $\zeta\in L^2(0,t;H^{\sigma+1})$ solves \eqref{eq:conj}, then
$\pa_\tau\zeta\in L^2(0,t;H^{\sigma-1})$, the quantities $D_j$ and $C_{ij}$ of
\eqref{eq:energy} are finite, and
$\tau\mapsto\|\langle\nabla\rangle^\sigma\zeta(\tau)\|^2_{L^2}$ is absolutely continuous with
$\frac{\mathrm d}{\mathrm d\tau}\|\langle\nabla\rangle^\sigma\zeta\|^2
=2\langle\pa_\tau\zeta,\langle\nabla\rangle^{2\sigma}\zeta\rangle_{H^{\sigma-1},H^{-\sigma+1}}$.
Consequently, for a.e.\ $\tau$,
\begin{equation}\label{eq:energyid}
 \tfrac12\tfrac{\mathrm d}{\mathrm d\tau}E=\sum_j D_j+\sum_{i,j}C_{ij}.
\end{equation}
\end{lemma}

\begin{proof}
For $\zeta\in H^{\sigma+1}$: $\mathcal D\zeta_j\in H^{\sigma-1}$;
$\widehat\Theta_j=e^{2\tau}\kappa_j+\tilde\theta_j$ with $\tilde\theta_j\in
H^{\sigma+1}$ and $H^{\sigma+1}\cdot H^{\sigma-1}\subset H^{\sigma-1}$ because
$\sigma-1>\frac d2$; hence $\widehat\Theta_j\mathcal D\zeta_j\in H^{\sigma-1}$,
and likewise $\mathcal K^\tau_{ij}(\widehat g_{ij}\mathcal D\zeta_i)\in
H^{\sigma-1}$ by Lemma~\ref{lem:kernelop} with $m=\sigma-1$, with bounds uniform in $\tau$ and $\eta$
by Lemma~\ref{lem:dilcost}. This gives \eqref{eq:Fbound}, and hence
$\|\pa_\tau\zeta\|_{H^{\sigma-1}}\le C\|\zeta\|_{H^{\sigma+1}}$. Since
$\langle\nabla\rangle^{2\sigma}\zeta\in H^{-\sigma+1}$, the
Lions--Magenes chain rule \cite{LiMa72} applies on the triple
$H^{\sigma+1}\subset H^\sigma\subset H^{\sigma-1}$. Each of $D_j$ and $C_{ij}$ is a
continuous bilinear form in $(\zeta,\zeta)$ on $H^{\sigma+1}\times H^{\sigma+1}$, so it is finite.
Since $\langle\pa_\tau\zeta,\langle\nabla\rangle^{2\sigma}\zeta\rangle
=\sum_j\big\langle\langle\nabla\rangle^\sigma(\mathcal F(\tau)\zeta)_j,v_j\big\rangle$, the
definition of $\mathcal F(\tau)$ in \eqref{eq:conj} gives \eqref{eq:energyid}.
\end{proof}

\begin{lemma}[Admissibility, and the forced Young step]\label{lem:identF}
Let $\zeta\in L^2(0,t;H^{\sigma+1})$ satisfy $\pa_\tau\zeta=\mathcal F(\tau)\zeta+\mathsf F$.
Then $\pa_\tau\zeta\in L^2(0,t;H^{\sigma-1})$, $E$ is absolutely continuous with
\[
 \tfrac{\dd}{\dtau}E=2\langle\pa_\tau\zeta,\langle\nabla\rangle^{2\sigma}\zeta\rangle_{H^{\sigma-1},H^{-\sigma+1}},
\]
and
\begin{equation}\label{eq:energyF}
 \tfrac12\tfrac{\dd}{\dtau}E
 =\sum_jD_j+\sum_{i,j}C_{ij}+\sum_j\br{\langle\nabla\rangle^\sigma\mathsf F_j,\langle\nabla\rangle^\sigma\zeta_j},
\end{equation}
$D_j,C_{ij}$ as in \eqref{eq:energy}. Moreover, for every
$\lambda>0$ and every $\tau$,
\begin{equation}\label{eq:young}
 2\Big|\sum_j\br{\langle\nabla\rangle^\sigma\mathsf F_j,\langle\nabla\rangle^\sigma\zeta_j}\Big|
 \ \le\ \lambda\,\theta(\tau)\,(E+D)+\tfrac1{\lambda\theta(\tau)}\|\mathsf F\|^2_{H^{\sigma-1}}.
\end{equation}
\end{lemma}

\begin{proof}
By \eqref{eq:Fbound}, $\mathcal F(\tau)\zeta\in H^{\sigma-1}$
with $\|\mathcal F(\tau)\zeta\|_{H^{\sigma-1}}\le C\|\zeta\|_{H^{\sigma+1}}$;
adding $\mathsf F\in H^{\sigma-1}$ changes nothing, and the Lions--Magenes chain rule on
$H^{\sigma+1}\subset H^\sigma\subset H^{\sigma-1}$ applies. In \eqref{eq:energyF}, the extra term is
$\br{\mathsf F,\langle\nabla\rangle^{2\sigma}\zeta}_{H^{\sigma-1},H^{-\sigma+1}}
=\br{\langle\nabla\rangle^\sigma\mathsf F,\langle\nabla\rangle^\sigma\zeta}$. For \eqref{eq:young}, by
Plancherel $\br{\langle\nabla\rangle^\sigma\mathsf F_j,\langle\nabla\rangle^\sigma\zeta_j}
=\br{\langle\nabla\rangle^{\sigma-1}\mathsf F_j,\langle\nabla\rangle^{\sigma+1}\zeta_j}$, so by
Cauchy--Schwarz, \eqref{eq:sob} and Young,
$2|\br{\cdot}|\le2\|\mathsf F_j\|_{H^{\sigma-1}}\|\zeta_j\|_{H^{\sigma+1}}
\le\lambda\theta\|\zeta_j\|^2_{H^{\sigma+1}}
+\frac1{\lambda\theta}\|\mathsf F_j\|^2_{H^{\sigma-1}}$, sum over $j$ and use
$\sum_j\|\zeta_j\|^2_{H^{\sigma+1}}=E+D$.
\end{proof}

\begin{proposition}[Energy estimate]\label{prop:apriori}
Let Condition~\ref{ass:budget} hold. There is a constant
$C_0=C_0(d,n,\sigma,a,T,A,\sigma_{\min},\max_j|\kappa_j|,C_\Theta,C_\sharp)$, independent of
$\eta$, $t$ and $\zeta$, such that the following holds.
\begin{enumerate}[label=\normalfont(\roman*)]
\item For every $\tau\in[0,t]$ and every $\zeta\in H^{\sigma+1}$, in the notation \eqref{eq:energy},
\begin{equation}\label{eq:Sbound}
 \sum_jD_j+\sum_{i,j}C_{ij}\ \le\ -\tfrac38\,\theta(\tau)\,D+\tfrac{C_0}2\,E.
\end{equation}
\item Let $\zeta\in C([0,t];H^\sigma)\cap L^2(0,t;H^{\sigma+1})$ be such that $E$ is absolutely
continuous and \eqref{eq:energyid} holds for a.e.\ $\tau$. By Lemma~\ref{lem:ident}, this is the
case if $\zeta$ solves \eqref{eq:conj}. Then
\begin{equation}\label{eq:diffineq}
 \frac{\mathrm d}{\mathrm d\tau}E+\tfrac34\,\theta(\tau)\,D\ \le\ C_0\,E
 \qquad\text{in }\mathcal D'(0,t),
\end{equation}
and consequently, since $\theta\ge\sigma_{\min}$,
\begin{equation}\label{eq:apriori}
 \sup_{0\le\tau\le t}\|\zeta(\tau)\|^2_{H^\sigma}
 +\sigma_{\min}\!\int_0^t\!\!D\dtau\ \le\ 3\,e^{C_0T}\,\|\zeta(0)\|^2_{H^\sigma}.
\end{equation}
\end{enumerate}
\end{proposition}

\begin{proof}
For (i) we fix $\tau$ and $\zeta\in H^{\sigma+1}$, and bound the diagonal terms $D_j$
(Step~1) and the coupling terms $C_{ij}$ (Step~2) with the Kato--Ponce commutator estimate.
The smallness \eqref{eq:epscond} controls the only coupling term with two gradients. Step~3
sums the bounds, which gives (i), and integrates \eqref{eq:energyid}, which gives (ii).
Throughout, $\delta_\ast:=\tfrac1{48}$ and every Young inequality is performed
against $\theta(\tau)$, not against $\sigma_{\min}$.

\smallskip\noindent\emph{Step 1 (diagonal).} Split
$\widehat\Theta\Delta\zeta=\diver(\widehat\Theta\na\zeta)-\na\widehat\Theta\cdot\na\zeta$,
so $D_j=D^{(1)}_j+D^{(2)}_j+D^{(3)}_j$, where $D^{(1)}_j$, $D^{(2)}_j$ and $D^{(3)}_j$ are
the pairings with $v_j$ of $\langle\nabla\rangle^\sigma\diver(\widehat\Theta_j\na\zeta_j)$, of
$-\langle\nabla\rangle^\sigma(\na\widehat\Theta_j\cdot\na\zeta_j)$ and of
$\langle\nabla\rangle^\sigma\big(\widehat\Theta_j(\mathbf b_w\cdot\na\zeta_j+c_w\zeta_j)\big)$. Since constants commute with $\langle\nabla\rangle^\sigma$,
$[\langle\nabla\rangle^\sigma,\widehat\Theta_j]=[\langle\nabla\rangle^\sigma,\widehat\Theta_j-e^{2\tau}\kappa_j]$
and the Kato--Ponce commutator estimate \cite{KaPo88} applies to the decaying function
$\widehat\Theta_j-e^{2\tau}\kappa_j\in H^{\sigma+1}$:
\[
 \big\|[\langle\nabla\rangle^\sigma,\widehat\Theta_j]\na\zeta_j\big\|_{L^2}
 \le C\Big(\|\na\widehat\Theta_j\|_{L^\infty}\|\langle\nabla\rangle^{\sigma-1}\na\zeta_j\|_{L^2}
 +\|\langle\nabla\rangle^\sigma(\widehat\Theta_j-e^{2\tau}\kappa_j)\|_{L^2}\|\na\zeta_j\|_{L^\infty}\Big)
 \le C\|\zeta_j\|_{H^\sigma},
\]
using $\|\langle\nabla\rangle^{\sigma-1}\na\zeta_j\|\le\|\zeta_j\|_{H^\sigma}$ and
$\|\na\zeta_j\|_{L^\infty}\le C\|\zeta_j\|_{H^\sigma}$ ($\sigma>\frac d2+1$).
Hence, with $\widehat\Theta_j\ge\theta(\tau)$,
\begin{equation}\label{eq:S1}
 D^{(1)}_j=-\!\int\!\langle\nabla\rangle^\sigma(\widehat\Theta_j\na\zeta_j)\cdot\na v_j
 \le-\theta\|\na v_j\|^2+C\|\zeta_j\|_{H^\sigma}\|\na v_j\|
 \le-\tfrac34\theta\|\na v_j\|^2+C\|v_j\|^2.
\end{equation}
By the algebra property of $H^\sigma$ and \eqref{eq:sob},
$|D^{(2)}_j|\le C\|\na\widehat\Theta_j\|_{H^\sigma}\|\zeta_j\|_{H^{\sigma+1}}\|v_j\|$,
and $\widehat\Theta_j\mathbf b_w,\widehat\Theta_jc_w$ are multipliers on $H^\sigma$
(a constant times $C^k_b$, plus $H^{\sigma+1}$ times $C^k_b$), whence
$|D^{(3)}_j|\le C\|\zeta_j\|_{H^{\sigma+1}}\|v_j\|$; both give
\begin{equation}\label{eq:S23}
 |D^{(2)}_j|+|D^{(3)}_j|\le 2\delta_\ast\theta\|\na v_j\|^2+C\|v_j\|^2.
\end{equation}

\smallskip\noindent\emph{Step 2 (coupling).} By \eqref{eq:krep} and
$\mathrm T_z\mathcal D=\mathcal D_z\mathrm T_z$ with
$\mathcal D_z:=\Delta+(\mathrm T_z\mathbf b_w)\cdot\na+\mathrm T_zc_w$ (translation does
not commute with $\mathcal D$. It translates the coefficients, whose
$C^k_b$ norms are translation invariant),
\[
 C_{ij}=\int_{\R^d}B^{\eta_\tau}_{ij}(z)\,
 \Big\langle\langle\nabla\rangle^\sigma\big(G_z\,\mathcal D_z \zeta^z_i\big),v_j\Big\rangle\dz,
 \qquad G_z:=\rho_z\,\mathrm T_z\widehat g_{ij},\quad \zeta^z_i:=\mathrm T_z\zeta_i.
\]
By Lemmas~\ref{lem:kernelop} and \ref{lem:dilcost},
$\|G_z\|_{L^\infty}\le\gold^{a}\|\widehat g_{ij}\|_{L^\infty}$ and
$\|G_z\|_{H^{\sigma+1}}\le C\|\widehat g_{ij}\|_{H^{\sigma+1}}\le C$. For fixed
$z$ split $G_z\Delta\zeta^z_i=\diver(G_z\na\zeta^z_i)-\na G_z\cdot\na\zeta^z_i$ and integrate by
parts:
\[
 \big\langle\langle\nabla\rangle^\sigma(G_z\Delta\zeta^z_i),v_j\big\rangle
 =-\!\int\!G_z\na\langle\nabla\rangle^\sigma\zeta^z_i\!\cdot\!\na v_j
 -\!\int\![\langle\nabla\rangle^\sigma,G_z]\na\zeta^z_i\!\cdot\!\na v_j
 -\!\int\!\langle\nabla\rangle^\sigma(\na G_z\!\cdot\!\na\zeta^z_i)v_j.
\]
Since $\mathrm T_z$ is an isometry, $\|\na\langle\nabla\rangle^\sigma\zeta^z_i\|=\|\na v_i\|$. The
remaining, first- and zeroth-order, part of $\mathcal D_z$ contributes
$\langle\langle\nabla\rangle^\sigma(G_z[(\mathrm T_z\mathbf b_w)\cdot\na\zeta^z_i+(\mathrm T_zc_w)\zeta^z_i]),v_j\rangle$,
bounded by $C\|\zeta^z_i\|_{H^{\sigma+1}}\|v_j\|$. For the commutator term
Kato--Ponce gives, as in Step~1,
$\|[\langle\nabla\rangle^\sigma,G_z]\na\zeta^z_i\|_{L^2}\le C\|\zeta^z_i\|_{H^\sigma}=C\|v_i\|$.
Collecting, and using $\int|B^{\eta_\tau}_{ij}|\le A$,
\begin{equation}\label{eq:Cij}
 |C_{ij}|\le A\gold^{a}\|\widehat g_{ij}\|_{L^\infty}\|\na v_i\|\|\na v_j\|
 +C A\Big[\big(\|\na v_i\|+\|v_i\|\big)\|v_j\|+\|v_i\|\,\|\na v_j\|\Big].
\end{equation}
Only the first term of \eqref{eq:Cij} carries a gradient on both factors. Every other term has at most one, because the commutator removes one derivative from $u_i$ and the $\na G_z$ and $\mathcal D_z$ pieces are of lower order. Such terms are $\lesssim\sqrt{ED}$ and are absorbed by an arbitrarily small multiple of $\theta D$.

Sum \eqref{eq:Cij} over $i,j$ \emph{first}, then apply Young once. Using
$\sum_{i,j}\|\na v_i\|\|\na v_j\|=\big(\sum_i\|\na v_i\|\big)^2\le nD$,
$\sum_{i,j}\|\na v_i\|\|v_j\|\le n\sqrt{DE}$ and
$\sum_{i,j}\|v_i\|\|v_j\|\le nE$, together with Lemma~\ref{lem:dilcost} and
\eqref{eq:epscond}, we obtain, with
$\|\widehat g\|_{L^\infty}:=\max_{ij}\|\widehat g_{ij}\|_{L^\infty}$,
\begin{equation}\label{eq:S4}
 \Big|\sum_{i,j}C_{ij}\Big|
 \le nA\gold^{a}\,\frac{\gamma}{\sigma_{\min}}\,\theta\,D
 \;+\;4\delta_\ast\theta D+C_0'E
 \;\le\;\tfrac14\theta D+4\delta_\ast\theta D+C_0'E,
\end{equation}
with $C_0'=C\,A^2n^2/(\delta_\ast\sigma_{\min})$. The factor $e^{2\tau}$ in
$\|\widehat g\|_{L^\infty}$ cancels against the one in $\theta(\tau)$, because Step~1 retained
it. Because Young's parameter is chosen after summing, the coefficient of $\theta D$ is $4\delta_\ast$ independently of $n$ and $A$, and the factor $n^2A^2$ enters only the coefficient of $E$, on which $C_0$ may depend.

\smallskip\noindent\emph{Step 3 (conclusion).} Summing \eqref{eq:S1}, \eqref{eq:S23},
\eqref{eq:S4} and using $6\delta_\ast=\frac18$ (two families from
\eqref{eq:S23}, four from \eqref{eq:S4}),
\[
 \sum_jD_j+\sum_{i,j}C_{ij}
 \le\Big(-\tfrac34+\tfrac14+\tfrac18\Big)\theta D+\tfrac{C_0}2E
 =-\tfrac38\theta D+\tfrac{C_0}2E,
\]
with $C_0$ twice the sum of the coefficients of $E$. This proves \eqref{eq:Sbound}. Under the
hypotheses of (ii), \eqref{eq:energyid} and \eqref{eq:Sbound} give \eqref{eq:diffineq}. Since $\theta\ge\sigma_{\min}$,
$E'+\frac34\sigma_{\min}D\le C_0E$; multiplying by $e^{-C_0\tau}$ and integrating
gives $\sup_\tau E\le e^{C_0T}E(0)$ and
$\frac34\sigma_{\min}\int_0^tD\le e^{C_0T}E(0)$, hence \eqref{eq:apriori} with
constant $1+\frac43\le3$.
\end{proof}

\begin{proposition}[Forced energy estimate]\label{prop:forced-forced}
Let $\zeta\in\mathcal X_{\rm cj}$ satisfy \eqref{eq:energyF} and
\eqref{eq:young} (in particular, let $\zeta$ solve \eqref{eq:conjF}, by Lemma~\ref{lem:identF}). Then
\begin{equation}\label{eq:forced-forced}
 \sup_{\tau\le t}\|\zeta(\tau)\|^2_{H^\sigma}
 +\int_0^t\|\zeta\|^2_{H^{\sigma+1}}\dtau
 \ \le\ C_F(T)\Big(\|\zeta(0)\|^2_{H^\sigma}
 +\int_0^t\|\mathsf F\|^2_{H^{\sigma-1}}\dtau\Big),
\end{equation}
\begin{equation}\label{eq:CF}
 C_F(T):=\Big[(1+T)+\tfrac2{\sigma_{\min}}\big(2+C_0^FT\big)\Big]\,
 e^{C_0^FT}\,\max\Big\{1,\tfrac4{\sigma_{\min}}\Big\},
 \qquad
 C_0^F:=C_0+\tfrac14e^{2T}\sigma_{\min},
\end{equation}
with $C_0$ the constant of Proposition~\ref{prop:apriori}. $C_F$ is independent of
$\eta,t,\zeta,\mathsf F$.
\end{proposition}

\begin{proof}
We insert \eqref{eq:Sbound} and \eqref{eq:young} into \eqref{eq:energyF} and apply
Gr\"onwall's lemma. By Proposition~\ref{prop:apriori}(i),
$\sum_jD_j+\sum_{ij}C_{ij}\le-\frac38\theta D+\frac{C_0}2E$. Doubling
\eqref{eq:energyF} and inserting \eqref{eq:young} with $\lambda=\frac14$, so
that $\frac34\theta D-\frac14\theta D=\frac12\theta D$,
\[
 \frac{\dd}{\dtau}E+\tfrac12\theta D
 \ \le\ \big(C_0+\tfrac14\theta\big)E+\tfrac4\theta\|\mathsf F\|^2_{H^{\sigma-1}}.
\]
Since $\sigma_{\min}\le\theta(\tau)\le e^{2T}\sigma_{\min}$, with
$g:=\frac4{\sigma_{\min}}\|\mathsf F\|^2_{H^{\sigma-1}}$,
\begin{equation}\label{eq:ode}
 \frac{\dd}{\dtau}E+\tfrac12\sigma_{\min}D\ \le\ C_0^FE+g
 \qquad\text{in }\mathcal D'(0,t).
\end{equation}
Put $R:=E(0)+\int_0^tg$. Then $\mathcal E:=\sup_{\tau\le t}E\le e^{C_0^FT}R$, and
integrating \eqref{eq:ode},
$\frac12\sigma_{\min}\int_0^tD\le R+C_0^FT\mathcal E\le(1+C_0^FTe^{C_0^FT})R
\le(2+C_0^FT)e^{C_0^FT}R$. Since
$\int_0^t\|\zeta\|^2_{H^{\sigma+1}}=\int_0^t(E+D)\le T\mathcal E+\int_0^tD$, the left side
of \eqref{eq:forced-forced} is $\le(1+T)\mathcal E+\int_0^tD
\le[(1+T)+\frac2{\sigma_{\min}}(2+C_0^FT)]e^{C_0^FT}R$, and
$R\le\max\{1,\frac4{\sigma_{\min}}\}(\|\zeta(0)\|^2_{H^\sigma}+\int_0^t\|\mathsf F\|^2_{H^{\sigma-1}})$.
\end{proof}

\begin{lemma}[The conjugated form]\label{lem:form}
Put $V:=H^{\sigma+1}$ and $H:=H^\sigma$ with the inner product
$(\zeta,\psi)_{H^\sigma}:=\br{\langle\nabla\rangle^\sigma\zeta,\langle\nabla\rangle^\sigma\psi}$. Identify
$V'$ with $H^{\sigma-1}$ through $\langle f,\psi\rangle_{V',V}:=\langle f,\langle\nabla\rangle^{2\sigma}\psi\rangle_{H^{\sigma-1},H^{-\sigma+1}}$.
This is an isometric isomorphism by Plancherel, and on $H$ it agrees with $(\cdot,\cdot)_{H^\sigma}$,
so $V\subset H\subset V'$ is a Gelfand triple and $\pa_\tau$ in $V'$ is $\pa_\tau$ in $H^{\sigma-1}$. For
$\tau\in[0,t]$ define
\[
 \mathsf a(\tau;\zeta,\psi):=-\big\langle\mathcal F(\tau)\zeta,\langle\nabla\rangle^{2\sigma}\psi\big\rangle_{H^{\sigma-1},H^{-\sigma+1}},
 \qquad\zeta,\psi\in V .
\]
Then, with the constants of \eqref{eq:Fbound} and Proposition~\ref{prop:apriori}:
\begin{enumerate}[label=\normalfont(\roman*)]
\item $|\mathsf a(\tau;\zeta,\psi)|\le C\|\zeta\|_V\|\psi\|_V$, uniformly in $\tau$ and $\eta$;
\item $\tau\mapsto\mathsf a(\tau;\zeta,\psi)$ is continuous on $[0,t]$ for all $\zeta,\psi\in V$;
\item $\mathsf a(\tau;\zeta,\zeta)\ge\frac38\sigma_{\min}\|\zeta\|^2_V-\lambda_0\|\zeta\|^2_H$ for all
$\tau$ and $\zeta\in V$, with $\lambda_0:=\frac38\sigma_{\min}+\frac{C_0}2$.
\end{enumerate}
The operator in $\mathcal L(V,V')$ associated with $\mathsf a(\tau;\cdot,\cdot)$ is
$-\mathcal F(\tau)$. So \eqref{eq:conjF} is the abstract Cauchy problem for this form.
\end{lemma}

\begin{proof}
\emph{(i)} is \eqref{eq:Fbound}, since
$\|\langle\nabla\rangle^{2\sigma}\psi\|_{H^{-\sigma+1}}=\|\psi\|_{H^{\sigma+1}}$.

\emph{(ii).} Since $\sigma\ge1$, $\mathcal D\zeta_j\in H^{\sigma-1}\subset L^2$, and $\mathcal D$ does not depend
on $\tau$. The multipliers $\widehat\Theta(\tau)$ and $\widehat g(\tau)$ are compositions of the
dilations $\mathrm{Dil}_\lambda$ at $\lambda=e^{\pm\tau}$, which are strongly continuous (by continuity
on $C^\infty_0$) and locally bounded (Lemma~\ref{lem:scaling-prop66}) on $L^2$, with multiplication
by $\Theta^\eta(t-\tau)$ and $g^\eta(t-\tau)$, which are bounded and strongly continuous on $L^2$ by
Condition~\ref{ass:budget}(a),(b). For $\eta>0$, \eqref{eq:krep} and the substitution $z=\eta_\tau z'$ give
$\mathcal K^\tau_{ij}h=\int B_{ij}(z)\,\rho_{\eta_\tau z}\,\mathrm T_{\eta_\tau z}h\dd z$ with
$|\rho_z|\le\gold^a$, so $\tau\mapsto\mathcal K^\tau h$ is continuous in $L^2$ by continuity of
translations and of $z\mapsto\rho_z$ and by dominated convergence; for $\eta=0$, $\mathcal K^\tau_{ij}$ is
multiplication by $a_{ij}$ and does not depend on $\tau$. A composition of strongly
continuous, locally bounded operator families is strongly continuous. Hence
$\tau\mapsto\mathcal F(\tau)\zeta$ is continuous into $L^2$, and
$\tau\mapsto\br{\mathcal F(\tau)\zeta,\chi}$ is continuous for $\chi\in L^2$. For
$\chi:=\langle\nabla\rangle^{2\sigma}\psi\in H^{-\sigma+1}$ take $\chi_k\in L^2$ with
$\chi_k\to\chi$ in $H^{-\sigma+1}$. By \eqref{eq:Fbound},
$|\br{\mathcal F(\tau)\zeta,\chi-\chi_k}|\le C\|\zeta\|_V\|\chi-\chi_k\|_{H^{-\sigma+1}}$
uniformly in $\tau$, so $\tau\mapsto\mathsf a(\tau;\zeta,\psi)$ is a uniform limit of continuous functions.

\emph{(iii).} By \eqref{eq:energy},
$-\mathsf a(\tau;\zeta,\zeta)=\sum_jD_j+\sum_{i,j}C_{ij}$, as in the proof of Lemma~\ref{lem:ident}.
By \eqref{eq:Sbound} and $\theta(\tau)\ge\sigma_{\min}$,
$\mathsf a(\tau;\zeta,\zeta)\ge\frac38\sigma_{\min}D-\frac{C_0}2E$, and $D=\|\zeta\|^2_V-\|\zeta\|^2_H$, $E=\|\zeta\|^2_H$.
\end{proof}

\begin{proposition}\label{prop:evolution-prop66}
Let Condition~\ref{ass:budget} hold, $\eta\in[0,\eta_0]$, $t\in[0,T]$. For every
$\varphi\in H^\sigma_a$ there is a unique solution of \eqref{eq:backward-prop66} with
$\zeta=w\,\Xi\in C([0,t];H^\sigma)\cap L^2(0,t;H^{\sigma+1})$, and
\begin{equation}\label{eq:evolbound-prop66}
 \sup_{0\le s\le t}\|\Phi(s)\|_{H^\sigma_a}\le C^{\rm ev}_{\sigma,a}(T)\|\varphi\|_{H^\sigma_a},
 \qquad
 C^{\rm ev}_{\sigma,a}(T):=C(a,\sigma,d)\,e^{(\sigma-\frac d2)T}\sqrt3\,e^{C_0T/2}.
\end{equation}
Setting $\Phi(s)=:S^\eta(t,s)^*\varphi$ and letting $S^\eta(t,s)$ be the adjoint
family on $H^{-\sigma}_{-a}$,
$\sup_{s\le t\le T}\|S^\eta(t,s)\|_{\mathcal L(H^{-\sigma}_{-a})}\le C^{\rm ev}_{\sigma,a}(T)$ and
$S^\eta(t,r)S^\eta(r,s)=S^\eta(t,s)$ for $s\le r\le t$.
\end{proposition}

\begin{proof}
Existence is J.-L.~Lions' theorem for the form of Lemma~\ref{lem:form}. Uniqueness and the bound
follow from the energy estimate \eqref{eq:apriori}, and uniqueness gives the evolution law.

\emph{Existence.} By Lemma~\ref{lem:form}, $\mathsf a$ is measurable in $\tau$, bounded on
$V\times V$ and quasi-coercive. J.-L.~Lions' theorem \cite[Ch.~XVIII, \S3, p.~513]{DaLi92},
\cite[p.~112]{Sho97} (stated for coercive forms; the shift $\lambda_0$ is removed by
passing to $e^{-\lambda_0\tau}\zeta$) gives a $\zeta\in L^2(0,t;H^{\sigma+1})$ with
$\pa_\tau\zeta\in L^2(0,t;H^{\sigma-1})$, $\pa_\tau\zeta=\mathcal F(\tau)\zeta$ and
$\zeta(0)=w\varphi\in H^\sigma$. Lions--Magenes \cite{LiMa72} gives $\zeta\in C([0,t];H^\sigma)$.
By Lemma~\ref{lem:ident} and Proposition~\ref{prop:apriori}(ii), $\zeta$ satisfies
\eqref{eq:apriori}.

\emph{Uniqueness.} The difference of two solutions in
$C([0,t];H^\sigma)\cap L^2(0,t;H^{\sigma+1})$ with the same initial condition solves \eqref{eq:conj} with zero
initial condition. By Lemma~\ref{lem:ident} and Proposition~\ref{prop:apriori}(ii), it satisfies
\eqref{eq:apriori}, so it vanishes.

\emph{The two-parameter law.} \eqref{eq:backward-prop66} has coefficients depending on
$s$ only, not on the terminal time $t$. Uniqueness on $[s,r]$ and $[r,t]$ then
gives $S^\eta(t,r)S^\eta(r,s)=S^\eta(t,s)$.

\emph{The bound.} $\|\Xi(\tau)\|_{H^\sigma_a}=\|\zeta(\tau)\|_{H^\sigma}$ and
$\Psi(\tau,\cdot)=\Xi(\tau,e^{\tau}\cdot)$, so Lemma~\ref{lem:scaling-prop66} with
$\mu=e^\tau\ge1$ gives
$\|\Psi(\tau)\|_{H^\sigma_a}\le Ce^{(\sigma-\frac d2)\tau}\|\zeta(\tau)\|_{H^\sigma}$.
Combine with \eqref{eq:apriori} and $\|\zeta(0)\|_{H^\sigma}=\|\varphi\|_{H^\sigma_a}$.
\end{proof}

\begin{proposition}[Persistence of regularity]\label{prop:persist}
Under Condition~\ref{ass:budget} with upper level,
Proposition~\ref{prop:evolution-prop66} holds simultaneously at the levels $(\sigma,a)$
and $(\sigma+2,a+1)$. Consequently, for $\varphi\in H^{\sigma+2}_{a+1}$ the
function $\Phi(s)=S^\eta(t,s)^*\varphi$ remains in $H^{\sigma+2}_{a+1}$, and
$x\cdot\na\Phi_j(s)\in H^{\sigma+1}_{a}\subset H^{\sigma}_{a}$.
\end{proposition}

\begin{proof}
Condition~\ref{ass:budget} holds at the level $(\sigma+2,a+1)$. Indeed, part (a) at that level
is the bound in $H^{(\sigma+2)+1}=H^{\sigma+3}$ of the upper level, part (b) does not depend on the
level, and part (c) at the weight $a+1$ is \eqref{eq:epscondup}. At that weight
Lemma~\ref{lem:kernelop} enters with $C_\rho(a+1)=\gold^{a+1}$. So
Proposition~\ref{prop:evolution-prop66} holds at $(\sigma+2,a+1)$. The two solutions coincide because
$H^{\sigma+2}_{a+1}\hookrightarrow H^\sigma_a$
($\langle x\rangle^{a}=\langle x\rangle^{-1}\langle x\rangle^{a+1}$ with
$\langle x\rangle^{-1}\in C^\infty_b$) and uniqueness holds at the lower level.
For the last claim, $\langle x\rangle^{a}x_k=\frac{x_k}{\langle
x\rangle}\langle x\rangle^{a+1}$ with $x/\langle x\rangle\in C^\infty_b$, so
multiplication by $x$ is bounded $H^{m}_{a+1}\to H^{m}_{a}$ for every $m$,
losing no derivative. Apply this with $m=\sigma+1$ to $\na\Phi_j\in H^{\sigma+1}_{a+1}$
(and $\na$ maps $H^{\sigma+2}_{a+1}\to H^{\sigma+1}_{a+1}$ because
$\na\langle x\rangle^{a+1}=\langle x\rangle^{a+1}(a+1)x\langle
x\rangle^{-2}$ with $x\langle x\rangle^{-2}\in C^\infty_b$).
\end{proof}

\begin{proposition}[Existence and uniqueness for the forced problem]
\label{prop:existF}
Under Condition~\ref{ass:budget}, for $\eta\in[0,\eta_0]$, $t\in[0,T]$,
$\zeta_0\in H^\sigma$ and $\mathsf F\in L^2(0,t;H^{\sigma-1})$, \eqref{eq:conjF}
has a unique solution, and it obeys \eqref{eq:forced-forced}.
\end{proposition}

\begin{proof}
By Lemma~\ref{lem:form}, \eqref{eq:conjF} is the abstract Cauchy problem for a measurable, bounded and
quasi-coercive form, with $\zeta_0\in H$ and $\mathsf F\in L^2(0,t;V')$. J.-L.~Lions' theorem
\cite[Ch.~XVIII, \S3, p.~513]{DaLi92}, \cite[p.~112]{Sho97} gives a unique solution
$\zeta\in L^2(0,t;H^{\sigma+1})$ with $\pa_\tau\zeta\in L^2(0,t;H^{\sigma-1})$, and Lions--Magenes
\cite{LiMa72} gives $\zeta\in C([0,t];H^\sigma)$, so $\zeta\in\mathcal X_{\rm cj}$. Every solution in
$\mathcal X_{\rm cj}$ has $\pa_\tau\zeta\in L^2(0,t;H^{\sigma-1})$ by \eqref{eq:Fbound}, so uniqueness
holds in $\mathcal X_{\rm cj}$. By Lemma~\ref{lem:identF}, $\zeta$ satisfies \eqref{eq:energyF} and
\eqref{eq:young}, and Proposition~\ref{prop:forced-forced} gives \eqref{eq:forced-forced}.
\end{proof}

\begin{corollary}[Existence at the physical level]\label{cor:physexist}
For every $\varphi\in H^\sigma_a$ and $G\in L^2(0,t;H^{\sigma-1}_a)$,
\eqref{eq:backF} has a unique solution in $\mathcal X_{\rm ph}$, namely
$\mathcal C^{-1}\zeta$ with $\zeta$ from Proposition~\ref{prop:existF}.
\end{corollary}

\begin{proof}
$\|\mathsf F(\tau)\|_{H^{\sigma-1}}=\|G(t-\tau,e^{-\tau}\cdot)\|_{H^{\sigma-1}_a}
\le Ce^{(\frac d2+a)T}\|G(t-\tau)\|_{H^{\sigma-1}_a}$, so
$\mathsf F\in L^2(0,t;H^{\sigma-1})$. Apply Proposition~\ref{prop:existF} and
Lemma~\ref{lem:bij}.
\end{proof}

\begin{corollary}[Dual maximal-regularity bound]\label{cor:maxreg}
For $F\in L^2(0,t;H^{-\sigma-1}_{-a})$ define $\int_0^tS^\eta(t,s)F(s)\,\mathrm ds
\in H^{-\sigma}_{-a}$ weakly by
$\varphi\mapsto\int_0^t\langle F(s),S^\eta(t,s)^*\varphi\rangle\,\mathrm ds$.
Then
\[
 \Big\|\int_0^tS^\eta(t,s)F(s)\,\mathrm ds\Big\|_{H^{-\sigma}_{-a}}
 \le C_{\rm mr}(T)\,\|F\|_{L^2(0,t;H^{-\sigma-1}_{-a})},
 \qquad C_{\rm mr}(T)=C\,e^{(\sigma+1-\frac d2)T}\big(T+\sigma_{\min}^{-1}\big)^{1/2}e^{C_0T/2}.
\]
\end{corollary}

\begin{proof}
Let $\zeta:=\mathcal C\big(S^\eta(t,\cdot)^*\varphi\big)$. By Lemma~\ref{lem:bij} and
Proposition~\ref{prop:evolution-prop66}, $\zeta$ solves \eqref{eq:conj} with $\zeta(0)=w\varphi$. Then
\eqref{eq:apriori} gives $\int_0^t\|\zeta\|^2_{H^{\sigma+1}}=\int_0^t(E+D)
\le(T+\sigma_{\min}^{-1})\,3e^{C_0T}\|\varphi\|^2_{H^\sigma_a}$. Undoing the
dilation at level $\sigma+1$ costs $e^{(\sigma+1-\frac d2)\tau}$ by
Lemma~\ref{lem:scaling-prop66}. Hence
$\|S^\eta(t,\cdot)^*\varphi\|_{L^2(0,t;H^{\sigma+1}_a)}\le C_{\rm mr}(T)\|\varphi\|_{H^\sigma_a}$.
Pair with $F$ and use Cauchy--Schwarz in time.
\end{proof}

\begin{proposition}[Physical forced estimate]\label{prop:physF}
For $G\in L^2(0,t;H^{\sigma-1}_a)$, the solution of \eqref{eq:backF} with
$\varphi=0$ given by Corollary~\ref{cor:physexist} satisfies
\begin{equation}\label{eq:physA-forced}
 \int_0^t\big\|\Phi(s)\big\|^2_{H^{\sigma+1}_a}\dd s
 \ \le\ C_3(T)^2\int_0^t\big\|G(s)\big\|^2_{H^{\sigma-1}_a}\dd s,
 \qquad
 C_3(T):=C(a,\sigma,d)\,e^{(\sigma+1+a)T}\,C_F(T)^{1/2}.
\end{equation}
\end{proposition}

\begin{proof}
By Lemma~\ref{lem:bij}, with $\mu=e^{-\tau}\le1$ at the forcing end and
$\mu=e^{\tau}\ge1$ at the output end,
\begin{equation}\label{eq:endF}
 \|\mathsf F(\tau)\|_{H^{\sigma-1}}\le C\,e^{(\frac d2+a)\tau}\|G(t-\tau)\|_{H^{\sigma-1}_a},
 \qquad
 \|\Phi(t-\tau)\|_{H^{\sigma+1}_a}\le C\,e^{(\sigma+1-\frac d2)\tau}\|\zeta(\tau)\|_{H^{\sigma+1}}.
\end{equation}
The change of variable $s=t-\tau$ has unit Jacobian in time. Hence, by
\eqref{eq:forced-forced} with $\zeta(0)=0$,
\[
 \int_0^t\|\Phi\|^2_{H^{\sigma+1}_a}\dd s
 \le C^2e^{2(\sigma+1-\frac d2)T}C_F(T)\int_0^t\|\mathsf F\|^2_{H^{\sigma-1}}\dtau
 \le C^4e^{2(\sigma+1-\frac d2)T+(d+2a)T}C_F(T)\int_0^t\|G\|^2_{H^{\sigma-1}_a}\dd s,
\]
and $2(\sigma+1-\frac d2)+(d+2a)=2\sigma+2+2a$.
\end{proof}

\begin{lemma}[Duhamel representation on a dense class]\label{lem:duhamel}
Assume Condition~\ref{ass:budget} with upper level, so that Proposition~\ref{prop:persist} holds at $(\sigma,a)$ and delivers
Proposition~\ref{prop:evolution-prop66} at $(\sigma+2,a+1)$. For
$G\in C^1_c\big((0,t);H^{\sigma+2}_{a+1}\big)$ the function
\begin{equation}\label{eq:duh}
 \Phi(s):=\int_s^tS^\eta(r,s)^*G(r)\,\dd r
\end{equation}
is the solution of \eqref{eq:backF} with $\varphi=0$ in the sense of
Corollary~\ref{cor:physexist}.
\end{lemma}

\begin{proof}
We check that $\mathcal A(s)$ maps $H^{\sigma+2}_{a+1}$ boundedly into $H^\sigma_a$, use the
integrated backward identity of each $S^\eta(r,\cdot)^*G(r)$, and integrate it in $r$ with
Fubini's theorem. Uniqueness in Corollary~\ref{cor:physexist} identifies the result.

\emph{The operator maps correctly.} For
$v\in H^{\sigma+2}_{a+1}$: $x\!\cdot\!\na v\in H^{\sigma+1}_a$, because
$\na:H^{\sigma+2}_{a+1}\to H^{\sigma+1}_{a+1}$ and multiplication by $x$ is
bounded $H^{m}_{a+1}\to H^{m}_a$ for every $m$, since
$\br x^ax_k=\frac{x_k}{\br x}\br x^{a+1}$ with $x/\br x\in C^\infty_b$. Hence,
$\Theta^\eta_j\Delta v_j\in H^\sigma_a$ since $\Delta v\in H^\sigma_{a+1}\subset
H^\sigma_a$ and $\Theta^\eta_j-\kappa_j\in H^{\sigma+1}$ multiplies $H^\sigma_a$
into itself ($\sigma>\frac d2$) and $B^\eta_{ij}*(g^\eta_{ij}\Delta v_i)\in
H^\sigma_a$ by Lemma~\ref{lem:kernelop} (whose representation
$\mathcal K h=\int B^{\eta}(z)\rho_z(\mathrm T_zh)\,\dd z$ with
$\operatorname{supp}B^\eta\subset B_1(0)$ is exactly the statement that
convolution is bounded on $H^m_a$ with constant $C_mA\gold^a$). Hence
\begin{equation}\label{eq:Amap}
 \mathcal A(s):H^{\sigma+2}_{a+1}\longrightarrow H^{\sigma}_a
 \quad\text{boundedly, uniformly in }s.
\end{equation}

\emph{Construction.} Fix $r$ and put $\Phi^{(r)}(s):=S^\eta(r,s)^*G(r)$ for
$s\le r$, $:=0$ for $s>r$. By Proposition~\ref{prop:evolution-prop66} at $(\sigma+2,a+1)$,
$\Phi^{(r)}\in C([0,r];H^{\sigma+2}_{a+1})$ with
$\sup_s\|\Phi^{(r)}(s)\|_{H^{\sigma+2}_{a+1}}\le C^{\rm ev}_{\sigma+2,a+1}(T)\|G(r)\|_{H^{\sigma+2}_{a+1}}$,
and it satisfies the integrated backward identity
$\Phi^{(r)}(s)=G(r)+\int_s^r\mathcal A(\varsigma)\Phi^{(r)}(\varsigma)\dd\varsigma$
in $H^\sigma_a$, the integrand being in $L^\infty(0,r;H^\sigma_a)$ by
\eqref{eq:Amap}.

\emph{Integration in $r$.} $r\mapsto\|G(r)\|_{H^{\sigma+2}_{a+1}}$ is continuous
with compact support and $(r,\varsigma)\mapsto\Phi^{(r)}(\varsigma)$ is jointly
measurable (on $\{\varsigma\le r\}$ it is continuous in $\varsigma$ and weakly continuous in $r$, shown below, so its pairings with a fixed vector are jointly Borel, and Pettis' theorem applies; it is $0$ off that set). Hence for
each fixed $\varsigma$ the integral
$\Phi(\varsigma)=\int_\varsigma^t\Phi^{(r)}(\varsigma)\dd r$ converges as a
Bochner integral in $H^{\sigma+2}_{a+1}$, and
$\mathcal A(\varsigma)$, being bounded $H^{\sigma+2}_{a+1}\to H^\sigma_a$ by
\eqref{eq:Amap}, commutes with it. The integrand is bounded by
$C^{\rm ev}_{\sigma+2,a+1}(T)\|G(r)\|_{H^{\sigma+2}_{a+1}}$. It is weakly continuous in $r\in[\varsigma,t]$ with values in
$H^{\sigma+2}_{a+1}$, hence strongly measurable by Pettis' theorem. Indeed, for $r\le r'$,
$S^\eta(r',\varsigma)^*G(r')-S^\eta(r,\varsigma)^*G(r)
=S^\eta(r,\varsigma)^*\big[S^\eta(r',r)^*G(r')-G(r')+G(r')-G(r)\big]$. The integrated backward
identity and \eqref{eq:Amap} give
$\|S^\eta(r',r)^*\psi-\psi\|_{H^\sigma_a}\le C(r'-r)\|\psi\|_{H^{\sigma+2}_{a+1}}$, and
$S^\eta(r,\varsigma)^*$ is bounded on $H^\sigma_a$. So the map is continuous into $H^\sigma_a$
and bounded in $H^{\sigma+2}_{a+1}$, hence weakly continuous in $H^{\sigma+2}_{a+1}$.

Fubini's theorem applies: $(r,\varsigma)\mapsto\mathcal A(\varsigma)\Phi^{(r)}(\varsigma)$ is
jointly measurable, as the composition of the jointly measurable $\Phi^{(r)}(\varsigma)$
with the strongly measurable, uniformly bounded $\mathcal A(\varsigma)$, and it is
dominated on the triangle by $C\,C^{\rm ev}_{\sigma+2,a+1}(T)\|G(r)\|_{H^{\sigma+2}_{a+1}}\in L^1(\dd r)$. It gives
\[
 \Phi(s)=\int_s^tG(r)\dd r
 +\int_s^t\!\!\int_s^r\!\mathcal A(\varsigma)\Phi^{(r)}(\varsigma)\dd\varsigma\dd r
 =\int_s^tG(r)\dd r+\int_s^t\mathcal A(\varsigma)\Phi(\varsigma)\dd\varsigma,
\]
using $\{s\le\varsigma\le r\le t\}$ and $\Phi(\varsigma)=\int_\varsigma^t\Phi^{(r)}(\varsigma)\dd r$
together with \eqref{eq:Amap}. The right-hand side is Lipschitz in $s$ with
values in $H^\sigma_a$ (both integrands lie in $L^\infty(0,t;H^\sigma_a)$), so
$\Phi$ is Lipschitz and $-\pa_s\Phi=\mathcal A(s)\Phi+G(s)$ holds for a.e.\ $s$,
equivalently in $\mathcal D'((0,t);H^\sigma_a)$, with $\Phi(t)=0$.
Since $H^{\sigma+2}_{a+1}\hookrightarrow H^{\sigma+1}_a$, $\Phi\in\mathcal
X_{\rm ph}$, and uniqueness (Corollary~\ref{cor:physexist}) identifies it.
\end{proof}

\begin{corollary}[Two-derivative maximal $L^2$ bound]
\label{cor:maxL2-forced}
Let $\sigma_0:=s_1$ and $a:=a_1$, so that, by \ref{H:budget},
\begin{equation}\label{eq:maxL2hyp}
 \sigma_0+4\le s,
\end{equation}
and put $\sigma:=\sigma_0+1$. For $F\in L^2(0,T;H^{-\sigma_0-2}_{-a})$ define
$\int_0^{\cdot}S^\eta(\cdot,s)F(s)\dd s$ weakly as the element of
$L^2(0,T;H^{-\sigma_0}_{-a})=\big(L^2(0,T;H^{\sigma_0}_{a})\big)^*$ acting by
\begin{equation}\label{eq:weakdef}
 G\ \longmapsto\ \int_0^T\Big\langle F(s),\ \int_s^TS^\eta(r,s)^*G(r)\,\dd r\Big\rangle\dd s.
\end{equation}
Then
$\big\|\int_0^{\cdot}S^\eta(\cdot,s)F(s)\dd s\big\|_{L^2(0,T;H^{-\sigma_0}_{-a})}
\le C_3(T)\|F\|_{L^2(0,T;H^{-\sigma_0-2}_{-a})}$,
with $C_3(T)$ as in \eqref{eq:physA-forced} at $\sigma=\sigma_0+1$.
\end{corollary}

\begin{proof}
We bound the operator $\mathcal E$ below on a dense class with Lemma~\ref{lem:duhamel} and
Proposition~\ref{prop:physF}, extend it by density, and identify the resulting functional with
the Duhamel integral of Corollary~\ref{cor:maxreg}.

By \eqref{eq:maxL2hyp}, $\sigma+3\le s$, so Lemma~\ref{lem:duhamel} applies at
level $\sigma$: for $G$ in the dense class $C^1_c((0,T);H^{\sigma+2}_{a+1})$ the
inner integral of \eqref{eq:weakdef} is the solution of \eqref{eq:backF} with
$\varphi=0$, and Proposition~\ref{prop:physF} bounds
\[
 \mathcal E:\ G\mapsto\Big(s\mapsto\int_s^TS^\eta(r,s)^*G(r)\dd r\Big),
 \qquad \|\mathcal EG\|_{L^2_sH^{\sigma+1}_a}\le C_3(T)\|G\|_{L^2_sH^{\sigma-1}_a}.
\]
$C^1_c((0,T);H^{\sigma+2}_{a+1})$ is dense in $L^2(0,T;H^{\sigma-1}_a)$
($C^\infty_0$ is dense in every $H^m_b$, and $\br\cdot^{b}$ is an isometric
isomorphism $H^m_b\to H^m$), so $\mathcal E$ extends with the same norm. With
$\sigma-1=\sigma_0$, $\sigma+1=\sigma_0+2$, and $(H^m_a)^*=H^{-m}_{-a}$ for
every real $m,a$, \eqref{eq:weakdef} is $G\mapsto\br{F,\mathcal EG}$, of norm
$\le C_3(T)\|F\|$ on $L^2_tH^{\sigma_0}_a$. Finally we identify this element with
$D(r):=\int_0^rS^\eta(r,s)F(s)\dd s$, the element of $H^{-\sigma_0-1}_{-a}$ given by
Corollary~\ref{cor:maxreg} at level $\sigma_0+1$, with $\|D(r)\|\le C_{\rm mr}(T)\|F\|$ for every
$r$. For fixed $f\in H^{-\sigma_0-2}_{-a}$ and $g\in\Test$, the map
$(s,r)\mapsto\br{f,S^\eta(r,s)^*g}$ is continuous in $s$ (Lemma~\ref{lem:strongcont}) and in
$r$ (Corollary~\ref{cor:fwdweak}\textup{(iii)}), hence jointly Borel. Approximating $F$ and $G$ by simple functions, the integrand
$\br{F(s),S^\eta(r,s)^*G(r)}$ is jointly measurable on $\{0\le s\le r\le T\}$ and
integrable. So $r\mapsto D(r)$ is weakly measurable, hence strongly measurable (Pettis),
and Fubini gives, for $G$ in the dense class,
$\eqref{eq:weakdef}=\int_0^T\br{D(r),G(r)}\dd r$. The dense class is dense in
$L^2(0,T;H^{\sigma_0}_a)$, and on it this functional is bounded by
$C_3(T)\|F\|\,\|G\|_{L^2_tH^{\sigma_0}_a}$. Hence $D\in L^2(0,T;H^{-\sigma_0}_{-a})$, it is the
element defined by \eqref{eq:weakdef}, and the bound holds for it.
\end{proof}

\begin{proposition}\label{prop:levels}
Under Assumption~\ref{ass:main}, the SKT pair of Lemma~\ref{lem:coeff} and the scalar pairs
$(\Theta^\eta_i,0)$ satisfy the following:
\begin{enumerate}[label=\normalfont(\roman*)]
\item \emph{$(s_1,a_1)$.} Condition~\ref{ass:budget} holds, so
Propositions~\ref{prop:evolution-prop66} and \ref{prop:existF} and Corollary~\ref{cor:physexist}
apply.
\item \emph{$(s_1+1,a_1)$.} Condition~\ref{ass:budget} holds with upper level, so in addition
Corollaries~\ref{cor:maxreg} and \ref{cor:maxL2-forced}, Propositions~\ref{prop:physF} and
\ref{prop:persist}, and Lemma~\ref{lem:duhamel} apply.
\item \emph{$(s_1+3,a_1+1)$.} Proposition~\ref{prop:evolution-prop66} applies, and its
solutions coincide with those at $(s_1+1,a_1)$.
\end{enumerate}
Moreover $\frac d2+1<\sigma$ at every level, and $s\ge s_1+4>\frac d2+5$.
\end{proposition}

\begin{proof}
Lemma~\ref{lem:coeff} at $\sigma=s_1$ requires $s_1+1\le s$, and gives (i). At $\sigma=s_1+1$,
the upper level requires $\sigma+3\le s$, that is $s_1+4\le s$, and part~(c) at the weight
$a_1+1$ is \eqref{eq:epscondup}, which is the second condition of \ref{H:small}. This gives (ii).
Both requirements on $s$ follow from \ref{H:budget}. Proposition~\ref{prop:persist} at
$(s_1+1,a_1)$ gives (iii). Finally, $\frac d2+1<s_1$ by \eqref{eq:choice}.
\end{proof}

\section{The fixed-time fourth moment}\label{sec:appk4}

\noindent
This appendix proves \eqref{eq:K4-mild}, the hypothesis of
Proposition~\ref{prop:close}, and, in Corollary~\ref{cor:CH2}, Lemma~\ref{lem:CH2}. The route is a coupling between the particle system
and the intermediate system, a fourth-moment estimate for the coupling defect, and
a law-of-large-numbers estimate for the intermediate ensemble, which is i.i.d.\
Lemma~\ref{lem:inputs} collects the inputs from
\cite{chen_holzinger_2021_rigorous_SKT}.

\begin{remark}[Where the logarithmic scale is used]\label{rem:betaceiling}
The scale \ref{H:scaling} is used in two ways. First, the coupling estimates
Proposition~\ref{prop:phi4} and \eqref{eq:phi2} carry the Gr\"onwall factor $e^{C\Lambda T}$,
which comes from the Lipschitz constant $\eta^{-(d+1+\alpha)}$ of the interaction and is
polynomial in $N$ only if $\Lambda\le\delta\log N$. Theorem~\ref{thm:K4},
Corollary~\ref{cor:CH2}, the bound \eqref{eq:cheap2} of Lemma~\ref{lem:taylor} (which
Lemma~\ref{lem:remainders}(b),(c) use) and Proposition~\ref{prop:bracketid}, through
\eqref{eq:phi2}, use them.
Second, Lemma~\ref{lem:polylog} is used for conditions $N^{-1/2}\eta^{-\kappa}\to0$ with fixed
$\kappa$: the self-interaction defect ($\kappa=d+\alpha$, in Lemma~\ref{lem:remainders}(a),
Corollary~\ref{cor:trace}, Lemma~\ref{lem:Mbounds} and Theorem~\ref{thm:noise}) and the smoothed
i.i.d.\ field in Proposition~\ref{prop:bracketid}
($\kappa=\alpha+\frac d2+\sigma$ for any $\sigma>\frac d2$). At $\eta=N^{-\varpi}$ these conditions hold for small
$\varpi>0$. The Gr\"onwall factor does not: $\Lambda=N^{2\varpi(d+1+\alpha)}$ and
$e^{C\Lambda T}$ exceeds every power of $N$. We do not know a coupling estimate without
this factor when the interaction sits in the diffusion coefficient. Relative-entropy
methods such as \cite{JaWa18} bound time marginals rather than pathwise coupling moments,
and Theorem~\ref{thm:K4} uses the latter.
\end{remark}

\begin{lemma}[Scaling of the mollifier]\label{lem:scaling-fourth}
For every $m\ge0$, every $p\in[1,\infty]$ and every $0<\eta\le1$,
\[
 \big\|B^\eta_{ij}\big\|_{H^{m}(\R^d)}\le \eta^{-(d/2+m)}\,\|B_{ij}\|_{H^m(\R^d)},
 \qquad
 \big\|B^\eta_{ij}\big\|_{L^{p}(\R^d)}=\eta^{-d/p'}\,\|B_{ij}\|_{L^p(\R^d)},
\]
$1/p+1/p'=1$. In particular $\|B^\eta_{ij}\|^2_{L^2}=\eta^{-d}\|B_{ij}\|^2_{L^2}$,
$\|B^\eta_{ij}\|^4_{L^4}=\eta^{-3d}\|B_{ij}\|^4_{L^4}$,
$\|B^\eta_{ij}\|_{L^1}=\|B_{ij}\|_{L^1}$ and
$\|\na B^\eta_{ij}\|_{L^\infty}=\eta^{-(d+1)}\|\na B_{ij}\|_{L^\infty}$.
\end{lemma}

\begin{proof}
For the Sobolev norm, with $z=\eta\xi$,
$\|B^\eta_{ij}\|^2_{H^m}=\eta^{-d}\int\br{z/\eta}^{2m}|\widehat B_{ij}(z)|^2\dd z
\le\eta^{-d-2m}\|B_{ij}\|^2_{H^m}$, using $\br{z/\eta}^2\le\eta^{-2}\br z^2$ for
$\eta\le1$. For the Lebesgue norm,
$\int|\eta^{-d}B(x/\eta)|^p\dd x=\eta^{-d(p-1)}\|B\|^p_{L^p}$, i.e.\
$\|B^\eta\|_{L^p}=\eta^{-d/p'}\|B\|_{L^p}$. The gradient case is
$\na B^\eta=\eta^{-d-1}(\na B)(\cdot/\eta)$.
\end{proof}

\begin{lemma}[Translation increments]\label{lem:transl}
For $f\in H^{m+1}(\R^d)$ and $a,b\in\R^d$,
$\ \|f(\cdot-a)-f(\cdot-b)\|_{H^m}\le|a-b|\,\|\na f\|_{H^m}\le|a-b|\,\|f\|_{H^{m+1}}$.
\end{lemma}

\begin{proof}
On the Fourier side the multiplier is $e^{-ia\cdot\xi}-e^{-ib\cdot\xi}$, of
modulus $|e^{-i(a-b)\cdot\xi}-1|\le|(a-b)\cdot\xi|\le|a-b|\,|\xi|$, and
$\int\br\xi^{2m}|\xi|^2|\hat f|^2=\|\na f\|^2_{H^m}$.
\end{proof}

\begin{lemma}[Fourth moment of a sum of independent centred vectors]
\label{lem:sum4}
Let $H$ be a separable real Hilbert space and let $\xi_1,\dots,\xi_M$ be
independent centred $H$-valued random variables with $\E\|\xi_\ell\|^4<\infty$.
Then
\[
 \E\Big\|\sum_{\ell=1}^M\xi_\ell\Big\|^4_H
 \ \le\ \sum_{\ell=1}^M\E\|\xi_\ell\|^4_H
 \ +\ 3\Big(\sum_{\ell=1}^M\E\|\xi_\ell\|^2_H\Big)^2;
\]
for i.i.d.\ $\xi_\ell\sim\xi$ this becomes
$\E\|\sum_\ell\xi_\ell\|^4\le M\,\E\|\xi\|^4+3M^2(\E\|\xi\|^2)^2$.
\end{lemma}

\begin{proof}
This is a Hilbert-space analogue, for $p=4$ and with explicit constants, of Rosenthal's inequality \cite{Ros70}. We give the short proof because we use these constants. Expand: $\E\|\sum\xi\|^4=\sum_{k,l,m,n}\E[\br{\xi_k,\xi_l}\br{\xi_m,\xi_n}]$.
Suppose some index value occurs exactly once in $(k,l,m,n)$, say
$k\notin\{l,m,n\}$. The random vector $Z:=\xi_l\br{\xi_m,\xi_n}$ is Bochner
integrable ($\E\|Z\|\le\E[\|\xi_l\|\|\xi_m\|\|\xi_n\|]<\infty$ by H\"older on
fourth moments) and independent of $\xi_k$, so
$\E[\br{\xi_k,Z}]=\br{\E\xi_k,\E Z}=0$. Only quadruples in which every value
occurs at least twice survive: all four positions equal, contributing
$\sum_\ell\E\|\xi_\ell\|^4$ or two pairs, in the three ways
$\{k{=}l\},\{m{=}n\}$; $\{k{=}m\},\{l{=}n\}$; $\{k{=}n\},\{l{=}m\}$, with
$\ell\ne m$ so no double count. The first contributes
$\sum_{\ell\ne m}\E\|\xi_\ell\|^2\E\|\xi_m\|^2$, the other two
$\sum_{\ell\ne m}\E\br{\xi_\ell,\xi_m}^2
\le\sum_{\ell\ne m}\E\|\xi_\ell\|^2\E\|\xi_m\|^2$ by Cauchy--Schwarz and
independence. Each of the three is $\le(\sum_\ell\E\|\xi_\ell\|^2)^2$.
\end{proof}

\begin{lemma}[Inputs for the coupling]\label{lem:inputs}\leavevmode
\begin{enumerate}[label=\normalfont(\roman*),ref=\thetheorem(\roman*)]
\item\label{H:W} For each fixed $N,\eta$ the particle system and the intermediate system have
unique global strong solutions with a.s.\ continuous paths on $[0,T]$.
\item\label{H:chaos} For every $t\in[0,T]$, the family $\{\bar X^\eta_{\ell,j}(t)\}_{\ell\le N,\,j\le n}$
is independent, and $\bar X^\eta_{\ell,j}(t)$ has density $u_{\eta,j}(t)$.
\item\label{H:U} $\sup_{\eta\le\eta_0}\sup_{t\le T}\max_j\|u_{\eta,j}(t)\|_{L^\infty}\le C_u$, with $C_u$
as in \eqref{eq:Cu},
$\mathrm{Lip}(f_\eta)\le\eta^{-\alpha}$, $\na U_i$ is globally Lipschitz with constant $L_U$ and
$\sigma_i\ge\sigma_{\min}>0$, $f_\eta\ge0$.
\end{enumerate}
\end{lemma}

\begin{proof}
(i) is \cite[Prop.~4(i)]{chen_holzinger_2021_rigorous_SKT}. (ii): the independence is
Lemma~\ref{lem:crossspecies}(ii), and the law of $\bar X^\eta_{\ell,j}(t)$ is $u_{\eta,j}(t)\dx$ at
every $t$ by Lemma~\ref{lem:lawid}(ii). (iii): the bound by $C_u$ is \eqref{eq:Cu}, which
holds at every $t$ for the weakly continuous representative of Lemma~\ref{lem:lawid}(i). The bound
on $\mathrm{Lip}(f_\eta)$ and $f_\eta\ge0$ are \ref{H:A4}, and $\sigma_i\ge\sigma_{\min}>0$ by
\ref{H:A1} and the definition of $\sigma_{\min}$. Finally $\na U_i(x)=-x$ by \ref{H:A2}, so
$L_U=1$.
\end{proof}

\begin{proposition}[Fourth moment, i.i.d.\ case]\label{prop:iid}
For every $t\le T$, every $i,j$ and every $N\ge1$,
\begin{equation}\label{eq:iid}
 \E\big\|B^\eta_{ij}*\bar\fluc^N_j(t)\big\|^4_{H^{s_1}}
 \ \le\ 19\,\|B_{ij}\|_{H^{s_1}}^4\ \eta^{-2d-4s_1}.
\end{equation}
\end{proposition}

\begin{proof}
The field $B^\eta_{ij}*\bar\fluc^N_j(t)$ is a normalised sum of $N$ i.i.d.\ centred
$H^{s_1}$-valued vectors, and Lemma~\ref{lem:sum4} bounds its fourth moment.
Fix $t,i,j$, put $G:=B^\eta_{ij}$ and
$Y_k:=G(\cdot-\bar X^\eta_{k,j}(t))-G*u_{\eta,j}(t)\in H^{s_1}(\R^d)$.
By Lemma~\ref{H:chaos} the $\bar X^\eta_{k,j}(t)$, $k\le N$, are i.i.d.\
with density $u_{\eta,j}(t)$, so the $Y_k$ are
i.i.d.\ and centred
($\E[G(\cdot-\bar X_{1,j})]=G*u_{\eta,j}$ as a Bochner integral), and
\begin{equation}\label{eq:Bconv}
 G*\bar\fluc^N_j=\sqrt N\Big(\tfrac1N\textstyle\sum_kG(\cdot-\bar X_{k,j})-G*u_{\eta,j}\Big)
 =\tfrac1{\sqrt N}\sum_{k=1}^NY_k.
\end{equation}
Translation is an isometry of $H^{s_1}$ and $\|G*u_{\eta,j}\|_{H^{s_1}}\le\|G\|_{H^{s_1}}$ by
Young's inequality with $\|u_{\eta,j}\|_{L^1}=1$, so by Lemma~\ref{lem:scaling-fourth}
\begin{equation}\label{eq:Ybound}
 \|Y_k\|_{H^{s_1}}\le2\|B^\eta_{ij}\|_{H^{s_1}}\le2\,C_B\,\eta^{-(d/2+s_1)}\ \text{a.s.},
 \qquad C_B:=\|B_{ij}\|_{H^{s_1}},
\end{equation}
whence $\E\|Y\|^4_{H^{s_1}}\le16C_B^4\eta^{-2d-4s_1}$; and since $H^{s_1}$ is a Hilbert
space, $\E\|Y\|^2_{H^{s_1}}=\E\|G(\cdot-\bar X)\|^2_{H^{s_1}}-\|G*u_{\eta,j}\|^2_{H^{s_1}}
\le\|G\|^2_{H^{s_1}}\le C_B^2\eta^{-d-2s_1}$. By \eqref{eq:Bconv} and
Lemma~\ref{lem:sum4} with $M=N$,
\[
 \E\Big\|\tfrac1{\sqrt N}\sum_kY_k\Big\|^4_{H^{s_1}}
 \le\frac{\E\|Y\|^4_{H^{s_1}}}{N}+3\big(\E\|Y\|^2_{H^{s_1}}\big)^2
 \le\Big(\frac{16}{N}+3\Big)C_B^4\eta^{-2d-4s_1}\le19\,C_B^4\,\eta^{-2d-4s_1}.\qedhere
\]
\end{proof}

For $p\in\{2,4\}$ and $t\le T$ we denote the moments of the coupling defect
$Y_{k,i}$ of \eqref{eq:dY} by
\begin{equation}\label{eq:phip}
 \phi_p(t):=\max_{k\le N}\max_{i\le n}\E|Y_{k,i}(t)|^p .
\end{equation}

\begin{lemma}[Three-term decomposition]\label{lem:decomp}
Let $L_{k,ij}:=f_\eta(G^N_{k,ij})-f_\eta(H_{k,ij})$, with $G^N_{k,ij}$ and $H_{k,ij}$ as in
\eqref{eq:Sig} and \eqref{eq:Sigbar}. Then $\Sigma^2_{k,i}-\bar\Sigma^2_{k,i}
=2\sum_{j}L_{k,ij}$ and, pathwise,
\begin{equation}\label{eq:Lbounds}
 |L_{k,ij}|\ \le\ \underbrace{C_L\Lambda^{1/2}|Y_{k,i}|}_{=:\ell^{(1)}_{k,i}}
 \ +\ \underbrace{C_L\Lambda^{1/2}\tfrac1N\!\sum_{\ell=1}^N|Y_{\ell,j}|}_{=:\ell^{(2)}_{j}}
 \ +\ \underbrace{\frac{\eta^{-\alpha}}{N}\Big(\Big|\sum_\ell{}'\mathsf D_{(k,i),(\ell,j)}\Big|
 +\delta_{ij}c_0\Big)}_{=:\ell^{(3)}_{k,ij}},
\end{equation}
where $C_L:=\max_{ij}\|\na B_{ij}\|_{L^\infty}$, $c_0:=C_u\max_{ij}\|B_{ij}\|_{L^1}$,
$\Lambda^{1/2}=\eta^{-(d+1+\alpha)}$, $\sum_\ell'$ omits $\ell=k$ when $j=i$ and
runs over all $\ell\le N$ otherwise, and
\[
 \mathsf D_{(k,i),(\ell,j)}:=B^\eta_{ij}\big(\bar X^\eta_{k,i}-\bar X^\eta_{\ell,j}\big)
 -\big(B^\eta_{ij}*u_{\eta,j}\big)\big(\bar X^\eta_{k,i}\big),
 \qquad(\ell,j)\ne(k,i).
\]
\end{lemma}

\begin{proof}
The Lipschitz bound on $f_\eta$ splits $L_{k,ij}$ into two differences. Replacing the particles
$X$ by $\bar X$ gives $\ell^{(1)}$ and $\ell^{(2)}$, and replacing the empirical average of the
intermediate particles by its mean gives $\ell^{(3)}$. By \eqref{eq:Sig} and \eqref{eq:Sigbar},
$\Sigma^2-\bar\Sigma^2=2\sum_jL_{k,ij}$. By
$\mathrm{Lip}(f_\eta)\le\eta^{-\alpha}$,
$|L_{k,ij}|\le\eta^{-\alpha}\big(|G^N_{k,ij}-\bar G^N_{k,ij}|
+|\bar G^N_{k,ij}-H_{k,ij}|\big)$, where $\bar G^N_{k,ij}$ is $G^N_{k,ij}$ with
every $X$ replaced by $\bar X$.

\emph{The first difference ($\ell^{(1)}+\ell^{(2)}$).} By the mean value theorem and
Lemma~\ref{lem:scaling-fourth},
\[
 \big|G^N_{k,ij}-\bar G^N_{k,ij}\big|
 \le\frac1N\!\!\sum_{(\ell,j)\ne(k,i)}\!\!\|\na B^\eta_{ij}\|_{L^\infty}
 \big(|Y_{k,i}|+|Y_{\ell,j}|\big)
 \le C_L\,\eta^{-(d+1)}\Big(|Y_{k,i}|+\frac1N\sum_{\ell=1}^N|Y_{\ell,j}|\Big),
\]
and multiplying by $\eta^{-\alpha}$ gives $\ell^{(1)}_{k,i}+\ell^{(2)}_j$ with
$\eta^{-(d+1+\alpha)}=\Lambda^{1/2}$.

\emph{The second difference ($\ell^{(3)}$).} If $j\ne i$ the sum $\sum_\ell'$ runs over
all $\ell$ and, subtracting and adding $H_{k,ij}$ inside each summand,
$\bar G^N_{k,ij}-H_{k,ij}=\frac1N\sum_\ell \mathsf D_{(k,i),(\ell,j)}$. If $j=i$
the summand $\ell=k$ is absent, so
\[
 \bar G^N_{k,ii}-H_{k,ii}
 =\frac1N\sum_{\ell\ne k}\mathsf D_{(k,i),(\ell,i)}-\frac1N H_{k,ii},
 \qquad |H_{k,ii}|\le\|B^\eta_{ii}*u_{\eta,i}\|_{L^\infty}\le C_u\|B_{ii}\|_{L^1}=c_0.
\]
Multiplying by $\eta^{-\alpha}$ gives $\ell^{(3)}_{k,ij}$.
\end{proof}

\begin{lemma}[Moments of the interaction defect]\label{lem:D}
Conditionally on $\bar X^\eta_{k,i}$
the variables $\{\mathsf D_{(k,i),(\ell,j)}\}_\ell$ are i.i.d.\ and centred, and, with
$M\in\{N-1,N\}$ the number of them,
\[
 \E\big[\mathsf D^2\,\big|\,\bar X^\eta_{k,i}\big]\le C_u\eta^{-d}\|B_{ij}\|^2_{L^2},
 \qquad
 \E\big[\mathsf D^4\,\big|\,\bar X^\eta_{k,i}\big]\le 8\big(C_u\eta^{-3d}\|B_{ij}\|^4_{L^4}+c_0^4\big),
\]
\[
 \E\Big|\sum_\ell{}'\mathsf D\Big|^2\le M\,C_u\eta^{-d}\|B_{ij}\|^2_{L^2},
 \qquad
 \E\Big|\sum_\ell{}'\mathsf D\Big|^4\le 8M\big(C_u\eta^{-3d}\|B_{ij}\|^4_{L^4}+c_0^4\big)
 +3M^2\big(C_u\eta^{-d}\|B_{ij}\|^2_{L^2}\big)^2.
\]
\end{lemma}

\begin{proof}
By Lemma~\ref{H:chaos} the $\bar X^\eta_{\ell,j}$, $(\ell,j)\ne(k,i)$, are
i.i.d.\ with density $u_{\eta,j}$ and independent of $\bar X^\eta_{k,i}$.
Conditioning on $\bar X^\eta_{k,i}=x$, the $\mathsf D$'s are i.i.d.\ and centred. For the
moments, $\E[B^\eta_{ij}(x-\bar X)^2]=\int B^\eta_{ij}(x-y)^2u_{\eta,j}(y)\dd y
\le C_u\|B^\eta_{ij}\|^2_{L^2}$, and centring only decreases the second moment.
Likewise, $\E[B^\eta_{ij}(x-\bar X)^4]\le C_u\|B^\eta_{ij}\|^4_{L^4}$ and
$|\mathsf D|^4\le8(|B^\eta_{ij}(x-\bar X)|^4+c_0^4)$. Lemma~\ref{lem:scaling-fourth} converts
the norms. The two bounds on $\sum_\ell'\mathsf D$ follow from Lemma~\ref{lem:sum4} with $H=\R$, applied
conditionally (the conditional bounds are uniform in $x$) and then integrated.
\end{proof}

\begin{lemma}\label{lem:transfer}
For every $t\le T$,
$\ \E\big\|B^\eta_{ij}*(\fluc^N_j-\bar\fluc^N_j)(t)\big\|^4_{H^{s_1}}
 \le C\,\eta^{-2d-4s_1-4}\,N^2\,\phi_4(t)$,
and the same with $\sup_{t\le T}$ inside both expectations.
\end{lemma}

\begin{proof}
$B^\eta_{ij}*(\fluc^N_j-\bar\fluc^N_j)=\frac1{\sqrt N}\sum_kZ_k$ with
$Z_k:=B^\eta_{ij}(\cdot-X^{N,\eta}_{k,j})-B^\eta_{ij}(\cdot-\bar X^\eta_{k,j})$,
and $\|Z_k\|_{H^{s_1}}\le|Y_{k,j}|\,\|B^\eta_{ij}\|_{H^{s_1+1}}\le
C|Y_{k,j}|\,\eta^{-(d/2+s_1+1)}$ by Lemmas~\ref{lem:transl} and \ref{lem:scaling-fourth}.
By the triangle inequality and Jensen applied to the average $\frac1N\sum_k$,
$\|\cdot\|^4_{H^{s_1}}\le\big(\sqrt N\cdot\frac1N\sum_k\|Z_k\|\big)^4
\le C\eta^{-2d-4s_1-4}N^2\cdot\frac1N\sum_k|Y_{k,j}|^4$; take expectations, or
$\sup_t$ then expectations.
\end{proof}

\begin{proposition}[Coupling estimate at orders $2$ and $4$]\label{prop:phi4}
Let $N\ge N_0$, with $N_0$ as in \eqref{eq:N0}, and let $\phi_p$ be as in \eqref{eq:phip}. Then there are constants
$C,C_2$ depending only on $d,n,\sigma_{\min},L_U,B,f,C_u$ (in particular
not on $T$, $N$ or $\eta$) such that
\begin{align}
 \sup_{t\le T}\phi_2(t)&\ \le\ C\,N^{-1}\,\eta^{-d-2\alpha}\,e^{\,C_2(T+1)\Lambda}
 \ \le\ C\,N^{-1+C_2(T+1)\delta}\,\eta^{-d-2\alpha},\label{eq:phi2}\\
 \sup_{t\le T}\phi_4(t)&\ \le\ C\,N^{-2}\,\eta^{-2d-4\alpha}\,e^{\,C_\star C_2(T+1)\Lambda}
 \ \le\ C\,N^{-2+\eps_\delta}\,\eta^{-2d-4\alpha},\label{eq:phi4}
\end{align}
where $C_\star$, $\eps_\delta$ and $\vartheta_0$ are as in \eqref{eq:smallness}. By
\eqref{eq:N0}, $\eta=\eta(N)\le\eta_0$, with $\eta_0$ as in \ref{H:cutoff}. The constant $C_2$ is $C_2:=2L_U+3dn^2C_L^2/\sigma_{\min}+1$, with $C_L$ as in
Lemma~\ref{lem:decomp}. Moreover $C_\star\le7$ for $\vartheta_0\le\frac16$ and all
$d\ge1$, and $C_\star>2$ always.
\end{proposition}

\begin{proof}
We apply It\^o's formula to $|Y_{k,i}|^2$ and $|Y_{k,i}|^4$ (Step~1), bound the diffusion
defect $\Sigma^{\rm d}_{k,i}$ by the three pieces of Lemma~\ref{lem:decomp} (Step~2), and close
with Gr\"onwall's lemma (Step~3). Step~0 localises, so that the stochastic integrals are
martingales.

\emph{Step 0 (localisation).} By Lemma~\ref{H:W} set
$\tau_R:=\inf\{t\ge0:\max_{k,i}|Y_{k,i}(t)|\ge R\}\wedge T$, a stopping time with
$\tau_R\uparrow T$ a.s. For fixed $\eta$ the diffusion coefficients are
bounded a.s.: $|G^N_{k,ij}|\le\|B^\eta_{ij}\|_{L^\infty}\le C\eta^{-d}$ (the
average of at most $N$ values), $|H_{k,ij}|\le c_0$, and $f_\eta$ is Lipschitz,
so $\Sigma_{k,i},\bar\Sigma_{k,i}\le c_\eta<\infty$. Hence on
$[\![0,\tau_R]\!]$ every stochastic integrand below is bounded (by a constant
depending on $R,\eta,N$), the stopped stochastic integrals
are true martingales, and $\phi^R_p(t):=\max_{k,i}\E|Y_{k,i}(t\wedge\tau_R)|^p\le
R^p<\infty$. All estimates are performed for the stopped processes. The resulting
Gr\"onwall bound is independent of $R$ (the forcing involves only the
intermediate particles), and Fatou's lemma with $\tau_R\uparrow T$ and continuity
of paths gives the statement. We suppress $R$.

\emph{Step 1 (It\^o).} For $y\in\R^d$ and every $d$,
\[
 \sum_{l}\pa_l\pa_l|y|^2=2d,\qquad
 \pa_l|y|^4=4|y|^2y_l,\quad \pa_l\pa_l|y|^4=8y_l^2+4|y|^2,\quad
 \sum_{l}\pa_l\pa_l|y|^4=(8+4d)|y|^2.
\]
With $Y=Y_{k,i}$, $\Sigma^{\rm d}=\Sigma^{\rm d}_{k,i}$ and the drift
$\mathrm b:=\na U_i(\bar X^\eta_{k,i})-\na U_i(X^{N,\eta}_{k,i})$ of \eqref{eq:dY}, which equals $Y$
since $\na U_i(x)=-x$, It\^o's formula on
\eqref{eq:dY} gives
\begin{align}
 \dd|Y|^2&=2Y\!\cdot\!\mathrm b\,\dd t+d\,(\Sigma^{\rm d})^2\dd t
 +2\,\Sigma^{\rm d}\,Y\!\cdot\!\dd W^k_i,\label{eq:ito2}\\
 \dd|Y|^4&=4|Y|^2Y\!\cdot\!\mathrm b\,\dd t+(4+2d)|Y|^2(\Sigma^{\rm d})^2\dd t
 +4|Y|^2\Sigma^{\rm d}\,Y\!\cdot\!\dd W^k_i,\label{eq:ito4}
\end{align}
and $|Y\cdot\mathrm b|\le L_U|Y|^2$ by Lemma~\ref{H:U}. The martingale terms have expectation zero, by Step~0.

\emph{Step 2 (the three pieces).} Since $\Sigma,\bar\Sigma\ge\sqrt{2\sigma_{\min}}$,
\[
 (\Sigma^{\rm d})^2=\frac{(\Sigma^2-\bar\Sigma^2)^2}{(\Sigma+\bar\Sigma)^2}
 \le\frac{(\Sigma^2-\bar\Sigma^2)^2}{8\sigma_{\min}}
 \le c_\Sigma\sum_jL_{k,ij}^2,\qquad c_\Sigma:=\frac{n}{2\sigma_{\min}},
\]
by Lemma~\ref{lem:decomp} and Cauchy--Schwarz over the $n$ species. By
\eqref{eq:Lbounds} and $(a+b+c)^2\le3(a^2+b^2+c^2)$, then Jensen
$\big(\frac1N\sum_\ell|Y_{\ell,j}|\big)^2\le\frac1N\sum_\ell|Y_{\ell,j}|^2$,
\[
 \sum_jL_{k,ij}^2\le3C_L^2\Lambda\Big(n|Y_{k,i}|^2
 +\sum_j\tfrac1N\sum_\ell|Y_{\ell,j}|^2\Big)+3\sum_j|\ell^{(3)}_{k,ij}|^2.
\]
Taking expectations and using $\E|Y_{\ell,j}|^2\le\phi_2$,
$\E[|Y_{k,i}|^2|Y_{\ell,j}|^2]\le\phi_4$ (Cauchy--Schwarz),
\begin{equation}\label{eq:pieces}
 \E\sum_jL_{k,ij}^2\le C_Y\Lambda\phi_2+3nF_2,\qquad
 \E(\Sigma^{\rm d})^2\le c_\Sigma\big(C_Y\Lambda\phi_2+3nF_2\big),\qquad
 \E\big[|Y_{k,i}|^2(\Sigma^{\rm d})^2\big]\le c_\Sigma\Big(C_Y\Lambda\phi_4
 +3\sum_j\E\big[|Y_{k,i}|^2|\ell^{(3)}_{k,ij}|^2\big]\Big),
\end{equation}
with $C_Y:=6nC_L^2$ and $F_p:=\max_{k,i,j}\E|\ell^{(3)}_{k,ij}|^p$. By Lemma~\ref{lem:D},
$\eqref{eq:Lbounds}$ and $\eta^{-d}\le N$,
\begin{equation}\label{eq:L34}
 F_2\le\frac{2\eta^{-2\alpha}}{N^2}\big(NC_u\eta^{-d}\|B\|^2_{L^2}+c_0^2\big)
 \le\frac{C}{N}\eta^{-d-2\alpha},\qquad
 F_4\le\frac{8\eta^{-4\alpha}}{N^4}\Big(3N^2C_u^2\eta^{-2d}\|B\|^4_{L^2}
 +8NC_u\eta^{-3d}\|B\|^4_{L^4}+C(1+N)c_0^4\Big)\le\frac{C}{N^2}\eta^{-2d-4\alpha},
\end{equation}
the three terms in $F_4$ arising respectively from the pair terms, the diagonal
of Lemma~\ref{lem:sum4}, and the deterministic self-correction. The last two are
dominated because $\eta^{-d}\le N$. Finally,
by Young with a weight $\varsigma>0$ at our disposal,
$\E[|Y|^2|\ell^{(3)}_{k,ij}|^2]\le\frac{\varsigma}{2}\phi_4+\frac1{2\varsigma}F_4$.

\emph{Step 3 (Gr\"onwall, and the ratio $C_\star$).} Inserting \eqref{eq:pieces}
into \eqref{eq:ito2}--\eqref{eq:ito4} and taking $\max_{k,i}$ (the bounds are
uniform in $(k,i)$),
\begin{align*}
 \phi_2(t)&\le\int_0^t\big[2L_U+dc_\Sigma C_Y\Lambda\big]\phi_2(r)\dd r
 +3ndc_\Sigma\,t\,F_2,\\
 \phi_4(t)&\le\int_0^t\Big[4L_U+(4+2d)c_\Sigma C_Y\Lambda(1+\vartheta)\Big]\phi_4(r)\dd r
 +\frac{3n(4+2d)c_\Sigma}{2\varsigma}\,t\,F_4.
\end{align*}
Choose $\varsigma:=C_Y\Lambda^{1/2}$. Then $\frac{3n\varsigma}{2}
=\vartheta C_Y\Lambda$ with $\vartheta:=\frac{3n}{2\Lambda^{1/2}}
=\frac{3n}{2}\eta^{\,d+1+\alpha}$, which is $\le\vartheta_0$
once $\eta\le\eta_0$; and the forcing coefficient is
$\frac{3n(4+2d)c_\Sigma}{2C_Y\Lambda^{1/2}}\le\frac{3n(4+2d)c_\Sigma}{2C_Y}$, bounded.
Comparing the two coefficients term by term,
\[
 \frac{4L_U+(4+2d)c_\Sigma C_Y\Lambda(1+\vartheta)}{2L_U+dc_\Sigma C_Y\Lambda}
 \ \le\ \max\Big\{2,\ \frac{4+2d}{d}(1+\vartheta_0)\Big\}
 \ =\ \Big(2+\frac4d\Big)(1+\vartheta_0)\ =\ C_\star,
\]
by the mediant inequality $\frac{p+q}{r+s}\le\max\{\frac pr,\frac qs\}$ for
positive $p,q,r,s$, and since $2+\frac4d>2$. Gr\"onwall's lemma in integral form now gives
$\phi_p(t)\le C\,t\,F_p\,e^{\lambda_pt}$\, with $\lambda_p$ the coefficient in brackets in the bound
on $\phi_p$. By the statement, $C_2=2L_U+dc_\Sigma C_Y+1$. Using
$\Lambda\ge1$ we have $\lambda_2\le(C_2-1)\Lambda$ and $\lambda_4\le C_\star(C_2-1)\Lambda$,
while the prefactor obeys $t\le T\le e^{T}\le e^{\Lambda(T+1)}$; hence the total
exponents are at most $(C_2-1)\Lambda(T+1)+\Lambda(T+1)=C_2(T+1)\Lambda$ and
$C_\star(C_2-1)\Lambda(T+1)+\Lambda(T+1)\le C_\star C_2(T+1)\Lambda$ (the last
step because $C_\star\ge1$). Inserting \eqref{eq:L34} and $\Lambda\le\delta\log
N$ gives \eqref{eq:phi2}--\eqref{eq:phi4}.
\end{proof}

\begin{corollary}[Proof of Lemma~\ref{lem:CH2}]\label{cor:CH2}
Let $L^{\rm full}_{k,ij}:=f_\eta\big((B^\eta_{ij}*\mu^N_j)(X^{N,\eta}_{k,i})\big)
-f_\eta\big((B^\eta_{ij}*u_{\eta,j})(\bar X^\eta_{k,i})\big)$ be the quantity of
Lemma~\ref{lem:CH2}. Then, for $N\ge N_0$ with $N_0$ as in \eqref{eq:N0},
\[
 \max_{k,i}\ \E\int_0^T\sum_j\big|L^{\rm full}_{k,ij}(t)\big|^2\dd t
 \ \le\ \mathsf e(N):=C\,T\,\delta\log N\cdot N^{-1+C_2(T+1)\delta}\,\eta^{-d-2\alpha}
 \ +\ C\,T\,N^{-2}\eta^{-2d-2\alpha},
\]
and $\mathsf e(N)\to0$ under \eqref{eq:smallness}. In particular
$\sup_N\max_{k,i}\E\int_0^T\sum_j|L^{\rm full}_{k,ij}|^2\dd t<\infty$, which is
Lemma~\ref{lem:CH2}.
\end{corollary}

\begin{proof}
$L^{\rm full}$ uses the full convolution
$B^\eta_{ij}*\mu^N_j$, whereas $L_{k,ij}$ of Lemma~\ref{lem:decomp} uses
$G^N_{k,ij}$, which for $j=i$ omits the pair $(\ell,j)=(k,i)$. They differ by
one application of $f_\eta$ at arguments differing by
$\frac1NB^\eta_{ii}(0)=\frac1N\eta^{-d}B_{ii}(0)$, that is, by the self-interaction defect $\varrho^N_{k,i}$ of \eqref{eq:self}. Hence
$|L^{\rm full}_{k,ij}-L_{k,ij}|\le\eta^{-\alpha}\eta^{-d}|B_{ii}(0)|/N
\le CN^{-1}\eta^{-(d+\alpha)}$ pathwise, contributing at most
$CTN^{-2}\eta^{-2(d+\alpha)}$ to the integral. The second term of $\mathsf e(N)$ is this bound,
and it tends to $0$ by Lemma~\ref{lem:polylog}. For $L_{k,ij}$ itself the bound is the first
inequality of \eqref{eq:pieces}, $\E\sum_j|L_{k,ij}|^2\le C_Y\Lambda\phi_2+3nF_2$ with $C_Y$ and $F_2$ as defined there,
together with \eqref{eq:L34}, \eqref{eq:phi2} and $\Lambda\le\delta\log N$.
Since $\mathsf e(N)\to0$, it is bounded on $\{N\ge N_0\}$, which is the range of $\sup_N$
by the convention after Assumption~\ref{ass:main}.
\end{proof}

\begin{theorem}[The fixed-time fourth moment]\label{thm:K4}
Let
$N\ge N_0$, with $N_0$ as in \eqref{eq:N0}. Then the quantity $\overline K_4$ of
\eqref{eq:K4def} satisfies
\begin{equation}\label{eq:main}
 \overline K_4
 \ \le\ C\,N^{\eps_\delta}\,\eta^{-(4d+4s_1+4+4\alpha)}
 \ =\ o\big(N\eta^{2\alpha}\big).
\end{equation}
\end{theorem}

\begin{proof}
Fix $t\le T$. By $\|a+b\|^4\le8(\|a\|^4+\|b\|^4)$, Proposition~\ref{prop:iid},
Lemma~\ref{lem:transfer} and Proposition~\ref{prop:phi4},
\[
 \E\|B^\eta_{ij}*\fluc^N_j(t)\|^4_{H^{s_1}}\ \le\ 152\,\|B_{ij}\|^4_{H^{s_1}}\eta^{-2d-4s_1}
 +C\eta^{-2d-4s_1-4}N^2\cdot N^{-2+\eps_\delta}\eta^{-2d-4\alpha}
 \ \le\ C\,N^{\eps_\delta}\,\eta^{-(4d+4s_1+4+4\alpha)},
\]
uniformly in $t$, the second term dominating. With
$\kappa_0:=4d+4s_1+4+6\alpha$, \eqref{eq:H3} gives
$\eta^{-\kappa_0}\le(\delta\log N)^{\kappa_0/(2(d+1+\alpha))}$, so
$\overline K_4/(N\eta^{2\alpha})\le CN^{\eps_\delta-1}
(\log N)^{\kappa_0/(2(d+1+\alpha))}\to0$ by \eqref{eq:smallness}.
\end{proof}
\bibliographystyle{amsplain}
\bibliography{references_v7}
\end{document}